\documentclass[amsmath,10pt,a4paper]{article}
\usepackage{amsfonts, amssymb, amsmath, amscd}
\usepackage[french]{babel}
\usepackage{amssymb}
\newtheorem{thm}{Th\'eor\`eme}[section]
\newtheorem{conj}[thm]{Conjecture}
\newtheorem{lem}[thm]{Lemme}
\newtheorem{prop}[thm]{Proposition}
\newtheorem{cor}[thm]{Corollaire}

\numberwithin{equation}{section}

\title{Propri\'et\'es de Lefschetz automorphes pour les groupes unitaires et orthogonaux}

\author{N. Bergeron}

\date{}

\begin{document}

\maketitle

\begin{footnotesize}
\begin{center}
{\bf Abstract}
\end{center}

Let $G$ be a connected semisimple group over ${\Bbb Q}$. Given a maximal compact subgroup $K\subset G({\Bbb R})$ - such that $X=G({\Bbb R})/K$ is a Riemannian symmetric space - and a convenient arithmetic subgroup $\Gamma\subset G({\Bbb Q})$, 
one constructs an arithmetic manifold $S=S(\Gamma)=\Gamma\backslash X$.
If $H\subset G$ is a connected semisimple subgroup such that $H({\Bbb R})\cap K$ is maximal compact, 
then $Y=H({\Bbb R})/H({\Bbb R})\cap K$ is a 
symmetric subspace of $X$. For each $g\in G({\Bbb Q})$ one can construct an arithmetic manifold  
$S(H,g)=(H({\Bbb Q})\cap g^{-1}\Gamma g)\backslash Y$ and a natural immersion $j_g \colon S(H,g)\rightarrow S$ induced by the map 
$H({\Bbb A})\rightarrow G({\Bbb A}), h\mapsto gh$. Let us assume that $G$ is anisotropic, which implies that $S$ and $S(H,g)$ 
are compact. Then, for each positive integer $k$, the map $j_g$ induces a restriction map
$$R_g \colon H^{k}(S , {\Bbb C})\rightarrow H^{k}(S(H,g) , {\Bbb C}).$$

In this paper we focus on symmetric spaces associated to the unitary and orthogonal groups, namely $O(p,q)$ and $U(p,q)$, 
and give explicit criterions for the injectivity of the product of the maps $R_g$ (for $g$ running through $G({\Bbb Q})$) 
when restricted to the strongly primitive (in the sense of Vogan and
Zuckerman) part of the cohomology. We also give explicit criterions for the injectivity of the map 
$$H^{k}(S(H) , {\Bbb C}) \rightarrow H^{k+{\rm dim} S - {\rm dim} S(H)} (S  , {\Bbb C})$$
dual to the restriction map $R_e$.

The results we obtain fit into a larger conjectural picture that we describe and which bare a strong analogy with the 
classical Lefschetz Theorems. This may sound quite surprising that such an analogy still exists in the case of the real 
arithmetic manifolds. 

We reduce the global problems mentioned above to
local ones by using Theorems of Burger and Sarnak and isolations properties of cohomological representations in the
automorphic dual. The methods used then are mainly representation-theoretic. 

We finally derive some applications concerning the non vanishing of some cohomology classes in arithmetic manifolds.
\end{footnotesize}

\tableofcontents

\section{Introduction}

Rappelons que si $M$ est une vari\'et\'e projective de dimension complexe $n$, l'\'etude des groupes de cohomologie complexes \footnote{Tous les groupes de cohomologie 
consid\'er\'es dans cet article sont \`a coefficients complexes.} $H^* (M)$ de degr\'e $* \neq n$ peut se r\'eduire \`a l'\'etude des groupes de cohomologie d'une vari\'et\'e
projective de dimension strictement plus petite. Plus pr\'ecisemment, Lefschetz a d\'emontr\'e le c\'el\`ebre r\'esultat suivant.

\medskip

\noindent
{\bf Th\'eor\`eme de Lefschetz} {\it
Soit $H$ un hyperplan g\'en\'erique de l'espace projectif ambiant. Alors, 
\begin{enumerate}
\item l'application naturelle de restriction 
$$H^i (M) \rightarrow H^i (M \cap H)$$
est  {\bf injective} pour $i \leq n-1$, et 
\item l'application naturelle ``cup-produit avec $[M \cap H]$''
$$H^i (M \cap H) \rightarrow H^{i+2} (M)$$
est {\bf injective} pour $i \leq n-2$.
\end{enumerate}
}

\medskip

Le but de cet article (qui est un prolongement de \cite{IRMN}) est de d\'ecrire un ph\'enom\`ene analogue au Th\'eor\`eme de Lefschetz dans le monde automorphe
et pour les groupes unitaires et orthogonaux ({\it i.e.} pour les vari\'et\'es arithm\'etiques associ\'ees aux groupes $U(p,q)$ ou $O(p,q)$). 

\subsubsection*{D\'efinitions des objets}

Dans tous le texte nous d\'esignerons par $G$ un groupe alg\'ebrique r\'eductif, connexe et {\bf anisotrope} sur ${\Bbb Q}$.
Les ad\`eles ${\Bbb A}$ de ${\Bbb Q}$
forment un anneau localement compact, dans lequel ${\Bbb Q}$ se plonge diagonalement comme un sous-anneau.
On peut consid\'erer le groupe $G({\Bbb A})$ des points ad\`eliques de $G$, qui contient $G ({\Bbb Q})$
comme sous-groupe discret.

Nous supposerons, pour simplifier et toujours dans tout le texte, que le groupe r\'eductif $G$ est presque simple sur ${\Bbb Q}$ modulo son
centre. Autrement dit, il n'a pas de sous-groupe distingu\'e, non central et connexe d\'efini sur ${\Bbb Q}$.
Il d\'ecoule de cette hypoth\`ese que tous les facteurs simples de l'alg\`ebre de Lie complexe $\mathfrak{g}$ de
$G$ (modulo son centre) sont isomorphes. Nous supposerons de plus que le groupe $G({\Bbb R})$ des points r\'eels est le
produit (avec intersection finie) d'un groupe compact et d'un groupe r\'eel non compact qui est {\it presque simple}
modulo son centre que l'on suppose compact. Nous notons ce dernier groupe $G^{{\rm nc}}$ (nc signifie ici non compact), et nous 
le supposerons g\'en\'eralement isomorphe soit au groupe $U(p,q)$ soit au groupe $O(p,q)$.

Un {\it sous-groupe de congruence} de $G({\Bbb Q})$ est un sous-groupe de la forme $\Gamma = G({\Bbb Q}) \cap K_f$,
o\`u $K_f$ est un sous-groupe compact ouvert du groupe $G({\Bbb A}_f)$ des points ad\`eliques finis de $G$.
Soient $X_G = G({\Bbb R}) / K_{\infty}$ l'espace sym\'etrique associ\'e au groupe $G$, o\`u $K_{\infty} \subset G({\Bbb R})$
est un sous-groupe compact maximal, et $d_G$ la dimension r\'eelle de $X_G$.

Dans cet article on \'etudie les quotients (compacts, puisque $G$ est anisotrope \footnote{Sauf dans la derni\`ere partie o\`u nous discuterons du cas isotrope.})
$\Gamma \backslash X_G$,
o\`u $\Gamma \subset G({\Bbb Q})$ est un sous-groupe de congruence; ces quotients s'identifient aux composantes
connexes de $G({\Bbb Q}) \backslash G({\Bbb A}) / K_{\infty} K_f = G({\Bbb Q}) \backslash (X_G \times G({\Bbb A}_f ) )/K_f $,
o\`u $K_f \subset G({\Bbb A}_f )$ est un sous-groupe compact ouvert.
Pr\'ecisemment, d\'esignons par $G_f$ l'adh\'erence de $G({\Bbb Q})$ dans le groupe $G({\Bbb A}_f )$. Un sous-groupe
de congruence $\Gamma \subset G({\Bbb Q})$ s'\'ecrit $\Gamma = G({\Bbb Q}) \cap K_f$ o\`u $K_f$ est l'adh\'erence de
$\Gamma$ dans $G_f$ et,
\begin{eqnarray} \label{quotient}
\Gamma \backslash X_G = G({\Bbb Q}) \backslash (X_G \times G_f ) / K_f .
\end{eqnarray}
Plus exactement, on s'interesse ici \`a la cohomologie
(\`a coefficients complexes) $H^* (\Gamma \backslash X_G )$ de ces quotients.

Une fois donn\'e deux sous-groupes de congruence $\Gamma ' \subset \Gamma \subset G({\Bbb Q})$, on obtient un
rev\^etement fini
$$\Gamma ' \backslash X_G \rightarrow \Gamma \backslash X_G$$
qui induit un morphisme injectif
$$H^* (\Gamma \backslash X_G ) \rightarrow H^* (\Gamma ' \backslash X_G )$$
en cohomologie. Les groupes de cohomologies $H^* (\Gamma \backslash X_G )$ (ou
$H^* (G({\Bbb Q}) \backslash (X_G \times G_f ) / K_f )$) forment donc un syst\`eme inductif index\'e par les
sous-groupes de congruence $\Gamma \subset G({\Bbb Q})$ (ou par les sous-groupes compacts ouverts $K_f \subset
G_f$).  En passant \`a la limite (inductive) on d\'efinit
\begin{eqnarray} \label{HSh}
H^* (Sh^0 G) = \lim_{
\begin{array}{c}
\rightarrow \\
\Gamma
\end{array}} H^* (\Gamma \backslash X_G )
= \lim_{
\begin{array}{c}
\rightarrow \\
K_f
\end{array}} H^* (G({\Bbb Q}) \backslash (X_G \times G_f ) / K_f ).
\end{eqnarray}
La notation ci-dessus provient de ce que lorsque l'espace $X_G$ est hermitien, on appelle {\it vari\'et\'e de Shimura}
l'espace topologique
\begin{eqnarray} \label{Sh}
Sh^0 G = \lim_{
\begin{array}{c}
\leftarrow \\
\Gamma
\end{array}} \Gamma \backslash X_G = G({\Bbb Q}) \backslash (X \times G_f ).
\end{eqnarray}
Cet espace est en tout cas toujours bien d\'efini et est un espace
topologique dont on peut consid\'erer la cohomologie de C\v{e}ch et il est de plus d\'emontr\'e dans \cite{Rohlfs} que sa
cohomologie co\"{\i}ncide avec (\ref{HSh}).
Pour ce qui nous concerne, il sera suffisant de consid\'erer que $H^* (Sh^0 G)$ n'est qu'une
notation pour la limite inductive (\ref{HSh}).

\subsubsection*{L'application de restriction et son application duale}

Soit $H \subset G$ un sous-groupe r\'eductif connexe d\'efini sur ${\Bbb Q}$. On suppose que
\begin{eqnarray} \label{cond1}
H({\Bbb R})  \cap K_{\infty} \mbox{ {\it est un sous-groupe compact maximal de }} H({\Bbb R}).
\end{eqnarray}
Alors la restriction \`a $H$ de l'involution de Cartan $\theta$ de $G$ est une involution de Cartan de $H({\Bbb R})$.
On a une d\'ecomposition correspondante
\begin{eqnarray}
\mathfrak{h} =  \mathfrak{k}_H \oplus \mathfrak{p}_H
\end{eqnarray}
avec $\mathfrak{p}_H = \mathfrak{p} \cap \mathfrak{h}$. 

Consid\'erons maintenant $\Gamma = G({\Bbb Q}) \cap K_f$ un sous-groupe de congruence sans torsion de $G({\Bbb Q})$.
Le quotient $S(\Gamma ) = \Gamma \backslash X_G$ est une vari\'et\'e compacte qui s'identifie \`a
$S(K_f )= G({\Bbb Q}) \backslash (X \times G_f ) / K_f$.

Soit $K_f^H \subset H({\Bbb A}_f )$ un sous-groupe compact ouvert. Si $K_f^H \subset K_f$, il existe une application
naturelle $j : S(K_f^H ) \rightarrow S(K_f )$. Puisque $\Gamma$ est sans torsion, l'application $j$ est finie et non ramifi\'ee.
On rappelle le r\'esultat suivant, d\^u \`a Deligne \cite{Deligne}.

\begin{lem}
\'Etant donn\'e $K_f^H \subset H({\Bbb A}_f )$, il existe un sous-groupe compact ouvert $K_f^1 \subset G({\Bbb A}_f )$
avec $K_f^H \subset K_f^1$, tel que l'application naturelle $j' : S(K_f^H ) \rightarrow S(K_f^1 )$ soit injective.
\end{lem}

En particulier, si l'on prend $K_f^H = K_f \cap H({\Bbb A}_f )$, on obtient une application naturelle finie $j : S(K_f^H )
\rightarrow S(K_f )$. Si l'on remplace $K_f$ par un sous-groupe suffisamment petit $K_f^1$, on obtient un diagramme
\begin{eqnarray}
\begin{array}{ccc}
S(K_f^H ) & \stackrel{j'}{\rightarrow} & S(K_f^1 ) \\
               & \stackrel{j}{\searrow}      & \downarrow \pi \\
               &                                         & S(K_f ) ,
\end{array}
\end{eqnarray}
o\`u $\pi$ est la projection naturelle de rev\^etement et $j'$ est injective.

En passant \`a la limite (inductive) sur les $K_f$, les applications $j$ induisent  l'application de restriction
\begin{eqnarray} \label{res}
\mbox{res}_H^G : H^* (Sh^0 G) \rightarrow H^* (Sh^0 H) .
\end{eqnarray}
Pour simplifier la lecture mous utiliserons la notation $Sh^0 H \subset Sh^0 G$ pour r\'esumer que $H$ est un sous-groupe alg\'ebrique r\'eductif,connexe, presque simple modulo son centre compact et v\'erifiant la condition
(\ref{cond1}).

L'application duale \`a l'application de restriction induite par $j$ 
$$H^* (S(K_f^H)) \rightarrow H^* (S(K_f))$$
induit, en passant \`a la limite (inductive) sur les $K_f$, l'application
\begin{eqnarray} \label{cup}
\bigwedge_H^G : H^* (Sh^0 H) \rightarrow H^{*+d_G -d_H} (Sh^0 G)
\end{eqnarray}
``cup-produit avec $[Sh^0 H]$'' (duale \`a (\ref{res})).

Nous aurons besoin de consid\'erer \'egalement une modification de l'application de restriction~: l'application de restriction
virtuelle. Expliquons sa construction. Soit $g\in G({\Bbb Q})$. Fixons $K_f^H$ un compact ouvert de $H({\Bbb A}_f )$, et consid\'erons
l'application $j_g : X_H \times H_f \rightarrow X_G \times G_f$ donn\'ee par $j_g (x,h) = (gx,gh)$. Il est facile de
v\'erifier que $j_g$ induit une application injective $H({\Bbb Q}) \backslash (X_H \times H_f ) /K_f^H \rightarrow
G({\Bbb Q}) \backslash (X_G \times G_f ) /K_f^H$. En supposant que $K_f^H = K_f \cap H({\Bbb A}) \subset K_f$ (o\`u
$K_f$ est un sous-groupe compact ouvert dans $G({\Bbb A}_f )$), on obtient alors une application naturelle
$j_g : S(K_f^H ) \rightarrow S(K_f )$, en utilisant les notations pr\'ec\'edentes. Cette application est finie et non ramifi\'ee.
On peut \'egalement d\'ecrire $j_g$ comme l'application naturelle $(H({\Bbb R}) \cap g^{-1} \Gamma g) \backslash X_H
\rightarrow  \Gamma \backslash X_G$. On obtient de cette mani\`ere toute une famille, param\`etr\'ee par $g\in G({\Bbb Q})$,
de sous-vari\'et\'es de $\Gamma \backslash X_G$ - les images des applications $j_g$. En cohomologie celles-ci
induisent l'application de {\it restriction virtuelle}
$$H^* (S(\Gamma )) \rightarrow \prod_{g\in G({\Bbb Q})} H^* (S_H (g)),$$
o\`u $S_H (g) = (H({\Bbb R}) \cap g^{-1} \Gamma g) \backslash X_H$, et l'application de restriction est d\'eduite de la famille
d'applications $(j_g )$.
En passant \`a la limite (inductive) sur les $\Gamma$, l'application de restriction virtuelle induit l'application de {\it restriction virtuelle} Res$_H^G$~:
\begin{eqnarray} \label{Res}
\mbox{Res}_H^G : H^* (Sh^0 G) \rightarrow \prod_{g\in G({\Bbb Q})} H^* (Sh^0 H).
\end{eqnarray}

\paragraph{Une Conjecture tir\'ee de \cite{IRMN}}

Dans \cite{IRMN} nous avons \'enonc\'e la conjecture suivante et montr\'e qu'elle pouvait \^etre d\'eduite de conjectures classiques sur le spectre automorphe des groupes semi-simples (Conjectures d'Arthur).

\begin{conj} \label{conj rg 1}
Supposons fix\'ees des donn\'ees $Sh^0 H \subset Sh^0 G$ avec $H^{{\rm nc}} = O(k,1)$ (resp. $U(k,1)$) et $G^{{\rm nc}} = O(n,1)$ (resp. $U(n,1)$), o\`u $n\geq k\geq1$ sont
des entiers. Alors,
\begin{enumerate}
\item pour tout entier $i \leq d_H /2$, l'application de restriction virtuelle
$${\rm Res}_H^G : H^i (Sh^0 G) \rightarrow \prod_{g\in G({\Bbb Q})} H^i (Sh^0 H)$$
est {\bf injective};
\item pour tout entier $i \leq d_H - d_G /2$, l'application ``cup-produit avec $[Sh^0 H]$''
$$\bigwedge_H^G : H^i (Sh^0 H) \rightarrow H^{i+d_G - d_H} (Sh^0 G)$$
est {\bf injective}.
\end{enumerate}
\end{conj}

Concernant les r\'esultats partiels inconditionnellement d\'emontr\'es nous renvoyons le lecteur \`a \cite{IRMN}, remarquons n\'eanmoins que dans le
cas {\bf unitaire} le point 1. est un Th\'eor\`eme de Venkataramana \cite{Venky}.

Le lien entre la Conjecture \ref{conj rg 1} et le Th\'eor\`eme de Lefschetz est particuli\`erement \'evident lorsque $k=n-1$, le but de cet article est d'explorer ces ``propri\'et\'es de Lefschetz'' en rang sup\'erieur. 
Dans ce cas la combinatoire est bien plus compliqu\'ee comme le montre d\'ej\`a \cite{BergeronTentative} qui g\'en\'eralise le Th\'eor\`eme de Venkataramana mentionn\'e ci-dessus 
au cas $G^{{\rm nc}} = U(p,q)$. Mais d'un autre c\^ot\'e la ``rigidit\'e" du rang sup\'erieur (plus pr\'ecisemment une g\'en\'eralisation de la propri\'et\'e (T) de Kazhdan d\'ecouverte
par Vogan) permet d'obtenir des r\'esultats inconditionnels en petit degr\'e.  Le but de cet article est alors double. D'abord d\'emontrer des propri\'et\'es de Lefschetz en petit degr\'e
puis d\'egager des conjectures concernant les autres degr\'es.

\subsubsection*{\'Enonc\'es des principaux r\'esultats}

Concernant les groupes unitaires nos m\'ethodes nous permettrons de d\'emontrer le th\'eor\`eme suivant.

\begin{thm} \label{upq}
Supposons fix\'ees des donn\'ees $Sh^0 H_i \subset Sh^0 G$ ($i=1,2$) avec $H_1^{\rm nc} = U(p,q-1)$, $H_2^{\rm nc} = U(p-1,q)$ et $G^{{\rm nc}} = U(p,q)$ avec $p,q \geq 2$. Alors,
\begin{enumerate}
\item pour tout entier $\mathbf{k < p+q-1}$, l'application de restriction virtuelle
$${\rm Res}_{H_1}^G \times {\rm Res}_{H_2}^G : H^k (Sh^0 G) \rightarrow \prod_{g\in G({\Bbb Q})} H^k (Sh^0 H_1) \times \prod_{g\in G({\Bbb Q})} H^k (Sh^0 H_2 )$$
est {\bf injective};
\item pour tout entier $\mathbf{k<q-p-1}$ (resp. $\mathbf{k<p-q-1}$), 
l'application ``cup-produit avec $[Sh^0 H_1]$ (resp. $[Sh^0 H_2]$)''
$$\bigwedge_{H_1}^G : H^k (Sh^0 H_1) \rightarrow H^{k+2p} (Sh^0 G)$$
$$(\mbox{resp. } \bigwedge_{H_2}^G : H^k (Sh^0 H_2) \rightarrow H^{k+2q} (Sh^0 G) )$$
est {\bf injective}.
\end{enumerate}
\end{thm}

De mani\`ere plus surprenante un r\'esultat analogue est \'egalement vrai pour les groupes orthogonaux.

\begin{thm} \label{opq}
Supposons fix\'ees des donn\'ees $Sh^0 H_i \subset Sh^0 G$ ($i=1,2$) avec $H_1^{\rm nc} = O(p,q-1)$, $H_2^{\rm nc} = O(p-1,q)$ et $G^{{\rm nc}} = O(p,q)$ avec $p,q \geq 3$. Alors,
\begin{enumerate}
\item pour tout entier $\mathbf{k \leq p+q-4}$, l'application de restriction virtuelle
$${\rm Res}_{H_1}^G \times {\rm Res}_{H_2}^G : H^k (Sh^0 G) \rightarrow \prod_{g\in G({\Bbb Q})} H^k (Sh^0 H_1) \times \prod_{g\in G({\Bbb Q})} H^k (Sh^0 H_2)$$
est {\bf injective};
\item pour tout entier $\mathbf{k\leq (q -p-3)/2}$ (resp. $\mathbf{k\leq (p-q-3)/2}$), l'application ``cup-produit avec $[Sh^0 H_1]$ (resp. $[Sh^0 H_2 ]$)''
$$\bigwedge_{H_1}^G : H^k (Sh^0 H_1 ) \rightarrow H^{k+p} (Sh^0 G)$$
$$(\mbox{resp. } \bigwedge_{H_2}^G : H^k (Sh^0 H_2) \rightarrow H^{k+q} (Sh^0 G) )$$
est {\bf injective}.
\end{enumerate}
\end{thm}

Le cas des groupes $O(2,n)$ ($n\geq 3$) est l\'eg\`erement diff\'erent nous obtenons le r\'esultat suivant.

\begin{thm} \label{o2n}
Supposons fix\'ees des donn\'ees $Sh^0 H \subset Sh^0 G$ avec $H^{\rm nc} = O(2,n-1)$ et $G^{{\rm nc}} = O(2,n)$ avec $n \geq 3$. Alors,
\begin{enumerate}
\item pour tout entier $\mathbf{k \leq n-1}$, l'application de restriction virtuelle
$${\rm Res}_H^G : H^k (Sh^0 G) \rightarrow \prod_{g\in G({\Bbb Q})} H^k (Sh^0 H)$$
est {\bf injective};
\item pour tout entier $\mathbf{k\leq \left[ \frac{n}{2} \right] -2}$, l'application ``cup-produit avec $[Sh^0 H]$''
$$\bigwedge_H^G : H^k (Sh^0 H) \rightarrow H^{k+2} (Sh^0 G)$$
est {\bf injective}.
\end{enumerate}
\end{thm}

En prenant $k=0$ dans le point 2. des Th\'eor\`emes \ref{opq} et \ref{o2n} on retrouve le r\'esultat suivant d\^u \`a Millson et Raghunathan \cite{MillsonRaghunathan} \footnote{Pour $q$
petit le Corollaire ne d\'ecoule pas directement des Th\'eor\`emes il faut travailler un petit peu plus, cf. \S 8.5.}.

\begin{cor} \label{C1}
Supposons fix\'ees les donn\'ees $Sh^0 H \subset Sh^0 G$ avec $H^{{\rm nc}} = O(p,q-1)$ et $G^{{\rm nc}} = O(p,q)$ avec $1\leq p\leq q$. Alors,
la classe fondamentale de $Sh^0 H$ est non triviale dans $H^p (Sh^0 G)$.
\end{cor}

Les Th\'eor\`emes \ref{opq} et \ref{o2n} sont en particulier une vaste g\'en\'eralisation de ce Corollaire. Le point
le plus surprenant est peut-\^etre que mises bout \`a bout les applications des points 1. et 2. impliquent une
sorte de d\'ecomposition de Lefschetz dans un cadre r\'eel.
Les Th\'eor\`emes \ref{upq} et \ref{opq} permettent en tout cas 
de comprendre g\'eom\'etriquement la fa\c{c}on dont certaines classes de cohomologie apparaissent dans l'esprit du Corollaire
ci-dessus. Il n'est pas facile en g\'en\'eral d'exhiber
des classes de cohomologie non triviales dans les vari\'et\'es arithm\'etiques associ\'ees aux groupes orthogonaux. 
Lorsque celles-ci proviennent de
groupes orthogonaux sur des corps de nombres, le Corollaire ci-dessus permet de telles constructions. 
C'est d'ailleurs historiquement le premier r\'esultat concernant ce
probl\`eme. Lorsqu'elles proviennent d'autres constructions Raghunathan et Venkataramana ont remarqu\'es qu'il pouvait \^etre 
utile de les plonger dans des
vari\'et\'es arithm\'etiques associ\'ees aux groupes unitaires. 
Le th\'eor\`eme suivant \'eclaire les relations entre leurs groupes de cohomologie.

\begin{thm} \label{u et o}
Supposons fix\'ees les donn\'ees $Sh^0 H \subset Sh^0 G$ avec $H^{{\rm nc }} = O(p,q)$, $G^{{\rm nc}} = U(p,q)$ et $p,q \geq 3$. Alors,
\begin{enumerate}
\item pour tout entier $\mathbf{k\leq p+q-3}$, l'application de restriction virtuelle
$$\mbox{Res}_H^G : H^{k,0} (Sh^0 G) \rightarrow \prod_{g\in G({\Bbb Q})} H^k (Sh^0 H)$$
est {\bf injective};
\item si $pq$ est {\bf pair}, la classe de $Sh^0 H$ dans $H^{pq} (Sh^0 G)$
est non triviale.
\end{enumerate}
\end{thm}

L'hypoth\`ese de parit\'e de $pq$ du point 2. est surement superflue mais nous ne savons pas comment s'en passer. 
L\`a encore le cas o\`u $p=2$ et $q\geq 3$
est l\'eg\`erement diff\'erent, nous montrons le th\'eor\`eme suivant.

\begin{thm} \label{u et o2}
Supposons fix\'ees les donn\'ees $Sh^0 H \subset Sh^0 G$ avec $H^{{\rm nc }} = O(2,n)$, $G^{{\rm nc}} = U(2,n)$ et $n \geq 3$. Alors,
\begin{enumerate}
\item pour tout entier $\mathbf{k\leq \left[ \frac{n}{2} \right]}$, l'application de restriction virtuelle
$$\mbox{Res}_H^G : H^{k,0} (Sh^0 G) \rightarrow \prod_{g\in G({\Bbb Q})} H^k (Sh^0 H)$$
est {\bf injective};
\item la classe de $Sh^0 H$ dans $H^{2n} (Sh^0 G)$
est non triviale.
\end{enumerate}
\end{thm}

En prenant $k=p$ dans les Th\'eor\`emes \ref{u et o} et \ref{u et o2} on obtient le
corollaire int\'eressant suivant \footnote{Pour $q$
petit le Corollaire ne d\'ecoule pas directement du Th\'eor\`eme il faut travailler un petit peu plus, cf. \S 8.2.}.

\begin{cor} \label{CORUO}
Supposons fix\'ees les donn\'ees $Sh^0 H \subset Sh^0 G$ avec $H^{{\rm nc }} = O(p,q)$, $G^{{\rm nc}} = U(p,q)$ et $1 \leq p\leq q $.
Alors, l'application de restriction virtuelle
$$\mbox{Res}_H^G : H^{p,0} (Sh^0 G) \rightarrow \prod_{g\in G({\Bbb Q})} H^p (Sh^0 H)$$
est {\bf injective}.
\end{cor}

De ce Corollaire et d'un Th\'eor\`eme classique de Borel et Wallach \cite{BorelWallach} 
sur la cohomologie des vari\'et\'es arithm\'etiques associ\'ees aux groupes unitaires, nous d\'eduisons
le Corollaire suivant.

\begin{cor} \label{application}
Soit $G$ un groupe alg\'ebrique sur ${\Bbb Q}$ obtenu par restriction des scalaires \`a partir d'un groupe de type $D_n$, $\neq {}^{3,6} D_4$ et $n > 2$, et tel que
$G^{{\rm nc}} = O(p,q)$, avec $1\leq p \leq q$. Alors,
$$H^p (Sh^0 G) \neq 0.$$
\end{cor}

Ce r\'esultat n'est nouveau que pour $n=3$ et $p >1$, les autres cas sont trait\'es dans diff\'erents articles de Li \cite{Li2}, 
Raghunathan et Venkataramana \cite{RaghunathanVenky} et Li-Millson \cite{LiMillson}.
La d\'emonstration que l'on en donne ici permet de traiter tous ces cas de mani\`ere uniforme, notons que ce Corollaire reste vrai
si $G$ est isotrope avec la m\^eme d\'emonstration (cf. \S 9.1).

\medskip

Nos m\'ethodes s'appliquent \'egalement au produit dans la cohomologie. Dans \cite{Produit}, nous avons pu d\'emontrer le Th\'eor\`eme suivant.

\begin{thm} \label{cup u}
Supposons fix\'ee une donn\'ee $Sh^0 G$ avec $G^{{\rm nc }} = U(p,q)$. Soient $\alpha$ et $\beta$ deux classes de cohomologie de degr\'es respectifs $k$ et $l$ dans
$H^* (Sh^0 G)$ avec $\mathbf{k+l \leq  q+p-1}$. Il existe alors un \'el\'ement $g \in G({\Bbb Q})$
tel que
$$g(\alpha ) \wedge \beta \neq 0 .$$
\end{thm}

Les m\'ethodes (diff\'erentes) de cet article nous permettrons de d\'emontrer (en rang sup\'erieur) le Th\'eor\`eme analogue suivant pour les groupes orthogonaux.

\begin{thm} \label{cup o}
Supposons fix\'ee une donn\'ee $Sh^0 G$ avec $G^{{\rm nc }} = O(p,q)$, avec $p,q \geq 2$. Soient $\alpha$ et $\beta$ deux classes de cohomologie de degr\'es respectifs $k$ et $l$ dans
$H^* (Sh^0 G)$ avec $\mathbf{k+l \leq  q+p-3}$. Il existe alors un \'el\'ement $g \in G({\Bbb Q})$
tel que
$$g(\alpha ) \wedge \beta \neq 0 .$$
\end{thm}

Concluons cette description des r\'esultats en remarquant qu'en rang $1$ la conjecture suivante d\'ecoule elle aussi des conjectures d'Arthur.

\begin{conj} \label{cup hyp}
Supposons fix\'ee une donn\'ee $Sh^0 G$ avec $G^{{\rm nc }} = O(1,n)$. Soient $\alpha$ et $\beta$ deux classes de cohomologie de degr\'es respectifs $k$ et $l$ dans
$H^* (Sh^0 G)$ avec $\mathbf{k+l \leq  n/2}$. Il existe alors un \'el\'ement $g \in G({\Bbb Q})$
tel que
$$g(\alpha ) \wedge \beta \neq 0 .$$
\end{conj}

\paragraph{Remerciements} Je tiens \`a remercier Patrick Delorme et Jacques Carmona qui m'ont signal\'e une erreur dans une premi\`ere
version du texte. J'ai contourn\'e le probl\`eme pos\'e par cette erreur \`a l'aide de la correspondance theta, cf. \S 5.1,
l'approche suivit dans ce paragraphe est sous-jacente dans \cite{HarrisLi}, Michael Harris et Jian-Shu Li 
m'ont d'ailleurs inform\'e qu'ils avaient eux aussi pens\'e \`a cette m\'ethode apr\`es l'\'ecriture de \cite{HarrisLi}. Enfin, merci
\`a Laurent Clozel, cet article s'inscrit \'evidemment dans le prolongement de notre travail en commun \cite{BergeronClozel}.

\section{Repr\'esentations cohomologiques}

Puisque $G$ est anisotrope sur ${\Bbb Q}$, un th\'eor\`eme de Borel et Harish-Chandra \cite{BorelHarishChandra}
affirme que si $\Gamma$ est un sous-groupe de congruence de $G$, la vari\'et\'e $S(\Gamma )$ est compacte.

Dans cette section les facteurs compact de $G({\Bbb R})$ ne nous int\'eresserons pas, nous noterons $\mathfrak{g}_0$ l'alg\`ebre de Lie de
$G^{{\rm nc}}$ et $\mathfrak{k}_0$ l'alg\`ebre de Lie de $K$ un sous-groupe compact maximal de $G^{{\rm nc}}$. Soit
$\mathfrak{g}_0 = \mathfrak{k}_0 \oplus \mathfrak{p}_0$ la d\'ecomposition de Cartan associ\'ee. Si $\mathfrak{l}_0$ est une alg\`ebre de Lie, nous
noterons $\mathfrak{l}=\mathfrak{l}_0 \otimes {\Bbb C}$ sa complexification.

Soit $\Gamma$ un sous-groupe de congruence de $G$.
Soit ${\cal E}^k (S(\Gamma ))$ l'espace des formes diff\'erentielles de degr\'e $k$ sur $S(\Gamma )$. Puisque le
fibr\'e cotangent $T^* S(\Gamma )$ est isomorphe au fibr\'e $\Gamma \backslash G^{{\rm nc}} \times_K \mathfrak{p}^*
\rightarrow \Gamma \backslash G^{{\rm nc}} /K = S(\Gamma )$, qui est associ\'e au $K$-fibr\'e principal $K \rightarrow \Gamma \backslash G^{{\rm nc}}
\rightarrow \Gamma \backslash G^{{\rm nc}} / K$ et \`a la repr\'esentation de $K$ dans $\mathfrak{p}^*$,
on a~:
\begin{eqnarray} \label{Ek}
{\cal E}^k (S(\Gamma )) \simeq (C^{\infty} (\Gamma \backslash G^{{\rm nc}} )\otimes \bigwedge {}^k \mathfrak{p}^* )
\simeq \mbox{Hom}_K (\bigwedge {}^k \mathfrak{p} , C^{\infty} (\Gamma \backslash G^{{\rm nc}}) ) \; \; (k \in {\Bbb N} ).
\end{eqnarray}

Notons $\Delta$ le laplacien de Hodge-de Rham sur la vari\'et\'e riemannienne (localement sym\'etrique) $S(\Gamma )$
(o\`u la m\'etrique est d\'eduite de la forme de Killing sur $\mathfrak{g}_0$). L'espace des formes harmoniques de
degr\'e $k$ est donn\'e par
$${\cal H}^k (S(\Gamma )) := \{ \omega \in {\cal E}^k (S(\Gamma )) \; : \; \Delta \omega = 0 \} .$$
La th\'eorie de Hodge fournit un isomorphisme naturel
$$H^* (S(\Gamma )) \simeq {\cal H}^* (S(\Gamma )) .$$
Soit $(\pi , V_{\pi} )$ un $(\mathfrak{g} , K)$-module irr\'eductible. \`A l'aide de (\ref{Ek}) on d\'efinit une application
lin\'eaire
\begin{eqnarray} \label{Tpi}
T_{\pi} : \left\{
\begin{array}{rcl}
\mbox{Hom}_K (\bigwedge^* \mathfrak{p}, \pi ) \otimes \mbox{Hom}_{\mathfrak{g} , K} (\pi , C^{\infty} (\Gamma \backslash G^{{\rm nc}} ))  & \rightarrow & {\cal E}^* (S(\Gamma )) , \\
\psi \otimes \varphi & \mapsto & \varphi \circ \psi .
\end{array} \right.
\end{eqnarray}

Soit $\widehat{G}^{{\rm nc}}$ l'ensemble des classes d'\'equivalence des $(\mathfrak{g} , K)$-modules irr\'eductibles
qui sont unitarisables. Rappelons qu'Harish-Chandra a d\'emontr\'e que $\widehat{G}^{{\rm nc}}$ s'identifie naturellement
au dual unitaire de la composante connexe de l'identit\'e $G_0^{{\rm nc}}$ de  $G^{{\rm nc}}$. Soient $U(\mathfrak{g} )$ l'alg\`ebre enveloppante de l'alg\`ebre de Lie complexe
$\mathfrak{g}$, $Z(\mathfrak{g} )$ son centre et $\Omega \in Z(\mathfrak{g} )$ le casimir d\'efinit par la forme de Killing sur $\mathfrak{g}_0$.
On d\'efinit le sous-ensemble ${}_0 \widehat{G}^{{\rm nc}}$ de $\widehat{G}^{{\rm nc}}$ par
$${}_0 \widehat{G}^{{\rm nc}} := \{ \pi \in \widehat{G}^{{\rm nc}} \; : \;  \pi (\Omega ) =0 \} ,$$
o\`u l'on a conserv\'e la m\^eme notation $\pi$ pour la repr\'esentation de $U(\mathfrak{g} )$.

L'action de $G_0^{{\rm nc}}$ sur $X_G$ induit la repr\'esentation de $U(\mathfrak{g} )$ sur l'espace des formes diff\'erentielles
sur $X_G$. En particulier, le casimir $\Omega ( \in Z (\mathfrak{g} ) \subset U(\mathfrak{g} ) )$ agit sur ${\cal E}^* (X_G )$
comme le laplacien de Hodge-de Rham, puisque la m\'etrique riemannienne sur $X_G$ est induite par la forme de Killing sur $\mathfrak{g}_0$.
Il d\'ecoule de tout ceci que
\begin{eqnarray} \label{im Tpi}
\mbox{Image} (T_{\pi}) \subset {\cal H}^* (S(\Gamma )) \simeq H^* (S(\Gamma ))
\end{eqnarray}
si et seulement si $\pi \in {}_0 \widehat{G}^{{\rm nc}}$. On dit dans ce cas que le sous-espace de $H^* (S(\Gamma ))$ correspondant
\`a l'image de $T_{\pi}$ est {\it la $\pi$-composante}, et on \'ecrit $H^* ( \pi : \Gamma )$. Autrement dit,
\begin{eqnarray} \label{pi composante}
H^k (\pi : \Gamma ) := \mbox{Image} (T_{\pi}) \cap H^k (S(\Gamma )) \; \;  (k \in {\Bbb N} ) ,
\end{eqnarray}
via l'isomorphisme (\ref{im Tpi}).

Un r\'esultat d\^u \`a Gel'fand et Piatetski-Shapiro \cite{GGPS} affirme que la repr\'esentation r\'eguli\`ere droite
dans $L^2 (\Gamma \backslash G_0^{{\rm nc}} )$ admet une d\'ecomposition en somme directe de Hilbert discr\`ete
$$L^2 (\Gamma \backslash G_0^{{\rm nc}} ) \simeq \sum {}^{\oplus} \mbox{Hom}_G (\pi , L^2 (\Gamma \backslash G_0^{{\rm nc}} )) \otimes \pi
= \sum {}^{\oplus} n_{\Gamma} (\pi ) \pi ,$$
o\`u $\pi$ parcourt cette fois le dual unitaire de $G_0^{{\rm nc}}$ et la multiplicit\'e
$$n_{\Gamma} (\pi ) := \mbox{dim}_{{\Bbb C}} \mbox{Hom}_G (\pi , L^2 (\Gamma \backslash G_0^{{\rm nc}} )) < \infty .$$
Alors la formule de Matsushima est r\'esum\'ee dans le lemme suivant.

\begin{lem}[\cite{BorelWallach}, \cite{Matsushima}]
Sous les notations pr\'ec\'edentes. On a
\begin{eqnarray} \label{lem1}
H^* (\pi : \Gamma ) \simeq n_{\Gamma} H^* (\mathfrak{g} , K ; \pi ) ,
\end{eqnarray}
\begin{eqnarray} \label{lem2}
H^* (S(\Gamma )) = \bigoplus_{\pi \in {}_0 \widehat{G}^{{\rm nc}}}  H^* (\pi : \Gamma ).
\end{eqnarray}
\end{lem}

En passant \`a la limite (inductive) sur les sous-groupes de congruence $\Gamma \subset G$, nous parlerons
de $\pi$-composante de la cohomologie $H^* (Sh^0 G)$ de la vari\'et\'e de Shimura $Sh^0 G$, ce que nous noterons
$H^* (\pi : Sh^0 G)$. La d\'ecomposition (\ref{lem2}) se traduit alors en
\begin{eqnarray} \label{dec cohom}
H^* (Sh^0 G) = \bigoplus_{\pi \in {}_0 \widehat{G}^{{\rm nc}}}  H^* (\pi : Sh^0 G ).
\end{eqnarray}

\medskip

D'apr\`es Parthasarathy \cite{Parthasarathy}, Kumaresan \cite{Kumaresan} et Vogan-Zuckerman \cite{VoganZuckerman},
les $(\mathfrak{g} , K)$-modules unitarisables ayant des groupes de $(\mathfrak{g}, K)$-cohomologie non triviaux peuvent
\^etre d\'ecrit comme suit. Notons $\mathfrak{t}_0 = $Lie$(T)$ une sous-alg\`ebre de Cartan de $\mathfrak{k}_0$.
On consid\`ere les sous-alg\`ebres paraboliques $\theta$-stable $\mathfrak{q} \subset \mathfrak{g}$~: $\mathfrak{q}  =
\mathfrak{l}  \oplus \mathfrak{u}$ \cite{VoganZuckerman}, o\`u $\mathfrak{l}$ est le centralisateur d'un \'el\'ement
$X\in i \mathfrak{t}_0$ et $\mathfrak{u}$ est le sous-espace engendr\'e par les racines positives de $X$ dans $\mathfrak{g}$.
Alors $\mathfrak{q} $ est stable sous $\theta$; on en d\'eduit une d\'ecomposition $\mathfrak{u}  = (\mathfrak{u}  \cap
\mathfrak{k} ) \oplus (\mathfrak{u}  \cap \mathfrak{p} )$. Soit $R = \mbox{dim} (\mathfrak{u} \cap \mathfrak{p} )$.

Associ\'e \`a $\mathfrak{q}$, se trouve un $(\mathfrak{g} , K)$-module irr\'eductible bien d\'efini $A_{\mathfrak{q}}$ caract\'eris\'e
par les propri\'et\'es suivantes. Supposons effectu\'e un choix de racines positives pour $(\mathfrak{k} , \mathfrak{t} )$ de fa\c{c}on
compatible avec $\mathfrak{u}$. Soit $e(\mathfrak{q})$ un g\'en\'erateur de la droite $\bigwedge^R (\mathfrak{u} \cap \mathfrak{p} )$.
Alors $e(\mathfrak{q})$ est le vecteur de plus haut poids d'une repr\'esentation irr\'eductible $V(\mathfrak{q})$ de $K$ contenue
dans $\bigwedge^R \mathfrak{p}$; et dont le plus haut poids est donc n\'ecessairement $2\rho (\mathfrak{u} \cap \mathfrak{p} )$.
La classe d'\'equivalence du $(\mathfrak{g} , K)$-module $A_{\mathfrak{q}}$ est alors uniquement caract\'eris\'ee
par les deux propri\'et\'es suivantes.
\begin{eqnarray} \label{Aq1}
\begin{array}{l}
A_{\mathfrak{q}} \mbox{ {\it est unitarisable avec le m\^eme caract\`ere infinit\'esimal que la}} \\
\mbox{{\it  repr\'esentation triviale}}
\end{array}
\end{eqnarray}
\begin{eqnarray} \label{Aq2}
\mbox{Hom}_K (V(\mathfrak{q} ), A_{\mathfrak{q}} ) \neq 0.
\end{eqnarray}
Remarquons que la classe du module $A_{\mathfrak{q}}$ ne d\'epend alors en fait que de l'intersection $\mathfrak{u} \cap \mathfrak{p}$,
autrement dit deux sous-alg\`ebres paraboliques $\mathfrak{q} = \mathfrak{l} \oplus \mathfrak{u}$ et $\mathfrak{q}' =\mathfrak{l}' \oplus
\mathfrak{u}'$ v\'erifiant $\mathfrak{u} \cap \mathfrak{p} = \mathfrak{u}' \cap \mathfrak{p}$ donnent lieu \`a une m\^eme classe de
module cohomologique.

De plus, $V(\mathfrak{q})$ intervient avec multiplici\'e $1$ dans $A_{\mathfrak{q}}$ et $\bigwedge^R (\mathfrak{p} )$, et
\begin{eqnarray} \label{gKcohom}
H^i (\mathfrak{g} , K , A_{\mathfrak{q}} ) \cong \mbox{Hom}_{L\cap K} ( \bigwedge {}^{i-R} (\mathfrak{l} \cap \mathfrak{p} ), {\Bbb C}).
\end{eqnarray}
Ici $L$ est un sous-groupe de $K$ d'alg\`ebre de Lie $\mathfrak{l}$.

Si $\Gamma$ est un sous-groupe de congruence de $G$, la $A_{\mathfrak{q}}$-composante $H^R (A_{\mathfrak{q}} : \Gamma )$ de $H^R (
S(\Gamma ))$ sera dite {\it fortement primitive}. D'apr\`es ce que nous avons rappel\'e ci-dessus la $A_{\mathfrak{q}}$-composante
fortement primitive est donc la somme sur une base
$\{ \varphi \}$ de Hom$_{\mathfrak{g} , K} (A_{\mathfrak{q}} , C^{\infty} (\Gamma \backslash
G_0^{nc} ))$ des formes diff\'erentielles $\omega_{\varphi}$ d\'efinies par
$$\omega_{\varphi} (g. \lambda ) = \varphi (\omega (\lambda )) (g) \; \; (\lambda \in \bigwedge {}^R \mathfrak{p} , \; g\in G_0^{{\rm nc}} ),$$
o\`u $\omega : \bigwedge^R \mathfrak{p} \rightarrow A_{\mathfrak{q}}$ est une $K$-application non nulle (uniquement d\'efinie \`a un
scalaire pr\`es) qui factorise n\'ecessairement via la composante isotypique $V(\mathfrak{q} )$.
De m\^eme nous parlerons de $H^R (A_{\mathfrak{q}} : Sh^0 G)$ comme de la {\it $A_{\mathfrak{q}}$-composante fortement primitive}, et notons
$H^*_{{\rm prim} +} (Sh^0 G)$ la composante fortement primitive de la cohomologie.

\subsection{Le cas des groupes unitaires}

Dans ce paragraphe $G^{{\rm nc}} = U(p,q)$, o\`u $p$ et $q$ sont des entiers strictement positifs avec $p\leq q$.
Le rang r\'eel de $G$ est donc $p$. On a
\begin{eqnarray} \label{Gnc}
G^{{\rm nc}} = \left\{ g= \left(
\begin{array}{cc}
A & B \\
C & D
\end{array} \right) \; : \; ^t \! \overline{g} \left(
\begin{array}{cc}
1_p & 0 \\
0 & -1_q
\end{array} \right) g = \left(
\begin{array}{cc}
1_p & 0 \\
0 & -1_q
\end{array} \right) \right\} ,
\end{eqnarray}
o\`u $A\in M_{p\times p} ({\Bbb C})$, $B\in M_{p\times q} ({\Bbb C})$, $C\in M_{q\times p} ({\Bbb C})$ et
$D\in M_{q\times q} ({\Bbb C})$. Soit
$$K = \left\{ g =  \left(
\begin{array}{cc}
A & 0 \\
0 & D
\end{array} \right) \in G^{{\rm nc}} \; : \;  A\in U(p) , \ D\in U(q) \right\} .$$
Le complexifi\'e $K_{{\Bbb C}}$ de $K$ est le groupe
$$K_{{\Bbb C}} = \left\{ g =  \left(
\begin{array}{cc}
A & 0 \\
0 & D
\end{array} \right) \in G^{{\rm nc}}_{{\Bbb C}} \; : \;  A\in GL_p , \ D\in GL_q \right\} .$$
L'involution de Cartan $\theta$ est donn\'ee par $x \mapsto - ^t \! \overline{x}$. Soit
$$T = \left\{ g\in K_{{\Bbb C}} \; : \; g = \left(
\begin{array}{cc}
A & 0 \\
0 & D
\end{array} \right) \mbox{ avec } A, \ D \mbox{ matrices diagonales} \right\} .$$

Rappelons que la multiplication par $i =\sqrt{-1}$ induit une d\'ecomposition
$$\mathfrak{p} = \mathfrak{p}^+ \oplus \mathfrak{p}^- .$$

Nous notons $(x_1 , \ldots , x_p ; y_1 , \ldots , y_q )$ les \'el\'ements de $T$ ou de son alg\`ebre de Lie.
L'alg\`ebre de Lie $\mathfrak{g}$ est bien \'evidemment $M_{(p+q) \times (p+q)} ({\Bbb C})$, et l'on voit ses
\'el\'ements sous forme de blocs comme dans (\ref{Gnc}). On a alors,
$$\mathfrak{p}^+ = \left\{ \left(
\begin{array}{cc}
0 & B \\
0 & 0
\end{array} \right) \mbox{ avec } B \in M_{p\times q} ({\Bbb C} ) \right\}$$
et
$$\mathfrak{p}^- = \left\{ \left(
\begin{array}{cc}
0 & 0 \\
C & 0
\end{array} \right) \mbox{ avec } C \in M_{q \times p} ({\Bbb C} )\right\} .$$
Soit $E= {\Bbb C}^p$ (resp. $F= {\Bbb C}^q$) la repr\'esentation standard de $U(p)$ (resp. $U(q)$). Alors, comme
repr\'esentation de $K_{{\Bbb C}}$, $\mathfrak{p}^+ = E \otimes F^*$.

Soient $(e_1 , \ldots , e_p )$ et $(f_1 , \ldots , f_q )$ les bases canoniques respectives de $E$ et $F$. Choisissons
comme sous-alg\`ebre de Borel $\mathfrak{b}_K$ dans $\mathfrak{k}$ l'alg\`ebre des matrices dans $\mathfrak{k}$,
qui sont triangulaires sup\'erieures sur $E$ et triangulaires inf\'erieures sur $F$ par rapport \`a ces bases.
Alors l'ensemble des racines simples compactes positives
\begin{eqnarray} \label{racines compactes}
\Phi (\mathfrak{b}_K , \mathfrak{t} ) = \{ x_i - x_j \; : \; 1\leq i < j \leq p \} \cup \{ y_j -y_i \; : \;  1 \leq i < j \leq q \} .
\end{eqnarray}
Les racines de $T$ apparaissant dans $\mathfrak{p}^+$ sont les formes lin\'eaires $x_i - y_j$ avec $1 \leq i \leq p$
et $1 \leq j \leq q$.

\subsubsection*{Modules cohomologiques et diagrammes de Young}

Nous avons vu comment associer une sous-alg\`ebre parabolique $\theta$-stable $\mathfrak{q}$ \`a un \'el\'ement $X
= (x_1 , \ldots , x_p ; y_1 , \ldots , y_q ) \in i \mathfrak{t}_0$ (les $x_i$, $y_j$ sont donc tous r\'eels).
Rappelons le choix fix\'e (\ref{racines compactes}) de racines simples compactes positives.
Apr\`es conjugaison par un \'el\'ement de $K$, on peut supposer, et nous le supposerons effectivement par la suite, que $X$ est dominant
par rapport \`a $\Phi (\mathfrak{b}_K , \mathfrak{t} )$, {\it i.e.} que $\alpha (X) \geq 0$ pour tout $\alpha \in \Phi (\mathfrak{b}_K ,
\mathfrak{t} )$; il satisfait alors aux in\'egalit\'es
$$x_1 \geq \ldots \geq x_p \: \mbox{ et } \; y_q \geq \ldots \geq y_1 .$$

Nous allons maintenant associer \`a $X$ un couple de partitions.
Rappelons qu'une {\it partition} est une suite d\'ecroissante $\lambda$ d'entiers naturels $\lambda_1 \geq \ldots \geq
\lambda_l \geq 0$. Les entiers $\lambda_1 , \ldots , \lambda_l$ sont des {\it parts}. La {\it longueur} $l(\lambda )$ d\'esigne le
nombre de parts non nulles, et le {\it poids} $|\lambda |$, la somme des parts. On se soucie peu, d'ordinaire, des parts nulles~: on se
permet en particulier, le cas \'ech\'eant, d'en rajouter ou d'en \^oter

Le {\it diagramme de Young} de $\lambda$, que l'on notera \'egalement $\lambda$, s'obtient en superposant, de haut en bas, des lignes
dont l'extr\'emit\'e gauche est sur une m\^eme colonne, et de longueurs donn\'ees par les parts de $\lambda$. Par sym\'etrie
diagonale, on obtient le diagramme de Young de la {\it partition conjugu\'ee}, que l'on notera $\lambda^*$.

Le diagramme de Young de la partition $\lambda = (5,3,3,2)$ et de sa conjugu\'e sont donc~:

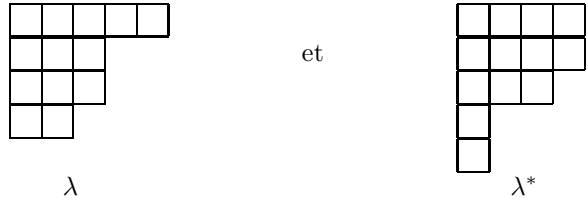
\bigskip

\begin{center}
\begin{tabular}{lcl}

\begin{tabular}{|c|c|c|c|c|} \hline
 & & & &\\ \hline
\end{tabular}
& \hspace{3cm}  &
\begin{tabular}{|c|c|c|c|} \hline
 & & & \\ \hline
\end{tabular}
\\

\begin{tabular}{|c|c|c|} \hline
 & & \\ \hline
\end{tabular}
& et  &
\begin{tabular}{|c|c|c|c|} \hline
 & & & \\ \hline
\end{tabular}
\\

\begin{tabular}{|c|c|c|} \hline
 & & \\ \hline
\end{tabular}
& \hspace{3cm} &
\begin{tabular}{|c|c|c|} \hline
 & & \\ \hline
\end{tabular}
\\

\begin{tabular}{|c|c|} \hline
 & \\ \hline
\end{tabular}
& \hspace{3cm} &
\begin{tabular}{|c|} \hline
 \\ \hline
\end{tabular}
\\

& \hspace{3cm} &
\begin{tabular}{|c|} \hline
 \\ \hline
\end{tabular}
\\

$\quad \quad \lambda$ & \hspace{3cm} & $\quad \quad \lambda^*$

\end{tabular}

\end{center}

\bigskip

Soient $\lambda$ et $\mu$ deux partitions telles que le diagramme de $\mu$ contienne $\lambda$, ce que nous noterons $\lambda \subset \mu$. Notons
$\mu / \lambda$ le compl\'ementaire du diagramme de $\lambda$ dans celui de $\mu$~: c'est une {\it partition gauche} son diagramme
est un {\it diagramme gauche}. Dans la pratique les partitions $\lambda$ que nous rencontrerons seront incluses dans la {\it partition
rectangulaire} $p\times q = (\underbrace{q, \ldots , q}_{p \; {\rm fois}}) =(q^p)$, le diagramme gauche $p\times q /\lambda$ est alors
le diagramme de Young d'une partition auquel on a appliqu\'e une rotation d'angle $\pi$; nous noterons $\hat{\lambda}$ cette partition,
la {\it partition compl\'ementaire} de $\lambda$ dans $p \times q$.
Par exemple, la partition $\lambda = (5,3,3,2)$ est incluse dans le rectangle $5\times 5$, et
dans ce rectangle, $\hat{\lambda} = (5 , 3 , 2 , 2)$.

\bigskip

Nous associons maintenant \`a notre \'el\'ement $X \in i \mathfrak{t}_0$ un couple $(\lambda , \mu )$ de partitions comme suit.
\begin{itemize}
\item La partition $\lambda \subset p\times q$ est associ\'ee au sous-diagramme de Young de $p\times q$ constitu\'e des cases de
coordonn\'ees $(i,j)$ telles que $x_i > y_j$.
\item La partition $\mu \subset p\times q$ est associ\'ee au sous-diagramme de Young de $p\times q$ constitu\'e des cases de
coordonn\'ees $(i,j)$ telles que $x_i \geq y_j$.
\end{itemize}

Le lemme suivant est absolument imm\'ediat.

\begin{lem} \label{triv}
Le couple de partitions $(\lambda , \mu)$ associ\'e \`a un \'el\'ement $X \in i \mathfrak{t}_0$ v\'erifie~:
\begin{enumerate}
\item la suite d'inclusion $\lambda \subset \mu \subset p\times q$, et
\item que le diagramme gauche $\mu / \lambda$ est une r\'eunion de diagrammes rectangulaires $p_i \times q_i$, $i=1, \ldots ,m$
ne s'intersectant qu'en des sommets.
\end{enumerate}
R\'eciproquement, \'etant donn\'e un couple de partitions $(\lambda , \mu )$ v\'erifiant 1 et 2, on peut toujours trouver
un \'el\'ement $X \in i \mathfrak{t}_0$ tel que $(\lambda , \mu )$ soit le couple de partitions associ\'e \`a $X$.
\end{lem}

Nous dirons d'un couple de partitions $(\lambda , \mu)$ qu'il est {\it compatible} (ou {\it compatible dans $p\times q$} en cas
d'ambiguit\'e) s'il v\'erifie les points 1 et 2 du Lemme \ref{triv}.

Remarquons maintenant que si $X$ et $X'$ sont deux \'el\'ements de $i \mathfrak{t}_0$ de m\^eme couple de partitions associ\'e $(\lambda
,\mu)$ et de sous-alg\`ebres paraboliques associ\'ees respectives $\mathfrak{q}$ et $\mathfrak{q}'$ alors
$\mathfrak{q} \cap \mathfrak{u} = \mathfrak{q}' \cap \mathfrak{u}'$. On d\'eduit donc de la remarque suivant la d\'efinition des modules
$A_{\mathfrak{q}}$ et du Lemme \ref{triv} que chaque couple compatible de partitions $(\lambda , \mu )$ d\'efinit sans ambiguit\'e une
classe d'\'equivalence de $(\mathfrak{g} , K)$-modules que nous notons $A(\lambda , \mu )$. Nous nous autoriserons \`a parler de
``la'' sous-alg\`ebre parabolique $\mathfrak{q}(\lambda , \mu ) = \mathfrak{l}(\lambda , \mu )\oplus \mathfrak{u}(\lambda , \mu )$ de
$(\mathfrak{g}, K)$-module associ\'e $A(\lambda , \mu )$, l'important
pour nous est qu'une telle sous-alg\`ebre existe (d'apr\`es le Lemme \ref{triv}). Nous supposerons de plus, ce que l'on peut
toujours faire, que le groupe $L(\lambda , \mu )$ associ\'e \`a la sous-alg\`ebre de Levi $\mathfrak{l} (\lambda , \mu )$
n'a pas de facteurs compacts non ab\'elien. Il est alors facile de voir que
\begin{eqnarray} \label{L}
L(\lambda , \mu )/(L(\lambda ,\mu )\cap K) = \prod_{i=1}^m U(p_i , q_i) /U(p_i ) \times U(q_i ).
\end{eqnarray}

Les r\'esultats de Parthasarathy, Kumaresan et Vogan-Zuckerman mentionn\'es plus haut affirment alors que
$$\widehat{G}^{{\rm nc}}_{{\rm VZ}} := \{ A(\lambda , \mu ) \; : \; (\lambda , \mu ) \mbox{ est un couple compatible de partitions} \} (
\subset {}_0 \widehat{G}^{{\rm nc}} \subset \widehat{G}^{{\rm nc}})$$
est l'ensemble des $(\mathfrak{g} , K)$-modules ayant des groupes de $(\mathfrak{g} , K)$-cohomologie non nuls.

\bigskip

Comme repr\'esentation de $K_{{\Bbb C}}$, $\mathfrak{p}^+ = E \otimes F^*$ et $\mathfrak{p} = (E \otimes F^*) \oplus (E \otimes F^* )^*$.
Il est bien connu (cf. \cite{Fulton}) qu'\`a chaque partition $\lambda$, il correspond une repr\'esentation irr\'eductible
$E^{\lambda}$ de $GL(E)$.

Consid\'erons la repr\'esentation de $K_{{\Bbb C}}$
\begin{eqnarray} \label{Vlambda}
V(\lambda ) := E^{\lambda} \otimes (F^{\lambda^*})^* .
\end{eqnarray}
C'est une sous-repr\'esentation irr\'eductible de $\bigwedge^{|\lambda |} (E \otimes F^* )$; son vecteur de plus haut poids est
\begin{eqnarray} \label{poids+}
v (\lambda ):= \bigwedge_{i=1}^p \bigwedge_{j=1}^{\lambda_i} e_i \otimes f_j^*
\end{eqnarray}
et son vecteur de plus bas poids est
\begin{eqnarray} \label{poids-}
w (\lambda ):= \bigwedge_{i=1}^p \bigwedge_{j=1}^{\lambda_i} e_{p-i+1} \otimes f_{q-j+1}^* .
\end{eqnarray}
On peut montrer, cf. \cite{Fulton}, que la repr\'esentation
\begin{eqnarray} \label{dec}
\bigwedge \mathfrak{p}^+ = \bigwedge (E \otimes F^* ) = \bigoplus_{\lambda  \subset p\times q} V(\lambda ),
\end{eqnarray}
o\`u chaque sous-espace irr\'eductible $V(\lambda )$ apparait avec multiplicit\'e 1.

Soit maintenant $(\lambda , \mu )$ un couple compatible de partitions. Le vecteur
\begin{eqnarray} \label{v(lambda,mu)}
v(\lambda ) \otimes w( \hat{\mu} )^* \in \bigwedge {}^{|\lambda|} (E\otimes F^* ) \otimes \bigwedge {}^{|\hat{\mu}|} (E\otimes F^* )^* = \bigwedge {}^{
|\lambda|,|\hat{\mu}|} \mathfrak{p} 
\end{eqnarray}
est un vecteur de plus haut poids $2\rho (\mathfrak{u}(\lambda , \mu ) \cap \mathfrak{p})$ et engendre donc sous l'action de
$K_{{\Bbb C}}$ un sous-module irr\'eductible que l'on note $V(\lambda , \mu )$. Ce module est isomorphe \`a $V(\mathfrak{q}(\lambda , \mu
))$ et appara\^{\i}t avec multiplicit\'e exactement $1$ dans $\bigwedge^{|\lambda| + |\hat{\mu}|} \mathfrak{p}$ (cf. \cite{VoganZuckerman}).

\medskip

Soit $(\lambda , \mu )$ un couple compatible de partitions.
Nous noterons $H^{\lambda , \mu} (Sh^0 G) = H^{|\lambda| +|\hat{\mu}|} (A(\lambda , \mu ) : Sh^0 G)$ la $A(\lambda , \mu )$-composante fortement primitive 
de la cohomologie de $Sh^0 G$.

\paragraph{Cohomologie holomorphe} 
Rappelons que le sous-espace $H^{\lambda , \mu} (Sh^0 G)$ apparait dans la cohomologie holomorphe si et seulement si $\mu = p\times q$. La partition $\lambda$
est alors naturellement param\'etr\'ee par un couple d'entier $(r,s)$ avec $0 \leq r \leq p$ et $0 \leq s \leq q$ tels que
$$\lambda = ( \underbrace{q, \ldots , q}_{r \ {\rm fois}} , \underbrace{s , \ldots , s}_{p-r \ {\rm fois}} )$$
de diagramme de Young~:
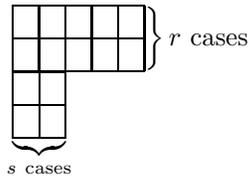
$$
\begin{array}{l}
\left. \hspace{0,035cm}
\begin{array}{|c|c|c|c|c|} \hline
 & & & &\\ \hline
 & & & &\\ \hline
\end{array}
\right\}  r \ {\rm cases} \\
\underbrace{
\begin{array}{|c|c|} \hline
 &  \\ \hline
 &  \\ \hline
\end{array}}_{s \ {\rm cases}}
\end{array}
$$
(Ici $p=4$ et $q=5$.)

Dans ce cas le sous-espace $H^{\lambda , \mu} (Sh^0 G)$ appara\^{\i}t dans la cohomologie holomorphe
de degr\'e $|\lambda | =  rq + s(p-r)$ (remarquons que $|\hat{\mu}|=0$), nous le noterons pour simplifier $H^{\lambda} (Sh^0 G)$.

\subsection{Le cas des groupes orthogonaux}

Dans ce paragraphe $G^{{\rm nc}} = O(p,q)$ et $G_0^{{\rm nc}} = SO(p,q)_0$, o\`u $p$ et $q$ sont des entiers strictement positifs. Le rang r\'eel de $G$ est donc $\min (p,q)$. On a 
\begin{eqnarray}
G^{{\rm nc}} = \left\{ g = \left(
\begin{array}{cc}
A & B \\
C & D 
\end{array} \right) \; : \; {}^t g \left(
\begin{array}{cc}
1_p & 0 \\
0 & -1_q 
\end{array} \right) g = \left(
\begin{array}{cc}
1_p & 0 \\
0 & -1_q 
\end{array} \right) \right\} ,
\end{eqnarray}
o\`u $A \in M_{p\times p} ({\Bbb R})$, $B \in M_{p\times q} ({\Bbb R})$, $C\in M_{q\times p} ({\Bbb R})$ et $D \in M_{q\times q } ({\Bbb R})$. Soit 
$$K = \left\{ g = \left( 
\begin{array}{cc}
A & 0 \\
0 & D 
\end{array} \right) \in G^{{\rm nc}} \; : \; A \in O(p) , D \in O(q) \right\} ,$$
et $K_0 = SO(p) \times SO(q)$ la composante connexe de l'identit\'e dans $K$.
L'involution de Cartan $\theta$ est donn\'ee par $x \mapsto - {}^t x$. On a alors, 
$$\mathfrak{p} = \left\{ \left( 
\begin{array}{cc}
0 & B \\
-{}^t B & 0
\end{array} \right) \; : \; B \in M_{p\times q} ({\Bbb C}) \right\} .$$
Soit $E= {\Bbb C}^p$ (resp. $F={\Bbb C}^q$) la repr\'esentation standard de $O(p,{\Bbb C})$ (resp. $O(q,{\Bbb C})$). Alors, comme repr\'esentation de 
$K_{{\Bbb C}}$, $\mathfrak{p} = E \otimes F^*$. 

Notons $r = \left[ \frac{p}{2} \right]$ et $s = \left[ \frac{q}{2} \right]$. Soit 
$i \mathfrak{t}_0$ la sous-alg\`ebre de $\mathfrak{k}$ constitu\'ee des \'el\'ements 
\begin{eqnarray*} 
\left(
\begin{array}{cccccc}
\begin{array}{cc}
0 & ix_1 \\
-ix_1 & 0 
\end{array} & & & & &\\
& \begin{array}{cc}
0 & ix_2 \\
-ix_2 & 0 
\end{array} &&&& \\
&& \ddots &&& \\
&&& \ddots && \\
&&&& \begin{array}{cc}
0 & iy_{s-1} \\
-iy_{s-1} & 0 
\end{array} & \\
&&&&& \begin{array}{cc}
0 & iy_s \\
-iy_s & 0 
\end{array}  
\end{array} \right) ,
\end{eqnarray*}
o\`u $x_1 , \ldots , x_r$ et $y_1 , \ldots , y_s$ sont r\'eels. 
L'alg\`ebre $\mathfrak{t}_0$ est une sous-alg\`ebre de Cartan de $\mathfrak{k}_0$. 
Nous noterons $(x_1 , \ldots , x_r ; y_1 , \ldots , y_s )$ les \'el\'ements de $\mathfrak{t}$.

\medskip

Distinguons trois cas.
 
\paragraph{p=2r et q=2s} 
Dans ces coordonn\'ees les syst\`emes de racines de $\mathfrak{k}$ et $\mathfrak{p}$ sont repr\'esent\'es par
\begin{eqnarray} \label{racines compactes o 1}
\Delta ( \mathfrak{k} , \mathfrak{t} ) = \{ \pm (x_i \pm x_j ) \; : \; 1 \leq i < j \leq r \} \cup \{ \pm (y_i \pm y_j ) \; : \; 1 \leq i < j \leq s \} , 
\end{eqnarray}
\begin{eqnarray} \label{racines o 1}
\Delta (\mathfrak{p} , \mathfrak{t} ) = \{ \pm (x_i \pm y_j ) \; : \; 1 \leq i \leq r,  \; 1 \leq j \leq s \}.
\end{eqnarray} 
 
Soit $j$ un entier $\in [1,r]$, nous notons $e_j$ (resp. $\overline{e}_j$) le vecteur de ${\Bbb C}^p$ 
$$(\ldots  , \underbrace{1 , -i }_{ j} ,  \ldots ) \  ({\rm resp.} \; ( \ldots  , \underbrace{1 , i }_{ j} , \ldots ) )$$ 
dont toutes les coordonn\'ees sont nulles sauf celles correspondant \`a la $j$-\`eme paire d'indices. 
De m\^eme si $j$ est un entier $\in [1,s]$, nous notons $f_j$ (resp. $\overline{f}_j$) le vecteur de ${\Bbb C}^q$ 
$$(\ldots  , \underbrace{1 , -i }_{ j} ,  \ldots ) \  ({\rm resp.} \; ( \ldots  , \underbrace{1 , i }_{ j} , \ldots ) )$$ 
dont toutes les coordonn\'ees sont nulles sauf celles correspondant \`a la $j$-\`eme paire d'indices. 

Les syst\`emes $(e_1 , \overline{e}_1 , \ldots , e_r , \overline{e}_r )$ et $(f_1 , \overline{f}_1 , \ldots , f_s , \overline{f}_s )$ 
sont des bases de $E$ et $F$ respectivement.
L'action adjointe d'un \'el\'ement $(x_1 , \ldots , x_r ; y_1 , \ldots , y_s ) \in \mathfrak{t}$ sur $\mathfrak{p} = E \otimes F^*$ est diagonalis\'ee par les vecteurs 
\begin{itemize}
\item $e_i \otimes f_j^*$ de valeur propre associ\'ee $x_i - y_j$;
\item $e_i \otimes \overline{f}_j^*$ de valeur propre associ\'ee $x_i +y_j$;
\item $\overline{e}_i \otimes f_j^*$ de valeur propre associ\'ee $-x_i-y_j$;
\item $\overline{e}_i \otimes \overline{f}_j^*$ de valeur propre associ\'ee $-x_i +y_j$.
\end{itemize}
Ici l'entier $i$ parcourt $[1,r]$ et l'entier $j$ parcourt $[1,s]$.

Pour la suite, nous fixons un syst\`eme positif de $\Delta (\mathfrak{k} , \mathfrak{t} )$ comme suit~:
\begin{eqnarray} \label{rac pos 1}
\Delta^+ ( \mathfrak{k} , \mathfrak{t} )  =  \{ x_i \pm x_j \; : \; 1 \leq i < j \leq r \} \cup \{ y_j \pm y_i \; : \; 1 \leq i < j \leq s \} .
\end{eqnarray} 
 
\paragraph{p=2r et q=2s+1}
Les  syst\`emes de racines de $\mathfrak{k}$ et $\mathfrak{p}$ sont alors repr\'esent\'es par 
\begin{eqnarray} \label{racines compactes o 2}
\begin{array}{ll}
\Delta ( \mathfrak{k} , \mathfrak{t} ) = & \{ \pm (x_i \pm x_j ) \; : \; 1 \leq i < j \leq r \} \\
& \cup \{ \pm (y_i \pm y_j ) \; : \; 1 \leq i < j \leq s \} \cup \{ \pm y_i \; : \; 1 \leq i \leq s \} , 
\end{array}
\end{eqnarray}
\begin{eqnarray} \label{racines o 2}
\begin{array}{ll}
\Delta (\mathfrak{p} , \mathfrak{t} ) = & \{ \pm (x_i \pm y_j ) \; : \; 1 \leq i \leq r,  \; 1 \leq j \leq s \} \\
& \cup \{ \pm x_i \; : \; 1 \leq i \leq r \}.
\end{array}
\end{eqnarray} 

Nous conservons les notations du paragraphe pr\'ec\'edent et notons $f_{s+1}$ le vecteur $(0, \ldots , 0 , 1) \in {\Bbb C}^q$. Les 
syst\`emes $(e_1 , \overline{e}_1 , \ldots , e_r , \overline{e}_r )$ et $(f_1 , \overline{f}_1 , \ldots , f_s , \overline{f}_s , f_{s+1})$ sont alors des bases de $E$ et $F$ respectivement.
L'action adjointe d'un \'el\'ement $(x_1 , \ldots , x_r ; y_1 , \ldots , y_s ) \in \mathfrak{t}$ sur $\mathfrak{p} = E \otimes F^*$ est diagonalis\'ee par les vecteurs
$e_i \otimes f_j^*$, $e_i \otimes \overline{f}_j^*$, $\overline{e}_i \otimes f_j^*$ et $\overline{e}_i \otimes \overline{f}_j^*$ comme au paragraphe pr\'ec\'edent ainsi que 
par les vecteurs $e_i \otimes f_{s+1}^*$ et $\overline{e}_i \otimes f_{s+1}^*$ ($1\leq i \leq r$) de valeur propres associ\'ees respectives $x_i$ et $- x_i$.

Pour la suite, nous fixons un syst\`eme positif de $\Delta (\mathfrak{k} , \mathfrak{t} )$ comme suit~:
\begin{eqnarray} \label{rac pos 2}
\begin{array}{ll}
\Delta^+ ( \mathfrak{k} , \mathfrak{t} )  = &  \{ x_i \pm x_j \; : \; 1 \leq i < j \leq r \} \cup \{ y_j \pm y_i \; : \; 1 \leq i < j \leq s \} \\
& \cup \{ y_i \; : \; 1 \leq i \leq s \}.
\end{array}
\end{eqnarray} 
 
\paragraph{p=2r+1 et q=2s+1}
Les  syst\`emes de racines de $\mathfrak{k}$ et $\mathfrak{p}$ sont alors repr\'esent\'es par 
\begin{eqnarray} \label{racines compactes o 3}
\begin{array}{ll}
\Delta ( \mathfrak{k} , \mathfrak{t} ) = & \{ \pm (x_i \pm x_j ) \; : \; 1 \leq i < j \leq r \} \cup \{ \pm x_i \; : \; 1 \leq i \leq r \} \\
& \cup \{ \pm (y_i \pm y_j ) \; : \; 1 \leq i < j \leq s \} \cup \{ \pm y_i \; : \; 1 \leq i \leq s \} , 
\end{array}
\end{eqnarray}
\begin{eqnarray} \label{racines o 3}
\begin{array}{ll}
\Delta (\mathfrak{p} , \mathfrak{t} ) = &  \{ \pm (x_i \pm y_j ) \; : \; 1 \leq i \leq r,  \; 1 \leq j \leq s \} \\
& \cup \{ \pm x_i \; : \; 1 \leq i \leq r \} \cup \{ \pm y_j \; : \; 1 \leq j \leq s \} \cup \{ 0 \}.
\end{array}
\end{eqnarray} 

Nous conservons les notations des paragraphes pr\'ec\'edents et notons $e_{r+1}$ le vecteur $(0, \ldots , 0 , 1) \in {\Bbb C}^p$. Les 
syst\`emes $(e_1 , \overline{e}_1 , \ldots , e_r , \overline{e}_r , e_{r+1})$ et $(f_1 , \overline{f}_1 , \ldots , f_s , \overline{f}_s , f_{s+1})$ sont alors des bases de $E$ et $F$ respectivement.
L'action adjointe d'un \'el\'ement $(x_1 , \ldots , x_r ; y_1 , \ldots , y_s ) \in \mathfrak{t}$ sur $\mathfrak{p} = E \otimes F^*$ est diagonalis\'ee par les vecteurs
$e_i \otimes f_j^*$, $e_i \otimes \overline{f}_j^*$, $\overline{e}_i \otimes f_j^*$, $\overline{e}_i \otimes \overline{f}_j^*$, $e_i \otimes f_{s+1}^*$ et $\overline{e}_i \otimes f_{s+1}^*$ 
comme au paragraphe pr\'ec\'edent ainsi que par les vecteurs $e_{r+1} \otimes f_{j}^*$, $e_{r+1} \otimes \overline{f}_{j}^*$ ($1\leq j \leq s$) et $e_{r+1} \otimes f_{s+1}^*$ de valeur propres associ\'ees respectives $-y_j$, $y_j$ et $0$.

Pour la suite, nous fixons un syst\`eme positif de $\Delta (\mathfrak{k} , \mathfrak{t} )$ comme suit~:
\begin{eqnarray} \label{rac pos 3}
\begin{array}{ll}
\Delta^+ ( \mathfrak{k} , \mathfrak{t} )  = &  \{ x_i \pm x_j \; : \; 1 \leq i < j \leq r \} \cup \{ x_i \; : \; 1 \leq i \leq r\}  \\
& \cup \{ y_j \pm y_i \; : \; 1 \leq i < j \leq s \} \cup \{ y_i \; : \; 1 \leq i \leq s \}.
\end{array}
\end{eqnarray}

\subsubsection*{Modules cohomologiques}

Nous avons vu comment associer une sous-alg\`ebre parabolique $\theta$-stable $\mathfrak{q}$ \`a un \'el\'ement $X= (x_1 , \ldots , x_r ; y_1 , \ldots , y_s ) \in i \mathfrak{t}_0$
(les $x_i$, $y_j$ sont donc tous r\'eels). Apr\`es conjugaison par un \'el\'ement de $K_0$, on peut supposer, et nous le supposerons effectivement par la suite, que $X$ est dominant
par rapport au choix (\ref{rac pos 1}), (\ref{rac pos 2}) ou (\ref{rac pos 3}) de $\Delta^+ (\mathfrak{k} , \mathfrak{t} )$, {\it i.e.} que $\alpha (X) \geq 0$ pour tout $\alpha \in \Delta^+ (\mathfrak{k} , \mathfrak{t} )$. Alors $X$ satisfait alors aux in\'egalit\'es
$$x_1 \geq \ldots \geq x_{r-1} \geq |x_r | \geq 0 \: \mbox{ et } \; y_s \geq  \ldots \geq y_2 \geq |y_1| \geq 0.$$
(Lorsque $p$ est impair (resp. $q$ est impair) on peut de plus supposer $x_r \geq 0$ (resp. $y_1 \geq 0$).)

Dans la suite nous notons
$$z_1 , \ldots , z_p  \  \  (\mbox{resp. } w_1 , \ldots , w_q )$$
les r\'eels $x_i$, $-x_i$ et $0$ (resp. $y_j$, $-y_j$ et $0$) rang\'es par ordre d\'ecroissant (resp. croissant).   

Nous associons alors \`a notre \'el\'ement $X \in i \mathfrak{t}_0$ la partition  $\lambda \subset p \times q$ dont le sous-diagramme de Young est
constitu\'e des cases de coordonn\'ees $(i,j)$ telles que $z_i > w_j$.
Il est facile de v\'erifier que le couple de partitions $(\lambda , \hat{\lambda} )$ est compatible. 

Nous dirons d'une partition $\lambda$ telle que le couple $(\lambda , \hat{\lambda} )$ est compatible qu'elle est {\it orthogonale}. 
Le diagramme gauche $\hat{\lambda} / \lambda$ est alors r\'eunion de diagrammes rectangulaires et est invariant par sym\'etrie centrale.
Lorsque $\hat{\lambda} / \lambda$ est constitu\'e d'un nombre impair de diagrammes rectangulaires nous dirons que la partition orthogonale 
$\lambda$ est {\it impaire} (et qu'elle est {\it paire} dans le cas contraire).

Comme dans le cas unitaire la donn\'ee d'une partition orthogonale {\bf impaire} $\lambda$ d\'efinit sans ambiguit\'e une classe d'\'equivalence de $(\mathfrak{g} , K)$-modules 
que nous notons $A(\lambda )$. Le $(\mathfrak{g} , K)$-module est cohomologique de degr\'e primitif $R= |\lambda |$. Nous notons 
$\mathfrak{q}(\lambda ) = \mathfrak{l}(\lambda ) \oplus \mathfrak{u} (\lambda )$ ``la'' sous-alg\`ebre parabolique correspondante (en supposant toujours que $\mathfrak{l}(\lambda )$
n'a pas de facteur compact non ab\'elien).  

\'Ecrivons le diagramme gauche $\hat{\lambda} / \lambda$ comme r\'eunion de diagrammes rectangulaires~: 
$(a_1 \times b_1 )* \ldots * (a_m \times b_m ) * (p_0 \times q_0 )* (a_m \times b_m )* \ldots *(a_1 \times b_1)$. 
Soit $L(\lambda )$ le sous-groupe de $G$ associ\'e \`a la sous-alg\`ebre de Levi $\mathfrak{l} (\lambda )$. Il est alors facile de v\'erifier que
\begin{eqnarray} \label{L}
\begin{array}{l}
L(\lambda ) /(L(\lambda ) \cap K ) \\
= \left( O(p_0 , q_0 )/O(p_0 ) \times O(q_0 ) \right) \times \prod_{i=1}^m U(a_i , b_i ) /U(a_i ) \times U(b_i ), 
\end{array}
\end{eqnarray}
o\`u le plongement du groupe $O(p_0 , q_0) \times \prod_{i=1}^m U(a_i , b_i)$ dans $O(p,q)$ est (\`a conjugaison dans $O(p,q)$ pr\`es) 
induit par le plongement (au niveau des groupes unitaires)
$$
\begin{array}{ccl}
U(p_0 , q_0 ) \times U(a_1 , b_1) \times \ldots \times U(a_m , b_m ) & \rightarrow & U(p,q) \\
(g_0 , g_1 , \ldots , g_m ) & \mapsto & (g_0 , g_1 , \overline{g}_1 , \ldots , g_m , \overline{g}_m ) .
\end{array}
$$

\medskip

Comme repr\'esentation de $(K_0)_{{\Bbb C}}$, $\mathfrak{p} = E \otimes F^*$. Soient 
$$(\varepsilon_1 , \ldots , \varepsilon_p ) = (e_1 ,  e_2 , \ldots \ldots , \overline{e}_2 , \overline{e}_1) $$ 
et 
$$(\gamma_1 , \ldots , \gamma_q ) = (\overline{f}_s , \overline{f}_{s-1} , \ldots \ldots , f_{s-1} , f_s )$$ 
les bases respectives de $E$ et $F$ obtenues en r\'eordonnant respectivement les vecteurs $e_i$, $\overline{e_i}$ et les vecteurs $f_j$, $\overline{f}_j$ de telle mani\`ere que 
$\varepsilon_i \otimes \gamma_j^*$ soit un vecteur propre pour l'action adjointe de $X$ sur $\mathfrak{p}$ de valeur propre $z_i - w_j$.

En supposant toujours que la partition orthogonale $\lambda$ est impaire, le vecteur
\begin{eqnarray}
v(\lambda ) := \bigwedge_{i=1}^p \bigwedge_{j=1}^{\lambda_i}  \varepsilon_i \otimes \gamma_j^*
\end{eqnarray}
appartient \`a $\bigwedge^R \mathfrak{p}$ et est un vecteur de plus haut poids $2\rho (\mathfrak{u} (\lambda ) \cap \mathfrak{p})$ pour l'action de $(K_0)_{{\Bbb C}}$. Il engendre
donc sous l'action de $(K_0)_{{\Bbb C}}$ un sous-module irr\'eductible $V(\lambda )$ qui est le plus bas $K$-type de $A(\lambda )$.

Lorsque la partition $\lambda$ est paire la situation est un petit peu plus compliqu\'ee car 
le diagramme $\lambda$ peut provenir d'un \'el\'ement $X$ avec $x_r$ ou $y_1$ non nul. 
Nous distinguons alors encore une fois trois cas.

\paragraph{p=2r et q=2s}
Soit $\lambda \subset p\times q$ une partition orthogonale {\bf paire}. Nous distinguons alors trois cas. 
\begin{enumerate}
\item $\lambda_r > \lambda_{r+1}$ et $\lambda^*_s = \lambda^*_{s+1}$~: On associe \`a $\lambda$ deux classes d'\'equivalence de 
$(\mathfrak{g} , K)$-modules que nous notons $A(\lambda )_+$ et $A(\lambda )_-$. 
Chacun de ces $(\mathfrak{g} , K)$-modules est cohomologique de degr\'e primitif $R= |\lambda |$. Nous notons 
$\mathfrak{q}(\lambda )_{\pm}  = \mathfrak{l}(\lambda )_{\pm} \oplus \mathfrak{u} (\lambda )_{\pm}$ ``la'' sous-alg\`ebre parabolique correspondante (en supposant toujours que $\mathfrak{l}(\lambda )_{\pm}$ n'a pas de facteur compact ab\'elien).  
L'espace sym\'etrique associ\'e au groupe $L(\lambda )_{\pm}$ est alors de la forme (\ref{L}) comme ci-dessus avec $p_0 = q_0 = 0$.
\item $\lambda_r = \lambda_{r+1}$ et $\lambda^*_s > \lambda^*_{s+1}$~: On associe \`a $\lambda$ deux classes d'\'equivalence de
$(\mathfrak{g} , K)$-modules que nous notons $A(\lambda )^+$ et $A(\lambda )^-$. 
Chacun de ces $(\mathfrak{g} , K)$-modules est cohomologique de degr\'e primitif $R= |\lambda |$. Nous notons 
$\mathfrak{q}(\lambda )^{\pm}  = \mathfrak{l}(\lambda )^{\pm} \oplus \mathfrak{u} (\lambda )^{\pm}$ ``la'' sous-alg\`ebre parabolique correspondante (en supposant toujours que $\mathfrak{l}(\lambda )^{\pm}$ n'a pas de facteur compact ab\'elien).  
L'espace sym\'etrique associ\'e au groupe $L(\lambda )^{\pm}$ est alors de la forme (\ref{L}) comme ci-dessus avec $p_0 = q_0 = 0$.
\item $\lambda_r > \lambda_{r+1}$ et $\lambda^*_s > \lambda^*_{s+1}$~:  On associe \`a $\lambda$ quatres classes d'\'equivalence de 
$(\mathfrak{g} , K)$-modules que nous notons $A(\lambda )^+_+$, $A(\lambda )^-_+$, $A(\lambda )^+_-$ et $A(\lambda )^-_-$. 
Chacun de ces $(\mathfrak{g} , K)$-modules est cohomologique de degr\'e primitif $R= |\lambda |$. Nous notons 
$\mathfrak{q}(\lambda )^{\pm_2}_{\pm_1}  = \mathfrak{l}(\lambda )^{\pm_2}_{\pm_1} \oplus \mathfrak{u} (\lambda )^{\pm_2}_{\pm_1}$ 
``la'' sous-alg\`ebre parabolique correspondante (en supposant toujours que $\mathfrak{l}(\lambda )^{\pm_2}_{\pm_1}$ n'a pas de facteur compact ab\'elien).  
L'espace sym\'etrique associ\'e au groupe $L(\lambda )^{\pm_2}_{\pm_1}$ est alors de la forme (\ref{L}) comme ci-dessus avec $p_0 = q_0 = 0$.
\end{enumerate}

D'apr\`es Vogan et Zuckerman \cite{VoganZuckerman}, l'ensemble 
$$
\begin{array}{l}
\left\{ A(\lambda )_{\pm_1}^{\pm_2} \; : \; 
\begin{array}{l}
\lambda \subset p\times q \mbox{ orthogonale et paire} \\
\lambda_r > \lambda_{r+1} \mbox{ et } \lambda^*_s > \lambda^*_{s+1}
\end{array}, \mbox{ et } \pm_1 , \pm_2 \mbox{ deux signes} \right\} \\
\bigcup \left\{ A(\lambda )_{\pm} \; : \; 
\begin{array}{l}
\lambda \subset p\times q \mbox{ orthogonale et paire} \\
\lambda_r > \lambda_{r+1} \mbox{ et } \lambda^*_s = \lambda^*_{s+1}
\end{array}, \mbox{ et } \pm  \mbox{ un signe} \right\} \\
\bigcup \left\{ A(\lambda )^{\pm} \; : \; 
\begin{array}{l}
\lambda \subset p\times q \mbox{ orthogonale et paire} \\
\lambda_r = \lambda_{r+1} \mbox{ et } \lambda^*_s > \lambda^*_{s+1}
\end{array}, \mbox{ et } \pm  \mbox{ un signe} \right\} \\
\bigcup \left\{ A(\lambda ) \; : \; \lambda \subset p\times q \mbox{ orthogonale et impaire} \right\} 
\end{array}
$$
exhauste l'ensemble des $(\mathfrak{g} , K)$-modules cohomologiques. 

\medskip

Faisons agir le signe $+$ trivialement sur $E$ et $F$ et le signe $-$ sur $E$ (resp. $F$) par 
l'application lin\'eaire qui fixe chacun des vecteurs $e_i$, $\overline{e}_i$ (resp. $f_j$, $\overline{f}_j$) pour $i=1 , \ldots , r-1$ 
(resp. $j=2 , \ldots , s$) et qui \'echange les vecteurs $e_r$ et $\overline{e}_r$  (resp. $f_1$ et $\overline{f}_1$) si (et seulement si) $\lambda_r > \lambda_{r+1}$ 
(resp. $\lambda^*_s > \lambda^*_{s+1}$). Nous notons de plus $v^{\pm}$  
l'action de $\pm$ sur un vecteur $v$ de $E$ ou de $F$. 

Si $\lambda$ est une partition orthogonale paire, le vecteur 
\begin{eqnarray} 
v(\lambda)_{\pm_1}^{\pm_2} := \bigwedge_{i=1}^p \bigwedge_{j=1}^{\lambda_i} \varepsilon_i^{\pm_1} \otimes (\gamma_j^{\pm_2} )^*
\end{eqnarray}
appartient \`a $\bigwedge^R \mathfrak{p}$. Et,
\begin{itemize}
\item si $\lambda_r > \lambda_{r+1}$ et $\lambda^*_s = \lambda^*_{s+1}$, $v(\lambda )_{\pm}$ est un vecteur de plus haut poids $2\rho (\mathfrak{u} (\lambda)_{\pm}  \cap \mathfrak{p})$ pour l'action de $(K_0)_{{\Bbb C}}$. Il engendre donc sous l'action de $(K_0)_{{\Bbb C}}$ un sous-module irr\'eductible $V(\lambda )_{\pm}$ qui est le plus bas $K$-type de $A(\lambda )_{\pm}$;
\item si $\lambda_r = \lambda_{r+1}$ et $\lambda^*_s > \lambda^*_{s+1}$, $v(\lambda )^{\pm}$ est un vecteur de plus haut poids $2\rho (\mathfrak{u} (\lambda)^{\pm}  \cap \mathfrak{p})$ pour l'action de $(K_0)_{{\Bbb C}}$. Il engendre donc sous l'action de $(K_0)_{{\Bbb C}}$ un sous-module irr\'eductible $V(\lambda )^{\pm}$ qui est le plus bas $K$-type de $A(\lambda )^{\pm}$;
\item si $\lambda_r > \lambda_{r+1}$ et $\lambda^*_s > \lambda^*_{s+1}$, $v(\lambda )_{\pm_1}^{\pm_2}$ est un vecteur de plus haut poids $2\rho (\mathfrak{u} (\lambda)_{\pm_1}^{\pm_2}  \cap \mathfrak{p})$ pour l'action de $(K_0)_{{\Bbb C}}$. Il engendre donc sous l'action de $(K_0)_{{\Bbb C}}$ un sous-module irr\'eductible $V(\lambda )_{\pm_1}^{\pm_2}$ qui est le plus bas $K$-type de $A(\lambda )_{\pm_1}^{\pm_2}$.
\end{itemize}
 
\paragraph{p=2r et q=2s+1} 
Soit $\lambda \subset p\times q$ une partition orthogonale {\bf paire}. On associe \`a $\lambda$ deux classes d'\'equivalence de
$(\mathfrak{g} , K)$-modules que nous notons $A(\lambda )_+$ et $A(\lambda )_-$.
Chacun de ces $(\mathfrak{g} , K)$-modules est cohomologique de degr\'e primitif $R= |\lambda |$. Nous notons 
$\mathfrak{q}(\lambda )_{\pm}  = \mathfrak{l}(\lambda )_{\pm} \oplus \mathfrak{u} (\lambda )_{\pm}  $ ``la'' sous-alg\`ebre parabolique correspondante (en supposant toujours que $\mathfrak{l}(\lambda )_{\pm}$ n'a pas de facteur compact ab\'elien).
L'espace sym\'etrique associ\'e au groupe $L(\lambda )_{\pm}$ est alors de la forme (\ref{L}) comme ci-dessus avec $p_0 = q_0 = 0$.

D'apr\`es Vogan et Zuckerman \cite{VoganZuckerman}, l'ensemble 
$$
\begin{array}{l}
\left\{ A(\lambda )_{\pm} \; : \; \lambda \subset p\times q \mbox{ orthogonale et paire et } \pm \mbox{ un signe} \right\} \\
\bigcup \left\{ A(\lambda ) \; : \; \lambda \subset p\times q \mbox{ orthogonale et impaire} \right\} 
\end{array}
$$
exhauste l'ensemble des $(\mathfrak{g} , K)$-modules cohomologiques. 

\medskip

Si $\lambda$ est une partition orthogonale paire, le vecteur 
\begin{eqnarray} 
v(\lambda)_{\pm} := \bigwedge_{i=1}^p \bigwedge_{j=1}^{\lambda_i} \varepsilon_i^{\pm} \otimes \gamma_j ^*
\end{eqnarray}
appartient \`a $\bigwedge^R \mathfrak{p}$ et est un vecteur de plus haut poids $2\rho (\mathfrak{u} (\lambda)_{\pm }  \cap \mathfrak{p})$ pour l'action de $(K_0)_{{\Bbb C}}$. 
Il engendre donc sous l'action de $(K_0)_{{\Bbb C}}$ un sous-module irr\'eductible $V(\lambda )_{\pm}$ qui est le plus bas $K$-type de $A(\lambda )_{\pm } $. 

\paragraph{p=2r+1 et q=2s+1}
Alors toute partition orthogonale dans $p\times q$ est impaire.
D'apr\`es Vogan et Zuckerman \cite{VoganZuckerman}, l'ensemble 
$$
\left\{ A(\lambda ) \; : \; \lambda \subset p\times q \mbox{ orthogonale} \right\}$$
exhauste l'ensemble des $(\mathfrak{g} , K)$-modules cohomologiques. 

\bigskip

 \`A toute partition $\lambda$ de longueur $l(\lambda ) \leq p$, nous associons la repr\'esentation irr\'eductible $\overline{\Gamma}_{\lambda}$ du groupe 
 $O(p,{\Bbb C})$ obtenue par la construction de Weyl (voir par exemple \cite{FultonHarris}) \`a partir de la partition 
 $(\lambda_1 -\lambda_p , \ldots , \lambda_{[p/2]} - \lambda_{p-[p/2]+1} )$. Rappelons alors que
\begin{itemize}
\item si $p$ impair $=2r +1$, la restriction $\Gamma_{\lambda}$ de la repr\'esentation 
$\overline{\Gamma}_{\lambda}$ au groupe  $SO(p,{\Bbb C})$ est irr\'eductible de plus haut poids $(\lambda_1 - \lambda_p , \ldots , \lambda_{r} - \lambda_{r +2} )$;
\item si $p$ est pair $=2 r$ et $\lambda_{r} = \lambda_{r +1}$, la restriction $\Gamma_{\lambda}$ de la repr\'esentation 
$\overline{\Gamma}_{\lambda}$ au groupe $SO(p,{\Bbb C})$ est irr\'eductible de plus haut poids $(\lambda_1 - \lambda_p , \ldots , \lambda_{r -1} - \lambda_{r +2})$;
\item si $p$ est pair $=2r$ et $\lambda_{r} > \lambda_{r +1}$, la restriction de la repr\'esentation 
$\overline{\Gamma}_{\lambda}$ au groupe  $SO(p,{\Bbb C})$ est la somme de deux repr\'esentations irr\'eductibles
$\Gamma_{\lambda}^+$ et $\Gamma_{\lambda}^-$ de plus haut poids respectifs $(\lambda_1 - \lambda_p , \ldots , \lambda_{r } - \lambda_{r +1})$
et $(\lambda_1 - \lambda_p , \ldots , \lambda_{r +1} - \lambda_{r})$.
\end{itemize}

\bigskip

Soit maintenant $\lambda \subset p\times q$ une partition orthogonale.
\begin{itemize}
\item Si $\lambda$ est impaire, la repr\'esentation $V(\lambda )$ de $K_0$ s'identifie \`a $\Gamma_{\lambda} \otimes \Gamma_{\lambda^*}^*$.
\item Si $\lambda$ est paire, $p=2r$, $q=2s$, $\lambda_{r} > \lambda_{r +1}$ et $\lambda_{s}^* > \lambda_{s +1}^*$, chaque repr\'esentation $V(\lambda )_{\pm_1}^{\pm_2}$
de $K_0$ s'identifie \`a $\Gamma_{\lambda}^{\pm_1} \otimes (\Gamma_{\lambda^*}^{\pm_2})^*$.
\item Si $\lambda$ est paire, $p=2r$, $q=2s$, $\lambda_{r} > \lambda_{r +1}$ et $\lambda_{s}^* = \lambda_{s +1}^*$, chaque repr\'esentation $V(\lambda )_{\pm}$
de $K_0$ s'identifie \`a $\Gamma_{\lambda}^{\pm} \otimes \Gamma_{\lambda^*}^*$.
\item Si $\lambda$ est paire, $p=2r$, $q=2s$, $\lambda_{r} = \lambda_{r +1}$ et $\lambda_{s}^* > \lambda_{s +1}^*$, chaque repr\'esentation $V(\lambda )^{\pm}$
de $K_0$ s'identifie \`a $\Gamma_{\lambda} \otimes (\Gamma_{\lambda^*}^{\pm} )^*$.
\item Si $\lambda$ est paire, $p=2r$ et $q=2s+1$, chaque repr\'esentation $V(\lambda )_{\pm}$
de $K_0$ s'identifie \`a $\Gamma_{\lambda}^{\pm} \otimes \Gamma_{\lambda^*}^*$.
\end{itemize}

\bigskip

On dit d'une partition paire $\lambda \subset p\times q$ qu'elle est de {\it type 1} (resp. {\it type 2}; {\it type 3}) si $p=2r$, $q=2s$, $\lambda_{r} > \lambda_{r +1}$ et $\lambda_{s}^* = \lambda_{s +1}^*$ ou $p=2r$ et $q=2s+1$ (resp. $p=2r$, $q=2s$, $\lambda_{r} = \lambda_{r +1}$ et $\lambda_{s}^* > \lambda_{s +1}^*$ ou $p=2r+1$ et
$q=2s$;  $p=2r$, $q=2s$, $\lambda_{r} > \lambda_{r +1}$ et $\lambda_{s}^* > \lambda_{s +1}^*$).

\bigskip

Soit $\lambda \subset p \times q$ une partition orthogonale. Nous notons
$H^{\lambda} (Sh^0 G)_{\pm_1}^{\pm_2}  = H^{|\lambda|} (A(\lambda )_{\pm_1}^{\pm_2} : Sh^0 G ) $ \footnote{Suivant
la parit\'e et le type 1, 2 ou 3 de $\lambda$ les signes $\pm_1$ et $\pm_2$ peuvent \^etre sans signification auquel cas
nous les ignorons.} la $A(\lambda )_{\pm_1}^{\pm_2}$-composante fortement primitive de la cohomologie de $Sh^0 G$.

\bigskip

D'apr\`es Harish-Chandra les modules $A(\lambda )_{\pm_1}^{\pm_2}$ correspondent \`a des repr\'esentations unitaire du groupe {\bf connexe} 
$G_0^{{\rm nc}} = SO_0 (p,q)$, nous aurons besoin d'\'etendre celles-ci \`a des repr\'esentations du groupe non connexe $G^{{\rm nc}} = O(p,q)$. Remarquons tout 
d'abord qu'il correspond \`a chaque repr\'esentation $\Gamma_{\lambda}^{\pm}$ la repr\'esentation $\overline{\Gamma}_{\lambda}$ du groupe 
compact non connexe $O(p)$. \`A toute partition orthogonale $\lambda \subset p\times q$ il correspond donc la repr\'esentation 
irr\'eductible $\overline{V(\lambda )} = \overline{\Gamma}_{\lambda} 
\otimes \overline{\Gamma}_{\lambda^*}^*$ du groupe $K$. Il existe alors une unique repr\'esentation unitaire irr\'eductible du groupe $G^{{\rm nc}} = O(p,q)$ de plus bas $K$-type
$\overline{\Gamma}_{\lambda} \otimes \overline{\Gamma}_{\lambda^*}^*$. Nous notons cette repr\'esentation $\overline{A(\lambda)}$.

\subsection{Restriction entre $K$-types}

Dans ce paragraphe nous rassemblons divers r\'esultats concernant la restriction des plus bas $K$-types des modules cohomologiques 
construits aux paragraphes pr\'ec\'edents. Dans la suite nous aurons \`a consid\'erer \`a la fois le cas du groupe $U(p,q)$ et du groupe $O(p,q)$. 
Pour distinguer ces deux cas nous utiliserons les notations \'evidentes $K_U$, $\mathfrak{p}_U$, $\mathfrak{p}_U^+$, $V_U (\lambda)$, $V_U (\lambda , \mu ) \ldots$ lorsque
nous parlerons du cas du groupe $U(p,q)$ et nous conserverons les notations du \S1.2 lorsque nous parlerons du groupe $O(p,q)$. 

\begin{lem} \label{resU}
Soient $(\lambda , \mu )$ un couple compatible de partitions dans $p\times q$ et $r$ un entier $\leq p\leq q$.
On plonge $GL_{q-r}$ dans $GL_q$ par
$$A \mapsto \left(
\begin{array}{cc}
A & 0 \\
0 & 1_{r}
\end{array} \right) .$$
Supposons que la partition $(r^p)$ s'inscrive \footnote{Cf. \cite{BergeronTentative} pour une d\'efinition g\'en\'erale.} dans le diagramme gauche $\mu / \lambda$, {\it i.e.}
que la partition $\mu -(r^p)$ contienne $\lambda$.
Alors, le $(GL_p \times  GL_q )$-module $V_U (\lambda , \mu )$ vu comme $(GL_p \times GL_{q-r} )$-module contient le module
$V_U (\lambda , \mu -(r^p) )$ avec multiplicit\'e 1.
\end{lem}
{\it D\'emonstration.} On peut facilement translat\'e $v(\lambda ) \otimes w(\hat{\mu} )^*$ par $1 \times GL_q$ de mani\`ere \`a obtenir
un vecteur de plus haut poids pour $GL_p \times GL_{q-r} $ qui engendre $V_U (\lambda , \mu -(r^p) )$. Enfin, on constate que la multiplicit\'e est triviale en projetant sur
$V_U (\lambda )$.~$\Box $

\bigskip

\noindent
{\bf Remarque.} Il d\'ecoule de plus de la d\'emonstration de \cite[Th\'eor\`eme 31]{BergeronTentative} que si $V_U (\alpha , \beta )$, avec 
$\alpha , \beta$ partitions dans $p\times (q-r)$, est un sous-$(GL_p \times GL_{q-r} )$-module
de $V_U (\lambda , \mu )$ avec $|\alpha | + |\hat{\beta}| = |\lambda | + |\hat{\mu}|$ alors $\alpha = \lambda$ et $\hat{\beta} = \hat{\mu}$. La partition $(r^p)$ doit
donc s'inscrire dans le diagramme gauche $\mu / \lambda$ et $V_U (\alpha , \beta ) = V_U (\lambda , \mu -(r^p) )$.

\bigskip

Par dualit\'e (pour la forme de Killing) l'inclusion $\mathfrak{p} \rightarrow \mathfrak{p}_U$ induit l'application $\mathfrak{p}_U \rightarrow \mathfrak{p}$ que nous dirons
``de restriction''. Cette derni\`ere application induit
\begin{eqnarray} \label{restr}
\bigwedge \mathfrak{p}_U \rightarrow \bigwedge \mathfrak{p} .
\end{eqnarray}

\begin{lem} \label{resOU}
Si $\lambda$ est une partition orthogonale dans $p\times q$, alors l'image de $V_U (\lambda )$ dans $\bigwedge \mathfrak{p}$ par l'application de restriction (\ref{restr})
est non triviale et contient avec multiplicit\'e 1 le $K$-module
$\overline{V(\lambda )}$.
\end{lem}
{\it D\'emonstration.}  Il est d\'emontr\'e dans \cite{Littlewood} (cf. aussi \cite{FultonHarris}) que si $\lambda$ est une partition de longueur $l(\lambda ) \leq p$, la repr\'esentation $E^{\lambda}$ (de $GL(E) = GL_p$), vue comme $O(p,{\Bbb C})$-module, contient avec multiplicit\'e un la repr\'esentation $\overline{\Gamma}_{\lambda}$. Le lemme \ref{resOU} d\'ecoule alors imm\'ediatement de la section pr\'ec\'edente.~$\Box $
 
\bigskip

\noindent
{\bf Remarque.} De m\^eme que dans le cas unitaire, si $\overline{V (\alpha )}$, o\`u 
$\alpha $ est une partition orthogonale dans $p\times q$, est un sous-$(O(p,{\Bbb C} )\times O(q,{\Bbb C}) )$-module
de $V_U (\lambda )$ avec $|\alpha | = |\lambda | $ alors $\alpha = \lambda$.

\bigskip

\begin{lem} \label{truc}
Soit $(\lambda , \mu )$ un couple compatible de partitions dans $p\times q$. Alors, l'image de $V_U (\lambda , \mu )$ dans $\bigwedge \mathfrak{p}$ par l'application de restriction (\ref{restr}) est triviale sauf si $\lambda= 0$ ou $\mu = p \times q$.
\end{lem}
{\it D\'emonstration.} Commen\c{c}ons par remarquer que l'application (\ref{restr}) envoie une matrice
$\left(
\begin{array}{cc}
0 & B \\
C & 0
\end{array} \right) \in \mathfrak{p}_U$ sur la matrice $\frac{B-{}^t C}{2}$. Quitte \`a conjuguer le plongement de $\mathfrak{p}$ dans $\mathfrak{p}_U$ par un \'el\'ement de $(K_U)_{{\Bbb C}}$, on peut donc supposer que les vecteur $e_i \otimes f_j^* \in E \otimes F^*$ et $(e_{p-i+1} \otimes f_{q-j+1}^* )^* \in (E \otimes F^*)^*$ (pour $1\leq i \leq p$ et $1 \leq j \leq q$)
s'envoient sur $\varepsilon_i \otimes \gamma_j^*$.

Il est donc imm\'ediat que l'image de $V_U (\lambda , \mu )$ dans $\bigwedge \mathfrak{p}$ par l'application de restriction (\ref{restr}) est triviale sauf si $\lambda= 0$ ou $\mu = p \times q$.~$\Box $

%Toujours d'apr\`es \cite{Littlewood}, vu comme $K$-module,
%$$E^{\lambda} = \overline{\Gamma}_{\lambda} \oplus \bigoplus_{\delta , \mu} \overline{\Gamma}_{\mu}^{\oplus c_{\delta \mu}^{\lambda}} ,$$
%o\`u $\delta$ parcourt l'ensemble des partitions dont toutes les parts sont paires, $\mu$ parcourt l'ensemble des partitions de longueur $l(\mu )\leq p$ et
%$c_{\delta \mu}^{\lambda}$ est le coefficient de Littlewood-Richardson. On peut de m\^eme d\'ecomposer le module $F^{\lambda^*}$.

\bigskip

\begin{lem} \label{resO}
Soit $\lambda$ une partition orthogonale dans $p\times q$ et $r$ un entier $\leq p\leq q$.
On plonge $O(q-r,{\Bbb C})$ dans $O(q,{\Bbb C})$ par
$$A \mapsto \left(
\begin{array}{cc}
A & 0 \\
0 & 1_{r}
\end{array} \right) .$$
Supposons que la partition $(r^p)$ s'inscrive dans le diagramme gauche $\hat{\lambda} / \lambda$, {\it i.e.}
que la partition $\hat{\lambda} - (r^p)$ contienne $\lambda$.
La partition $\lambda$ est alors contenue et orthogonale dans $p\times (q-r)$ et le $(O(p,{\Bbb C}) \times  O(q,{\Bbb C}) )$-module $\overline{V (\lambda )}$
vu comme $(O(p,{\Bbb C} )\times O(q-r,{\Bbb C}) )$-module contient le module $\overline{V (\lambda )}$ avec multiplicit\'e 1.
\end{lem}
{\it D\'emonstration.} Il est clair que si $(r^p)$ s'inscrit dans $\hat{\lambda} / \lambda$ alors la partition $\lambda$ est contenue et orthogonale dans $p\times (q-r)$.
Le lemme d\'ecoule alors des d\'ecompositions (25.34) et (25.35) de \cite[\S 25.3]{FultonHarris}.~$\Box$

\bigskip

\noindent
{\bf Remarque.} De m\^eme que dans le cas unitaire, si $\overline{V (\alpha )}$, o\`u 
$\alpha $ est une partition orthogonale dans $p\times (q-r)$, est un sous-$(O(p,{\Bbb C} )\times O(q-r,{\Bbb C}) )$-module
de $\overline{V (\lambda )}$ avec $|\alpha | = |\lambda | $ alors (et toujours d'apr\`es les d\'ecompositions (25.34) et (25.35) de \cite[\S 25.3]{FultonHarris}) 
$\alpha = \lambda$. La partition $(r^p)$ doit
donc s'inscrire dans le diagramme gauche $\hat{\lambda} / \lambda$ et $\overline{V(\alpha )} = \overline{V (\lambda )}$.

\bigskip

\section{Th\'eor\`eme \ref{upq} et cohomologie $L^2$}

\subsection{D\'emonstration du Th\'eor\`eme \ref{upq}}

\paragraph{Premier point}
Le premier point du Th\'eor\`eme \ref{upq} est d\'emontr\'e dans \cite{BergeronTentative}. Plus g\'en\'eralement, le Th\'eor\`eme suivant d\'ecoule de 
\cite[Th\'eor\`emes 25 et 31]{BergeronTentative}.

\begin{thm} \label{analogue}
Soit $G$ un groupe alg\'ebrique r\'eductif, connexe et anisotrope sur ${\Bbb Q}$ tel que $G^{{\rm nc}} = U(p,q)$. 
Soit $Sh^0 H \subset Sh^0 G$ avec $H^{{\rm nc}} = U(p , q-r)$ plong\'e de mani\`ere standard dans $G^{{\rm nc}}$.
Soient $\lambda$ et $\mu$ deux partitions incluses dans $p\times q$ formant un couple compatible avec 
$\mu / \lambda = (p_1 \times q_1) * \ldots * (p_m \times q_m )$. 
Alors, l'application
$$H^{\lambda , \mu} (Sh^0 G) \rightarrow \prod_{g\in G({\Bbb Q})} H^{|\lambda| +|\hat{\mu}|}_{{\rm prim} +} (Sh^0 H)$$
obtenue en composant l'application ${\rm Res}_H^G$ et la projection sur la composante fortement primitive de la cohomologie de $Sh^0 H$
est {\bf injective} si et seulement si la partition $(r^p)$ s'inscrit dans le diagramme gauche $\mu / \lambda$ ({\it i.e.} si $p_1 + \ldots + p_m = p$ et $r\leq q_i$ pour $i=1 , \ldots ,m$).
Son image est alors contenue dans $\prod_{g\in G({\Bbb Q})} H^{\lambda , \mu -(r^p)} (Sh^0 H)$.
\end{thm}

\medskip

\noindent
(L'assertion sur la composante fortement primitive, qui n'est pas explicit\'ee dans \cite{BergeronTentative}, d\'ecoule de la Remarque suivant le Lemme \ref{resU}.)

\bigskip

\paragraph{Deuxi\`eme point}
Le deuxi\`eme point du Th\'eor\`eme \ref{upq} est quant \`a lui d\'emontr\'e dans \cite{BergeronClozel}. Plus g\'en\'eralement nous conjecturons le r\'esultat suivant.

\begin{conj} \label{conj2}
Soit $G$ un groupe alg\'ebrique r\'eductif, connexe et anisotrope sur ${\Bbb Q}$ tel que $G^{{\rm nc}} = U(p,q+r)$. 
Soit $Sh^0 H \subset Sh^0 G$ avec $H^{{\rm nc}} = U(p , q)$ plong\'e de mani\`ere standard dans $G^{{\rm nc}}$.
Soient $\lambda$ et $\mu$ deux partitions incluses dans $p\times q$ formant un couple compatible avec 
$\mu / \lambda = (p_1 \times q_1) * \ldots * (p_m \times q_m )$. 
Alors, l'application
$$H^{\lambda , \mu} (Sh^0 H) \rightarrow H^{|\lambda| +pr +|\hat{\mu}|}_{{\rm prim} +} (Sh^0 G)$$
obtenue en composant l'application $\bigwedge_H^G$ et la projection sur la composante fortement primitive de la cohomologie de $Sh^0 G$ 
est {\bf injective} si et seulement si la partition $(r^p)$ s'inscrit dans le diagramme gauche $\mu /\lambda$ ({\it i.e.} si $p_1 + \ldots + p_m = p$ et $r\leq q_i$ pour $i=1 , \ldots ,m$).
Son image est alors contenue dans $H^{\lambda +(r^p), \mu} (Sh^0 G)$.
\end{conj}

Nous allons maintenant chercher \`a motiver cette Conjecture.

\medskip

Dans \cite{BergeronClozel}, le deuxi\`eme point du  Th\'eor\`eme \ref{upq} est principalement d\'eduit de deux th\'eor\`emes, l'un concernant l'isolation des repr\'esentations 
cohomologiques de petit degr\'e dans le dual unitaire de $G^{{\rm nc}}$ et l'autre calculant la cohomologie $L^2$ (r\'eduite) en petit degr\'e des vari\'et\'es $\Lambda \backslash X_G$, o\`u $\Lambda$ est un sous-groupe de congruence de $H$.

Le r\'esultat concernant l'isolation des repr\'esentations cohomologiques est d\'eduit d'un th\'eor\`eme g\'en\'eral de Vogan sur lequel nous revenons dans la section suivante. Nous
n'avons en fait besoin que de l'isolation dans le dual automorphe, or on montre dans \cite{BergeronClozel} que celle-ci d\'ecoule de 
la Conjecture de changement de base (version faible des Conjectures d'Arthur) formul\'ee dans \cite{BergeronClozel}.

D'apr\`es \cite{BergeronClozel} et compte tenu des Conjectures d'Arthur, la Conjecture \ref{conj2} est principalement motiv\'ee par les calculs de cohomologie $L^2$ que nous
pr\'esentons dans les paragraphes suivants.

\subsection{Un calcul de cohomologie $L^2$}

Dans ce paragraphe, nous fixons $G$ un groupe de Lie r\'eductif {\bf r\'eel} connexe de type non compact et \`a centre compact. Soit $\tau$ une involution sur $G$ et soit 
$H=G^{\tau}$ la composante connexe de l'identit\'e du groupe des points fixes de $\tau$. Nous supposons que $G$ est la forme r\'eelle d'un groupe de Lie 
complexe et notons $K$ un sous-groupe compact maximal $\tau$-stable de $G$. L'espace sym\'etrique $X_G = G/K$ est de courbure n\'egative. Soit $\Lambda$ un
r\'eseau cocompact de $H$. Nous nous int\'eressons \`a la cohomologie $L^2$ du quotient $\Lambda \backslash X_G$. Nous notons $M= \Lambda \backslash X_G$ 
et $F = \Lambda \backslash X_H$.

\medskip

La vari\'et\'e $M$ est riemannienne et compl\`ete. On note $C_0^{\infty} (\Lambda^k T^* M)$ (respectivement 
$L^2 (\Lambda^k T^* M)$, etc...) l'ensemble des $k$-formes lisses \`a support compact (respectivement de carr\'e 
int\'egrable, etc...) dans $M$. Le $k$-i\`eme espace de cohomologie $L^2$ (r\'eduite) de $M$ est d\'efini par 
$$H_2^k (M) = \{ \alpha \in L^2 (\Lambda^k T^* M) \; : \; d\alpha =0 \} / \overline{dC_0^{\infty} (\Lambda^{k-1} T^* M)} ^{L^2} .$$
Un autre espace tr\`es proche souvent consid\'er\'e est l'espace de cohomologie $L^2$ non r\'eduite, qui, en degr\'e $k$, est le 
quotient de $\{ \alpha \in L^2 (\Lambda^k T^* M) \; : \; d\alpha =0 \}$ par  $\{ d\alpha \; : \; \alpha \in L^2 (\Lambda^{k-1} T^* M), \; d\alpha \in L^2 \}$,
sans prendre d'adh\'erence. En g\'en\'eral, cohomologie $L^2$ r\'eduite et non r\'eduite sont diff\'erentes.
Il y a n\'eanmoins \'egalit\'e en degr\'e $k$ lorsque $0$ n'est pas dans le spectre essentiel du laplacien $\Delta$ sur 
les formes diff\'erentielles de degr\'e $k$. Dans la suite, ``cohomologie $L^2$'' voudra dire ``cohomologie $L^2$ r\'eduite''.

Il y a une interpr\'etation de la cohomologie $L^2$ en termes de formes harmoniques. En effet, notons ${\cal H}^k_{2}$ l'espace
des $k$-formes harmoniques $L^2$ de $M$~:
$${\cal H}^k_2 (M) = \{ \alpha \in L^2 (\Lambda^k T^* M) \; : \; d\alpha = \delta \alpha =0 \}$$
o\`u $\delta$ est l'op\'erateur d\'efini initialement sur les formes lisses \`a support compact comme l'adjoint de $d$. 
Comme $M$ est compl\`ete, ${\cal H}_2^* (M)$ est aussi le noyau $L^2$ du laplacien $\Delta = d\delta +\delta d$. Un fait important est la d\'ecomposition
de Hodge-de Rham-Kodaira~:
$$L^2 (\Lambda^k T^*M )={\cal H}^k_2 (M) \oplus \overline{dC_0^{\infty} (\Lambda^{k-1} T^*M )} \oplus \overline{\delta C_0^{\infty} (\Lambda^{k+1} T^* M)} ,$$
et de plus,
$$\{ \alpha \in L^2 (\Lambda^k T^*M ) \; : \; d\alpha =0 \} = {\cal H}^k_2 (M) \oplus \overline{dC_0^{\infty} (\Lambda^{k-1} T^*M )}.$$
On en d\'eduit que 
$$H_2^k (M) \cong {\cal H}_2^k (M) .$$
Nous noterons $C_0^{\infty} (\Lambda \backslash G , \bigwedge^* \mathfrak{p}^* )$ l'image de $C^{\infty} (M, \bigwedge^* T^* M)$ par l'application naturelle consistant
\`a tir\'ee en arri\`ere les formes diff\'erentielles. Toute forme harmonique $L^2$ $\varphi$ sur $M$ d\'efinit donc un \'el\'ement
$\tilde{\varphi}$ de 
$$C_0^{\infty} (\Lambda \backslash G , \bigwedge ^* \mathfrak{p}^* ) \cong {\rm Hom}_K (\bigwedge ^* \mathfrak{p} , C^{\infty} (\Lambda \backslash G ; {\Bbb C} )) .$$

La formule de Matsushima se g\'en\'eralise \`a ce cadre (cf. \cite{BorelWallach}). L'espace $H_2^* (M)$ se d\'ecompose donc en somme directe~:
\begin{eqnarray}
H_2^* (M) & = & \bigoplus_{\pi \in \widehat{G}_0}
H_2^* (\pi : M),
\end{eqnarray}
o\`u nous avons not\'e $H_2^* (\pi : M)$ la $\pi$-composante de la
cohomologie $L^2$ de $M$. Si $R$ est le degr\'e fortement primitif d'une repr\'esentation $\pi \in {}_0 \widehat{G}$, l'espace $H_2^R (\pi : M)$ est aussi la partie de 
la cohomologie $L^2$ de $M$ qui est repr\'esent\'ee par des formes harmoniques dans 
\begin{eqnarray} \label{3}
C_0^{\infty} (\Lambda \backslash G , \bigwedge {}^* \mathfrak{p}^* )_{\delta} 
\end{eqnarray}
(avec $*=R$) 
le sous-ensemble de $C_0^{\infty} (\Lambda \backslash G , \bigwedge^* \mathfrak{p}^* )$ constitu\'e des \'el\'ements de la forme $\tilde{\varphi} = \sum_i  \tilde{\varphi}_i X_i$ avec 
$\tilde{\varphi}_i$ dans la composante isotypique $C^{\infty} ( \Lambda \backslash G)_{\delta}$ de type $\delta$, l'unique $K$-type minimal de $\pi$ (cf. section 1).
Plus g\'en\'eralement, nous notons $H^*_2 (M)_{\delta}$ la partie de la cohomologie $L^2$ de $M$ repr\'esent\'ee par des formes harmoniques dans (\ref{3}). Puisque le
laplacien de Hodge-de Rham commute \`a la projection sur les $K$-types, on a alors la d\'ecomposition
\begin{eqnarray} \label{dec en K types}
H_2^* (M) & = & \bigoplus_{\delta} H_2^* (M)_{\delta} .
\end{eqnarray}

\bigskip
 
Nous notons $\delta^*$ le $K$-type dual d'un $K$-type donn\'e $\delta$. La dualit\'e 
\begin{eqnarray} \label{dualite}
H_2^* (M)_{\delta} \times H_2^* (M)_{\delta^* } \rightarrow {\Bbb C}
\end{eqnarray}
est alors donn\'ee par 
$$(\varphi , \psi ) \mapsto \int_M \varphi \wedge \psi .$$
D'un autre c\^ot\'e, l'op\'erateur $*$ de Hodge induit un isomorphisme lin\'eaire 
\begin{eqnarray} \label{isom hodge}
H^*_2 (M)_{\delta } \stackrel{*}{\rightarrow} H^{d_G -*}_2 (M)_{\delta^*}.
\end{eqnarray}
On en d\'eduit un produit scalaire sur $H^*_2 (M)_{\delta}$
$$(\varphi_1 , \varphi_2 ) \mapsto \int_M \varphi_1 \wedge * \varphi_2 .$$

\medskip

Notre but est de comprendre la cohomologie $L^2$ de $M$ en termes de la cohomologie (usuelle) de $F = \Lambda \backslash X_H$. Remarquons que $F$ est naturellement
plong\'ee dans $M$, par dualit\'e elle d\'efinit une classe ``($L^2$-)duale'' $[F] \in H_2^{d_G -d_H} (M)$ (qui peut \^etre nulle a priori). 
Le th\'eor\`eme principal de cette section se d\'eduit presque imm\'ediatement des travaux de Tong et Wang \cite{TongWang}.

\begin{thm} \label{cohoml2}
La classe $[F] \in H_2^{d_G - d_H} (M)$ est non nulle si et seulement si ${\rm rang}_{{\Bbb C}} (G/H) = {\rm rang}_{{\Bbb C}} (K/(K\cap H))$.
\end{thm}

Nous l'avons dit, la d\'emonstration du Th\'eor\`eme \ref{cohoml2} repose sur des travaux de Tong et Wang \cite{TongWang}. Ceux-ci permettent de r\'ealiser 
g\'eom\'etriquement certaines s\'eries discr\`etes d'espaces sym\'etriques (non n\'ecessairement riemanniens). 
Nous commen\c{c}ons par des rappels sur les s\'eries discr\`etes. Puis nous suivons Tong et Wang \cite{TongWang}
en incluant les d\'emonstrations de ceux de leurs lemmes et propositions que nous utilisons.

\subsection{Rappels sur les s\'eries discr\`etes de $G/H$}

Commen\c{c}ons par rappeler qu'une repr\'esentation $\pi$ de $G$ est dite {\it membre de la s\'erie discr\`ete de $G/H$} si elle est isomorphe \`a une sous-repr\'esentation irr\'eductible
de la repr\'esentation r\'eguli\`ere gauche ${\cal L}$ sur $L^2 (G/H)$. Nous notons $L^2_d (G/H)$ l'{\it espace des  s\'eries dicr\`etes de $G/H$} {\it i.e.} le sous-espace lineaire ferm\'e
dans $L^2 (G/H)$ engendr\'e par les sous-repr\'esentations irr\'eductibles de ${\cal L}$.  

Si ${\rm rang}_{{\Bbb C}} (G/H) = {\rm rang}_{{\Bbb C}} (K/K\cap H)$ et d'apr\`es Flensted-Jensen \cite{FlenstedJensen}, l'espace $L_d^2 (G/H)$ des s\'eries discr\`etes est non nul.
Les repr\'esentations de la s\'erie discr\`ete sont classifi\'ees dans \cite{OshimaMatsuki} (cf. aussi \cite{FlenstedJensen2}). Il s'av\`ere que la condition 
${\rm rang}_{{\Bbb C}} (G/H) = {\rm rang}_{{\Bbb C}} (K/K\cap H)$ est en fait n\'ecessaire \`a l'existence des s\'eries discr\`etes. 
Nous supposerons dor\'enavant ${\rm rang}_{{\Bbb C}} (G/H) = {\rm rang}_{{\Bbb C}} (K/K\cap H)$ et extrayons de leurs r\'esultats ceux dont nous aurons besoin. 

Soit $\mathfrak{q}_0$ le suppl\'ementaire orthogonal (pour la forme de Killing) de $\mathfrak{h}_0$ dans $\mathfrak{g}_0$, et soit $\mathfrak{t}_0$ un sous-espace de 
Cartan compact de $\mathfrak{q}_0$. Soit $\Sigma$ le syst\`eme de racines de $\mathfrak{t}$ dans $\mathfrak{g}$ et $\Sigma_c$ le sous-syst\`eme des racines de $\mathfrak{t}$
dans $\mathfrak{k}$. Soient $W$ et $W_c$ les groupes de Weyl correspondant. Fixons un sous-syst\`eme positif $\Sigma_c^+$ dans $\Sigma_c$. Le choix d'un sous-syst\`eme 
positif $\Sigma^+$ dans $\Sigma$ compatible avec $\Sigma_c^+$ d\'etermine une bijection entre le quotient
$W_c \backslash W $ et l'ensemble des sous-syst\`emes positifs de $\Sigma$ compatibles avec $\Sigma_c^+$ : dans chaque classe de $W_c \backslash W$, il 
existe un unique repr\'esentant $w \in W$ tel que 
\begin{eqnarray} \label{w}
w(\Sigma^+ ) \cap \Sigma_c = \Sigma_c^+.
\end{eqnarray} 
Nous notons $W^c$ l'ensemble des $w \in W$ v\'erifiant (\ref{w}).
Soient $\rho$ et $\rho_c$ les demi-sommes respectives des racines dans $\Sigma^+$ et $\Sigma_c^+$, compt\'ees avec multiplicit\'ees.
Soit $\Lambda \subset i \mathfrak{t}^*_0$ l'ensemble des param\`etres $\lambda \in i \mathfrak{t}^*_0$ tels que
\begin{enumerate}
\item $\langle \lambda , \alpha \rangle >0$ pour toute racine $\alpha \in \Sigma^+$,
\item $\lambda + \rho$ est un poids pour $T_H$, {\it i.e.} $e^{\lambda + \rho}$ est bien d\'efini sur $T_H$ le tore dans $G/H$ correspondant \`a $\mathfrak{t}_0$, et
\item $\langle \lambda + \rho , \beta \rangle \geq 0$ pour toute racine simple compacte $\beta \in \Sigma^+$.
\end{enumerate} 

On \'etend $\mathfrak{t}$ en une sous-alg\`ebre de Cartan $\widetilde{\mathfrak{t}}$, $\tau$ et $\theta$-invariante dans $\mathfrak{g}$. Notons
$$\Delta_c = \Delta (\widetilde{\mathfrak{t}} , \mathfrak{k}) \; \mbox{ et } \Delta = \Delta (\widetilde{\mathfrak{t}} , \mathfrak{g})$$
les syst\`emes de racines correspondants et soient $\Delta_c^+$ et $\Delta^+$ deux choix fix\'es de sous-syst\`emes positifs compatibles avec $\Sigma^+$.
Les repr\'esentations irr\'eductibles de dimension finie de $\mathfrak{k}$ sont classifi\'ees par leur plus haut (resp. plus bas) poids par rapport au choix $\Delta_c^+$.
Sauf mention du  contraire nous identifierons dans la suite un $K$-type avec son plus haut par rapport \`a $\Delta_c^+$. Remarquons que si $\delta$ est un $K$-type de
plus haut poids $\mu$ par rapport \`a $\Delta_c^+$, la repr\'esentation contragr\'ediente ou duale $\delta^*$ a pour plus haut poids $-\mu$ par rapport \`a
$\Delta_c^+$.

On peut associer \`a chaque couple $(w , \lambda ) \in W^c \times \Lambda$ une sous-repr\'esentation $\pi_{w \lambda}$ dans ${\cal L}$ qui est soit irr\'eductible 
soit nulle et dont les sous-espaces de $L^2 (G/H)$ correspondant sont deux \`a deux disjoints. Voici quelques propri\'et\'es de ces repr\'esentations.
\begin{enumerate}
\item Le caract\`ere infinit\'esimal de $\pi_{w  \lambda}$ est $\chi_{\lambda}$.
\item Si 
\begin{eqnarray} \label{cond}
\langle w  (\lambda + \rho ) - 2 \rho_c , \alpha \rangle \geq 0,  \    \  \alpha \in \Sigma_c^+,
\end{eqnarray} 
alors $\pi_{w  \lambda}$ est non nulle et $-w  (\lambda + \rho ) + 2 \rho_c$ est (le plus haut poids par rapport \`a $\Delta_c^+$ de) l'unique $K$-type minimal de $\pi_{w  \lambda}$.
\item Si $(w_0 , \lambda_0 )$ et $(w_1 , \lambda_1 ) \in W^c \times \Lambda$ v\'erifient (\ref{cond}), les repr\'esentations $\pi_{w_0 \lambda_0}$ et $\pi_{w_1 \lambda_1}$
sont isomorphes si et seulement si $(w_0 , \lambda_0 ) = (w_1 , \lambda_1 )$.
\end{enumerate} 
Le r\'esultat principal de \cite{OshimaMatsuki} est alors que l'espace $L_d^2 (G/H)$ est engendr\'e par les sous-repr\'esentations $\pi_{w  \lambda}$ pour $(w, \lambda ) \in 
W^c \times \Lambda$.

Pour les ``grands'' $\lambda$, c'est \`a dire ceux v\'erifiant (\ref{cond}), la non trivialit\'e de $\pi_{w \lambda}$ provient de la construction explicite de la {\it fonction de 
Flensted-Jensen} $\psi_{w,\lambda}$ qui est une fonction $K$-finie dans $L^2 (G/H) \cap C^{\infty} (G/H)$ telle que 
\begin{enumerate}
\item $\psi_{w, \lambda}$ engendre le $K$-type $-w (\lambda + \rho)+2\rho_c$.
\item $\psi_{w, \lambda}$ engendre la sous-repr\'esentation $\pi_{w \lambda}$ de $L^2 (G/H)$.
\end{enumerate}

\medskip

\'Etant donn\'e deux choix $\Sigma^+$ et $\Sigma_c^+$ 
de sous-syst\`emes positifs compatibles dans $\Sigma$ et $\Sigma_c$ respectivement, Flensted-Jensen associe
\`a tout $\lambda \in \Lambda$, une fonction $\psi_{\lambda}$, la fonction not\'ee $\psi_{e , \lambda}$ ci-dessus. 
Si $w \in W^c$, les sous-ensembles de
racines $w(\Sigma^+ )$ et $\Sigma_c^+$ forment deux sous-syst\`emes positifs compatibles dans $\Sigma$ et $\Sigma_c$ respectivement. 
L'ensemble $\Lambda$ 
correspondant \`a ce nouveau choix est $w\Lambda$ et la fonction de Flensted-Jensen $\psi_{w\lambda}=\psi_{w, \lambda}$. 
Dans la nous ne nous pr\'eoccuperons que des fonctions $\psi_{\lambda}$. 

\medskip

Notons $\delta_{\mu}$ la (classe d'\'equivalence d'une) repr\'esentation irr\'eductible de dimension finie de $\mathfrak{k}$ de plus haut poids $\mu$
par rapport au sous-syst\`eme positif $\Delta_c^+$.
D'apr\`es un th\'eor\`eme d'Helgason \cite[Chap. III, \S 3]{Helgason} la repr\'esentation $\delta_{\mu}$ a un vecteur non nulle $\mathfrak{h} \cap \mathfrak{k}$-invariant si et
seulement si
\begin{eqnarray} \label{helgason}
\mu \in \mathfrak{t}^* \; \mbox{ et } \;  \frac{\langle \mu , \alpha \rangle}{\langle \alpha , \alpha \rangle } \in {\Bbb N} \; \; \mbox{ pour tout } \alpha \in \Sigma_c^+ .
\end{eqnarray}
Dans ce cas le vecteur invariant est unique \`a un multiple scalaire pr\`es.

Soit maintenant $\lambda \in i \mathfrak{t}^*_0$. Supposer
\begin{itemize}
\item $\lambda + \rho$ est un poids pour $T_H$, {\it i.e.} $e^{\lambda + \rho}$ est bien d\'efini sur $T_H$ le tore dans $G/H$ correspondant \`a $\mathfrak{t}_0$, et
\item $\langle \lambda + \rho  - 2 \rho_c , \alpha \rangle \geq 0$,  pour tout $\alpha \in \Sigma_c^+$,
\end{itemize}
revient \`a supposer que $\mu = \mu_{\lambda} := \lambda + \rho -2 \rho_c \in \mathfrak{t}^*$ (\'etendue trivialement \`a une sous-alg\`ebre de Cartan de $\mathfrak{k}$ contenant
$\mathfrak{t}_0$) est le plus haut poids d'une repr\'esentation de dimension finie de $K$ avec un vecteur $L$-invariant. 
C'est en particulier le cas pour les ``grands'' $\lambda  \in \Lambda$.
Rappelons que la fonction $\psi_{\lambda} \in C^{\infty} (G/H) \cap L^2 (G/H )$ engendre alors 
le $K$-type $\delta_{\mu}^*$ dual de $\delta_{\mu}$ et de plus haut poids $-\mu$ par rapport \`a $\Delta_c^+$.

\subsubsection*{Un raffinement du Th\'eor\`eme \ref{cohoml2}}

Nous allons maintenant pouvoir commencer la d\'emonstration du Th\'eor\`eme \ref{cohoml2}. Celui-ci d\'ecoulera du r\'esultat plus g\'en\'eral suivant.

\begin{thm} \label{cohoml2+}
Supposons ${\rm rang}_{{\Bbb C}} (G/H) = {\rm rang}_{{\Bbb C}} (K/(K\cap H))$. Si $w \in W^c$, v\'erifie
\begin{eqnarray} \label{condw}
\langle 2w\rho -2\rho_c , \alpha \rangle \geq 0 , \  \  \alpha \in \Sigma_c^+ ,
\end{eqnarray}
alors la projection de la classe $[F]$ dans $H_2^{d_G - d_H} (M)_{2w\rho -2\rho_c}$ est non nulle.
\end{thm}

Commen\c{c}ons par d\'eduire le Th\'eor\`eme \ref{cohoml2} du Th\'eor\`eme \ref{cohoml2+}. Supposons ${\rm rang}_{{\Bbb C}} (G/H) = {\rm rang}_{{\Bbb C}} (K/(K\cap H))$.
La condition (\ref{condw}) est toujours v\'erifi\'ee pour $w=e$ puisque $\rho-\rho_c$ et $\alpha$ sont tous deux dans une chambre de Weyl positive (au sens large). Le Th\'eor\`eme
\ref{cohoml2+} implique donc l'implication r\'eciproque dans l'\'enonc\'e du Th\'eor\`eme \ref{cohoml2}. D\'emontrons maintenant l'implication directe~: la classe de cohomologie 
$L^2$, $[F]$, admet clairement un repr\'esentant $H$-invariant. Si $[F]\neq 0$ il doit donc exister des s\'eries discr\`etes dans $L^2 (G/H)$ et nous avons rappel\'e qu'alors
on a n\'ecessairement  ${\rm rang}_{{\Bbb C}} (G/H) = {\rm rang}_{{\Bbb C}} (K/(K\cap H))$.~$\Box$

\bigskip

Il nous reste donc \`a d\'emontrer le Th\'eor\`eme \ref{cohoml2+}. Il nous suffit de traiter le cas $w=e$. Remarquons que $\rho$ appartient toujours \`a $\Lambda$ et donc que 
la fonction de Flensted-Jensen $\psi_{\rho}$ est bien d\'efinie.

\subsection{Construction de la forme harmonique}

En suivant Tong et Wang \cite{TongWang} nous allons maintenant associer \`a $\psi_{\rho}$ une forme diff\'erentielle harmonique $\tilde{\omega}$ sur  
$G$ et \`a valeurs dans ${\Bbb C}$. Les translat\'es \`a droite de (coefficients de) $\tilde{\omega}$ r\'ealiseront la repr\'esentation contragr\'ediente $\pi_{\rho}^*$ de
$\pi_{\rho}$ et $\tilde{\omega}$ sera invariante \`a gauche sous $H$, et sous l'action \`a droite de $K$ elle se transformera selon le plus bas $K$-type 
$\delta$ (de plus haut poids $2(\rho - \rho_c)$) de $\pi_{\rho}^*$.
La forme $\tilde{\omega}$ sera tir\'ee arri\`ere d'une forme diff\'erentielle sur $X_G$ dont nous v\'erifierons plus loin qu'elle d\'efinit une classe de cohomologie $L^2$ non nulle
dans $H_2^{d_G -d_H} (M)_{\delta}$ et \'egale \`a un multiple du projet\'e de $[F]$.

\medskip

Soit $\mathfrak{q}= \mathfrak{l} \oplus \mathfrak{u}$ la sous-alg\`ebre parabolique $\theta$-stable de $\mathfrak{g}$ telle que 
$\mathfrak{l}$ soit \'egale au centralisateur de $\mathfrak{t}$
dans $\mathfrak{g}$ et $\mathfrak{u} = \oplus_{\alpha \in \Sigma^+} \mathfrak{g}_{\alpha}$. 
Puisque les repr\'esentations $\pi_{\rho}^*$ et $A_{\mathfrak{q}}$ ont toutes deux le m\^eme caract\`ere infinit\'esimal, \'egal \`a celui de la repr\'esentation triviale $\chi_{\rho}$, de
telle mani\`ere que le casimir $\Omega$ agit trivialement sur chacun de leurs $(\mathfrak{g},K)$-modules associ\'ees, et puisqu'elles contiennent toutes les 
deux le $K$-type $\delta$ de plus haut poids $2(\rho -\rho_c )$, elles doivent \^etre isomorphe d'apr\`es \cite[Proposition 6.1]{VoganZuckerman}.

Notons $R=\dim_{{\Bbb C}} (\mathfrak{u} \cap \mathfrak{p})$ et $\rho_n = \rho - \rho_c$. 
Rappelons que le vecteur $e(\mathfrak{q}) \in \bigwedge^R (\mathfrak{u} \cap \mathfrak{p})$ est alors un vecteur de plus 
haut poids $2 \rho (\mathfrak{u} \cap \mathfrak{p})=2\rho_n$ d'une repr\'esentation irr\'eductible de $K$, isomorphe \`a $\delta$, contenue dans
$\bigwedge^R \mathfrak{p}$ et qui appara\^{\i}t avec multiplicit\'e $1$.

La fonction $\psi_{\rho} \in L^2 (G/H)$ est $K$-finie et engendre le $K$-type $\delta^*$.
Pour un choix convenable de $k_i \in K$, les fonctions
\begin{eqnarray} \label{base}
f_i (g) = \psi_{\rho} (k_i^{-1} g^{-1} ) \in L^2 (H \backslash G)
\end{eqnarray}
forment donc une base d'un $K$-type isomorphe \`a $\delta^*$ et apparaissant dans (la repr\'esentation r\'egul\`ere droite dans) $L^2 (H \backslash G)$ de plus haut poids $2\rho_n$. 
Les fonctions $\overline{f}_i$ forment alors une base d'un $K$-type isomorphe \`a $\delta$. 

Nous venons de rappeler que la composante isotypique $\left( \bigwedge^R \mathfrak{p}^* \right)_{\delta^*}$ dans $ \bigwedge^R \mathfrak{p}^* $ appara\^{i}t avec multiplicit\'e 1. 
Choisissons alors une base $\{ X_i \}$ de $\left( \bigwedge^R \mathfrak{p}^* \right)_{\delta^*}$ duale \`a 
$\{ \overline{f}_i \}$ et posons
\begin{eqnarray} \label{omega tilde}
\tilde{\omega} (g) = \sum_i \overline{f_i (g)} X_i . 
\end{eqnarray}

\begin{prop}[Tong et Wang] {\rm \cite[Proposition 3.5]{TongWang}} \label{omega}
\begin{enumerate}
\item La forme $\tilde{\omega}$ appartient \`a  $C^{\infty}_0 (G , \bigwedge^R \mathfrak{p}^* )$.
\item La forme $\tilde{\omega}$ est invariante \`a gauche sous l'action de $H$. En restriction \`a $H$, $\tilde{\omega}$ prend ses valeurs dans le sous-espace
de dimension 1 invariant sous l'action de $H\cap K$ : $\left( \bigwedge^R \mathfrak{p}^* \right)_{\delta^*}^{H \cap K}$.
\end{enumerate}
\end{prop}
{\it D\'emonstration.} Il est imm\'ediat par construction que $\tilde{\omega} \in C^{\infty} (G , \bigwedge^R \mathfrak{p}^* )$. L'invariance \`a gauche de $\tilde{\omega}$ sous
l'action de $H$ d\'ecoule imm\'ediatement de l'invariance \`a droite de $\psi_{\rho}$ sous l'action de $H$. L'application $\tilde{\omega}_{|H}$ est donc constante \'egale
\`a $\tilde{\omega} (e)$. D'apr\`es le premier point, ce vecteur est invariant sous l'action de $H\cap K$. On conclut alors la d\'emonstration de la Proposition \ref{omega}
gr\^ace au th\'eor\`eme d'Helgason (\ref{helgason}).~$\Box$

\bigskip

D'apr\`es la Proposition \ref{omega}, $\tilde{\omega}$ est tir\'ee arri\`ere d'une forme diff\'erentielle
\begin{eqnarray} \label{fd}
\omega \in C^{\infty} ( X_G , \bigwedge {}^R T^* X_G ) . 
\end{eqnarray}

Celle-ci descend en une forme diff\'erentielle $\omega \in C^{\infty} (M, \bigwedge^R T^* M)$.

\begin{lem} \label{harmo}
La forme diff\'erentielle $\omega$ est $L^2$ et harmonique.
\end{lem}
{\it D\'emonstration.} Puisque $\psi_{\rho} \in L^2 (G/H)$, $\tilde{\omega} \in L^2 (\Lambda \backslash G)$. Les coefficients de $\tilde{\omega}$ appartiennent donc 
\`a un $(\mathfrak{g},K)$-module unitaire et nous pouvons appliquer le lemme de Kuga \cite{BorelWallach}. Soit $\Delta$ le laplacien sur 
$C_0^{\infty} (G, \bigwedge^* (\mathfrak{p}^* ))$. Alors
$$\Delta \tilde{\omega} = - \pi_{\rho}^* (\Omega ) \tilde{\omega} = 0.$$
O\`u la derni\`ere \'egalit\'e provient de ce que la repr\'esentation $\pi_{\rho}$ a la m\^eme caract\`ere infinit\'esimal que la repr\'esentation triviale.
On en d\'eduit imm\'ediatement le Lemme \ref{harmo}.~$\Box$

\subsubsection*{Lien avec la sous-vari\'et\'e $F$}

Soit $\omega_H$ la forme volume invariante sur $X_H$. Nous voyons cette forme diff\'erentielle comme une forme diff\'erentielle sur $X_H \subset X_G$ \`a valeurs dans
$\left( \bigwedge^* T^* X_G \right)_{|X_H}$. Notons toujours $*$ l'op\'erateur de Hodge-de Rham. Alors $*\omega_H$ prend elle aussi ses valeurs dans
$\left( \bigwedge^* T^* X_G \right)_{|X_H}$. Il est clair que $*\omega_H$ prend plus pr\'ecisemment ses valeurs dans la puissance ext\'erieure maximale de l'espace des 
vecteurs cotangents normaux \`a $X_H$. Notons $\tilde{\omega}_H$ et $*\tilde{\omega}_H$ les tir\'ees arri\`ere respectifs de $\omega_H$ et 
$* \omega_H$. Il est imm\'ediat que $\tilde{\omega}_H$ (resp. $*\tilde{\omega}_H$) n'est autre que le produit ext\'erieur d'une base orthonorm\'ee de 
$(\mathfrak{p} \cap \mathfrak{h} )^*$ (resp. $(\mathfrak{p} \cap \mathfrak{q} )^*)$.

Soit
\begin{eqnarray} \label{proj}
P_0 : \bigwedge {}^R \mathfrak{p}^* \rightarrow \left( \bigwedge {}^R \mathfrak{p}^* \right)_{\delta^*}
\end{eqnarray}
la projection orthogonale sur le $K$-type $\delta^*$ de plus haut poids $2\rho_n$. Soit $e(\mathfrak{q})^* \in \bigwedge^R \mathfrak{p}^*$ le vecteur dual \`a $e(\mathfrak{q})$,
c'est un vecteur de plus haut poids $-2\rho_n$.
Et, d'apr\`es le th\'eor\`eme d'Helgason (\ref{helgason}), 
$${\rm dim} \left( \bigwedge {}^R \mathfrak{p}^* \right)^{H\cap K}_{\delta^* } = 1.$$

\begin{prop}[Tong et Wang] {\rm \cite[Proposition 4.6]{TongWang}} \label{geom}
$P_0 (* \tilde{\omega}_H )$ est un vecteur non nul dans $\left( \bigwedge^R \mathfrak{p}^* \right)^{H \cap K}_{\delta^*}$.
\end{prop}
{\it D\'emonstration.} Puisque $*\tilde{\omega}_H$ est tir\'ee arri\`ere d'une forme diff\'erentielle $H$-invariante sur $X_H$, elle est $(H\cap K)$-invariante.
Le projet\'e $P_0 (* \tilde{\omega}_H )$ est donc $(H \cap K)$-invariant \`a droite et donc invariant sous Ad$^* (H \cap K)$~:
$$P_0 (* \tilde{\omega}_H ) \in \left( \bigwedge {}^R \mathfrak{p}^* \right)^{H\cap K}_{\delta^*} .$$
Il nous reste \`a v\'erifier que ce vecteur est non nul. Remarquons que le produit ext\'erieur maximal $\bigwedge^{d_G} \mathfrak{p}^*$ est un module trivial et que l'op\'erateur 
$*$ envoie un $K$-module sur son dual. Il suffit donc de montrer que 
$$*\tilde{\omega}_H \wedge * e(\mathfrak{q})^* \neq 0$$
autrement dit que 
$$\tilde{\omega}_H \wedge e(\mathfrak{q})^* \neq 0.$$
Ce dernier fait d\'ecoule imm\'ediatement du fait que $\tilde{\omega}_H$ et $e(\mathfrak{q})^*$ sont des vecteurs non nuls dans les puissances ext\'erieures maximales respectives des
espaces $(\mathfrak{p} \cap \mathfrak{h} )^*$ et $(\mathfrak{p} \cap \mathfrak{u})^*$.~$\Box$

\bigskip

\subsection{D\'emonstration du Th\'eor\`eme \ref{cohoml2+}}

Nous voulons montrer que la projection $[F]_{\delta}$ de la classe de cohomologie $[F] \in H^{R}_{2} (M)$ dans $H_2^R (M)_{\delta}$ est non nulle. Pour ce faire nous allons montrer que 
\begin{eqnarray} \label{bidule}
[\omega ] = c [F]_{\delta}
\end{eqnarray}
dans $H^R_2 (M)_{\delta}$, o\`u $c$ est une constante non nulle.

\medskip

Si $[\varphi ] \in H^R_2 (M)_{\delta}$, son tir\'e arri\`ere $\tilde{\varphi} = \sum_i \tilde{\varphi}_i X_i $ est harmonique et le $K$-espace engendr\'e par
$\tilde{\varphi}_i (g)$ est soit nul soit irr\'eductible isomorphe \`a $\delta$.
Chaque fonction $\overline{\tilde{\varphi}_i} \in L^2 (\Lambda \backslash G)$ engendre alors, sous l'action \`a droite de $K$, 
soit le sous-espace nul soit un sous-espace irr\'eductible isomorphe \`a $\delta^*$. D'apr\`es le lemme de Kuga \cite{BorelWallach}, le casimir $\Omega$
appliqu\'e \`a $\overline{\tilde{\varphi}_i }$ donne $0$. D'apr\`es le lemme suivant, les fonctions $\overline{\tilde{\varphi}_i (g)}$ engendrent un sous-module de $L^2 (\Lambda \backslash G)$ isomorphe \`a $\pi_{\rho}$.

\begin{lem}[Tong et Wang] {\rm \cite[Lemma 5.3]{TongWang}}
Soit $\pi$ une repr\'esentation unitaire de $G$ dans un espace de Hilbert $V$ et $A$ un module cohomologique de $K$-type minimal $\tau$.
Supposons que $W$ soit un $K$-type de type $\tau$ dans $V$ et qu'en restriction \`a $W$, le casimir $\pi (\Omega )$ soit nul. Alors, le $G$-sous-module ferm\'e $M$ engendr\'e 
par $W$ est isomorphe \`a $A$.
\end{lem}
{\it D\'emonstration.} Consid\'erons la restriction $\pi_{|M}$ de la repr\'esentation unitaire $\pi$ \`a $M$, celle-ci se d\'ecompose en une int\'egrale directe
$$\pi_{|M} = \int_{\widehat{G}} \sigma d\mu (\sigma ),$$
o\`u $\mu$ est une mesure sur $\widehat{G}$. Puisque $M$ est engendr\'e par $W$, pour presque tout $\sigma$, $\sigma (C) =0$ et $\tau \subset \sigma_{|K}$.
Mais son $K$-type minimal $\tau$ et son carat\`ere infinit\'esimal d\'etermine compl\`etement le module cohomologique $A$.
On en d\'eduit donc que presque tout $\sigma$ est isomorphe \`a $A$ et donc que la restriction de $\pi$ \`a $M$ est un multiple de $A$. Mais puisque $\tau$ 
intervient avec multiplicit\'e 1 dans $A$, tout $K$-sous-module irr\'eductible de type $\tau$ engendre un $G$-sous-module irr\'eductible. La multiplicit\'e 
de $A$ dans $\pi_{|M}$ est donc \'egale \`a $1$. Ce qui conclut la d\'emonstration du Lemme.~$\Box$

\bigskip

D'apr\`es (\ref{dualite}) et (\ref{isom hodge}), l'\'egalit\'e (\ref{bidule}) est \'equivalente \`a~: \'etant donn\'ee $[\varphi ] \in H^{R}_2 (M)_{\delta}$ telle que son tir\'e arri\`ere
$\tilde{\varphi} = \sum_i \tilde{\varphi}_i X_i$ soit harmonique et $\{ \tilde{\varphi}_i \}$ engendre une composante isotypique de plus haut poids $2\rho_n$, alors
\begin{eqnarray} \label{dualite2}
\int_M \omega \wedge * \varphi = c \int_F * \varphi .
\end{eqnarray}

\medskip

Soient $dg$, $dh$ et $dHg$ les mesures invariantes respectives sur $G$, $H$ et $H\backslash G$. On a alors~:
\begin{eqnarray*}
\int_M \omega \wedge * \varphi & = & \int_{H\backslash G} \left\{ \int_{\Lambda \backslash H} \langle \tilde{\omega} (hg) , \tilde{\varphi} (hg) \rangle dh \right\} dHg \\
                                                         & = & \int_{H \backslash G} \langle \tilde{\omega} (g) , \int_{\Lambda \backslash H} \tilde{\varphi} (hg) dh \rangle dHg ,
\end{eqnarray*}
o\`u $\langle . , . \rangle$ d\'esigne le produit scalaire standard sur $\bigwedge^R \mathfrak{p}^*$ et o\`u la derni\`ere \'egalit\'e d\'ecoule de l'invariance \`a gauche
de $\tilde{\omega}$ sous l'action de $H$. Soit
\begin{eqnarray} \label{d}
\tilde{\Phi}_H (g) = \int_{\Lambda \backslash H} \tilde{\varphi} (hg) dh
\end{eqnarray}
et soit $\tilde{\Phi}_H (g) = \sum_i \tilde{\Phi}_i X_i$ son d\'eveloppement. 

\begin{lem} \label{5.4}
\begin{enumerate}
\item La forme $\Phi_H$ est \`a valeurs complexes et harmonique sur $G$.
\item Sous l'action \`a droite de $K$, $\overline{\tilde{\Phi}_i (g)}$ engendre soit le sous-espace nulle soit  un sous-espace irr\'eductible isomorphe \`a $\delta^*$.
\item Si $\overline{\tilde{\Phi}_i (g)}$ est non nul, le casimir agit trivialement sur celui-ci~: $\Omega. \overline{\tilde{\Phi}_i (g)} =0$.
\item La forme $\tilde{\Phi}_H$ est invariante \`a gauche sous l'action de $H$.
\item La forme $\tilde{\Phi}_H$ est born\'ee sur $G$.
\end{enumerate}
\end{lem}
{\it D\'emonstration.} Les trois premiers points d\'ecoulent des propri\'et\'es analogues pour $\tilde{\varphi}$ d\'ecrites plus haut.
Le point 4. d\'ecoule de la d\'efinition (\ref{d}) de $\tilde{\Phi}_H$. Il nous reste \`a montrer que la forme $\tilde{\varphi}$ est born\'ee, ce qui d\'ecoule du lemme suivant.

\begin{lem} \label{borne}
Soit $\varphi$ une forme diff\'erentielle harmonique $L^2$ sur une vari\'et\'e riemannienne \`a coubure partout n\'egative (ou nulle) et uniform\'ement minor\'ee.
Il existe alors une constante $c>0$ ne d\'ependant que de la borne sur la courbure de $M$ telle que
\begin{eqnarray} \label{moser}
||\varphi ||_{L^{\infty} (M)} \leq c ||\varphi ||_{L^2 (M)} .
\end{eqnarray}
\end{lem}
{\it D\'emonstration.} La d\'emonstration repose sur le proc\'ed\'e classique d'it\'eration \`a la Nash-Moser. L'\'enonc\'e ci-dessus et sa d\'emonstration sont une l\'eg\`ere modification d'un lemme de Yeung, cf.
\cite{Yeung}.

Soit $\varphi$ un forme diff\'erentielle de degr\'e $i$.
Une formule de Bochner-Weitzenbock que l'on obtient par un calcul direct permet de comparer le laplacien usuel et le {\it laplacien brut} $\Delta^b= - \nabla^* \nabla$~:
$$\langle \Delta \varphi , \varphi \rangle = \langle \Delta^b \varphi , \varphi \rangle + i \langle {\cal R} \varphi , \varphi \rangle , $$
o\`u ${\cal R}$ est un op\'erateur ne d\'ependant que du tenseur de courbure tel que
$$\langle {\cal R} \varphi , \varphi \rangle \geq -k \langle \varphi , \varphi \rangle ,$$
o\`u $k$ est une constante $\geq 0$ ne d\'ependant que de la minoration sur la courbure.
Supposons maintenant $\varphi$ harmonique, alors $\Delta \varphi = 0$ et puisque
$$\Delta^b |\varphi |^2 = \langle \Delta^b \varphi , \varphi \rangle + |\nabla \varphi |^2,$$
on obtient
\begin{eqnarray} \label{aa}
\Delta^b |\varphi |^2 + k |\varphi |^2 \geq 0 .
\end{eqnarray}

Nous allons maintenant pouvoir appliquer le proc\'ed\'e d'it\'eration \`a la Nash-Moser. Soit $f =|\varphi |^2$.

Commen\c{c}ons par remarquer que puisque
$$ \int g \Delta^b f    = -\int \langle \nabla f , \nabla g \rangle ,$$
en prenant $g= \eta^2 f^{\beta}$ (o\`u $\eta$ est une fonction et $\beta$ un r\'eel $\geq 1$), on obtient l'identit\'e suivante~:
$$\beta \int \eta^2 f^{\beta-1} |\nabla f |^2 = - 2 \int_M \langle \eta f^{\frac{\beta -1}{2}} \nabla f , f^{\frac{\beta +1}{2}} \nabla \eta \rangle - \int \eta^2 f^{\beta} \Delta^b f .$$
Or,
\begin{eqnarray*}
\int |\nabla (f^{\frac{\beta +1}{2}} \eta )|^2 & = & \int \langle \eta \nabla (f^{\frac{\beta +1}{2}} ) + f^{\frac{\beta +1}{2}} \nabla \eta , \eta \nabla (f^{\frac{\beta +1}{2}} ) + f^{\frac{\beta +1}{2}} \nabla \eta \rangle \\
& = & \left( \frac{\beta+1}{2} \right)^2 \int \eta^2 f^{\beta -1} |\nabla f|^2 + \int f^{\beta+1} |\nabla \eta|^2 \\
& & + (\beta +1) \int \langle  \eta f^{\frac{\beta -1}{2}} \nabla f , f^{\frac{\beta +1}{2}} \nabla \eta \rangle .
\end{eqnarray*}
D'o\`u l'on d\'eduit~:
\begin{eqnarray*}
\frac{4\beta}{(\beta +1)^2} \int |\nabla (f^{\frac{\beta+1}{2}} \eta ) |^2  & = &
 2 \frac{\beta -1}{\beta +1} \int \langle \eta f^{\frac{\beta -1}{2}} \nabla f , f^{\frac{\beta +1}{2}} \nabla \eta \rangle \\
& & - \int \eta^2 f^{\beta} \Delta^b f  + \frac{4\beta}{(\beta +1)^2}  \int f^{\beta+1} |\nabla \eta|^2 .
\end{eqnarray*}
En appliquant l'in\'egalit\'e de Cauchy-Schwarz au premier terme \`a droite, on obtient alors que pour tout $\varepsilon >0$~:
\begin{eqnarray} \label{fond}
\begin{array}{cl}
\frac{4\beta}{(\beta +1)^2} \int |\nabla (f^{\frac{\beta+1}{2}} \eta ) |^2 & \leq \varepsilon \int \eta^2 f^{\beta -1} |\nabla f|^2 + \frac{1}{\varepsilon} \int f^{\beta +1} |\nabla \eta|^2 \\
& - \int \eta^2 f^{\beta} \Delta^b f 
+ \frac{4\beta}{(\beta +1)^2}  \int f^{\beta+1} |\nabla \eta|^2 .
\end{array}
\end{eqnarray}

Soit $C$ une constante de Sobolev valable pour toutes les boules isom\'etriquement plong\'ees dans $M$ (une telle constante existe en vertu des restrictions
faites sur la courbure). Le fait suivant traduit alors l'in\'egalit\'e de Sobolev standard.

\medskip
\noindent
{\bf Fait 1.} Pour toute boule $B$ isom\'etriquement plong\'ee dans $M$, et toute fonction $u \in H^1 (B)$,
$$||u^{\frac{n}{n-2}} ||_{L^2 (B)}^2 \leq C \left[ ||\nabla u||^{\frac{2n}{n-2}}_{L^2 (B)} + || u||^{\frac{2n}{n-2}}_{L^2 (B)} \right],$$
o\`u $n$ est la dimension de $M$.

\medskip

Soient maintenant $B_b \subset B_a$ deux boules de m\^eme centre et de rayons respectifs $a$ et $b$ dans $M$ et soit $\eta$ une fonction \`a support dans $B_a$, partout $\leq 1$,
constante \'egale \`a $1$ sur $B_b$ et telle que $|\nabla \eta | \leq 2/(a-b)$. Les expressions (\ref{aa}) et (\ref{fond}) impliquent le fait suivant.

\medskip 
\noindent
{\bf Fait 2.} Il existe une constante $k_1$ ne d\'ependant que de $k$ et du rayon d'injectivit\'e de $M$ (mais ni de $a$ ni de $\beta$) telle que~:
$$||\nabla f^{\frac{\beta +1}{2}} ||_{L^2 (B_b)}^2 \leq \frac{k_1 (\beta +1 )}{(a-b)^2} ||f^{\beta +1}||^2_{L^2 (B_a )}.$$

\medskip

En effet, d'apr\`es (\ref{aa}) et (\ref{fond}),
\begin{eqnarray*}
\frac{4\beta}{(\beta +1)^2} \int |\nabla (f^{\frac{\beta+1}{2}} \eta ) |^2 & \leq \varepsilon \int \eta^2 f^{\beta -1} |\nabla f|^2 + \frac{1}{\varepsilon} \int f^{\beta +1} |\nabla \eta|^2  \\
& +k \int \eta^2 f^{\beta+1}   + \frac{4\beta}{(\beta +1)^2}  \int f^{\beta+1} |\nabla \eta|^2 .
\end{eqnarray*}
Or,
\begin{eqnarray*}
\int \eta^2 f^{\beta -1} |\nabla f|^2 & = & \frac{4}{(\beta +1)^2} \int | \nabla (f^{\frac{\beta +1}{2}} \eta ) - f^{\frac{\beta +1}{2}} \nabla \eta |^2 \\
                                                           & \leq & \frac{8}{(\beta +1)^2} \left[ \int |\nabla (f^{\frac{\beta +1}{2}} \eta )|^2 + \int f^{\beta +1} |\nabla \eta |^2 \right].
\end{eqnarray*}
D'o\`u l'on d\'eduit, en prenant $\varepsilon = \beta /4$,
$$\frac{2\beta}{(\beta +1)^2} \int |\nabla (f^{\frac{\beta +1}{2}} \eta ) |^2 \leq \frac{6\beta^2 +4(\beta +1 )^2}{\beta (\beta +1)^2} \int f^{\beta+1} |\nabla \eta |^2 + k \int \eta^2 f^{\beta +1}.$$
Et, puisque $\eta \leq 1$ et $|\nabla \eta | \leq \frac{2}{a-b}$,
$$\frac{2\beta}{(\beta +1)^2} \int_{B_a} |\nabla (f^{\frac{\beta +1}{2}} \eta ) |^2 \leq \left(\frac{10}{\beta (a-b)^2} + k \right) \int_{B_a} f^{\beta +1} .$$
On conclut imm\'ediatement la d\'emonstration du fait 2 en utilisant que la fonction $\eta$ est constante \'egale \`a $1$ sur $B_b$.

\medskip

Posons maintenant $A(l,r ) = (\int_{B_r} |f|^l )^{1/l}$ et $\alpha = n/(n-2)$. En prenant $\beta = l-1$, les faits 1 et 2 impliquent l'in\'egalit\'e suivante~:
\begin{eqnarray} \label{z}
A( \alpha l , b ) \leq \frac{k_2^{1/l} l^{1/l}}{(a-b)^{2/l}} A(l,a) ,
\end{eqnarray}
pour tout $l\geq 2$ et $k_2$ une constante ne d\'ependant que de $k$, de $C$ et du rayon d'injectivit\'e de $M$.

Nous pouvons maintenant it\'erer l'in\'egalit\'e (\ref{z}). Soient donc $r_0$ un r\'eel strictement positif et inf\'erieur au rayon d'injectivit\'e de $M$ et 
$r_m =(r_0 /2) (1+ 1/2^m)$ pour $m=0 , 1 , 2 , \ldots$. L'in\'egalit\'e (\ref{z}) implique alors~:
$$A(2 \alpha^{m+1} , r_m ) \leq c A(2, r_0 ),$$
o\`u $c$ est une constante ne d\'ependant que de $k_2$ et $n$ (et donc ind\'ependante de $m$). En passant \`a la limite $m \rightarrow + \infty$, ceci implique 
$$\sup_{B_{r_0 /2}} f \leq c \left( \int_{B_{r_0}} f^2 \right)^{1/2} .$$
Puisque 
$$\int_{B_{r_0}} f^2 \leq (\sup_{B_{r_0}} f ) . || \varphi ||_{L^2 (M)} .$$
On en d\'eduit que pour tout nombre r\'eel $r_0$ inf\'erieur au rayon d'injectivit\'e et pour toute paire de boules $B_{r_0 /2} \subset B_{r_0}$ dans $M$ centr\'ees en un
m\^eme point et de rayons respectifs $r_0 /2$ et $r_0$, 
$$\sup_{B_{r_0 /2}} f \leq c (\sup_{B_{r_0}} f)^{1/2} || \varphi ||_{L^2 (M)}.$$
En utilisant la continuit\'e de $f$, on en d\'eduit imm\'ediatement que $f$ est n\'ecessairement born\'ee sur $M$ tout entier par
$\sqrt{c} || \varphi ||_{L^2 (M)}^{1/2}$ ce qui conclut la d\'emonstration du Lemme \ref{borne}.~$\Box$

\bigskip

\begin{prop} \label{5.5}
Il existe une constante $c_{\varphi}$ (d\'ependant de $\varphi$) telle que 
$$\tilde{\Phi}_H  = c_{\varphi} \tilde{\omega} .$$
\end{prop}
{\it D\'emonstration.} On peut supposer $\tilde{\Phi}_H$ non nulle. D'apr\`es le troisi\`eme point du Lemme \ref{5.4}, le casimir agit trivialement sur 
$\overline{\tilde{\Phi}_i (g)}$ or d'apr\`es le cinqui\`eme  point cette fonction est born\'ee et appartient donc \`a $L^2 (H \backslash G)$. Enfin, puisque sous l'action 
\`a droite de $K$ elle engendre le $K$-type $\delta^*$, la fonction $\overline{\tilde{\Phi}_i (g)}$ engendre un sous-module de $L^2 (H \backslash G)$
isomorphe au module $\pi_{\rho}$. Un tel module appara\^{\i}t avec multiplicit\'e 1. La fonction 
$\overline{\tilde{\Phi}_i (g)}$ appartient donc au $K$-type $\delta^*$ dans $\pi_{\rho}$, celui-ci est irr\'eductible, la fonction $\tilde{\Phi}_H$ est donc
n\'ecessairement \'egale \`a $ c_{\varphi} \tilde{\omega}$ pour une certaine constante $c_{\varphi}$.

\bigskip

Dans la suite, nous aurons \`a introduire de nombreuses constantes non nulles ind\'ependantes de $\varphi$ (et de $\Lambda$). Nous les appellerons
{\it constantes universelles} et nous les noterons $c_1 , c_2 , \ldots$.
Le lemme suivant va nous permettre de mieux comprendre la constante $c_{\varphi}$.

\begin{lem} \label{5.6}
L'int\'egrale
$$\int_{\Lambda \backslash H} \langle \tilde{\omega} (h) , \tilde{\varphi} (h) \rangle dh = c_1  \int_F * \varphi$$
o\`u $c_1$ est une constante universelle.
\end{lem}
{\it D\'emonstration.} D'apr\`es les Propositions \ref{omega} et \ref{geom},
$$\tilde{\omega} (h) = c_1 P_0 (*\tilde{\omega}_H ) .$$
Puisque $\tilde{\varphi}(g) \in \left( \bigwedge^R \mathfrak{p}^* \right)_{\delta^*}$ on a~:
\begin{eqnarray} \label{i}
\langle \tilde{\omega} (h),\tilde{\varphi} (h) \rangle = c_1 \langle *\tilde{\omega}_H , \tilde{\varphi} (h) \rangle .
\end{eqnarray}
Soit $\psi = * \varphi$. Alors, 
\begin{eqnarray} \label{iia}
\langle *\tilde{\omega}_H , \tilde{\varphi} (h) \rangle d{\rm vol}_{X_G} = * \omega_H \wedge \psi (hK).
\end{eqnarray}
Soit $i$ le plongement de $F \rightarrow M$, alors
\begin{eqnarray} \label{iii}
*\omega_H \wedge \psi (hK) = *\omega_H \wedge (i^* \psi (hK)).
\end{eqnarray}
Cette derni\`ere \'egalit\'e r\'esulte du fait que $*\omega_H$ prend ses valeurs dans le produit ext\'erieur maximal de l'espace des vecteurs cotangents normaux \`a $X_H$
dans $X_G$. Les composantes de $\psi (hK)$ selon des vecteurs cotangents normaux \`a $X_H$ sont donc tu\'ees par le produit ext\'erieur avec $*\omega_H$.

Puisque $(d{\rm vol}_{X_G})_{|X_H} = \omega_H \wedge * \omega_H$ avec $\omega_H$ et $*\omega_H$ des sections inversibles d'un fibr\'e en droites, on 
d\'eduit de (\ref{i}), (\ref{iia}) et (\ref{iii})~:
$$\langle  \tilde{\omega} (h),\tilde{\varphi} (h) \rangle \omega_H = c_1 i^* \psi (hK).$$
Ce qui conclut la d\'emonstration du Lemme \ref{5.6}.~$\Box$

\bigskip

On peut maintenant d\'eterminer la constante $c_{\varphi}$. D'apr\`es la Proposition \ref{5.5} et le Lemme \ref{5.6},
\begin{eqnarray*}
c_1 \int_F * \varphi & = & \int_{\Lambda \backslash H} \langle  \tilde{\omega} (h),\tilde{\varphi} (h) \rangle dh \\
& = & \langle \tilde{\omega} (e) , \int_{\Lambda \backslash H} \tilde{\varphi} (h) dh \rangle \\
& = & c_{\varphi} \langle  \tilde{\omega} (e) , \tilde{\omega} (e) \rangle .
\end{eqnarray*}
On a finalement montr\'e que
\begin{eqnarray} \label{a}
c_{\varphi} = c_2 \int_F * \varphi , 
\end{eqnarray}
o\`u $c_2$ est une constante universelle. D'un autre c\^ot\'e et d'apr\`es la Proposition \ref{5.5}, on a~:
\begin{eqnarray} \label{b}
\begin{array}{lcl}
\int_M \omega \wedge * \varphi  & = & c_{\varphi} \int_{H\backslash G} \langle \tilde{\omega} (g) , \tilde{\omega} (g) \rangle dHg \\
& = & c_{\varphi} c_3 . 
\end{array}
\end{eqnarray}
Le Th\'eor\`eme \ref{cohoml2+} d\'ecoule maintenant imm\'ediatement de (\ref{a}) et (\ref{b}).~$\Box$

\subsection{Application aux groupes unitaires et orthogonaux}

\subsubsection*{Groupes unitaires}

Dans ce paragraphe $G$ est un groupe alg\'ebrique r\'eductif, connexe et anisotrope sur ${\Bbb Q}$ tel que $G^{{\rm nc}} = U(p,q+r)$. On suppose
fix\'ee une donn\'ee $Sh^0 H \subset Sh^0 G$ avec $H^{{\rm nc}} = U(p,q)$ plong\'e de mani\`ere standard dans $G^{{\rm nc}}$ et $\Lambda$ un sous-groupe de
congruence sans torsion de $H$. Notons enfin $M=\Lambda \backslash X_G$ et $F= \Lambda \backslash X_H$. 

Soient $(\lambda , \mu )$ un couple compatible de partitions $\subset p\times (q+r)$.
Conform\'ement aux notations des sections pr\'ec\'edentes nous notons $H^*_2 (M)_{\lambda , \mu} = H_2^* (A(\lambda , \mu ) : M)$ la $A(\lambda , \mu )$-composante
de la cohomologie $L^2$ de $M$, enfin nous notons $H_2^{\lambda , \mu} (M)$ la partie {\it fortement primitive} $H_2^{|\lambda | + |\hat{\mu}|} (M)_{\lambda , \mu}$.

Commen\c{c}ons par montrer que le Th\'eor\`eme \ref{cohoml2+} implique le corollaire suivant.

\begin{cor} \label{cohoml2U}
La classe $[F] \in H_2^{d_G -d_H} (M)_{(r^p) , (q^p)}$ est non nulle si et seulement si $d_H \geq d_G /2$ ({\it i.e.} si et seulement si $q\geq r$).
\end{cor}
{\it D\'emonstration.} Nous appliquons le Th\'eor\`eme \ref{cohoml2} au groupe $G=U(p,q+r)$ muni de l'involution $\tau$ standard telle que
$H=G^{\tau} = U(p,q) \times U(r)$. Commen\c{c}ons par remarquer que rang$_{\Bbb C} (G/H) =$rang$_{\Bbb C} (K/(K\cap H))$ si et seulement si
$q \geq r$. Et dans ce cas rang$_{\Bbb C} (G/H) = \dim_{{\Bbb C}} \mathfrak{t} = r$.

Consid\'erons maintenant la repr\'esentation $\pi_{\rho}^*$. Nous avons vu au cours de la d\'emonstration du Th\'eor\`eme \ref{cohoml2+} que celle-ci est isomorphe
\`a la repr\'esentation cohomologique de $K$-type minimal $2\rho_n$. Il nous suffit donc de calculer $2\rho_n$.
Comme au \S 1.1, nous consid\'erons la sous-alg\`ebre de Cartan $\widetilde{\mathfrak{t}}$ dans $\mathfrak{g}$ constitu\'ee des matrices diagonales.
Nous pouvons prendre la sous-alg\`ebre de Cartan $\mathfrak{t}$ dans $\mathfrak{q}$ telle que 
$i\mathfrak{t}_0$ soit \'egale au sous-ensemble
$$\{ (\underbrace{0,\ldots , 0}_p ; \underbrace{-u_r , \ldots , - u_1 , 0 , \ldots , 0, u_1 , \ldots , u_r }_{q+r} ) \; : \; u_i \in {\Bbb R}  \mbox{ pour tout } i= 1, \ldots ,r \} $$
de $i\widetilde{\mathfrak{t}}_0$.
Repr\'esentons les syst\`emes de racines respectifs de $\mathfrak{k}$ et $\mathfrak{g}$ comme $\Delta_c = \{ \pm (x_i - x_j ) \; : \; 1 \leq i< j \leq p \} \cup
\{ \pm (y_i -y_j ) \; : \; 1 \leq i < j \leq q+r \}$ et $\Delta = \Delta_c \cup \{ \pm (x_i - y_j ) \; : \; 1 \leq i \leq p \; {\rm et } \; 1 \leq j \leq q+r \}$.
On peut alors consid\'erer les syst\`emes positifs compatibles $\Delta_c^+ = \{ x_i - x_j  \; : \; 1 \leq i< j \leq p \} \cup \{ y_j -y_i \; : \;  1\leq  i < j \leq q+r \}$,
$\Sigma_c^+ = \{ y_j - y_i \; : \; i< j \mbox{ et soit } 1 \leq i \leq r \mbox{ soit } q+1 \leq j \leq q+r \} $ et
$\Sigma^+ = \Sigma_c^+ \cup \{x_i -y_j \; : \; 1\leq i \leq p \mbox{ et } 1 \leq j \leq r \} \cup \{ y_j -x_i \; : \; 1 \leq i \leq p \mbox{ et } q+1 \leq j \leq q+r \}$.
Un calcul simple montre alors que 
$$2 \rho_n = p \sum_{j=1}^r (y_{q+j} -y_j) $$
{\it i.e.}, le plus haut poids par rapport \`a $\Delta_c^+$ de la repr\'esentation $V((r^p), (q^p ))$. Ce qui conclut la d\'emonstration du Corrolaire \ref{cohoml2U}.~$\Box$

\bigskip

\noindent
{\bf Remarque.} Dans \cite{Schlichtkrull} Schlichtkrull identifie les param\`etres de Langlands des repr\'esentations de Flensted-Jensen.
Avec ses notations les syst\`emes de racines respectifs de $\mathfrak{g}$ et 
$\mathfrak{k}$ sont $\Delta = \{ \pm (e_i - e_j ) \; : \; 1 \leq i < j \leq p+q+r \} $ et $\Delta_c = \{ \pm (e_i - e_j ) \; : \; 1 \leq i< j \leq p \; {\rm ou} \; p < i < j \leq p+q+r \}$. 
Il fixe comme sous-syst\`eme positif $\Delta_c^+ = \{ e_i - e_j  \; : \; 1 \leq i< j \leq p \; {\rm ou} \; p < i < j \leq p+q+r \}$. Par rapport \`a ce syst\`eme positif, le $K$-type
$V_r = V((r^p) , (q^p))$ a pour plus haut poids :
$$2 \rho_r := 2\rho ( \mathfrak{u} ((r^p) ,(q^p)) \cap \mathfrak{p}) = - p \sum_{j=1}^{r} (e_{p+j} - e_{p+q+r+1-j} )$$
\'egal \`a la somme des racines de $\mathfrak{u} ((r^p) , (q^p)) \cap \mathfrak{p}$ par rapport \`a $\mathfrak{k}$.
Soit $P_r = MAN$ le parabolique cuspidal et $\delta_r \in \hat{M}$ la repr\'esentation de la s\'erie discr\`ete associ\'es au $K$-type $2\rho_r$ par Vogan \cite{Vogan3}.
Ceux-ci co\"{\i}ncident avec le parabolique et la repr\'esentation de la s\'erie discr\`ete que l'on associe \`a une sous-alg\`ebre parabolique, dans notre cas $\mathfrak{q} (
(r^p), (q^p))$, dans \cite[\S 5.2]{BergeronClozel}. Le groupe $A$ est alors la composante neutre d'un tore maximal d\'eploy\'e du sous-groupe de Levi
$L((r^p) , (q^p)) \cong U(p,q-r)$. Soit $k = \min (p, q-r)$ la dimension de $\mathfrak{a}$. 
Dans \cite[\S 5.2]{BergeronClozel}, nous associons \'egalement \`a la sous-alg\`ebre parabolique $\mathfrak{q} ((r^p) , (q^p))$ un \'el\'ement $\nu \in \mathfrak{a}^*$. 
Si l'on choisit une base $(f_1 , \ldots , f_k)$ de $\mathfrak{a}$ telle que les racines de $\mathfrak{a}$ dans $N$  soient
$$\{ \frac12 (f_i \pm f_j ) , \frac12 f_l , f_l \; : \; 1 \leq i < j \leq k, \; 1 \leq l \leq k \},$$
alors  
$$\nu = \frac12 (p+q-r-1 , \ldots , p+q-r+1-2k ) \in \mathfrak{a}^* .$$
Nous avons donc d\'etermin\'e une  repr\'esentation $\delta_{r} \otimes \nu \otimes 1$ de $P_{r} = MAN$. Comme rappel\'e dans \cite[\S 5.2]{BergeronClozel}, il d\'ecoule alors de
\cite{VoganZuckerman} que la repr\'esentation unitairement induite
$$I(\delta_{r} , \nu ) = {\rm ind}_{P_{r}}^G (\delta_{r} \otimes \nu \otimes 1)$$
admet un unique quotient irr\'eductible $J(\delta_{r} , \nu )$ qui est isomorphe \`a la repr\'esentation cohomologique $A((r^p) ,(q^p))$.
Enfin et puisque $q \geq r$, il d\'ecoule de \cite[Theorem 7.7]{Schlichtkrull} que la repr\'esentation
$J(\delta_{r} , \nu )^*$ peut-\^etre r\'ealis\'ee sur un sous-espace ferm\'e de $L^2 (G / H)$. 

\medskip

La Conjecture \ref{conj2} est principalement motiv\'ee par la conjecture suivante.

\begin{conj} \label{conjl2} (Sous les hypoth\`eses du Corollaire \ref{cohoml2U}.)
Si $\lambda$ et $\mu$ sont deux partitions incluses dans $p\times q$ formant un couple compatible avec 
$\mu / \lambda = (p_1 \times q_1) * \ldots * (p_m \times q_m )$, 
alors l'application
Alors, l'application
$$H^{\lambda , \mu} (F) \rightarrow H^{|\lambda| +pr +|\hat{\mu}|}_{2, \; {\rm prim} +} (M)$$
obtenue en composant l'application ``cup-produit avec $[F]$'' et la projection sur la composante fortement primitive de la cohomologie $L^2$ de $M$ 
est {\bf injective} si et seulement si la partition $(r^p)$ s'inscrit dans le diagramme gauche $\mu /\lambda$ ({\it i.e.} si $p_1 + \ldots + p_m = p$ et $r\leq q_i$ pour $i=1 , \ldots ,m$).
Son image est alors contenue dans $H^{\lambda +(r^p), \mu}_2 (M)$.
\end{conj}

Le Corollaire \ref{cohoml2U} est un cas particulier de cette Conjecture (le cas $\lambda = \emptyset$, $\mu = p\times q$).

La combinatoire pr\'edite par la conjecture provient du Lemme suivant. Soit $(\lambda , \mu)$ un couple compatible de partitions dans $p\times q$. 
Le groupe $K_{\Bbb C} = GL_p \times GL_{q+r}$, il contient le groupe $GL_p \times( GL_q \times GL_r)$, o\`u le groupe $GL_q \times GL_r$ est
plong\'e dans $GL_{q+r}$ par
$$(A,B) \mapsto \left(
\begin{array}{cc}
A & 0 \\
0 & B
\end{array}
\right) .$$
Le groupe $GL_p \times (GL_q \times GL_r)$ pr\'eserve la d\'ecomposition orthogonale
$\mathfrak{p} = (\mathfrak{p} \cap \mathfrak{h}) \oplus (\mathfrak{p} \cap \mathfrak{q})$. 
Dans la suite de ce paragraphe, nous noterons $V_H (\lambda , \mu)$ le sous-$(GL_p \times GL_q)$-module de $\bigwedge (\mathfrak{p} \cap \mathfrak{h})$
associ\'e au couple compatible $(\lambda , \mu)$ de partitions $\subset p\times q$. Remarquons que le $GL_p \times (GL_q \times GL_r)$-module
$V_H (\lambda , \mu) \otimes \bigwedge^{R} (\mathfrak{p} \cap \mathfrak{q})$ est un sous-module de $\bigwedge^{|\lambda|+|\hat{\mu}|} (\mathfrak{p} \cap \mathfrak{h})
\otimes \bigwedge^R (\mathfrak{p} \cap \mathfrak{q})  \subset \bigwedge^{|\lambda|+|\hat{\mu}|+R} \mathfrak{p}$.

\begin{lem} \label{j}
Supposons que la partition $(r^p)$ s'inscrive dans le diagramme gauche $\mu / \lambda$. Alors, 
le $K$-module $V(\lambda + (r^p) , \mu)$, vu comme $GL_p \times (GL_q \times GL_r)$-module, contient (avec multiplicit\'e 1) le 
module $V_H (\lambda , \mu) \otimes \bigwedge^R (\mathfrak{p} \cap \mathfrak{q})$.
\end{lem}
{\it D\'emonstration.} Comme au \S 1.1, notons $v_H (\lambda ) \otimes w_H (\hat{\mu})^*$ le vecteur de plus haut poids de $V_H (\lambda , \mu )$.
Soit $\xi$ un g\'en\'erateur de la droite $\bigwedge^{rp} (\mathfrak{p}^+ \cap \mathfrak{q})$. 

Le vecteur $w_H (\hat{\mu}) \otimes \xi \in \bigwedge^{|\hat{\mu}|} (\mathfrak{p}^+ \cap \mathfrak{h})  \otimes 
\bigwedge^{rp} (\mathfrak{p}^+ \cap \mathfrak{q}) \subset \bigwedge^{|\hat{\mu}| +rp} \mathfrak{p}^+$ est alors un vecteur de plus bas poids qui engendre le sous-$K$-module 
$V(\hat{\mu})$ de $ \bigwedge^{|\hat{\mu}| +rp} \mathfrak{p}^+$ puisqu'il est colin\'eaire au vecteur $w(\hat{\mu})$ \footnote{Ici $\hat{\mu}$ d\'esigne le compl\'ementaire de 
$\mu$ dans $p\times (q+r)$ alors que plus haut (lorsque l'on se place dans $\mathfrak{p} \cap \mathfrak{h}$) $\hat{\mu}$ d\'esigne le compl\'ementaire de $\mu$ dans 
$p\times q$.}.
On peut de plus facilement v\'erifier que le vecteur $v(\lambda + (r^p))$ est dans l'orbite du vecteur $v_H (\lambda) \otimes \xi$ sous l'action du sous-groupe de Borel de $K_{{\Bbb C}}$
fix\'e au \S 1.1. 

Le vecteur $w_H (\hat{\mu})^* \otimes \xi^* \in \bigwedge^{|\hat{\mu}|+rp} \mathfrak{p^-}$ est un vecteur de plus haut poids, il engendre donc une droite sous l'action du sous-groupe de
Borel de $K_{{\Bbb C}}$, le vecteur $v(\lambda + (r^p)) \otimes w(\hat{\mu})^* $ appartient donc \`a l'orbite sous l'action de $K$ du vecteur
$(v_H (\lambda ) \otimes w_H (\hat{\mu})^* ) \otimes (\xi \otimes \xi^*) \in V_H (\lambda , \mu ) \otimes \bigwedge^R (\mathfrak{p} \cap \mathfrak{q}) \subset \bigwedge^{|\lambda| +
|\hat{\mu}| + R} \mathfrak{p}$. Puisque $V(\lambda + (r^p), \mu  )$ appara\^{\i}t avec multiplicit\'e 1 dans $\bigwedge^{|\lambda| +
|\hat{\mu}| + R} \mathfrak{p}$, le Lemme \ref{j} est d\'emontr\'e.~$\Box$

\bigskip

\`A l'aide de ce Lemme et du Corollaire \ref{cohoml2U}, il semble raisonnable de penser pouvoir construire une application
$$H^{\lambda , \mu} (F) \rightarrow H_2^{\lambda +(r^p) , \mu} (M)$$
qui co\"{\i}ncide avec celle consid\'er\'ee dans la Conjecture \ref{conjl2} et qui soit injective si et seulement si
la partition $(r^p)$ s'inscrit dans le diagramme gauche $\mu / \lambda$. Malheureusement nous ne sommes parvenus \`a le v\'erifier que dans le cas $p=1$ (d\'ej\`a trait\'e
par une autre m\'ethode dans \cite{BergeronClozel}). 

\bigskip

Concluons ce paragraphe par quelques remarques et cons\'equences de la Conjecture \ref{conjl2}. 

L'application ``cup-produit avec $[F]$''
$$H^{\lambda , \mu} (F ) \rightarrow H^{|\lambda | + |\hat{\mu}| + 2pr}_2 (M)$$
devrait \^etre injective d\`es que $(r^p)$ s'inscrit dans $\mu / \lambda$. Cette condition est-elle n\'ecessaire ?

Si la Conjecture \ref{conjl2} est vraie, 
pour tout entier $k\leq q-r$, l'application ``cup-produit avec $[F]$'' 
\begin{eqnarray} \label{c}
H^{k} (F ) \rightarrow H^{k+ 2pr}_2 (M)
\end{eqnarray}
devrait \^etre {\bf injective}.

Il nous suffit en effet de d\'emontrer que si $(\lambda , \mu )$ est un couple compatible de partitions dans $p\times q$ telles
que $|\lambda | + |\hat{\mu}| \leq q- r$ (si $p \geq 2$ on peut en fait remplacer $q-r$ par $q$) alors la partition $(r^p)$ s'inscrit dans le diagramme gauche $\mu / \lambda$. 
Ce qui d\'ecoule du fait suivant.

\medskip

\noindent
{\bf Fait.} Soient $p_1 , q_1, \ldots , p_m , q_m $ des entiers $\geq 1$. Supposons soit que $p_1 + \ldots + p_m \leq p-1$ soit que l'un des $q_i$ soit $<r$, alors 
$p_1 q_1 + \ldots + p_m q_m < pq - q + r $.

\medskip

En effet, 
$$p_1 q_1 + \ldots + p_m q_m \leq  \sum_i p_i (q_i -r ) + r \sum_i p_i .$$
Donc si l'un des $q_i$, par exemple $q_{1}$, est $<r$, 
$$p_1 q_1 + \ldots + p_m q_m < (q -r ) \sum_{i\neq 1} p_i + r \sum_i p_i \leq (q-r-1)(p-1) +rp= pq  -q-p +r .$$
Alors que si $p_1 + \ldots + p_m \leq p-1$,
$$p_1 q_1 + \ldots + p_m q_m \leq (q-r) (p-1) + r(p-1) = pq -q < pq -q +r .$$

\bigskip

Dans \cite{BergeronClozel} l'injectivit\'e de l'application ``cup-produit'' (\ref{c}) est obtenue pour 
$k< q- pr$, et dans ce cas nous montrons  \cite[Th\'eor\`eme 4.0.6]{BergeronClozel} que l'application est en fait un isomorphisme. Est-ce encore le cas pour tout
$k \leq q-r$ ?

\medskip

\subsubsection*{Groupes orthogonaux}

Dans ce paragraphe $G$ est un groupe alg\'ebrique r\'eductif, connexe et anisotrope sur ${\Bbb Q}$ tel que $G^{{\rm nc}} = O(p,q+r)$. On suppose fix\'ee une donn\'ee
$Sh^0 H \subset Sh^0 G$ avec $H^{{\rm nc}} = O(p,q)$ plong\'e de mani\`ere standard dans $G^{{\rm nc}}$ et $\Lambda$ un sous-groupe de congruence sans torsion de $H$.
Notons enfin $M= \Lambda \backslash X_G$ et $F= \Lambda \backslash X_H$.

Soit $\lambda$ une partition orthogonale $\subset p \times (q+r)$. Nous notons $H_2^* (M)_{\lambda}^{\pm_1 , \pm_2} = H_2^* (A(\lambda )_{\pm_1}^{\pm_2} \; : \; M)$ la
$A(\lambda )_{\pm_1}^{\pm_2}$-composante de la cohomologie $L^2$ de $M$, Nous notons plus g\'en\'eralement $H_2^* (M)_{\lambda}$ la somme directe de tous les
$H_2^* (M)_{\lambda}^{\pm_1 , \pm_2}$ lorsque les signes $\pm_1$ et $\pm_2$ varient.

Le Th\'eor\`eme \ref{cohoml2+} implique le corollaire suivant.

\begin{cor} \label{cohoml2O}
La classe de $[F] \in H_2^{d_G -d_H} (M)_{(r^p)}$ est non nulle si et seulement si $d_H \geq d_G /2$ ({\it i.e.} si et seulement si $q \geq r$).
\end{cor}
{\it D\'emonstration.} Nous appliquons le Th\'eor\`eme \ref{cohoml2} au groupe $G=O(p,q+r)$ muni de l'involution $\tau$ standard telle que
$H=G^{\tau} = O(p,q) \times O(r)$. Commen\c{c}ons par remarquer que rang$_{\Bbb C} (G/H) =$rang$_{\Bbb C} (K/(K\cap H))$ si et seulement si
$q \geq r$. Et dans ce cas rang$_{\Bbb C} (G/H) = \dim_{{\Bbb C}} \mathfrak{t} = r$. Dans la suite de la d\'emonstration
nous supposons $p$ et $q+r$ pairs respectivement \'egaux \`a $2\alpha$ et $2\beta$. (Les autres cas se traitent de mani\`ere
similaire.)

Consid\'erons maintenant la repr\'esentation $\pi_{\rho}^*$. Nous avons vu au cours de la d\'emonstration du Th\'eor\`eme \ref{cohoml2+} que celle-ci est isomorphe
\`a la repr\'esentation cohomologique de $K$-type minimal $2\rho_n$. Il nous suffit donc de calculer $2\rho_n$.
Comme au \S 1.2, nous consid\'erons la sous-alg\`ebre de Cartan $\widetilde{\mathfrak{t}}$ dans $\mathfrak{g}$ constitu\'ee des matrices diagonales.
Nous pouvons prendre la sous-alg\`ebre de Cartan $\mathfrak{t}$ dans $\mathfrak{q}$ telle que
$$i\mathfrak{t}_0 = \{ (\underbrace{0,\ldots , 0}_{\alpha} ;
\underbrace{0 , \ldots , 0, u_1 , \ldots , u_r }_{\beta} ) \; : \; u_i \in {\Bbb R}  \mbox{ pour tout } i= 1, \ldots ,r \} \subset i\widetilde{\mathfrak{t}}_0 .$$
Repr\'esentons les syst\`emes de racines respectifs de $\mathfrak{k}$ et $\mathfrak{g}$ comme
$\Delta_c = \{ \pm (x_i \pm x_j ) \; : \; 1 \leq i< j \leq \alpha \ \} \cup
\{ \pm (y_i \pm y_j ) \; : \; 1 \leq i < j \leq \beta \}$ et $\Delta = \Delta_c \cup \{ \pm (x_i \pm y_j ) \; : \; 1 \leq i \leq \alpha \; {\rm et } \; 1 \leq j \leq \beta \}$.
On peut alors consid\'erer les syst\`emes positifs compatibles $\Delta_c^+ = \{ x_i \pm x_j  \; : \; 1 \leq i< j \leq \alpha \} \cup \{ y_j \pm y_i \; : \;  1\leq  i < j \leq \beta \}$,
$\Sigma_c^+ = \{ y_j \pm y_i \; : \; 1 \leq i< j \mbox{ et } \beta -r+1 \leq j \leq \beta \}$ et
$\Sigma^+ = \Sigma_c^+ \cup \{ y_j \pm x_i \; : \; 1 \leq i \leq \alpha \mbox{ et } \beta -r+1 \leq j \leq \beta \}$.
Un calcul simple montre alors que
$$2 \rho_n = p \sum_{j=0}^{r-1} y_{\beta-j}$$
{\it i.e.}, le plus haut poids par rapport \`a $\Delta_c^+$ de la repr\'esentation $V((r^p))$. Ce qui conclut la d\'emonstration du Corrolaire \ref{cohoml2U}.~$\Box$

\bigskip

\noindent
{\bf Remarque.} Lorsque $r=q$ la partition orthogonale $(r^p) \subset p \times (q+r)$ est paire, il ya donc
lieu de consid\'erer les projections respectives $[F]_{\pm}$ de la classe $[F]$ dans les groupes
$H_2^{d_G -d_H} (M)_{(r^p)}^{\pm}$. La d\'emonstration du Corollaire \ref{cohoml2O} montre que ces deux projections sont
non triviales. Ceci correspond au fait que l'on peut remplacer le sous-syst\`eme positif $\Sigma^+$ consid\'erer dans la
d\'emonstration par celui o\`u les \'el\'ements $y_1 \pm x_i$ ($1 \leq i \leq \alpha$) sont rempla\c{c}\'es par leurs oppos\'es.

\medskip

De mani\`ere analogue \`a la Conjecture \ref{conjl2} nous conjecturons~:

\begin{conj} \label{conjl2O}
Si $\lambda$ est une partition orthogonale incluse dans $p\times q$, alors l'application
$$H^{\lambda} (F) \rightarrow H^{|\lambda| +pr}_{2, \; {\rm prim} +} (M)$$
obtenue en composant l'application ``cup-produit avec $[F]$'' et la projection sur la composante fortement primitive de la cohomologie $L^2$ de $M$ 
est {\bf injective} si et seulement si la partition $(r^p)$ s'inscrit dans le diagramme gauche $\hat{\lambda} /\lambda$.
Son image est alors contenue dans $H_2^{\lambda +(r^p)} (M)$.
\end{conj}

En particulier, pour tout entier $k\leq (q-   r)/2$, l'application ``cup-produit avec $[F]$'' 
\begin{eqnarray} \label{c}
H^{k} (F ) \rightarrow H^{k+ pr}_2 (M)
\end{eqnarray}
devrait \^etre {\bf injective}. Nous reviendrons sur cette application au \S 6.

\section{Isolation des repr\'esentations cohomologiques}

\subsection{Isolation dans le dual unitaire}

Dans ce paragraphe nous explicitons pour les groupes unitaires et orthogonaux un th\'eor\`eme de Vogan \cite[Theorem A.10]{Vogan} caract\'erisant les repr\'esentations cohomologiques {\bf isol\'ees}.
Le cas du groupe $U(p,q)$ est d\'ej\`a trait\'e dans \cite[\S 5.4]{BergeronClozel}.
Dans les termes de ce texte, on obtient la Proposition suivante.

\begin{prop} \label{Uisol}
Soit $(\lambda , \mu)$ un couple compatible de partitions dans $p\times q$ avec $\mu/\lambda = (p_1 \times q_1 )* \ldots * (p_m \times q_m)$.
Alors la repr\'esentation $A(\lambda , \mu )$ est isol\'ee dans le dual unitaire de $SU(p,q)$ si
et seulement si
\begin{enumerate}
\item $\min_i (p_i , q_i ) \geq 2$, et
\item si $\lambda_i = \mu_i > \lambda_{i+1}$ ($i=1, \ldots , p$) alors $\mu_{i+1} = \mu_i$ (o\`u nous adoptons exceptionnellement ici
la convention que $\lambda_{p+1} = \mu_{p+1} =-1$).
\end{enumerate}
\end{prop}

On peut visualiser le point 2.~: il signifie que $\lambda$ et $\mu$ ne sont pas tous les deux triviaux et n'ont aucun angle
$\begin{array}{c|}
\\ \hline
\end{array}$ en commun. En particulier, les repr\'esentations $A(\lambda , \mu)$ telles que $\lambda = \mu$ (qui sont exactement les repr\'esentations cohomologique de la s\'erie discr\`ete)
ne sont {\bf jamais isol\'ees}.

Dans cette section nous prouvons la Proposition analogue suivante pour les groupes orthogonaux $SO_0 (p,q)$.

\begin{prop} \label{Oisol}
Soit $\lambda \subset p\times q$ une partition orthogonale avec $\hat{\lambda} / \lambda =
(a_1 \times b_1 ) * \ldots * (a_m \times b_m ) * (p_0 \times q_0 ) * (a_m \times b_m ) * \ldots *(a_1 \times b_1)$. Alors
les diff\'erentes repr\'esentations $A(\lambda )_{\pm_1}^{\pm_2}$ associ\'ees \`a $\lambda$ sont
soit toutes isol\'ees dans le dual unitaire de $SO_0(p,q)$ soit toutes non isol\'ees. Elles sont effectivement isol\'ees si et seulement si
\begin{enumerate}
\item $\min_i (a_i, b_i ) \geq 2$, et
\item soit $p_0 , q_0 \geq 2$ et $p_0 +q_0 \geq 5$, soit $p_0 q_0 =0$, et
\item $\lambda$ et $\hat{\lambda}$ ne sont pas tous les deux triviaux et n'ont aucun angle $\begin{array}{c|}
\\ \hline
\end{array}$ en commun.
\end{enumerate}
\end{prop}
{\it D\'emonstration.} Nous allons d\'eduire la Proposition \ref{Oisol} des r\'esultats de Vogan. D'apr\`es le Th\'eor\`eme A.10 de \cite{Vogan}, une repr\'esentation cohomologique
$A_{\mathfrak{q}}$ est isol\'ee si et seulement si $\mathfrak{q}$ v\'erifie certaines conditions $(0)-(3)$. Vogan montre que l'on peut toujours, quitte \`a changer de repr\'esentant 
pour la classe d'\'equivalence de $A_{\mathfrak{q}}$, choisir $\mathfrak{q}$ de fa\c{c}on \`a ce que la condition $(0)$ soit satisfaite. Lorsque $\lambda \subset p\times q$ est 
une partition orthogonale et $\pm_1$, $\pm_2$ deux signes \'eventuels, notre choix ``canonique'' $\mathfrak{q}(\lambda )_{\pm_1}^{\pm_2}$ v\'erifie toujours la condition $(0)$
($\mathfrak{l} (\lambda)_{\pm_1}^{\pm_2}$ n'a pas de facteur compact non ab\'elien). La condition $(1)$ (le centre du groupe $L (\lambda)$ est compact) est automatiquement v\'erifi\'ee si
$p$ ou $q$ est pair, elle d\'ecoule en toute g\'en\'eralit\'e du fait que $(p_0 , q_0 ) \neq (1,1)$ (d'apr\`es le point 3. de la Proposition). La condition $(2)$ (le groupe $L(\lambda)$
n'a pas de facteurs simples locallement isomorphe \`a $SO(n,1)$ ($n\geq 2$), ou \`a $SU(n,1)$ ($n\geq 1$)) correspond exactement aux points 1. et 2. de la Proposition. Il
nous reste \`a exprimer la condition $(3)$.

La construction-classification de Vogan-Zuckerman \cite{VoganZuckerman} pour les $A_{\mathfrak{q}}$ est la suivante. Tout d'abord on a $A_{\mathfrak{q}}= R_{\mathfrak{q}} ({\Bbb C})$ o\`u ${\Bbb C}$ est la repr\'esentation triviale de $\mathfrak{l}$ et le foncteur $R_{\mathfrak{q}} =R_{\mathfrak{q}}^{{\rm dim}(\mathfrak{u} \cap \mathfrak{k})}$ est d\'efini
dans les r\'ef\'erences cit\'ees par \cite{VoganZuckerman}. Soit $\mathfrak{h}$ une sous-alg\`ebre de Cartan $\theta$-stable de $\mathfrak{l}$ contenant $\mathfrak{t}$ \footnote{La
sous-alg\`ebre $\mathfrak{h}$ est donc \'egale \`a $\mathfrak{t}$ sauf si $p$ et $q$ sont tous les deux impairs respectivement \'egaux \`a $2r+1$ et $2s+1$,
auquel cas on peut ajouter \`a $\mathfrak{t}$ l'alg\`ebre $\left\{
\left(
\begin{array}{cc}
0 & z \\
z & 0 
\end{array} \right) \; : \; z \in {\Bbb C} \right\}$ dans le $r+1$-\`eme bloc diagonal.}. Supposons donn\'e un syst\`eme de racines positives $\Delta^+ (\mathfrak{g} , \mathfrak{h})$ pour $(\mathfrak{g} ,
\mathfrak{h})$ tel que les racines de $\mathfrak{u}$ soient positives. Alors, avec les notations usuelles~:
$$\rho_{\mathfrak{g}} = \rho_{\mathfrak{u}} + \rho_{\mathfrak{l}}.$$
Le caract\`ere infinit\'esimal de ${\Bbb C}_{\mathfrak{l}}$ est $\rho_{\mathfrak{l}}$; $\lambda = \rho_{\mathfrak{g}}$ v\'erifie les hypoth\`eses du Th\'eor\`eme A.10 de \cite{Vogan}
pour $\Delta^+ (\mathfrak{g} , \mathfrak{h})$.

Soit $\Pi \subset \Delta^+ (\mathfrak{g} , \mathfrak{h})$ l'ensemble des racines simples et $\Pi (\mathfrak{l})$ le sous-ensemble form\'e des racines simples de $\mathfrak{l}$. 
Alors la condition $(3)$ de Vogan s'\'ecrit~:
\begin{eqnarray} \label{cond3}
\langle \beta^{\vee} , \lambda \rangle = \langle \beta^{\vee} , \rho_{\mathfrak{g}} \rangle \neq 1
\end{eqnarray}
{\it pour toute racine imaginaire non compacte $\beta \in \Pi$ orthogonale \`a $\Pi (\mathfrak{l})$.}

Il s'agit d'expliciter (\ref{cond3}). Nous allons encore une fois distinguer trois cas suivant les parit\'es des entiers $p$ et $q$.

\paragraph{p=2r et q=2s} Dans ce cas $\mathfrak{h} = \mathfrak{t}$.
Rappelons (\S 1.2) que nous consid\'erons une sous-alg\`ebre parabolique $\mathfrak{q}(\lambda )_{\pm_1}^{\pm_2}$
associ\'ee \`a un \'el\'ement $X=(x_1 , \ldots , x_r ; y_1 , \ldots , y_s ) \in i \mathfrak{t}_0$ avec 
$$x_1 \geq \ldots \geq x_{r-1} \geq |x_r | \geq 0 \mbox{ et } y_s \geq \ldots \geq y_2 \geq |y_1 | \geq 0.$$

Nous rassemblons les valeurs $x_1 , \ldots , x_{r-1} , |x_r | , y_s , \ldots , y_2 , |y_1 |$ en une suite strictement d\'ecroissante
$$u_1 > u_2 > \ldots > u_l \geq 0$$
et notons 
$$u_{j_1} > u_{j_2} > \ldots > u_{j_m}$$
les valeurs non nulles de multiplicit\'e $>1$. Alors,
\begin{eqnarray} \label{x}
(x_1 , \ldots , x_{r-1}, |x_r| ) = ( x_1 > x_2 > \ldots > x_{\alpha_1} > \underbrace{u_{j_1}}_{a_1} > \ldots > x_{\alpha_2} > \underbrace{u_{j_2}}_{a_2} > \ldots )
\end{eqnarray}
et
\begin{eqnarray} \label{y}
(y_s , \ldots , y_2, |y_1| ) = ( y_s > y_{s-1} > \ldots > y_{\beta_1} > \underbrace{u_{j_1}}_{b_1} > \ldots > y_{\alpha_2} > \underbrace{u_{j_2}}_{b_2} > \ldots ) .
 \end{eqnarray}
Si $u_l =0$ nous notons $r_0$ (resp. $s_0$) la multiplicit\'e avec laquelle $0$ intervient dans $(x_1 , \ldots , x_{r-1}, |x_r| )$ (resp. $(y_s , \ldots , y_2, |y_1| )$).

Rappelons (cf. \S 1.2) que nous avons suppos\'e $\mathfrak{l}(\lambda )_{\pm_1}^{\pm_2}$ sans facteur compact non ab\'elien. Alors pour $i=1, \ldots ,m$,
$a_i$ et $b_i$ sont tous deux $>0$ et si $r_0$ (resp. $s_0$) est nul alors $s_0$ (resp. $r_0$) est $\leq 1$. 

Avec ces notations, le diagramme gauche 
$$\hat{\lambda} / \lambda = (a_1 \times b_1 ) * \ldots * (a_m \times b_m ) * (p_0 \times q_0 ) * (a_m \times b_m ) * \ldots *(a_1 \times b_1) ,$$
o\`u $p_0 = 2r_0$, $q_0 = 2s_0$ et le diagramme rectangulaire $p_0 \times q_0$ peut \^etre trivial (si $x_r$ ou $y_1$ est non nul).

Il s'agit maintenant de montrer que la condition (\ref{cond3}) est \'equivalente au point 3. de la Proposition \ref{Oisol}.

Soit $W=W_G$ le groupe de Weyl de $G$ (associ\'e au syst\`eme $\Delta (\mathfrak{g} , \mathfrak{h})$). Le groupe $W$ est isomorphe au groupe 
$\Sigma_{r+s} \ltimes \{ -1 \}^{r+s-1}$, $\{ -1 \}^{r+s-1}$ \'etant le sous-groupe de $\{ -1 \}^{r+s}$, op\'erant diagonalement, d\'efini par $\prod s_i = 1$.

Il existe un \'el\'ement $w \in W$ tel que
$$w X = (v_1 , \ldots , v_{r+s} ),$$
o\`u 
$$(v_1 , \ldots v_{r+s-1} ,|v_{r+s}| ) = (u_1 > u_2 > \ldots > \underbrace{u_{j_1}}_{a_1 + b_1} > \ldots > \underbrace{u_{j_2}}_{a_2 + b_2} > \ldots > u_k ) $$
et $u_k$ appara\^{\i}t avec multiplicit\'e $r_0 +s_0$ s'il est nul. Nous notons $X' = wX$.

La base $\Delta^+ (\mathfrak{g} , \mathfrak{h} )$ est une ensemble de racines positives $\alpha$ telles que $\langle \alpha , X \rangle \geq 0$; $w$ l'envoie sur un ensemble
de racines telles que $\langle \alpha , X' \rangle \geq 0$, que l'on prendra \'egal \`a l'ensemble usuel puisque $X'$ est dominant. Alors
\begin{eqnarray*}
w \Pi & = &  \{ \epsilon_i^{\pm} = (0 , \ldots ,0 , 1 , \pm 1 , 0 , \ldots ,0) , \; i= 1 , \ldots , r+s-1 \} \\
w \Pi & = & \{ \epsilon_i^{\pm} \; : \; X_i ' = \mp X_{i+1} ' \} , 
\end{eqnarray*}
et $w$ envoie l'orthogonal de $\Pi (\mathfrak{l} )$ sur l'ensemble des racines $\epsilon_i$ telles que
\begin{eqnarray} \label{ii}
\{ i , i+1 \} \subset \{ 1 , \ldots , r+s \} - \cup_{i=0}^m I_i ,
\end{eqnarray}
$I_i$ ($i\geq 1$) \'etant le support de $u_{j_i}$ dans l'expression de $X'$ et $I_0$ \'etant le support de $0$ dans $X'$. 

Par ailleurs, le syst\`eme de racines \'etant de type $D_{r+s}$, $\beta^{\vee} = \beta$ si $\beta \in \Pi$ et donc $\langle \beta^{\vee} , \rho_{\mathfrak{g}} \rangle = 1$. La condition
(\ref{cond3}) est donc \'equivalente \`a
\begin{eqnarray} \label{cond3'}
\mbox{Il n'y a pas de racine imaginaire non compacte dans }\Pi \mbox{ orthogonale \`a } \Pi (\mathfrak{l}).
\end{eqnarray}

Soit donc $\{ i, i+1 \}$ v\'erifiant (\ref{ii}). Ceci implique donc que $i$ et $i+1$ appartiennent \`a une composante connexe de cardinal $\geq 2$ de $\{ 1 , \ldots , r+s \}
- \cup I_i $. Supposons par exemple que $\alpha_1$ ou $\beta_1 >1$ et que $a_1 + b_1 \geq 2$. Alors, dans l'expression de $X'$, $i$ et $i+1$ sont deux indices associ\'es 
\`a $u_i > u_{i+1} > u_{j_1}$. Alors $w^{-1} \{ i, i+1 \}$ est associ\'ee dans les expressions (\ref{x}) et (\ref{y}) \`a deux couples de la forme $(x_i ,| x_{i+1}|)$, $(x_j ,| y_{j'}|)$ ou $(y_j , |y_{j+1} | )$ \`a gauche de $u_{j_1}$. Dans le second cas, la racine de valeur $x_j - | y_{j'}|$ est non compacte ce qui contredit (\ref{cond3'}). 

Donc $u_1$ et $u_2$ sont tous deux (par exemple) de la forme $(x_1 , x_2 )$; il en est de m\^eme pour $x_2$ et $x_3$, etc. Ceci veut dire que pour tout $j<j_1$ on a $u_j = x_j$
et donc $\alpha _1 = j_1 -1$ et $\beta_1 = 0$. Le m\^eme argument s'applique \`a toute composante connexe, et il est clair que (\ref{cond3'}) (et donc (\ref{cond3})) est en fait \'equivalente au point 3. de la Proposition \ref{Oisol}.

\paragraph{p=2r et q=2s+1} Cette fois encore $\mathfrak{h} = \mathfrak{t}$.
Nous consid\'erons maintenant une sous-alg\`ebre parabolique $\mathfrak{q}(\lambda )_{\pm}$ 
associ\'ee \`a un \'el\'ement $X=(x_1 , \ldots , x_r ; y_1 , \ldots , y_s ) \in i \mathfrak{t}_0$ avec 
$$x_1 \geq \ldots \geq x_{r-1} \geq |x_r | \geq 0 \mbox{ et } y_s \geq \ldots \geq y_2 \geq y_1  \geq 0.$$

Nous rassemblons l\`a encore les valeurs $x_1 , \ldots , x_{r-1} , |x_r | , y_s , \ldots , y_2 , y_1 $ en une suite strictement d\'ecroissante
$$u_1 > u_2 > \ldots > u_l \geq 0$$
et notons 
$$u_{j_1} > u_{j_2} > \ldots > u_{j_m}$$
les valeurs non nulles de multiplicit\'e $>1$. Les identit\'es (\ref{x}) et (\ref{y}) sont encore v\'erifi\'ees, nous conservons les m\^emes notations $r_0, s_0 , \ldots$.

Il s'agit maintenant de montrer que la condition (\ref{cond3}) est \'equivalente au point 3. de la Proposition \ref{Oisol}.

Soit $W=W_G$ le groupe de Weyl de $G$ (associ\'e au syst\`eme $\Delta (\mathfrak{g} , \mathfrak{h})$). Le groupe $W$ est maintenant isomorphe au groupe
$\Sigma_{r+s} \ltimes \{ -1 \}^{r+s}$.

Il existe un \'el\'ement $w \in W$ tel que  
$$w X = (v_1 , \ldots , v_{r+s} ),$$
o\`u 
$$(v_1 , \ldots v_{r+s-1} ,v_{r+s} ) = (u_1 > u_2 > \ldots > \underbrace{u_{j_1}}_{a_1 + b_1} > \ldots > \underbrace{u_{j_2}}_{a_2 + b_2} > \ldots > u_k ) $$
et $u_k$ appara\^{\i}t avec multiplicit\'e $r_0 +s_0$ s'il est nul. Nous notons $X' = wX$.

La base $\Delta^+ (\mathfrak{g} , \mathfrak{h} )$ est une ensemble de racines positives $\alpha$ telles que $\langle \alpha , X \rangle \geq 0$; $w$ l'envoie sur un ensemble 
de racines telles que $\langle \alpha , X' \rangle \geq 0$, que l'on prendra \'egal \`a l'ensemble usuel puisque $X'$ est dominant. Alors
\begin{eqnarray*}
w \Pi & = &  \{ \epsilon_i = (0 , \ldots , 0 , 1 , - 1 , 0 , \ldots  ,0) , \; i= 1 , \ldots , r+s-1 \} \\ 
& & \cup \{ \eta_i = (0 , \ldots , 0, 1 , 0 , \ldots , 0) , \; i= 1, \ldots ,r+s \} , \\  
\end{eqnarray*}
et $w$ envoie l'orthogonal de $\Pi (\mathfrak{l} )$ sur l'ensemble des racines $\epsilon_i$ et $\eta_j$ telles que 
\begin{eqnarray} \label{ii2}
\{ i , i+1 \} \subset \{ 1 , \ldots , r+s \} - \cup_{i=0}^m I_i ,
\end{eqnarray}
et
\begin{eqnarray} \label{ii22}
j \notin \cup_{i=0}^m I_i .
\end{eqnarray}

Par ailleurs, le syst\`eme de racines \'etant de type $B_{r+s}$, $\epsilon_i^{\vee} = \epsilon_i$ et $\eta_i^{\vee} = 2 \eta_i$. On a donc
$\langle \epsilon_i^{\vee} , \rho_{\mathfrak{g}} \rangle = 1$ pour $i=1, \ldots , r+s-1$, $\langle \eta_{r+s}^{\vee} , \rho_{\mathfrak{g}} \rangle = 1$ et
$\langle \eta_i^{\vee} , \rho_{\mathfrak{g}} \rangle \neq 1$ pour $i=1, \ldots ,r+s-1$. La condition
(\ref{cond3}) est donc \'equivalente \`a 
\begin{eqnarray} \label{cond32}
\begin{array}{l}
\mbox{Aucune des racines } \{ \epsilon_i \; : \; i=1, \ldots , r+s-1 \} \cup \{ \eta_{r+s} \}\\
\mbox{ n'est orthogonale \`a } \Pi (\mathfrak{l}).
\end{array}
\end{eqnarray}

Soit donc $\epsilon_i$ v\'erifiant (\ref{ii2}). Ceci implique donc que $i$ et $i+1$ appartiennent \`a une composante connexe de cardinal $\geq 2$ de $\{ 1 , \ldots , r+s \}
- \cup I_i \}$. Supposons par exemple que $\alpha_1$ ou $\beta_1 >1$ et que $a_1 + b_1 \geq 2$. Alors, dans l'expression de $X'$ $i$ et $i+1$ sont deux indices associ\'es
\`a $u_i > u_{i+1} > u_{j_1}$. Alors $w^{-1} \epsilon_i$ est associ\'ee dans les expressions (\ref{x}) et (\ref{y}) \`a deux couples de la forme $(x_i , x_{i+1})$, $(x_j , y_{j'})$ ou $(y_j , y_{j+1} )$ \`a gauche de $u_{j_1}$. Dans le second cas, la racine de valeur $x_j - y_{j'}$ est non compacte.

Donc $u_1$ et $u_2$ sont tous deux (par exemple) de la forme $(x_1 , x_2 )$; il en est de m\^eme pour $x_2$ et $x_3$, etc. Ceci veut dire que pour tout $j<j_1$ on a $u_j = x_j$
et donc $\alpha _1 = j_1 -1$ et $\beta_1 = 0$. Le m\^eme argument s'applique \`a toute composante connexe, et il est clair que l'on obtient que (\ref{cond32}) pour les 
racines $\epsilon_i$ implique que $\lambda$ et $\hat{\lambda}$ ne sont pas tous les deux triviaux et n'ont aucun angle $\begin{array}{c|}
\\ \hline 
\end{array}$ en commun sauf peut-\^etre celui correspondant \`a la case $(r,s+1)$. Ce dernier cas est \'equivalent \`a ce que $v_{r+s} = x_r > 0$, alors la racine de valeur $x_r$ est non compacte ce qui contredit (\ref{cond32}) pour la racine $\eta_{r+s}$. On obtient bien finalement que (\ref{cond32}) (et donc (\ref{cond3})) est en fait \'equivalente au point 3. de la Proposition \ref{Oisol}.

\paragraph{p=2r+1 et q=2s+1} Ce cas se traite de la m\^eme mani\`ere que les deux pr\'ec\'edents sans aucune difficult\'e suppl\'ementaire.
Ce qui conclut la d\'emonstration de la Proposition \ref{Oisol}.~$\Box$

\bigskip

Remarquons l\`a encore que la condition 3. de la Proposition \ref{Oisol} est viol\'ee lorsque la repr\'esentation cohomologique est une s\'erie discr\`ete~: celles-ci ne sont pas 
isol\'ees.

\medskip 

Rappelons (\S 1.1) que la cohomologie de $A_{\mathfrak{q}}$ n'appara\^{\i}t qu'en degr\'es $\geq R=R(\mathfrak{q} )$.

\begin{cor} \label{mino}
Lorsque $G$ est de rang $1$ aucune repr\'esentation cohomologique n'est isol\'ee. Si $G$ est de rang $2$, il est du type $SO_0 (2,n)$ avec $n\geq 3$ et
si une repr\'esentation cohomologique $A_{\mathfrak{q}}$ n'est pas isol\'ee, la cohomologie de $A_{\mathfrak{q}}$ n'appara\^{\i}t qu'en degr\'es $k \geq \left[
\frac{n}{2} \right]$. Enfin si $G$ est de rang $\geq 3$, il est du type $SO_0 (p,q)$ avec $p,q \geq 3$ et si une repr\'esentation cohomologique $A_{\mathfrak{q}}$ n'est pas isol\'ee, la cohomologie de $A_{\mathfrak{q}}$ n'appara\^{\i}t qu'en degr\'es $k \geq p+q-3$.
\end{cor}
{\it D\'emonstration.} Tout d'abord un calcul simple montre que si $\mathfrak{q} = \mathfrak{q}(\lambda )_{\pm_1}^{\pm_2}$,
$$R= \frac12 \left( pq - 2 \sum_{j=1}^{m} a_j b_j -p_0 q_0 \right) .$$

Pour simplifier nous ne v\'erifions le Corollaire \ref{mino} que dans le cas $p=2r$ et $q=2s$ (les deux autres cas sont plus facile \`a traiter).
Nous distinguons diff\'erents cas.

\begin{itemize}
\item Supposons tout d'abord $p_0 q_0 = 0$. 

Supposons que $\mathfrak{q}$ viole le point 1. de la Proposition \ref{Oisol}. Alors (\`a l'ordre pr\`es) on peut supposer $a_1 = 1$, $b_1 \geq 1$. Alors
$$\sum_{j=2}^m a_j b_j \leq (\sum a_i ) (\sum b_i ) \leq (r-1)(s-b_1 )$$
donc
\begin{eqnarray*}
R & \geq & 2rs -b_1 -(r-1)(s-b_1 ) \\
    & \geq & rs+s+(r-2)b_1 .
\end{eqnarray*}
Nous devons alors distinguer le cas $r=1$ du cas $r\geq 2$. Si $r=1$ et puisque $b_1 \leq s$, 
\begin{eqnarray} \label{regal1}
R \geq s = \left[ \frac{q}{2} \right] .
\end{eqnarray}
Alors que si $r\geq 2$ (et puisque $b_1 \geq 1$) ,
\begin{eqnarray} \label{rgeq2}
R \geq  rs + r + s -2 .
\end{eqnarray}

Supposons maintenant que $\mathfrak{q}$ viole le point 3. de la Proposition \ref{Oisol}. Alors, en particulier, $\sum a_i \leq r-1$ et $\sum b_i \leq s-1$ et dans ce cas
\begin{eqnarray} \label{1.2}
R \geq 2rs - (r-1)(s-1) = rs +r+s -1. 
\end{eqnarray}

\item Supposons maintenant $p_0$ et $q_0 >0$. Nous continuons de noter $p_0 =2r_0$ et $q_0 =2s_0$.

Supposons que $\mathfrak{q}$ viole le point 2. de la Proposition \ref{Oisol}. Alors $p_0 q_0 = 4r_0 s_0 = 4$
et $\sum a_i b_i \leq (r-1)(s-1)$ et dans ce cas
\begin{eqnarray} \label{2.1}
R \geq 2rs - (r-1)(s-1) - 2 = rs +r+s -3. 
\end{eqnarray}

Supposons que $\mathfrak{q}$ viole le point 1. de la Proposition \ref{Oisol}. Alors (\`a l'ordre pr\`es) on peut supposer $a_1 = 1$, $b_1 \geq 1$ . Alors
$$ \sum_{j=1}^m a_j b_j + 2r_0 s_0 \leq 2 (r-1) (s-1) +1$$
donc
\begin{eqnarray} \label{2.2}
R & \geq & 2rs -2(r-1)(s-1)-1 = 2r +2s -3.
\end{eqnarray}

Supposons maintenant que $\mathfrak{q}$ viole le point 3. de la Proposition \ref{Oisol}. Alors, en particulier, $\sum a_i + r_0 \leq r-1$ et $\sum b_i +s_0 \leq s-1$ et dans ce cas
\begin{eqnarray} \label{2.3}
R \geq 2rs - 2(r-1)(s-1) = 2r +2s -2. 
\end{eqnarray}
\end{itemize}

Les diff\'erentes in\'egalit\'es (\ref{regal1}), (\ref{rgeq2}), (\ref{1.2}), (\ref{2.1}), (\ref{2.2}) et (\ref{2.3}) impliquent le Corollaire \ref{mino}.~$\Box$

\bigskip

\noindent
{\bf Remarque.} Du point de vue de la th\'eorie des repr\'esentations de $G$, ce r\'esultat ne peut \^etre am\'elior\'e. Ainsi, lorsque $p=2$ la repr\'esentation $A((\left[ \frac{q}{2} \right] ))$
est cohomologique de degr\'e $R=\left[ \frac{q}{2} \right]$ et n'est pas isol\'ee d'apr\`es la Proposition \ref{Oisol}. Et, lorsque $p,q \geq 3$, la repr\'esentation 
$A((q-1, 1^{p-2}))$ est cohomologique de degr\'e $R=p+q-3$ et n'est pas isol\'ee d'apr\`es la Proposition \ref{Oisol}.

\subsection{Isolation sous la condition $d=0$}

Soient $G$ un groupe simple r\'eel, $K$ un sous-groupe compact maximal de $G$ et $(\pi , V_{\pi})$ une repr\'esentation irr\'eductible unitaire de $G$.
Soit $V_{\pi}^K$ l'espace des vecteurs $K$-finis de la repr\'esentation $\pi$ et consid\'erons le complexe
calculant la $(\mathfrak{g} , K)$-cohomologie de $\pi$~:
\begin{eqnarray} \label{cpx}
\ldots \rightarrow C^i (\pi ) = {\rm Hom}_K (\bigwedge {}^i \mathfrak{p} , V_{\pi}^K) \stackrel{d_i}{\rightarrow} C^{i+1} (\pi ) \rightarrow \ldots .
\end{eqnarray}

Nous dirons d'une repr\'esentation cohomologique $A_{\mathfrak{q}}$ de $G$ de
degr\'e primitif $R=R(\mathfrak{q})$ qu'elle est {\it isol\'ee sous la condition $d=0$} si elle est isol\'ee de l'ensemble des repr\'esentations 
irr\'eductibles unitaires telles que Im$(d_{R} )=0$.

Soit ${\cal P}$ l'ensemble de toutes les sous-alg\`ebres paraboliques $\theta$-stables de $\mathfrak{g}$~: $\mathfrak{q}= \mathfrak{l}(\mathfrak{q}) + \mathfrak{u}(\mathfrak{q})$.
Notons
\begin{eqnarray} \label{ccccc}
r_G = \min \{ \dim (\mathfrak{u} (\mathfrak{q}) \cap \mathfrak{p}) \; : \; \mathfrak{q} \in {\cal P} \} .
\end{eqnarray}
D'apr\`es Parthasarathy \cite{Parthasarathy}, Kumaresan \cite{Kumaresan} et Vogan-Zuckerman \cite{VoganZuckerman}, toute repr\'esentation cohomologique non triviale de $G$ est 
de degr\'e primitif $\geq r_G$. 

Remarquons que $r_G \geq {\rm rang}_{{\Bbb R}} (G)$ et que les cas o\`u l'in\'egalit\'e est stricte sont rassembl\'es dans \cite[\S 10.3]{BorelWallach}. En ce qui concerne les 
groupes $O(p,q)$ et $U(p,q)$, $r_G= {\rm rang}_{{\Bbb R}} (G) = \min (p,q)$.

Nous aurons besoin de la proposition suivante que nous d\'eduisons facilement des travaux de Parthasarathy, Kumaresan et Vogan-Zuckerman.

\begin{prop} \label{isold=0}
Si $G$ n'est localement isomorphe ni \`a $SO(n,1)$ ($n\geq 2$), ni \`a $SU(n,1)$ ($n\geq1$), alors toute repr\'esentation cohomologique $A_{\mathfrak{q}}$ de degr\'e primitif $R=r_G$ est isol\'ee sous la condition $d=0$ dans $\widehat{G}$.
\end{prop}
{\it D\'emonstration.} Soit $\mathfrak{t} \subset \mathfrak{k}$ une sous-alg\`ebre de Cartan comme dans la premi\`ere section. Soit $\mathfrak{q} = \mathfrak{l} + \mathfrak{u} \in {\cal P}$ telle que $A_{\mathfrak{q}}$ soit de degr\'e primitif $R=r_G$. Fixons un sous-syst\`eme positif $\Delta^+ (\mathfrak{l} \cap \mathfrak{k})$ du syst\`eme de racines $\Delta (\mathfrak{l}
\cap \mathfrak{k} , \mathfrak{t} )$. L'ensemble 
$$\Delta^+ (\mathfrak{k}) = \Delta^+ (\mathfrak{l} \cap \mathfrak{k}) \cup \Delta (\mathfrak{u} \cap \mathfrak{k}) $$
est alors un syst\`eme de racines positives de $\mathfrak{t}$ dans $\mathfrak{k}$.

Rappelons que ($\Delta^+ (\mathfrak{k})$ \'etant fix\'e) les repr\'esentations irr\'eductibles de groupe compact $K$ peuvent \^etre param\'etr\'ees par leurs plus haut poids, qui 
sont des \'el\'ements de $\mathfrak{t}^*$. Notons
$$\mu (\mathfrak{q}) = \mbox{ repr\'esentation de } K \mbox{ de plus haut poids } 2 \rho (\mathfrak{u}\cap \mathfrak{p}) .$$
C'est le $K$-type minimal de la repr\'esentation $A_{\mathfrak{q}}$.

Rappelons maintenant la c\'el\`ebre in\'egalit\'e de Dirac de Parthasarathy (cf. \cite[(2.26)]{Parthasarathy}, \cite[II.6.11]{BorelWallach}, \cite[Lemma 4.2]{VoganZuckerman}).

\begin{lem} \label{dirac}
Soit $(\pi , V_{\pi})$ une repr\'esentation unitaire irr\'eductible de $G$ et $V_{\pi}^K$ son $(\mathfrak{g} , K)$-module associ\'e.
Fixons une repr\'esentation de $\mathfrak{k}$ de plus haut poids $\chi \in \mathfrak{t}^*$ et apparaissant dans $V_{\pi}^K$; 
et un sous-syst\`eme positif de racines $\Delta^+ (\mathfrak{g})$
de $\mathfrak{t}$ dans $\mathfrak{g}$. Notons $\rho$ (resp. $\rho_c$, $\rho_n$) dans $\mathfrak{t}^*$ la demi-somme des racines dans $\Delta^+ (\mathfrak{g})$ (resp. $\Delta^+ (\mathfrak{k})$, $\Delta^+ (\mathfrak{p})$), de sorte que $\rho = \rho_c + \rho_n$. Soit $w$ un \'el\'ement du groupe de Weyl $W_K = W( \mathfrak{k} , \mathfrak{t})$ de $\mathfrak{t}$
dans $\mathfrak{k}$, tel que $w(\chi - \rho_n)$ soit dominant pour $\Delta^+ (\mathfrak{k})$. Alors,
$$-\pi (\Omega ) \geq || \rho ||^2 - || w(\chi - \rho_n ) + \rho_c ||^2,$$
o\`u $\Omega$ d\'esigne le casimir de $\mathfrak{g}$ et la norme $||.||$ est d\'eduite de la forme de Killing sur $\mathfrak{g}$.
\end{lem}

Pour toute repr\'esentation unitaire 
irr\'eductible de $G$ non triviale $(\pi , V_{\pi})$, le groupe $H^{r_G -1} (\mathfrak{g} , K , V_{\pi}^K ) = \{ 0 \}$. Il en est de m\^eme si $\pi$ est triviale puisque, $G$ n'\'etant localement isomorphe ni \`a $SO(n,1)$ ($n\geq 2$) ni \`a $SU(n,1)$ ($n\geq1$), $0< r_G -1$ et donc, Hom$_K (\bigwedge^{r_G -1} \mathfrak{p} , {\Bbb C} ) = \{0 \}$.
La suite (\ref{cpx}) est donc exacte en $i=r_G -1$. Ceci reste vrai en $i= r_G$ en dehors des repr\'esentations $\pi$ telles que $H^{r_G} (\mathfrak{g} , K , V_{\pi}^K )=0$. Ces
diff\'erentes repr\'esentations ont des $K$-types diff\'erents et sont donc isol\'ees les unes des autres.

Si $A_{\mathfrak{q}}$ n'est pas isol\'ee sous la condition $d=0$, il existe donc une suite $\{ \pi_i \}$ de repr\'esentations unitaires irr\'eductibles de $G$ telle que
\begin{enumerate}
\item ${\rm Hom}_K (\bigwedge^{R-1} \mathfrak{p} , V_{\pi_i}^K ) \neq \{ 0 \} $ et,
\item la suite $\pi_i (\Omega )$ tende vers $0$ lorsque $i$ tend vers l'infini.
\end{enumerate}

Mais, l'ensemble des $K$-types de $\bigwedge^{R-1} \mathfrak{p}$ est fini et d'apr\`es Kumaresan \cite{Kumaresan}, si $\chi \in \mathfrak{t}^*$ est le plus haut poids d'un tel $K$-type,
$$||\rho ||^2 -  || w(\chi - \rho_n ) + \rho_c ||^2 >0 , $$
o\`u $w$ est comme dans le Lemme \ref{dirac}. (Kumaresan montre plus pr\'ecisemment qu'un $K$-type de $\bigwedge^* \mathfrak{p}$ v\'erifie l'\'egalit\'e 
$$||\rho ||^2 -  || w(\chi - \rho_n ) + \rho_c ||^2 =0 $$
est n\'ecessairement de la forme $\mu (\mathfrak{q})$ pour une certaine sous-alg\`ebre parabolique $\mathfrak{q} \subset \mathfrak{g}$.)
D'apr\`es le Lemme \ref{dirac}, il existe donc une borne inf\'erieure strictement positive uniforme de l'ensemble des nombres r\'eels  
$-\pi (\Omega )$ tels que $\pi$ soit une repr\'esentation unitaire irr\'eductible de $G$ v\'erifiant Hom$_K (\bigwedge^{R-1} \mathfrak{p} , V_{\pi}^K ) \neq \{ 0\}$. Ce qui
contredit l'existence de la suite $\{ \pi_i \}$ et conclut la d\'emonstration de la Proposition \ref{isold=0}.~$\Box$

\bigskip

\begin{cor} \label{isolq}
La repr\'esentation cohomologique $A((1^p))$ du groupe $SO_0 (p,q)$ ($2\leq p < q$) est isol\'ee sous la condition $d=0$.
\end{cor}
{\it D\'emonstration.} Puisque si $G= SO_0 (p,q)$ ($1\leq p < q$), le degr\'e $r_G = p$, le Corollaire \ref{isolq} d\'ecoule 
trivialement de la Proposition \ref{isold=0}, tant que $p>1$.~$\Box$

\bigskip

\subsection{Isolation dans le dual automorphe}

Dans deux articles fondamentaux, Arthur a donn\'e une description conjecturale
des repr\'esentations des groupes r\'eductifs qui peuvent appara\^{\i}tre dans $L^2 (\Gamma \backslash
G)$ pour un sous-groupe de congruence. Avec Clozel dans \cite{BergeronClozel}, nous avons d\'eduit de la th\'eorie d'Arthur {\bf a minima} des limitations s\'ev\`eres sur les
caract\`eres infinit\'esimaux des repr\'esentations pouvant appara\^{\i}tre dans $L^2 (\Gamma \backslash G)$ lorsque $G^{{\rm nc}} = U(n,1)$ ou $O(n,1)$.

\bigskip

Soit $G$ un groupe alg\'ebrique semisimple et connexe sur ${\Bbb Q}$. Comme dans \cite{BurgerLiSarnak} nous notons
$$\sigma (\Gamma \backslash G) = \{ \pi \in \widehat{G({\Bbb R})} \; : \; \pi \propto L^2 (\Gamma \backslash G({\Bbb R})) \}$$
le {\it spectre de $L^2 (\Gamma \backslash G)$}, o\`u $\Gamma$ est un sous-groupe de congruence de $G$ et $\propto$ signifie ``\^etre
faiblement contenue''. Nous notons de plus
$$\widehat{G}_{{\rm Aut}} = \overline{ \bigcup_{\Gamma} \sigma (\Gamma \backslash G)}$$
le {\it dual automorphe de $G$}, o\`u la r\'eunion est prise sur l'ensemble des sous-groupes de congruence de $G$ et l'adh\'erence est prise dans le dual unitaire
$\widehat{G({\Bbb R})}$ de $G({\Bbb R})$.

Dans \cite{BergeronClozel} avec Clozel nous avons montr\'e comment d\'eduire la conjecture suivante d'une Conjecture 
de changement de base pour les groupes unitaires, cas faible des Conjectures d'Arthur.

\begin{conj} \label{isolautomU}
Supposons $G^{\rm nc} = U(p,q)$. Soit $\pi$ une repr\'esentation cohomologique de $G({\Bbb R})$, alors la repr\'esentation $\pi$ est isol\'ee dans 
$$\{ \pi \} \cup \widehat{G}_{\rm Aut}.$$
\end{conj}

\bigskip

Le cas du groupe $O(p,q)$ est plus d\'elicat. L'analogue de la Conjecture \ref{isolautomU} est d'ailleurs fausse 
en g\'en\'eral pour le groupe $O(p,q)$.

La Conjecture suivante pourrait surement \^etre d\'eduite des conjectures g\'en\'erales d'Arthur (lorsque $p=1$ nous montrons dans \cite{BergeronClozel} qu'elle
d\'ecoule d'une version faible des Conjectures d'Arthur).

\begin{conj} \label{isolautomO}
Supposons $G^{\rm nc} = O(p,q)$. Chaque repr\'esentation cohomologique $A((1^i))$, pour $0 \leq i \leq q/2 -1$ est isol\'ee dans le dual automorphe
$\widehat{G}_{\rm Aut}$ de $G$.
\end{conj}

On ne connait que tr\`es peu de cas des Conjectures \ref{isolautomU} et \ref{isolautomO}. Dans \cite{Clozel} Clozel d\'emontre que 
la repr\'esentation trivial est toujours isol\'ee dans le dual automorphe. Dans \cite{BergeronClozel} avec Clozel, nous d\'emontrons
la Conjecture \ref{isolautomU} pour $p=1$ et $\pi$ de degr\'e fortement primitif \'egal \`a $1$. Dans les deux cas les d\'emonstrations
utilisent un argument de r\'eduction d\^u \`a Burger et Sarnak \cite{BurgerSarnak}.

\bigskip

De mani\`ere similaire au \S 3.2 on peut \'evidemment \'etudier l'isolation d'une repr\'esentation cohomologiques sous la condition 
$d=0$. La d\'emonstration de \cite[Th\'eor\`eme 2.3.3]{BergeronClozel} et les r\'esultats mentionn\'es dans le paragraphe pr\'ec\'edent
impliquent alors les deux analogues automorphes suivant du Corollaire \ref{isolq}.

\begin{prop} \label{isolq1}
La repr\'esentation cohomologique $A((1))$ du groupe $SO_0 (1,q)$ ($2 \leq q$) est isol\'ee dans
le dual automorphe sous la condition $d=0$.
\end{prop}

\begin{prop} \label{isolq2}
Les repr\'esentations cohomologiques de degr\'e $2$ du groupe $SU(1,q)$ ($1 \leq q$) sont toutes isol\'ees dans
le dual automorphe sous la condition $d=0$.
\end{prop}

\section{Restriction des repr\'esentations cohomologiques}

Dans cette section nous \'etudions la possibilt\'e pour la restriction d'une repr\'esentation cohomologique de $G$ \`a un sous-groupe $H$ 
de contenir (discr\`etement) une repr\'esentation cohomologique. Ce probl\`eme a \'et\'e \'etudi\'e par diff\'erents auteurs citons notamment les
travaux de Kobayashi \cite{Kobayashi}, Harris et Li \cite{HarrisLi} et mon article \cite{IRMN}. Nous nous pla\c{c}ons tout d'abord dans une cadre
g\'en\'eral.

\subsection{Restriction et s\'eries discr\`etes}

Soit $G$ un groupe de Lie (r\'eel) r\'eductif connexe \`a centre compact et avec un sous-groupe de Cartan compact. 
Le groupe $G$ poss\`ede alors une s\'erie
discr\`ete. Commen\c{c}ons par \'etudier le probl\`eme de la restriction des repr\'esentations de la s\'erie discr\`ete de $G$ \`a un 
sous-groupe de $G$. Soit donc $H$ un sous-groupe r\'eductif connexe ferm\'e dans $G$. Supposons que l'intersection $K^H = K \cap H$
d'un sous-groupe compact maximal $K$ de $G$ soit un sous-groupe compact maximal dans $H$. Le th\'eor\`eme suivant se d\'eduit 
imm\'ediatement des travaux de Li \cite{Li} et de Harris et Li, en particulier de \cite[Proposition 1.2.3]{HarrisLi}.

\begin{thm} \label{res disc}
Soit $\rho$ une repr\'esentation unitaire irr\'eductible de la s\'erie discr\`ete de $G$ de plus bas $K$-type $\tau$.
Soit $\pi$ une repr\'esentation de la s\'erie discr\`ete de $H$ de plus bas $K^H$-type $\sigma$. Supposons que le
$K^H$-type $\sigma$ intervienne dans la restriction de $\tau$ \`a $K^H$. Alors, la repr\'esentation $\pi$
est \'equivalente \`a une sous-repr\'esentation irr\'eductible de $\rho_{|H}$.
\end{thm}
{\it D\'emonstration.} Il est bien connu, cf. \cite{Knapp}, que la repr\'esentation $\rho$ admet un unique plus bas $K$-type qui est 
donc $\tau$.
Notons $P_{\tau}$ la projection orthogonale sur la $\tau$-composante de $\rho$ et posons
$$\psi_{\rho} (x) = \frac{1}{{\rm dim} (\tau)} {\rm tr} (P_{\tau} \rho (x) P_{\tau} ), \  x \in G.$$
Pour appliquer \cite[Proposition 1.2.3]{HarrisLi} \`a la repr\'esentation $\rho$, il nous suffit donc de v\'erifier que
\begin{enumerate}
\item la restriction de $\rho$ \`a $H$ est fortement $L^{2+\varepsilon}$;
\item la fonction $\psi_{\rho}$ v\'erifie la formule de Flensted-Jensen.
\end{enumerate}
Nous renvoyons \`a \cite{Li} ou \cite{HarrisLi} pour la d\'efinition de la formule de Flensted-Jensen. Il suffit ici de noter
que Flensted-Jensen d\'emontre pr\'ecisemment dans \cite{FlenstedJensen} que celle-ci est v\'erifi\'ee par $\psi_{\rho}$
d\`es que $\rho$ appartient \`a la s\'erie discr\`ete de $G$. Le point 2. est donc v\'erifi\'e.

Rappelons maintenant qu'une repr\'esentation unitaire $\rho$ de $H$ dans un espace de Hilbert ${\cal H}_{\rho}$ est fortement $L^{2+\varepsilon}$ si elle 
est fortement $L^p$ pour tout $p>2$. Et qu'elle est fortement $L^p$ si, pour un ensemble dense de vecteurs dans ${\cal H}_{\rho}$, les coefficients matriciels 
associ\'es sont tous dans $L^p (H)$.

Ici la repr\'esentation $\rho$ appartient \`a la s\'erie discr\`ete de $G$, elle est en particulier temp\'er\'ee. Notons ${\cal H}_{\rho}^K$ l'espace des vecteurs $K$-finis dans 
${\cal H}_{\rho}$. Soit $\Xi_G$ la fonction sph\'erique d'Harish-Chandra \cite{HarishChandra} (dans \cite{Knapp}, la fonction $\Xi_G$ est appel\'ee $\varphi_0^G$). 
Puisque $\rho$ est temp\'er\'ee, pour tous $u,v \in {\cal H}_{\rho}^K$ on a  
\begin{eqnarray} \label{coefft}
| (\rho (g) u , v)| \leq C_{u,v} \Xi_G (g) , 
\end{eqnarray}
o\`u $C_{u,v}$ est une constante qui d\'epend de $u$ et $v$ mais pas de $g$.
Fixons deux d\'ecompositions compatibles $G=KAN$ et $H=K^H A^H N^H$ et notons $\rho$ (resp. $\rho^H$) la demi-somme
des racines de $(\mathfrak{g} , \mathfrak{a})$ (resp. $(\mathfrak{h} , \mathfrak{a}^H )$) positives pour $N$ (resp.  $N^H$). Il est imm\'ediat que 
$\rho \geq \rho^H$. La Proposition 7.15 de \cite{Knapp} implique alors 
$$\Xi_G (a) \leq ({\rm const}) e^{-\rho^H \log a} (1+||a||)^d , $$
pour $a \in  A^H_+$ et o\`u $d$ est une certaine constante $\geq 0$. La d\'emonstration de \cite[(7.52)]{Knapp} implique alors imm\'ediatement que 
la fonction $\Xi_G$ restreinte \`a $H$ est dans $L^p (H)$ pour tout $p>2$. Puisque le sous-espace ${\cal H}_{\rho}^K$ est dense dans ${\cal H}_{\rho}$ l'in\'egalit\'e
(\ref{coefft}) implique finalement que la restriction de $\rho$ \`a $H$ est fortement $L^{2+\varepsilon}$. Ce qui permet de conclure la d\'emonstration du Th\'eor\`eme \ref{res
disc} en appliquant \cite[Proposition 1.2.3]{HarrisLi}.~$\Box$

\bigskip

Nous allons maintenant chercher \`a \'etudier un probl\`eme analogue dans le cas des repr\'esentations cohomologiques des groupe unitaires et orthogonaux. L'id\'ee est 
de se ramener au Th\'eor\`eme \ref{res disc} \`a l'aide de la correspondance th\'eta $L^2$ locale. Commen\c{c}ons par quelques rappels \`a ce sujet.

Soit $(G,G')$ une paire r\'eductive duale irr\'eductible de type I dans le groupe symplectique $Sp=Sp_{2n} ({\Bbb R})$ (nous renvoyons \`a l'article \cite{Howe} de Howe 
pour plus de pr\'ecisions concernant cette terminologie). Soit $\widetilde{Sp}$ le rev\^etement m\'etaplectique \`a deux feuillet de $Sp$. \'Etant donn\'e un sous-groupe $E$ de $Sp$ nous notons
$\tilde{E}$ son image inverse dans $\widetilde{Sp}$. Dans \cite{Li}, Li \'etudie la correspondance th\'eta locale entre les repr\'esentations
de la s\'erie discr\`ete de $\tilde{G}'$ et les repr\'esentations cohomologiques unitaires de $\tilde{G}$. Nous exploitons maintenant les r\'esultats (et m\'ethodes) de Li pour faire
correspondre au Th\'eor\`eme \ref{res disc} un th\'eor\`eme sur la restriction des repr\'esentations cohomologiques. 

Rappelons (cf. \cite{Howe}) qu'une paire duale irr\'eductible de type I est construite comme suit. Soit $D$ l'une des trois alg\`ebres \`a division sur ${\Bbb R}$ ($D$ est donc \'egal
\`a ${\Bbb R}$, ${\Bbb C}$ ou ${\Bbb H}$, l'alg\`ebre des quaternions), munie de son involution standard $*$. (L'involution $*$ est donc triviale dans le premier cas et est la 
conjugaison complexe (resp. quaternionique) dans les deux derniers cas.) Soient $V$ et $V'$ deux espaces vectoriels de dimension finie sur $D$ \'equip\'e de deux formes 
$*$-sesquilin\'eaires non d\'eg\'en\'er\'ees $(.,.)$ et $(.,.)'$, l'une $*$-hermitienne et l'autre $*$-anti-hermitienne. Soient $G$ et $G'$ les groupes d'isom\'etries respectifs de 
$(.,.)$ et $(.,.)'$. Alors $(G,G')$ est une paire duale irr\'eductible dans $Sp=Sp_{2n} ({\Bbb R})$, o\`u 
$$2n = \dim_{{\Bbb R}} (D) (\dim_D V)(\dim_D V').$$
Nous supposerons toujours que la ``taille'' de $G'$ est inf\'erieure \`a celle de $G$, \`a savoir que 
\begin{eqnarray} \label{Li1}
\dim_D V \geq \dim_D V' .
\end{eqnarray}
Nous notons enfin 
\begin{eqnarray} \label{Li70}
G_1  = \left\{
\begin{array}{ll}
SO(p,q) & \mbox{ si } G  = O(p,q) \\
SU(p,q) & \mbox{ si } G  = U(p,q) \\
G  & \mbox{ sinon}.
\end{array} \right.
\end{eqnarray}

Consid\'erons maintenant $A_{\mathfrak{q}}$ une repr\'esentation cohomologique de $G_1$ associ\'ee \`a une sous-alg\`ebre parabolique 
$\theta$-stable $\mathfrak{q}=\mathfrak{l} \oplus \mathfrak{u}$ du complexifi\'e $\mathfrak{g}$
de l'alg\`ebre de Lie $\mathfrak{g}_0$ de $G_1$. Posons $\mathfrak{l}^0 = \mathfrak{l} \cap \mathfrak{g}_0$. Nous consid\'erons dans cette section les repr\'esentations
cohomologiques  $A_{\mathfrak{q}}$ v\'erifiant la condition 
\begin{eqnarray} \label{condsurcohom}
\mathfrak{l}^0 \cong \mathfrak{l}_0 ' \oplus \mathfrak{g}_0^1,
\end{eqnarray} 
o\`u $\mathfrak{l}_0 '$ est une alg\`ebre de Lie compacte et $\mathfrak{g}_0^1$ est du ``m\^eme type''
que $\mathfrak{g}_0$, c'est \`a dire isomorphe \`a l'alg\`ebre de Lie du groupe des isom\'etries de $(.,.)_{|V^1}$, o\`u $V^1$ est un sous-espace non d\'eg\'en\'er\'e de $V$.

D'apr\`es \cite[Theorem 6.2]{Li}, il existe une repr\'esentation $\pi '$ de la s\'erie discr\`ete de $\tilde{G}'$ telle que $\pi '$ admette un 
relev\'e th\'eta non trivial au groupe $\tilde{G}$ $\pi$ dont la restriction au sous-groupe $G_1 \subset  \tilde{G}$ soit pr\'ecisemment la 
repr\'esentation $A_{\mathfrak{q}}$.

Pr\'ecisons un peu ce r\'esultat en supposant $G$ non compact. Soit $(\omega , {\cal Y})$ la repr\'esentation de Weil du groupe
$\widetilde{Sp}$ munie de sa structure unitaire. L'un des deux groupe $G$, $G'$ est de type hermitien nous supposerons que c'est le cas de $G'$.
Soient $K$ et $K'$ deux sous-groupes compacts maximaux respectifs de $G$ et $G'$ et $M' \supset G'$ le centralisateur de $K$ dans $Sp$.
Puisque $G$ n'est pas compact, $M' \cong G' \times G'$ et $(K,M')$ forme une paire duale dans $Sp$. Et puisque $K$ est compact, il est bien connu que comme repr\'esentation
unitaire dans l'espace de Hilbert ${\cal Y}$, la restriction de $\omega$ \`a $\tilde{K} \cdot \tilde{M}'$ se d\'ecompose en une somme directe 
$\sum_i \sigma_i \otimes \rho_i$ de repr\'esentations unitaires irr\'eductibles de $\tilde{K} \cdot \tilde{M}'$. La correspondance th\'eta est alors $\sigma_i \leftrightarrow \rho_i$.
Cette correspondance est connue et explicit\'ee dans \cite{Li}, chaque $\rho_i$ ainsi obtenue est une repr\'esentation unitaire de plus haut poids de $M'$.
Soit $\sigma$ le plus bas $\tilde{K}$-type de la repr\'esentation $\pi$. La repr\'esentation $\sigma$ intervient
dans la correspondance duale avec $M'$. Notons $\sigma \otimes \rho '$ le facteur direct correspondant dans la d\'ecomposition en irr\'eductibles de la 
la restriction de $\omega$ \`a $\tilde{K} \cdot \tilde{M}'$. Soit $\tau'$ le plus bas $\tilde{L}'$-type de $\rho'$, vue comme repr\'esentation de $\tilde{K}' \times \tilde{K}'$, 
$\tau ' = \sigma_1 \otimes \sigma_2$. 
La restriction de $\tau '$ \`a la diagonale $\tilde{K}' \subset \tilde{K}' \times \tilde{K}'$ contient un facteur irr\'eductible $\sigma '$ dont le plus haut poids est \'egal \`a
la somme des plus hauts poids de $\sigma_1$ et $\sigma_2$. La repr\'esentation $\sigma '$ est pr\'ecisemment le plus bas $\tilde{K}'$-type de $\pi '$.   

Suivant Kudla \cite{Kudla}, nous dirons que deux paires r\'eductives duales irr\'eductibles et de type I $(H,H')$ et $(G,G')$ dans le groupe symplectique $Sp$
sont en {\it balance} (``see-saw'') si $H\subset G$ et (donc) $G' \subset H'$. Consid\'erons donc deux telles paires en balance, ce que nous repr\'esentons par le diagramme
$$\begin{array}{ccc}
G & & H' \\
| & \times & | \\
H & & G' 
\end{array}$$

Mettons ici en balance la paire $(G,G')$ et la paire $(K,M')$. Puisque $K$ est compact, la correspondance de Howe $L^2$ est classique pour la paire $(K,M')$, autrement dit
la repr\'esentation $\sigma \otimes \rho '$ de $\tilde{K} \cdot \tilde{M}'$ est \'equivalente \`a une sous-repr\'esentation irr\'eductible de la restriction de $\omega$ 
\`a $\tilde{K} \cdot \tilde{M}'$. En particulier la repr\'esentation $\rho '$ de $\tilde{M}'$ est \'equivalente \`a une sous-repr\'esentation irr\'eductible de la restriction de $\omega$ 
\`a $\tilde{M}'$. Une g\'en\'eralisation due \`a Li \cite[Theorem 4.1]{Li} du Th\'eor\`eme \ref{res disc} implique (puisque 1) $\rho'$ est une repr\'esentation unitaire de plus haut
poids et 2) la restriction de $\omega$ \`a $\tilde{G}'$ est fortement $L^{2+\varepsilon}$) que la repr\'esentation $\pi '$ de $\tilde{G}'$ est \'equivalente \`a une sous-repr\'esentation irr\'eductible de la restriction de $\rho'$ et donc de $\omega$ \`a $\tilde{G}'$. La repr\'esentation $\pi'$ intervient donc dans la correspondance de Howe $L^2$. 
D'apr\`es Howe \cite[Theorem 6.1]{Howe2}, la repr\'esentation $\pi \otimes \pi'$ de $\tilde{G} \cdot \tilde{G}'$ est alors \'equivalente \`a une sous-repr\'esentation 
irr\'eductible de la restriction de $\omega$ \`a $\tilde{G} \cdot \tilde{G}'$ et elle intervient avec multiplicit\'e un. La composante $\pi$-isotypique ${\cal Y} (\pi )$ de la repr\'esentation 
unitaire $(\omega , {\cal Y})$ est isomorphe \`a $\pi \otimes \pi '$ et co\"{\i}ncide donc avec la composante $\pi '$-isotypique ${\cal Y} (\pi ')$.

On peut alors appliquer dans ce contexte une id\'ee due \`a Howe \cite{Howe3}. \'Etant donn\'ees deux paires r\'eductives duales irr\'eductibles de type I $(H,H')$ et $(G,G')$ 
en balance dans le groupe symplectique $Sp$ et $\pi \leftrightarrow \pi '$ (resp. $\sigma \leftrightarrow \sigma '$) des repr\'esentations intervenant dans la correspondance
de Howe $L^2$ pour la paire $(G,G')$ (resp. $(H,H')$). Supposons que la repr\'esentation $\pi '$ de $\tilde{G}'$ soit \'equivalente \`a une sous-repr\'esentation irr\'eductible de 
la restriction de $\sigma '$ \`a $\tilde{G}'$. La composante $\pi '$-isotypique de ${\cal Y} (\sigma )$ est alors non triviale $={\cal Y} (\sigma \otimes \pi ' )$. Mais 
$G'$ et $H$ commutent dans $Sp$ et donc
$${\cal Y} (\pi ' ) (\sigma ) = {\cal Y} (\sigma \otimes \pi ' ) = {\cal Y} (\sigma ) (\pi ') .$$
Puisqu'enfin ${\cal Y} (\pi ') = \pi \otimes \pi '$, la repr\'esentation $\sigma$ de $\tilde{H}$ est n\'ecessairement \'equivalente \`a une sous-repr\'esentation irr\'eductible de 
la restriction de $\pi$ \`a $\tilde{H}$ et intervient avec la m\^eme multiplicit\'e que $\pi '$ dans $\sigma '$.

Dans la prochaine section nous appliquons ce principe pour d\'eduire du Th\'eor\`eme \ref{res disc} des r\'esultats analogues pour les repr\'esentations cohomologiques
des groupes unitaires et orthogonaux.

\subsection{Restriction de repr\'esentations cohomologiques}

Nous revenons maintenant au cas $G^{{\rm nc}} = U(p,q)$ ou $O(p,q)$ et supposons que $G$ contient un sous-groupe $H$ tel que 
$H^{{\rm nc}} = U(p,q-r)$ ou $O(p,q-r)$, plong\'e de mani\`ere usuelle (stable par l'involution de Cartan) et avec $1 \leq p , q$ et $1 \leq r < q$.

Soient $\pi$ une repr\'esentation unitaire irr\'eductible de $G$ et $\pi '$ une repr\'esentation unitaire irr\'eductible de $H$. Par d\'efinition la multiplicit\'e de $\pi '$ dans 
$\pi_{|H}$ est le plus grand entier $m$ tel que $\pi_{|H}$ contienne $m$ sous-espaces irr\'eductibles $2$ \`a $2$ orthogonaux  sur chacun desquels l'action de $H$ est 
\'equivalente \`a la repr\'esentation $\pi '$. Soit $V$ l'espace de la repr\'esentation $\pi$ (c'est \'egalement l'espace de la repr\'esentation $\pi_{|H}$). Soit $V' \subset V$
le plus grand sous-espace sur lequel l'action de $H$ est \'equivalente \`a un multiple de $\pi '$. Soit $\beta : V \rightarrow V'$ la projection orthogonale correspondante.
Soit $K$ le sous-groupe compact maximal de $G^{\rm nc}$ et $K^H = K \cap H^{\rm nc}$.
Notons $V_0$ (resp. $V_0 '$) le $(\mathfrak{g} , K)$-module ($(\mathfrak{h} , K^H)$-module) constitu\'e des vecteurs $K$-finis de $V$ (resp. $K^H$-finis de $V'$). On a alors
une application naturelle 
\begin{eqnarray} \label{147}
H^* (\mathfrak{g}, K; V_0 ) \rightarrow  H^* (\mathfrak{h} , K^H ; V_0 ') ,
\end{eqnarray}
obtenue en composant 
\begin{eqnarray} \label{148}
\begin{array}{ccl}
H^* (\mathfrak{g}, K; V_0 ) & \rightarrow & H^* (\mathfrak{h} , K^H ; V_0 ) \\
& \rightarrow &  H^* (\mathfrak{h} , K^H ; V_0 ')
\end{array}
\end{eqnarray}
o\`u la premi\`ere application est obtenue en restreignant de $(\mathfrak{g} , K)$ \`a $(\mathfrak{h} , K^H)$ et la seconde est induite par la projection
$\beta : V_0 \rightarrow V_0 '$.

Nous d\'emontrons maintenant les deux th\'eor\`emes suivant concernant la restriction de  certaines repr\'esentations cohomologiques respectivement
dans le cas unitaire et dans le cas orthogonal.

\begin{thm} \label{res cohomU}
Soient $H=U (p,q-r) \subset U (p,q)= G$ o\`u l'inclusion est l'inclusion standard, $1\leq p ,q$ et $1 \leq r <q$. Alors, pour
tout couple d'entiers naturels $(i,j)$ tel que $i +j\leq q-r$,
\begin{enumerate}
\item la repr\'esentation cohomologique $A((i^p) , ((q-r-j)^p ))_H$ de $H$ appara\^{\i}t avec multiplicit\'e un dans la restriction \`a $H$ de la repr\'esentation cohomologique $A((i^p) , ((q-j)^p))$ de $G$, et
\item l'application naturelle en cohomologie
\begin{eqnarray} \label{res cohomU'}
\begin{array}{rl}
H^{pi , pj} (\mathfrak{g} , K ; A((i^p) , ((q-j)^p ))) & \\
 \rightarrow &  H^{pi, pj} (\mathfrak{h} , K^H ;  A((i^p) , ((q-r-j)^p ))_H )
\end{array}
\end{eqnarray}
est un isomorphisme d'espaces de dimension un.
\end{enumerate}
\end{thm}

\begin{thm} \label{res cohomO}
Soient $H=SO_0 (p,q-r) \subset SO_0 (p,q)= G$ o\`u l'inclusion est l'inclusion standard, $1\leq p ,q$ et $1 \leq r <q$. Alors, pour
tout entier naturel $i\leq (q-r)/2$,
\begin{enumerate}
\item la repr\'esentation cohomologique $A((i^p))^{\pm}_H$ de $H$ appara\^{\i}t avec multiplicit\'e un dans la restriction \`a $H$ de la repr\'esentation cohomologique
$A((i^p))$ de $G$, et
\item l'application naturelle en cohomologie
\begin{eqnarray} \label{res cohomO'}
H^{pi} (\mathfrak{g} , K ; A((i^p))) \rightarrow H^{pi} (\mathfrak{h} , K^H ;  A((i^p))^{\pm}_H)
\end{eqnarray}
est un isomorphisme d'espaces de dimension un.
\end{enumerate}
\end{thm}

\noindent
{\it D\'emonstrations.} Ces deux Th\'eor\`emes se d\'emontrent de la m\^eme mani\`ere. 
Nous commen\c{c}ons par d\'emontrer le point dans le cas des groupes orthogonaux. 
Nous allons travailler avec les groupes $O(p,q)$ et $O(p,q-r) \times O(r)$ et donc avec les repr\'esentations
$\overline{A((i^p))}$ et $\overline{A((i^p))}_H$; cette derni\`ere repr\'esentation \'etant triviale sur le facteur $O(r)$ 
(nous la consid\`ererons tour \`a tour comme une
repr\'esentation de $O(p,q-r) \times O(r)$ ou de $O(p,q-r)$). En regardant les $K$-types, il est facile de v\'erifier que la restriction de $\overline{A((i^p))}$ au groupe
$SO_0 (p,q)$ reste irr\'eductible et est \'egale \`a $A((i^p))$, tant que $i< q/2$. Quand $q-r=2i$, la restriction de $\overline{A((i^p))}_H$ au groupe
$SO_0 (p,q-r) \times SO_0 (r)$ est somme directe des deux repr\'esentations irr\'eductibles $A((i^p))^+_H$ et $A((i^p))^-_H$ de $SO_0 (p,q-r)$. Dans tous les cas la
multiplicit\'e de la repr\'esentation $A((i^p))^{\pm}_H$ dans la restriction \`a $H$ de la repr\'esentation $A((i^p))$ de $G$ est donc \'egale \`a la multiplicit\'e de
$\overline{A((i^p))}_H$ dans $\overline{A((i^p))}$. C'est cette derni\`ere multiplicit\'e que nous calculons maintenant.

Ce calcul repose sur la correspondance th\'eta locale. Nous consid\'erons les inclusions
$$O(r) \times O(p, q-r) \subset O(p,q)$$
et
$$Sp (2i , {\Bbb R}) \subset Sp (2i ,{\Bbb R}) \times Sp (2i , {\Bbb R})$$
(plongement diagonal). Les paires $(O(p,q ), Sp (2i ,{\Bbb R}))$ et
$(O(r) \times O(p, q-r) ,Sp (2i ,{\Bbb R}) \times Sp (2i , {\Bbb R}) )$ sont des paires r\'eductives duales en balance dans 
$Sp_{2i(p+q)}$~:
\begin{eqnarray} \label{seesaw}
\begin{array}{ccc}
O(p,q) &  & Sp (2i ,{\Bbb R}) \times Sp (2i , {\Bbb R}) \\
| & \times & | \\
 O(r) \times O(p, q-r) & & Sp (2i , {\Bbb R})
\end{array}
\end{eqnarray}
Soient $\pi$, $\pi_1$ et $\pi_2$ les relev\'es respectifs au groupe $\widetilde{Sp}(2i , {\Bbb R})$,
via la correspondance th\'eta locale, de la repr\'esentation $\overline{A((i^p))}$ de $O(p,q)$, de la repr\'esentation
triviale de $O(r)$ et de la repr\'esentation $\overline{A((i^p))}_H$ de $O(p, q-r)$.
Dans \cite{Li}, Li montre par la m\'ethode que nous avons rappel\'e au paragraphe pr\'ec\'edent que ces repr\'esentations 
interviennent dans la correspondance $L^2$, il d\'ecrit de plus explicitement les
relev\'es th\'eta locaux de ces repr\'esentations cohomologiques. Sous l'hypoth\`ese $2i\leq q-r$, 
Li montre que les trois repr\'esentations $\pi$, $\pi_1$ et $\pi_2$
sont toutes des s\'eries discr\`etes, de plus bas $K$-types (ayant pour plus haut poids) respectifs
$(-(p+q)/2 , \ldots , -(p+q)/2)$, $(-r/2 , \ldots , -r/2 )$ et $(-(p+q-r)/2 , \ldots , -(p+q-r)/2 )$. Remarquons que
la somme des deux derniers poids est \'egale au premier poids. Puisque $\pi_1$ et $\pi_2$ sont des s\'eries discr\`etes,
elles contiennent leur plus bas $K$-type avec multiplicit\'e un, le $K$-type $(-(p+q)/2 , \ldots , -(p+q)/2)$ appara\^it
donc avec multiplicit\'e exactement un dans la restriction du produit tensoriel $\pi_1 \otimes \pi_2$. Le Th\'eor\`eme \ref{res disc}
implique alors que la repr\'esentation $\pi$ appara\^{i}t avec multiplicit\'e exactement \'egale \`a $1$ dans le
produit tensoriel $\pi_1 \otimes \pi_2$. Le r\'esultat de Howe rappel\'e au pr\'ec\'edent paragraphe implique alors que la
multiplicit\'e de la repr\'esentation $A((i^p))^{\pm}_H$ dans la restriction \`a $H$ de la repr\'esentation $A((i^p))$ de $G$
est exactement \'egale \`a $1$. Ce qui conclut la d\'emonstration du premier point dans le cas orthogonal.

La d\'emonstration dans le cas unitaire est identique (elle est m\^eme un peu plus simple puisque les groupes sont
connexes) en utilisant les groupes en balance suivants~:
\begin{eqnarray*}
\begin{array}{ccc}
U(p,q) &  & U (i ,j) \times U (i , j) \\
| & \times & | \\
U(r) \times U(p, q-r) & & U (i ,j)
\end{array}
\end{eqnarray*}

\medskip

Il nous reste \`a montrer les deuxi\`emes points des Th\'eor\`emes \ref{res cohomU} et \ref{res cohomO}. Nous
voulons donc d\'emontrer que les applications naturelles (\ref{res cohomU'}) et (\ref{res cohomO'}) sont des isomorphismes
entre espaces de dimension un. La d\'emonstration est similaire dans les cas orthogonaux et unitaires. Pour changer
nous allons maintenant traiter le cas unitaire.

Dans \cite{VoganZuckerman} Vogan et Zuckerman d\'emontrent que les espaces des deux c\^ot\'es de
(\ref{res cohomU'}) sont de dimension un. Rappelons que $V((i^p) , ((q-j)^p))$ (resp. $V ((i^p) , ((q-r-j)^p))_H$)
est le plus bas $K$-type de la repr\'esentation $A((i^p) , ((q-j)^p))$ (resp. $A((i^p) , ((q-r-j)^p))_H$).
Le membre de gauche (resp. droite) de (\ref{res cohomU'}) est isomorphe \`a
\begin{eqnarray} \label{1410}
\begin{array}{l}
{\rm Hom}_{K} (\bigwedge^{pi,pj} \mathfrak{p} , V((i^p) , ((q-j)^p)) ) \\   
(\mbox{resp. } \ {\rm Hom}_K ( \bigwedge^{pi,pj} \mathfrak{p}_H , V((i^p) , ((q-r-j)^p))_H ) ).
\end{array}
\end{eqnarray}
La d\'emonstration de la premi\`ere partie du Th\'eor\`eme et le Lemme \ref{resU} (ce serait le Lemme \ref{resO} dans le cas
orthogonal) impliquent que l'application compos\'ee
\begin{eqnarray} \label{1411}
V((i^p) , ((q-j)^p)) \rightarrow A((i^p) , ((q-r-j)^p))_H \rightarrow V ((i^p) , ((q-r-j)^p))_H ,
\end{eqnarray}
o\`u la premi\`ere application est la restriction de la projection orthogonale $\beta$ et la deuxi\`eme
application est la projection sur le plus bas $K^H$-type, est non nulle. Elle est donc surjective, puisque
le membre de droite de (\ref{1411}) est irr\'eductible.

De plus, le plus haut poids appara\^{\i}t avec multiplicit\'e un dans chacune des extr\'emit\'es de (\ref{1411}).
L'application (\ref{1411}) envoie donc un vecteur de plus haut poids non nul vers un vecteur de plus haut poids non nul.
Compte tenu de (\ref{1410}) ceci implique que l'application (\ref{res cohomU'}) est bien un isomorphisme.~$\Box$

\bigskip

Aux vus des Lemmes \ref{resU} et \ref{resO}, nous conjecturons plus g\'en\'eralement les \'enonc\'es suivants.

\begin{conj} \label{CU1}
Soient $H=U (p,q-r) \subset U (p,q)= G$ o\`u l'inclusion est l'inclusion standard, $1\leq p ,q$ et $1 \leq r <q$.
Soient $\lambda$ et $\mu$ deux partitions incluses dans $p \times q$ formant un couple compatible.
\begin{enumerate}
\item La restriction \`a $H$ de la repr\'esentation cohomologique $A(\lambda , \mu )$ de $G$ contient (discr\`etement)
une repr\'esentation cohomologique de degr\'e fortement primitif $|\lambda |+ |\hat{\mu}|$ si et seulement si
la partition $(r^p)$ s'inscrit dans le diagramme gauche $\mu/ \lambda$. Elle contient dans ce cas la repr\'esentation
cohomologique $A(\lambda , \mu -(r^p))$ de $H$ avec multiplicit\'e exactement \'egale \`a $1$.
\item Supposons que la partition $(r^p)$ s'inscrive dans le diagramme gauche $\mu/ \lambda$. Alors,
l'application naturelle en cohomologie
\begin{eqnarray*}
H^{|\lambda | , |\hat{\mu}|} (\mathfrak{g} , K ; A(\lambda , \mu)) \rightarrow H^{|\lambda |, |\hat{\mu}|} (\mathfrak{h} , K^H ;  A(\lambda , \mu -(r^p) )_H )
\end{eqnarray*}
est un isomorphisme d'espaces de dimension un.
\end{enumerate}
\end{conj}

\begin{conj} \label{CU2}
Soient $H=SO_0 (p,q-r) \subset SO_0 (p,q)= G$ o\`u l'inclusion est l'inclusion standard, $1\leq p ,q$ et $1 \leq r <q$. Soit $\lambda$ une partition orthogonale dans $p\times q$.
\begin{enumerate}
\item La restriction \`a $H$ de la repr\'esentation cohomologique $A(\lambda)^{\pm_2}_{\pm_1}$ de $G$
contient (discr\`etement) une repr\'esentation cohomologique de degr\'e fortement primitif $|\lambda |$ si et
seulement si la partition $(r^p)$ s'inscrit dans le diagramme gauche $\hat{\lambda} / \lambda$.
Elle contient dans ce cas la repr\'esentation cohomologique $A(\lambda)_{\pm_1}^{\pm_2}$ de $H$ avec multiplicit\'e exactement \'egale \`a $1$.
\item Supposons que la partition $(r^p)$ s'inscrive dans le diagramme gauche $\hat{\lambda}/ \lambda$. Alors,
l'application naturelle en cohomologie
\begin{eqnarray*}
H^{|\lambda|} (\mathfrak{g} , K ; A(\lambda)_{\pm_1}^{\pm_2}) \rightarrow H^{|\lambda|} (\mathfrak{h} , K^H ;  A(\lambda)^{\pm_2}_{\pm_1})
\end{eqnarray*}
est un isomorphisme d'espaces de dimension un.
\end{enumerate}
\end{conj}

\bigskip

Remarquons enfin qu'en consid\'erant les paires $(U(p,q) , U(i))$ et $(O(p,q) , Sp(2i, {\Bbb R}))$ en balance dans $Sp_{2i(p+q)}$~:
\begin{eqnarray*}
\begin{array}{ccc}
U(p,q) &  & U (i ) \\
| & \times & | \\
O(p,q) \times Sp(2i , {\Bbb R}),
\end{array}
\end{eqnarray*}
la d\'emonstration des Th\'eor\`emes \ref{res cohomU} et \ref{res cohomO} implique le th\'eor\`eme suivant.

\begin{thm} \label{rescohomUO}
Soient $H=O (p,q) \subset U (p,q)= G$ o\`u l'inclusion est l'inclusion standard, $1\leq p ,q$. Alors, pour
tout entier naturel $i\leq q/2$,
\begin{enumerate}
\item la repr\'esentation cohomologique $\overline{A((i^p))}_H$ de $H$ appara\^{\i}t avec multiplicit\'e un dans la restriction \`a 
$H$ de la repr\'esentation cohomologique $A((i^p))$ de $G$, et
\item l'application naturelle en cohomologie
\begin{eqnarray} \label{res cohomO'}
H^{pi,0} (\mathfrak{g} , K ; A((i^p))) \rightarrow H^{pi} (\mathfrak{h} , K^H ;  \overline{A((i^p))}_H)
\end{eqnarray}
est un isomorphisme d'espaces de dimension un.
\end{enumerate}
\end{thm}
 
Nous conjecturons plus g\'en\'eralement l'\'enonc\'e suivant.

\begin{conj} \label{CUO}
Soient $H=O (p,q) \subset U (p,q)= G$ o\`u l'inclusion est l'inclusion standard, $1\leq p ,q$. Soient $\lambda$ et $\mu$ deux partitions
incluses dans $p\times q$ formant un couple compatible.
\begin{enumerate}
\item Si la restriction \`a $H$ de la repr\'esentation cohomologique $A(\lambda , \mu )$ de $G$ contient 
(discr\`etement) une repr\'esentation de degr\'e fortement primitif $|\lambda | + |\hat{\mu}|$ alors $\lambda=0$ ou 
$\mu = p\times q$.
\item Supposons par exemple $\mu = p \times q$. Alors la restriction \`a $H$ de la repr\'esentation cohomologique 
$A(\lambda )$ de $G$ contient (discr\`etement) une repr\'esentation de degr\'e fortement primitif $|\lambda |$ si et seulment si 
la diagramme est $\lambda$ est orthogonal. Elle contient dans ce cas la repr\'esentation cohomologique 
$\overline{A(\lambda)}_H$ de $H$ avec multiplicit\'e exactement \'egale \`a $1$.
\item Supposons que la partition $\lambda$ est orthogonale. Alors, l'application naturelle en cohomologie 
\begin{eqnarray} \label{res cohomO'}
H^{|\lambda |,0} (\mathfrak{g} , K ; A(\lambda)) \rightarrow H^{|\lambda|} (\mathfrak{h} , K^H ;  \overline{A(\lambda)}_H)
\end{eqnarray}
est un isomorphisme d'espaces de dimension un.
\end{enumerate}
\end{conj}
 
\subsection{Produits tensoriels de repr\'esentations cohomologiques}

La d\'emonstration des Th\'eor\`emes \ref{res cohomU} et \ref{res cohomO} peut \'egalement 
s'appliquer aux paires d'espaces sym\'etriques r\'eels respectifs
\begin{eqnarray*}
\begin{array}{l}
U(p,q) /(U(p, q-l-k) \times U(l) \times U(k)) \\
\subset (U(p,q) \times U(p,q)) / ((U(p,q-l) \times U(l)) \times (U(p,q-k) \times U(k)))
\end{array}
\end{eqnarray*}
et
\begin{eqnarray*}
\begin{array}{l}
SO_0 (p,q) / (SO_0 (p,q-l-k) \times SO (l) \times SO (k)) \\
\subset (SO_0 (p,q) \times SO_0 (p,q)) / ((SO_0 (p,q-l) \times SO (l)) \times (SO_0 (p,q-k) \times SO(k))) .
\end{array}
\end{eqnarray*}

La consid\'eration des paires r\'eductives duales en balance~:
\begin{eqnarray*}
\begin{array}{ccc}
U(p,q) \times U(p,q) &  & U(i+k,j+l)  \\
| & \times & | \\
U(p, q) & & U(i,j) \times U(k,l)
\end{array}
\end{eqnarray*}
et
\begin{eqnarray*}
\begin{array}{ccc}
O (p,q) \times O (p,q) &  & Sp (2 (k+l) , {\Bbb R})  \\
| & \times & | \\
O (p, q) & & Sp (2k ,{\Bbb R}) \times Sp (2l , {\Bbb R})
\end{array}
\end{eqnarray*}
permet alors de d\'eriver les deux th\'eor\`emes suivants.

\begin{thm} \label{pdt cohomU}
Soit $G=U (p,q)$, $1\leq p ,q$. Alors, pour tout quadruplet $(i,j,k,l)$ d'entiers $\geq 0$ de somme $i+j+k+l \leq q$,
\begin{enumerate}
\item la repr\'esentation cohomologique $A(((i+k)^p) , ((q-j-l)^p ))$ de $G$ appara\^{\i}t avec multiplicit\'e un dans le produit tensoriel
des repr\'esentations cohomologiques $A((i^p) , ((q-j)^p))$ et $A((k^p) , ((q-l)^p))$ de $G$, et
\item l'application ``cup-produit''
\begin{eqnarray} \label{pdt cohomU'}
\begin{array}{l}
H^{pi , pj} (\mathfrak{g} , K ; A((i^p) , ((q-j)^p ))) \otimes  H^{pk,pl} (\mathfrak{g} , K ; A((k^p) , ((q-l)^p ))) \\
\rightarrow H^{p(i+k), p(j+l)} (\mathfrak{g} , K ;  A(((i+k)^p) , ((q-j-l)^p )) )
\end{array}
\end{eqnarray}
est un isomorphisme d'espaces de dimension un.
\end{enumerate}
\end{thm}

\begin{thm} \label{pdt cohomO}
Soit $G=SO_0 (p,q)$, $1\leq p ,q$. Alors, pour tout couple d'entiers $\geq 0$ de somme $k+l \leq q/2$,
\begin{enumerate}
\item la repr\'esentation cohomologique $A(((k+l)^p))^{\pm}$ de $G$ appara\^{\i}t avec multiplicit\'e un dans  le produit tensoriel
des repr\'esentations cohomologiques $A((k^p))^{\pm}$ et $A((l^p))^{\pm}$ de $G$, et
\item l'application ``cup-produit''
\begin{eqnarray} \label{pdt cohomO'}
\begin{array}{l}
H^{pk} (\mathfrak{g} , K ; A((k^p))^{\pm}) \otimes H^{pl} (\mathfrak{g} , K ; A((l^p))^{\pm}) \\
 \rightarrow H^{p(k+l)} (\mathfrak{g} , K ;  A(((k+l)^p))^{\pm})
\end{array}
\end{eqnarray}
est un isomorphisme d'espaces de dimension un.
\end{enumerate}
\end{thm}

De mani\`ere similaire \`a la Conjecture \ref{CU1}, et au vu de \cite{Produit} nous conjecturons de plus le r\'esultat suivant.

\begin{conj} \label{CP1}
Soit $G= U (p,q)$, $1\leq p ,q$.
Soient $(\lambda ,\mu )$ et $(\alpha , \beta )$ deux couples compatibles de partitions incluses dans $p \times q$.
\begin{enumerate}
\item Le produit tensoriel des repr\'esentations cohomologiques $A(\lambda , \mu )$ et $A(\alpha , \beta )$ de $G$ contient (discr\`etement)
une repr\'esentation cohomologique de degr\'e fortement primitif $|\lambda |+ |\hat{\mu}|+|\alpha |+ |\hat{\beta}|$ si et seulement si
il existe une partition $\nu \subset p\times q$ telle que $\nu$ (resp. $\hat{\nu}$) s'inscrive dans le diagramme gauche $\mu/ \lambda$ (resp. $\beta /\alpha$).
\item S'il existe une partition $\nu \subset p\times q$ telle que $\nu$ (resp. $\hat{\nu}$) s'inscrive dans le diagramme gauche $\mu/ \lambda$ (resp. $\beta /\alpha$), on
peut choisir $\nu$ telle que $|\nu | = |\alpha |+ |\hat{\beta}|$ et telle que l'image de $\nu$ dans $\mu/\lambda$ r\'eunit \`a $\lambda$ d\'efinisse
une partition $\lambda \uplus \nu \subset \mu$ formant un couple compatible avec $\mu$
\footnote{On peut effectivement d\'eduire cela dans le cas de \cite{Produit} \`a l'aide de \cite[Lemme 27]{BergeronTentative} par r\'ecurrence et en distingant les cas $p_1+ \ldots +p_m <$ ou $= p$.}.
Le produit tensoriel de $A(\lambda , \mu )$ et $A(\alpha , \beta )$ de $G$ contient
alors (discr\`etement) la repr\'esentation cohomologique $A (\lambda \uplus \nu , \mu )$.
\item La  repr\'esentation cohomologique $A (\lambda \uplus \nu , \mu )$ intervient avec multiplicit\'e exactement \'egale \`a $1$, et
l'application naturelle en cohomologie
\begin{eqnarray*}
H^{|\lambda | , |\hat{\mu}|} (\mathfrak{g} , K ; A(\lambda , \mu)) \otimes H^{|\alpha | , |\hat{\beta}|} (\mathfrak{g} , K ; A(\alpha , \beta)) \\
\rightarrow H^{|\lambda |+ |\hat{\mu}| + |\alpha |+ |\beta|} (\mathfrak{g} , K ;  A(\lambda \uplus \nu , \mu ) )
\end{eqnarray*}
est un isomorphisme d'espaces de dimension un.
\end{enumerate}
\end{conj}

Finalement les Conjectures \ref{CUO} et \ref{CP1} implique une conjecture analogue dans le cas des groupes orthogonaux, nous laissons
le soin au lecteur d'\'eventuellement l'\'enonc\'e.

\subsection{Repr\'esentations cohomologiques discr\`etement d\'ecomposables}

Soit toujours $G$ un groupe de Lie r\'eel r\'eductif et connexe. Consid\'erons $H$ un sous-groupe ferm\'e r\'eductif tel
que l'intersection $K^H = K \cap H$ du sous-groupe compact maximal $K$ de $G$ avec $H$ soit un sous-groupe compact
maximal de $H$. On dit (voir \cite[Definition 1.1, 1.2]{Kobayashi}) qu'un $(\mathfrak{g} , K)$-module unitarisable
et irr\'eductible est {\it discr\`etement d\'ecomposable} si, vu comme $(\mathfrak{h} , K^H)$-module, il est isomorphe
\`a une somme directe de $(\mathfrak{h} , K^H)$-modules irr\'eductibles.

Soit $\pi$ une repr\'esentation unitaire irr\'eductible de $G$. Remarquons que si la restriction $\pi_{|K^H}$ est
{\it $K^H$-admissible}, autrement dit si chaque repr\'esentation irr\'eductible de $K^H$ n'intervient qu'avec une
multiplicit\'e finie (peut-\^etre nulle) dans $\pi_{|K^H}$, alors le $(\mathfrak{g}, K)$-module associ\'e \`a $\pi$ est
discr\`etement d\'ecomposable comme $(\mathfrak{h} , K^H)$-module et chaque sous-$(\mathfrak{h} , K^H)$-module irr\'eductible
est de multiplicit\'e finie (voir \cite[Proposition 1.6(2)]{Kobayashi}).

Supposons maintenant que $H$ soit un sous-groupe ouvert du groupe des points fixes d'une involution $\tau$ sur $G$
qui commute \`a l'involution de Cartan de $G$. Nous notons $\mathfrak{k}_{0 \pm} = \{ X \in \mathfrak{k}_0 \; : \;
\tau (X) = \pm X \}$. Fixons $\mathfrak{t}_0$ une sous-alg\`ebre de Cartan $\tau$-stable de $\mathfrak{k}_0$ telle
que $\mathfrak{t}_{0 -} = \mathfrak{t}_0 \cap \mathfrak{k}_{0-}$ soit un sous-espace ab\'elien maximal de
$\mathfrak{k}_{0-}$. Soit $\Delta (\mathfrak{k} , \mathfrak{t} )$ (resp. $\Sigma (\mathfrak{k} , \mathfrak{t}_- )$)
le syst\`eme de racine (restreint) de $\mathfrak{k}$ par rapport \`a $\mathfrak{t}$ (resp. $\mathfrak{t}_-$). Fixons
alors deux sous-syst\`emes positifs $\Delta^+ (\mathfrak{k} , \mathfrak{t} )$ $\Sigma^+ (\mathfrak{k} , \mathfrak{t}_- )$
compatibles.

Soit maintenant $\mathfrak{q} = \mathfrak{q} (\lambda )$ une sous-alg\`ebre parabolique $\theta$-stable de $\mathfrak{g}$
d\'efinie par $\lambda \in i \mathfrak{t}_0$ dominant par rapport \`a $\Delta^+ (\mathfrak{k} , \mathfrak{t} )$. L'ensemble
$${\Bbb R}^+ \langle \mathfrak{u} \cap \mathfrak{p} \rangle = \left\{ \sum_{\beta \in \Delta (\mathfrak{u} \cap
\mathfrak{p} ; \mathfrak{t} )} n_{\beta} \beta \; : \; n_{\beta} \geq 0 \right\} $$
d\'efinit alors un c\^one ferm\'e dans $it_{0}$.

Le th\'eor\`eme suivant, d\^u \`a Kobayashi, se d\'eduit alors de \cite[Theorem 3.2]{Kobayashi2} et \cite[Theorem 4.2]{Kobayashi}.

\begin{thm} \label{kobayashi}
Conservons les notations ci-dessus. Alors, les trois conditions suivantes sont \'equivalentes.
\begin{enumerate}
\item ${\Bbb R}^+ \langle \mathfrak{u} \cap \mathfrak{p} \rangle \cap i \mathfrak{t}_{0-} = \{ 0 \}$.
\item la restriction $A_{\mathfrak{q} \; |K^H}$ de la repr\'esentation cohomologique $A_{\mathfrak{q}}$ de $G$ est
$K^H$-admissible.
\item le $(\mathfrak{g} , K)$-module $A_{\mathfrak{q}}$ est discr\`etement d\'ecomposable comme
$(\mathfrak{h} , K^H)$-module.
\end{enumerate}
\end{thm}

Nous allons maintenant v\'erifier le crit\`ere 1. du Th\'eor\`eme pour $(G,H) = (U(p,q) , U(p,q-r) \times U(r))$ (resp.
$=(SO_0 (p,q) , SO_0 (p,q-r) \times SO(r))$). Nous en d\'eduirons finalement de nouveaux cas de la Conjecture \ref{CU1}.

\medskip

\paragraph{Le cas unitaire}
Nous supposons donc maintenant $(G,H) = (U(p,q) , U(p,q-r) \times U(r))$ avec $2r \leq q$, o\`u le plongement $U(p,q-r) \times U(r) \rightarrow
U(p,q)$ est standard donn\'e par~:
$$(A,B) \mapsto \left(
\begin{array}{cc}
A & 0 \\ 
0 & B 
\end{array} \right) .$$
Apr\`es conjugaison par la matrice 
$$\left(
\begin{array}{cccc}
1_p & & & \\
& \frac{1}{\sqrt{2}} 1_r & 0 & \frac{i}{\sqrt{2}} 1_r  \\
& 0 & 1_{q-2r} & 0 \\
& \frac{i}{\sqrt{2}} J_r & 0 & \frac{1}{\sqrt{2}} J_r 
\end{array} \right) \in G $$
(ici $1_r$ d\'esigne la matrice identit\'e de taille $r$ et $J_r$ d\'esigne la matrice carr\'ee $r\times r$ avec que des $1$ sur l'anti-diagonale),
on obtient un plongement de $H$ dans $G$ pour lequel le sous-espace $\mathfrak{t}_0$ de $\mathfrak{g}_0$ constitu\'e des matrices diagonales (\`a coefficients imaginaires
purs) soit une une sous-alg\`ebre de Cartan $\tau$-stable de $\mathfrak{k}_0$ telle que $\mathfrak{t}_{0 -} = \mathfrak{t}_0 \cap \mathfrak{k}_{0-}$ 
soit un sous-espace ab\'elien maximal de $\mathfrak{k}_{0-}$. Si nous notons comme d'habitude $(x_1 , \ldots , x_p ; y_1 , \ldots , y_q)$ les \'el\'ements de $i\mathfrak{t}_0$, le
sous-espace $i\mathfrak{t}_{0-}$ est constitu\'e des \'el\'ements
$$(0, \ldots , 0 ; -u_1 , \ldots , -u_r , 0 , \ldots , 0 , u_r , \ldots , u_1 ).$$ 

Consid\'erons maintenant $\mathfrak{q} = \mathfrak{q} (\lambda , \mu )$ la sous-alg\`ebre parabolique $\theta$-stable associ\'ee \`a un couple compatible de partitions
$(\lambda , \mu)$. Les racines de $\mathfrak{t}$ dans $\mathfrak{p} \cap \mathfrak{u}$ sont alors les $x_i -y_j$ lorsque la case de coordonn\'ees $(i,j)$ dans $p\times q$
appartient \`a $\lambda$ et les $y_j - x_i$ lorsque la case de coordonn\'ees $(i,j)$ n'appartient pas \`a $\mu$. Celles-ci correspondent respectivement aux vecteurs
$(e_i ; -f_j )$ et $(-e_i ; f_j)$ dans $i\mathfrak{t}_0$. Le c\^one ferm\'e ${\Bbb R}^+ \langle \mathfrak{u} \cap \mathfrak{p} \rangle$ intersecte donc non trivialement le 
sous-espace $i\mathfrak{t}_{0-}$ si et seulement s'il existe un couple d'entiers $(i,j)$ dans $[1,p] \times [1,r]$ tel que 
\begin{itemize}
\item la case de coordonn\'ees $(i,j)$ dans $p\times q$ appartienne au diagramme $\lambda$, et 
\item la case de coordonn\'ees $(i, q-j+1)$ dans $p\times q$ n'appartienne pas au diagramme $\mu$.
\end{itemize}
Remarquons que si c'est le cas pour un couple $(i,j)$, c'est \'egalement le cas pour le couple $(i,1)$. Finalement, le Th\'eor\`eme \ref{kobayashi} implique 
le th\'eor\`eme suivant.

\begin{thm} \label{kobaU}
Soient $p,q,r$ des entiers naturels avec $2r \leq q$ et $H= U(p,q-r) \times U(r)$ plong\'e de mani\`ere standard dans $G=U(p,q)$. 
Soit $(\lambda , \mu )$ un couple compatible de partitions dans $p\times q$. Alors, la restriction $A(\lambda , \mu )_{|K^H}$ de la
repr\'esentation cohomologique $A(\lambda , \mu)$ de $G$ au compact maximal $K^H= U(p) \times U(q-r) \times U(r)$ de $H$ est
$K^H$-admissible si et seulement si deux \'el\'ements de $\lambda$ et $p\times q /\mu $ dans $p\times q$ ne sont jamais align\'es
(autrement dit, $\lambda_i (q-\mu_i)=0$ pour tout $i=1,\ldots ,p$).
\end{thm}

\medskip

\paragraph{Le cas orthogonal}
Nous supposons donc maintenant $(G,H) = (SO_0 (p,q) , SO_0 (p,q-r) \times SO(r))$ avec $2r \leq q$, o\`u le plongement $SO_0 (p,q-r) \times SO(r) \rightarrow
SO_0 (p,q)$ est standard donn\'e par~:
$$(A,B) \mapsto \left(
\begin{array}{cc}
A & 0 \\ 
0 & B 
\end{array} \right) .$$
Comme dans le cas unitaire, quitte \`a conjuguer $H$ dans $G$, on peut supposer que le sous-espace $\mathfrak{t}_0$ de $\mathfrak{g}_0$ est constitu\'e des matrices
\begin{eqnarray*}
\left(
\begin{array}{cccccc}
\begin{array}{cc}
0 & x_1 \\
-x_1 & 0 
\end{array} & & & & &\\
& \begin{array}{cc}
0 & x_2 \\
-x_2 & 0 
\end{array} &&&& \\
&& \ddots &&& \\
&&& \ddots && \\
&&&& \begin{array}{cc}
0 & y_{s-1} \\
-y_{s-1} & 0 
\end{array} & \\
&&&&& \begin{array}{cc}
0 & y_s \\
-y_s & 0 
\end{array}  
\end{array} \right)
\end{eqnarray*}
(o\`u $t = \left[ \frac{p}{2} \right]$, $s= \left[ \frac{q}{2} \right]$ et $x_1 , \ldots , x_t$ et $y_1 , \ldots , y_s$ sont r\'eels)
est une une sous-alg\`ebre de Cartan $\tau$-stable de $\mathfrak{k}_0$ telle que $\mathfrak{t}_{0 -} = \mathfrak{t}_0 \cap \mathfrak{k}_{0-}$
soit un sous-espace ab\'elien maximal de $\mathfrak{k}_{0-}$. Si nous notons comme d'habitude $(x_1 , \ldots , x_t ; y_1 , \ldots , y_s)$ les \'el\'ements de $i\mathfrak{t}_0$, le
sous-espace $i\mathfrak{t}_{0-}$ est constitu\'e des \'el\'ements 
$$(0, \ldots , 0 ; 0 , \ldots , 0 , u_r , \ldots , u_1 ).$$

Consid\'erons maintenant $\mathfrak{q} = \mathfrak{q} (\lambda )_{\pm_1}^{\pm_2}$ la sous-alg\`ebre parabolique $\theta$-stable associ\'ee \`a une partition orthogonale
$\lambda $. Notons toujours
$$z_1 , \ldots , z_p \  (\mbox{resp. } w_1 , \ldots ,w_q )$$
les r\'eels $x_i$, $-x_i$ et $0$ (resp. $y_j$, $-y_j$ et $0$) rang\'es par ordre d\'ecroissant (resp. croissant). 
Les racines de $\mathfrak{t}$ dans $\mathfrak{p} \cap \mathfrak{u}$ sont alors les $z_i - w_j$ et les $w_{q-j+1} - z_{p-j+1}$ o\`u
la case de coordonn\'ees $(i,j)$ dans $p\times q$
appartient \`a $\lambda$. Comme dans le cas unitaire, il d\'ecoule de tout ceci
que le c\^one ferm\'e ${\Bbb R}^+ \langle \mathfrak{u} \cap \mathfrak{p} \rangle$ intersecte donc non trivialement le
sous-espace $i\mathfrak{t}_{0-}$ si et seulement s'il existe un couple d'entiers $(i,j)$ dans $[1,p] \times [1,r]$ tel que 
les cases de coordonn\'ees $(i,j)$ et $(p-i+1 , j)$ dans $p\times q$ appartiennent toutes deux au diagramme $\lambda$. 
Remarquons que si c'est le cas pour un couple $(i,j)$, c'est \'egalement le cas pour le couple $(i,1)$ et donc pour le couple $(t,1)$. 
Finalement, le Th\'eor\`eme \ref{kobayashi} implique le th\'eor\`eme suivant.

\begin{thm} \label{kobaO}
Soient $p,q,r$ des entiers naturels avec $2r \leq q$ et $H= SO_0 (p,q-r) \times SO(r)$ plong\'e de mani\`ere standard dans $G=SO_0 (p,q)$. 
Soit $\lambda$ une partition orthogonale dans $p\times q$. Alors, la restriction $A(\lambda )_{\pm_1 \; |K^H}^{\pm_2}$ de la 
repr\'esentation cohomologique $A(\lambda )_{\pm_1}^{\pm_2}$ de $G$ au compact maximal $K^H= SO(p) \times SO(q-r) \times SO(r)$ de $H$ est
$K^H$-admissible si et seulement si $\lambda \subset [p/2] \times q$.
\end{thm}

\medskip

\paragraph{Applications}
Nous allons maintenant d\'eduire du Th\'eor\`eme \ref{kobaU} le nouveau cas suivant de la Conjecture \ref{CU1}.

\begin{thm}\label{res cohomholomU}
Soient $H=U (p,q-r) \subset U (p,q)= G$ o\`u l'inclusion est l'inclusion standard, $1\leq p ,q$ et $1 \leq r <q$. Alors, pour  
$i, j$ entiers naturels $ \leq q-r$ tels que $i+j \geq q$ et pour tout entier naturel $k \leq p$,
la repr\'esentation cohomologique $A(((i^k) , ((q-r)^k , (q-r-j)^{p-k})))_H$ de $H$ est \'equivalente \`a une sous-repr\'esentation
irr\'eductible de la repr\'esentation cohomologique $A(((i^k), (q^k , (q-j)^{p-k})))$ de $G$.
\end{thm}
{\it D\'emonstration.} Fixons $i,j,k$ comme dans l'\'enonc\'e du Th\'eor\`eme et consid\'erons la repr\'esentation cohomologique
$A(((i^k), (q^k , (q-j)^{p-k})))$ de $G$. D'apr\`es le Th\'eor\`eme \ref{kobaU}, 
le $(\mathfrak{g} , K)$-module associ\'e (toujours not\'e)
$A(((i^k), (q^k , (q-j)^{p-k})))$ est discr\`etement d\'ecomposable comme $(\mathfrak{h} , K^H)$-module. \'Ecrivons donc
\begin{eqnarray} \label{decoA}
A(((i^k), (q^k , (q-j)^{p-k}))) \cong \bigoplus_{\pi \in \widehat H} \underbrace{\pi \oplus \ldots \oplus \pi}_{{\rm m} (\pi )}
\end{eqnarray}
comme $(\mathfrak{h},K^H)$-module. Pour tout $\pi \in \widehat H$ et tout entier $l$ ($1\leq l \leq {\rm m} (\pi)$), notons
$${\rm pr}_{\pi}^{(l)} : A(((i^k), (q^k , (q-j)^{p-k}))) \rightarrow \pi $$
la projection sur la $l$-i\`eme composante de $\pi$ dans la somme directe (\ref{decoA}), et
$$e_{\pi}^{(l)} : \pi \rightarrow A(((i^k), (q^k , (q-j)^{p-k})))$$
l'injection dans la $l$-i\`eme composante de $\pi$. Les applications ${\rm pr}_{\pi}^{(l)}$ et $e_{\pi}^{(l)}$ sont toutes
deux des $(\mathfrak{h}, K^H )$-morphismes. Remarquons que si $U \subset A(((i^k), (q^k , (q-j)^{p-k})))$ 
est un sous-espace $K$-invariant de
dimension finie, alors ${\rm pr}_{\pi}^{(l)} (U)=0$ sauf peut \^etre pour un nombre fini de $\pi \in \widehat H$ (voir
\cite[Proposition 1.6(1)]{Kobayashi}). En particulier, si $\psi \in {\rm Hom}_K (\bigwedge^*\mathfrak{p} , 
A(((i^k), (q^k , (q-j)^{p-k}))) )$, alors le membre de droite de
$$\psi = \sum_{\pi \in \widehat H} \sum_{l=1}^{{\rm m} (\pi)} e_{\pi}^{(l)} \circ {\rm pr}_{\pi}^{(l)} \circ \psi$$
est une somme finie puisque ${\rm dim}  \bigwedge^* \mathfrak{p} < \infty$. On a donc un isomorphisme d'alg\`ebres
gradu\'ees
$$H^* (\mathfrak{h} , K^H ; A(((i^k), (q^k , (q-j)^{p-k})))) 
\cong \bigoplus_{\pi \in \widehat H} {\rm m} (\pi ) H^* (\mathfrak{h} , K^H ; \pi ).$$
D'un autre c\^ot\'e les m\'ethodes de \cite{BergeronTentative} (notamment la d\'emonstration de la Proposition 11) impliquent 
que l'application naturelle de restriction
$$H^* (\mathfrak{g} , K ; A(((i^k), (q^k , (q-j)^{p-k})))) \rightarrow H^* (\mathfrak{h} , K^H ; A(((i^k), (q^k , (q-j)^{p-k}))))$$
est injective. 
L'un des $(\mathfrak{h} , K^H )$-modules $\pi$ de multiplicit\'e ${\rm m}(\pi ) \neq 0$ dans $A(((i^k), ((q-r)^k , (q-r-j)^{p-k})))$
doit donc \^etre cohomologique. Un examen des $K$-types montre finalement que $\pi = A((i^k), ((q-r)^k , 
(q-r-j)^{p-k} )))_H$.~$\Box$

\bigskip

De la m\^eme mani\`ere on peut d\'eduire du Th\'eor\`eme \ref{kobaO} le nouveau cas suivant de la Conjecture \ref{CU2}.

\begin{thm} \label{res cohomO}
Soient $H=SO_0 (2p,q-r) \subset SO_0 (2p,q)= G$ o\`u l'inclusion est l'inclusion standard, $1\leq p ,q$ et $1 \leq r <q$. 
Alors, pour tout entier $q/2 \leq i \leq q-r$, la repr\'esentation cohomologique $A((i^p))^{\pm}_H$ de $H$ est \'equivalente \`a une 
sous-repr\'esentation irr\'eductible de la repr\'esentation $A((i^p))$ de $G$.
\end{thm}

Les r\'esultats correspondants aux Th\'eor\`emes \ref{kobaO} et \ref{kobaU} pour le cup-produit pourraient \'evidemment
\^etre obtenus de la m\^eme mani\`ere.

\section{G\'eom\'etrie de l'espace sym\'etrique associ\'e au groupe $O(p,q)$}

Dans cette section nous \'etudions la g\'eom\'etrie de l'espace sym\'etrique associ\'e au groupe $O(p,q)$.
Nous suivons essentiellement le m\^eme plan que dans \cite{BergeronClozel} pour l'\'etude du cas du groupe $U(p,q)$.
Nous avons grandement profit\'e des travaux de Wang \cite{Wang} sur des probl\`emes proches.

\subsection{Pr\'eliminaires}

Soient $p$, $q$ et $r$ trois entiers strictement positifs. Dans cette section $G= O(q+r,p)$, $K=O(q+r) \times O(p)$ et
$X_{p,q+r} = G/K$, l'espace sym\'etrique associ\'e. Nous aurons besoin d'un mod\`ele pour $X_{p,q+r}$. Posons donc
$$X_{p,q+r} = \left\{ Z \in M_{q+r, p} ({\Bbb R}) \; : \; {}^t ZZ < 1_{p} \right\}.$$
\'Etant donn\'e un \'el\'ement $g \in G$, on peut \'ecrire
$$g= \left(
\begin{array}{cc}
A & B \\
C & D
\end{array} \right)$$
o\`u $A \in M_{q+r,q+r} ({\Bbb R})$, $B\in M_{q+r,p} ({\Bbb R})$, $C \in M_{p,q+r} ({\Bbb R})$ et $D \in M_{p,p} ({\Bbb R})$.
Remarquons que le fait que $g \in G$ \'equivaut \`a ce que
\begin{eqnarray} \label{ginv}
g^{-1} = \left(
\begin{array}{cc}
{}^t A & -{}^t C \\
-{}^t B & {}^t D
\end{array} \right).
\end{eqnarray}
L'action de $g$ sur $X_{q+r,p}$ est donn\'ee par
\begin{eqnarray} \label{gaction}
gZ = (AZ+B)(CZ+D)^{-1} .
\end{eqnarray}
Le groupe $G$ agit transitivement sur $X_{p,q+r}$ et le groupe d'isotropie du point $Z=0$ est $K$. Sur $X_{p,q+r}$ on a
une m\'etrique riemannienne $G$-invariante d\'efinie par
\begin{eqnarray} \label{metrique}
{\rm tr} \left( (1_{q+r} -Z {}^tZ)^{-1} dZ (1_{p} - {}^t ZZ)^{-1} d {}^t Z \right) .
\end{eqnarray}
Cette m\'etrique est la m\'etrique sym\'etrique induite par la forme de Killing de $G$. Nous en donnons \'egalement
la description suivante.

Notons toujours $\mathfrak{g}_0$, $\mathfrak{k}_0$ les alg\`ebres de Lie respectives de $G$ et $K$ et
$\mathfrak{p}_0$ le suppl\'ementaire orthogonal de $\mathfrak{k}_0$ dans $\mathfrak{g}_0$
par rapport \`a la forme de Killing.
\'Etant donn\'ee une matrice $Z \in M_{q+r,p} ({\Bbb R})$, notons
$$\xi (Z) = \left(
\begin{array}{cc}
0 & Z \\
{}^t Z & 0
\end{array}
\right) .$$
Rappelons alors que $\mathfrak{p}_0 = \{ \xi (Z) \; : \; Z\in M_{q+r,p} ({\Bbb R}) \}$. La forme de
Killing induit sur $\mathfrak{p}_0$ le produit scalaire
${\rm tr}(Z {}^t W )$.

Nous identifions $\mathfrak{p}_0$ avec l'espace tangent
$T_0 (X_{p,q+r} )$ \`a $X_{p,q+r}$ au point $Z=0$. Pour $Z \in M_{q+r,p} ({\Bbb R})$, soit $\tau_t$ la courbe
$\tau_t = (\exp t\xi (Z) ) 0$. L'image de $\xi (Z)$ dans $T_0 (X_{p,q+r} )$ est le vecteur tangent $\dot{\tau}_0$
\`a la courbe $\tau_t$ en $t=0$. Sous cette identification, la m\'etrique riemannienne $g$ de $X_{p,q+r}$ est
induite par la forme de Killing~:
$$g_0 (\xi (Z) , \xi (W)) = {\rm tr}(Z {}^t W )
 .$$

\subsection{Sous-espaces totalement g\'eod\'esiques}

Si $v\in {\Bbb R}^n$, o\`u $n=p+q+r$, nous d\'ecomposons $v$ en
$$v= \left(
\begin{array}{c}
v_+ \\
v_-
\end{array} \right), \; v_+ \in {\Bbb R}^{q+r} , \; v_- \in {\Bbb R}^{p} .$$
Soit $g = \left(
\begin{array}{cc}
A & B \\
C & D 
\end{array} \right) \in G$ et $Z \in X_{p,q+r}$, on introduit les facteurs d'automorphie~:
\begin{eqnarray}
J(g,Z) & = & \left( 
\begin{array}{cc}
l(g,Z) & 0 \\
0 & j(g,Z ) 
\end{array}
\right), \\
j(g,Z)  & = & CZ+D, \\
l(g,Z) & = & A -(gZ)C .
\end{eqnarray}
L'action de $g$ sur $X_{p,q+r}$ peut alors prendre la forme suivante~:
$$g \left( \begin{array}{c} 
Z \\
1_{p}
\end{array} \right) = \left( \begin{array}{c}
gZ \\
1_{p} 
\end{array} \right) j(g,Z ) , \; Z\in X_{p,q+r} .$$

\medskip

Dans la suite, $n=p+q+r$ et $Q$ est la forme quadratique sur ${\Bbb R}^n$ de matrice (elle aussi not\'ee $Q$)~: 
$\left( \begin{array}{cc}
1_{q+r} & 0 \\
0 & -1_{p} 
\end{array} \right)$.

Soit $V$ un sous-espace de ${\Bbb R}^n$ de dimension $r$ et positif par rapport \`a $Q$. On associe \`a un tel espace 
un sous-groupe $G_V$ de $G$ et une sous-vari\'et\'e $X_V$ de $X=X_{p,q+r}$ d\'efinis par~:
$$G_V = \{ g \in G \; : \; g \mbox{ laisse invariant le sous-espace } V \},$$
$$X_V = \{ Z \in X \; : \; {}^t Z v_+ = v_- \mbox{ pour tout } v\in V \}.$$

\begin{lem} \label{sous-espace}
La sous-vari\'et\'e $X_V$ et le sous-groupe $G_V$ ont les propri\'et\'es suivantes.
\begin{enumerate}
\item Pour tout $g\in G$, $g X_V = X_{gV}$.
\item Le groupe $G_V$ agit transitivement sur $X_V$.
\item La sous-vari\'et\'e $X_V$ est un sous-espace sym\'etrique totalement g\'eod\'esique de dimension 
$pq$. En tant qu'espace sym\'etrique $X_V$ est isomorphe \`a $X_{p,q}$.
\end{enumerate}
\end{lem}
{\it D\'emonstration.} Il d\'ecoule facilement des d\'efinitions que 
$$X_V = \left\{ Z \in X \; : \; {}^t v Q \left( 
\begin{array}{c}
Z \\
1_p 
\end{array}
\right) = 0, \mbox{ pour tout } v \in V \right\}.$$
Alors si $g\in G$, $Z \in X_V$ et $v \in V$, on v\'erifie facilement que 
$$0 = {}^t v  Q \left(
\begin{array}{c}
Z \\
1_p 
\end{array} \right) = {}^t v  {}^t g  Q g \left(
\begin{array}{c}
Z \\
1_p 
\end{array} \right) = {}^t gv Q \left( 
\begin{array}{c}
gZ \\
1_p
\end{array} \right) j(g,Z).$$
Donc ${}^t v  {}^t g  Q g \left(
\begin{array}{c}
Z \\
1_p 
\end{array} \right) = {}^t gv Q \left( 
\begin{array}{c}
gZ \\
1_p
\end{array} \right) =0$ et le premier point est d\'emontr\'e.

Puis, il d\'ecoule du Th\'eor\`eme de Witt et du premier point que l'on peut supposer $V={\Bbb R}^r$. Dans 
ce cas,
$$X_V = \left\{ \left( 
\begin{array}{c}
0 \\
W
\end{array} \right) \; : \; W \in M_{q,p} ({\Bbb R}), \; {}^t W W < 1_p \right\}$$
et 
$$G_V = \left\{ \left( 
\begin{array}{cc}
u & 0 \\
0 & h 
\end{array} \right) \in G \; : \; h \in O(q,p), \; u\in O(r) \right\} .$$
Et les points 2. et 3. du Lemme \ref{sous-espace} s'en d\'eduisent facilement.~$\Box$

\bigskip

Soit $e_1 , \ldots , e_n$ la base canonique de ${\Bbb R}^n$. Et soit $V$ le sous-espace 
engendr\'e par $e_{q+1} , \ldots , e_{q+r}$. Nous \'etudions maintenant la fonction distance $d(Z, X_V )$ d'un \'el\'ement $Z \in X$
\`a $X_V$. \'Etant donn\'e $Z \in X$, nous d\'ecomposons $Z$ en 
$$Z = \left( 
\begin{array}{c}
Z_1 \\
Z_2
\end{array} \right),$$
o\`u $Z_1 \in M_{q,p} ({\Bbb R})$ et $Z_2 \in M_{r,p} ({\Bbb R})$. Le sous-espace $X_V$ est alors 
donn\'e par 
$$X_V = \{ Z \in X \; : \; Z_2 =0 \}.$$
Un \'el\'ement $g\in G_V$ s'\'ecrit comme matrice par blocs
$$g = \left( 
\begin{array}{ccc}
A_1 & 0 & B_1 \\
0 & u & 0  \\
C_1 & 0 & D_1  
\end{array} \right).$$
Notons $X_1$ l'espace 
$$X_1 = \{ W \in M_{q,p} ({\Bbb R}) \; : \; {}^t W W < 1_p \} .$$
Le groupe $G_V$ agit transitivement sur $X_1$ par~:
$$gW
= (A_1 W + B_1 ) (C_1 W + D_1 )^{-1}, \; g\in G_V, \; W \in X_1.$$
L'action de $G_V$ sur $X$ s'\'ecrit~:
\begin{eqnarray} \label{action}
gZ = \left(
\begin{array}{c}
gZ_1 \\
uZ_2 j(g,Z)^{-1} 
\end{array} \right) , \; g \in G_V , \; Z \in X .
\end{eqnarray}
D'apr\`es (\ref{action}), il existe un \'el\'ement $g\in G_V$ tel que $gZ=Z'$ avec $Z_1 ' =0$. La m\'etrique 
riemannienne de $X$ \'etant $G$-invariante, $d(Z, X_V )= d(Z' , X_V) = d(0,Z')$. Il nous 
suffit donc d'\'etudier la fonction distance $d(0,Z)$ de $0$ \`a $Z$.

\begin{lem} \label{distance}
Soit $Z \in X$, $d= d(0,Z)$ et $m$ le rang de $Z {}^t Z$. Alors,
\begin{enumerate}
\item si $m=1$, $\cosh^2 d = (\det (1_{q+r} -Z{}^t Z ))^{-1}$,
\item et en g\'en\'eral,
$$\frac{1}{2^m} e^d \leq (\det (1_{q+r} -Z {}^t Z ))^{-1} \leq e^{\sqrt{m} d}.$$
\end{enumerate}
\end{lem}
{\it D\'emonstration.} Il existe un unique $Y \in M_{q+r,p} ({\Bbb R})$ v\'erifiant 
\begin{eqnarray} \label{expo}
\exp (\xi (Y)) 0 = Z .
\end{eqnarray}
La courbe $\exp (t \xi (Y))0$, $0 \leq t \leq 1$, est une g\'eod\'esique joignant $0$ \`a $Z$. On a donc
$$d^2 = {\rm tr} (Y {}^t Y  ).$$
Soit $A$ (resp. $B$) une matrice hermitienne positive v\'erifiant 
$$A^2 = Y {}^t Y  \  ({\rm resp.} \  B^2 = {}^t Y Y ) .$$
Il d\'ecoule des d\'efinitions que 
$$\exp (\xi (Y)) = \left(
\begin{array}{cc}
\cosh A & \sum_{k=0}^{\infty} \frac{A^{2k}}{(2k+1)!} Y \\
\sum_{k=0}^{\infty} \frac{B^{2k}}{(2k+1)!} {}^t Y & \cosh B
\end{array} \right) .$$
Puisque $A^{2k} Y= Y B^{2k}$, on d\'eduit de l'expression ci-dessus et de (\ref{expo}) que 
$Z {}^t Z = \tanh^2 (A)$ et donc que~:
\begin{eqnarray} \label{expA}
e^A = \frac{1_{q+r} + \sqrt{Z {}^t Z }}{(1_{q+r} -Z{}^t Z )^{1/2}} .
\end{eqnarray}
Il d\'ecoule facilement du fait que $A$ est de rang $m$ que
$$d \leq {\rm tr} (A) \leq \sqrt{m} d .$$
Puisque $Z{}^t Z < 1_{q+r}$, si l'on applique le d\'eterminant \`a (\ref{expA}), on obtient le point 2. du Lemme. Si
$m=1$, $d={\rm tr} (A)$ et le point 1. d\'ecoule encore de (\ref{expA}).~$\Box$

\bigskip

\begin{lem} \label{h}
Soient $h_Z = (1_p - {}^t Z_1 Z_1 )^{-1}$ et $\tilde{h}_Z =(1_p - {}^t Z Z )^{-1}$. Alors, 
$$\tilde{h}_{gZ} = j(g,Z) \tilde{h}_Z  {}^t j(g,Z) , \mbox{ pour tout } g\in G,$$
$$h_{gZ} = j(g,Z) h_Z {}^t j(g,Z) , \mbox{ pour tout } g\in G_V .$$
\end{lem}
{\it D\'emonstration.} Soit $l_Z =(1_p - {}^t Z Z )$. On a~:
\begin{eqnarray*}
l_Z & = & -( {}^t Z \   1_p) Q \left( 
\begin{array}{c}
Z \\
1_{p} 
\end{array} \right) \\
& = & -({}^t Z \   1_p) {}^t g Q g \left( 
\begin{array}{c}
Z \\
1_{p} 
\end{array} \right) \\
& = & {}^t j(g,Z) l_{gZ} j(g,Z) .
\end{eqnarray*}
Puisque $\tilde{h}_Z = l_Z^{-1}$, on obtient la premi\`ere propri\'et\'e annonc\'ee. Mais 
$h_Z = \tilde{h}_{\left(
\begin{array}{c}
Z_1 \\
0
\end{array} \right) }$ et pour tout $g\in G_V$, $j(g,Z) = j \left( g, \left( 
\begin{array}{c}
Z_1 \\
0
\end{array} \right) \right)$, d'o\`u la seconde propri\'et\'e annonc\'ee.~$\Box$

\bigskip

Pour $Z \in X$, on introduit les fonctions $A$ et $B$ sur $X$ d\'efinies par 
\begin{eqnarray} \label{A et B}
A & = & \mbox{det} (1_p - {}^t Z Z ) \\
B & = & \mbox{det} (1_p - {}^t Z_1 Z_1 ).
\end{eqnarray}
La fonction $B$ est obtenue en restreignant la fonction $A$ \`a $X_V$, puis en l'\'etendant \`a $X$
tout entier de fa\c{c}on constante dans la direction de $Z_2$.

\begin{lem} \label{B sur A}
La fonction $\frac{B}{A}$ est $G_V$-invariante.
\end{lem}
{\it D\'emonstration.} Cela d\'ecoule imm\'ediatement du Lemme \ref{h}.~$\Box$

\bigskip

Nous pouvons maintenant estimer la fonction $d(Z, X_V )$.

\begin{prop} \label{distance2}
Soient $Z \in X$ et $m$ le rang de la matrice $Z_2 {}^t Z_2$.
\begin{enumerate}
\item Si $m=1$, $(\cosh d(Z, X_V ))^2 = \frac{B}{A}$.
\item En g\'en\'eral, on a~: 
$$4^m \left( \frac{B}{A} \right) \geq e^{2 d(Z, X_V )} , $$
et 
$$e^{2 \sqrt{m} d(Z, X_V )} \geq \frac{B}{A} .$$
\end{enumerate}
\end{prop}
{\it D\'emonstration.} Les fonctions $\frac{B}{A}$ et $d(., X_V )$ sont toutes deux $G_V$-invariantes. On a 
vu, cf. (\ref{action}), que l'on pouvait se ramener \`a ce que $Z_1 =0$ et donc $d(Z, X_V )=d(0,Z_2) $.
Mais alors, $\frac{B}{A} = (\mbox{det}(1_p - Z_2 {}^t Z_2 ))^{-1}$. Et la Proposition d\'ecoule alors 
du Lemme \ref{distance}.~$\Box$

\bigskip

Nous aurons \'egalement besoin dans la suite des expressions suivantes.

\begin{lem} \label{outils}
On a l'\'egalit\'e 
$$\frac{B}{A} = \det \{ 1_r + Z_2 (1_p - {}^t Z Z )^{-1} {}^t Z_2 \} $$
et l'in\'egalit\'e
$$1+ r^{-1} {\rm tr} \left( Z_2 (1_p - {}^t Z Z^{-1} {}^t Z_2 \right) \geq \left( \frac{B}{A} \right)^{1/r} .$$
\end{lem}
{\it D\'emonstration.} Remarquons d'abord que 
\begin{eqnarray*}
1_p - {}^t Z_1 Z_1 & = & 1_p - {}^t Z Z + {}^t Z_2 Z_2 \\
& = & (1_p - {}^t Z Z )^{1/2} \left\{ 1_p + (1_p - {}^t Z Z )^{-1/2} {}^t Z_2 Z_2 (1_p - {}^t Z Z )^{-1/2} 
\right\} (1_p - {}^t Z Z )^{1/2} .
\end{eqnarray*}
On en d\'eduit que 
\begin{eqnarray*}
\frac{B}{A} & = & \det \left\{ 1_p + (1_p - {}^t Z Z )^{-1/2} {}^t Z_2 Z_2 (1_p - {}^t Z Z )^{-1/2} \right\} \\
                 & = & \det \left\{ 1_r + Z_2 (1_p -{}^t Z Z )^{-1} {}^t Z_2 \right\} .
\end{eqnarray*}
Ce qui d\'emontre la premi\`ere partie du Lemme. Remarquons maintenant que la matrice $Z_2 (1_p - {}^t Z Z )^{-1} {}^t Z_2$
est positive. Notons $\lambda_1 , \ldots , \lambda_r$ ses valeurs propres (r\'eelles positives). Alors,
\begin{eqnarray*}
1+\frac{1}{r} {\rm tr} \left( Z_2 (1_p - {}^t Z Z )^{-1} {}^t Z_2 \right) & = &  
\frac{(1+\lambda_1 ) + \ldots + (1+\lambda_r )}{r} \\
& \geq & \left\{ (1+\lambda_1 ) \ldots (1+\lambda_r ) \right\}^{1/r} \\
& = & \det \left\{ 1_r + Z_2 (1_p - {}^t Z Z )^{-1} {}^t Z_2 \right\}^{1/r} \\
& = & \left( \frac{B}{A} \right)^{1/r} .
\end{eqnarray*}
$\Box$

\subsection{Croissance du volume}

\'Examinons maintenant les champs de Jacobi \'emanant de $X_V$. (Une r\'ef\'erence g\'en\'erale pour les champs
de Jacobi est \cite{Sakai}.)

\begin{lem} \label{champs de Jacobi}
Soient $Z\in X_V$, $T_Z (X_V )$ l'espace tangent \`a $X_V$ en $Z$ et $T_Z (X_V )^{\perp}$
le suppl\'ementaire orthogonal de $T_Z (X_V )$ dans $T_Z (X)$. Soit $Y$ un vecteur dans 
$T_Z (X_V )^{\perp}$ avec $g_Z (Y,Y)=1$. Alors, 
\begin{enumerate}
\item les espaces $T_Z (X_V )$ et $T_Z (X_V )^{\perp}$ sont invariants sous 
l'application $R(.,Y)Y$ (o\`u $R$ d\'esigne le tenseur de courbure de $X$), 
\item il existe $\lambda_1 \geq \lambda_2 \geq \ldots \geq \lambda_l \geq 0$, $l=\max \{ r,p \}$ tels que 
\begin{itemize}
\item $\lambda_1^2 + \ldots + \lambda_l^2 =1$,
\item $\lambda_i =0$, si $i>\min \{ r, p \}$,
\item l'op\'erateur $R(.,Y)Y_{|{T_Z (X_V )}}$ a pour valeurs propres 
$$\begin{array}{cccc}
\underbrace{-\lambda_1^2 , \ldots , -\lambda_1^2} , & \underbrace{-\lambda_2^2 , \ldots , -\lambda_2^2 }, &  \ldots  &, \underbrace{-\lambda_p^2 , \ldots , - \lambda_p^2 } ,\\
q & q & & q 
\end{array}$$
\item l'op\'erateur $R(.,Y)Y_{|{T_Z (X_V )^{\perp} }}$ a pour valeurs propres 
$$-(\lambda_i -\lambda_j )^2 , \; 1\leq i \leq r , \; 1\leq j \leq p .$$
\end{itemize}
\end{enumerate}
\end{lem}
{\it D\'emonstration}. D'apr\`es (\ref{action}), on peut supposer que $Z=0$ et
$Y=\xi \left(
         \begin{array}{c}
          0 \\
          M
          \end{array}  \right) $. 
D'apr\`es \cite[Theorem 3.2, Chap. XI]{KobayashiNomizu}, \'etant donn\'e $X \in T_0 (X )$, le tenseur de 
courbure est donn\'e par~:
\begin{eqnarray} \label{courbure}
R(X,Y)Y= -[[X,Y],Y].
\end{eqnarray}
Si $X= \xi \left( 
\begin{array}{c}
N \\
0
\end{array} \right) \in T_0 (X_V )$,
un calcul simple donne alors
\begin{eqnarray} \label{cas1}
R(X,Y)Y = \xi \left(
                  \begin{array}{c}
                       -N {}^t M M \\
                        0
                  \end{array} \right) .
\end{eqnarray}
Si maintenant $X= \xi \left(
                          \begin{array}{c}
                           0 \\
                           L
                           \end{array} \right) \in T_0 (X_V )^{\perp} $, un autre calcul simple 
\`a l'aide de (\ref{courbure}) donne
\begin{eqnarray} \label{cas2}
R(X,Y)Y = \xi \left(
                  \begin{array}{c}
                  0 \\
                  -L {}^t M M +2M {}^t L M - M {}^t M L
                   \end{array}  \right) .
\end{eqnarray}
Il d\'ecoule de (\ref{cas1}) et (\ref{cas2}) que les espaces $T_0 (X_V )$ et $T_0 (X_V )^{\perp}$ sont invariants
sous l'application $R(.,Y)Y$. Remarquons que l'on peut toujours supposer que $M$ est 
de la forme 
\begin{eqnarray*}
M & = & \left( 
\begin{array}{c}
\begin{array}{ccc}
\lambda_1 & \ldots & 0 \\
\vdots & \ddots & \vdots \\
0 & \ldots & \lambda_p 
\end{array} \\
0 
\end{array} \right) \; \mbox{ si } r>p, \\
M & = & \left( 
\begin{array}{ccc}
\lambda_1 & \ldots & 0 \\
\vdots & \ddots & \vdots \\
0 & \ldots & \lambda_p 
\end{array} \right) \; \mbox{ si } r=p, \\
M & = & \left( 
\begin{array}{cc}
\begin{array}{ccc}
\lambda_1 & \ldots & 0 \\
\vdots & \ddots & \vdots \\
0 & \ldots & \lambda_r 
\end{array} & 0 
\end{array} \right) \; \mbox{ si } r<p,
\end{eqnarray*}
avec $\lambda_1 \geq \ldots \geq \lambda_l \geq 0$, $l=\max \{r,p \}$ et $\lambda_i =0$ si $i> \min \{ r,p \}$.
Puisque $\lambda_1^2 + \ldots + \lambda_l^2 = {\rm tr} (M {}^t M ) =g_0 (Y,Y) = 1$, le deuxi\`eme 
point du Lemme \ref{champs de Jacobi} d\'ecoule alors des formules (\ref{cas1}) et (\ref{cas2}).~$\Box$

\bigskip

Soit $\tau$ une g\'eod\'esique perpendiculaire \`a $X_V$. Nous pouvons maintenant \'etudier les 
champs de Jacobi le long de $\tau$. Soit $\tau =\tau_t$, o\`u $t$ est la longueur d'arc de $X_V$ \`a 
$\tau_t$ et $Y= \dot{\tau}_0$. Dans la suite nous d\'ecrivons les champs de Jacobi $X=X(t)$ le long de $\tau$ v\'erifiant 
\begin{eqnarray} \label{conditions initiales}
X(0) \in T_{\tau_0} (X_V ) \mbox{  et  } \nabla_Y X \in T_{\tau_0} (X_V )^{\perp} .
\end{eqnarray}

L'\'equation de Jacobi est donn\'ee par 
$$\nabla_{\tau_t}^2 X +R(X,\dot{\tau}_t )\dot{\tau}_t =0 .$$

D'apr\`es le Lemme \ref{champs de Jacobi}, les valeurs propres de $R(.,Y)Y$ sont n\'egatives.
Soit $X_0 \in T_{\tau_0} (X_V )$ un vecteur propre de $R(.,Y)Y$ pour la valeur propre 
$-\lambda^2$ ($\lambda \geq 0$). Soit $X_t$ le transport parall\`ele de $X_0$ le long de $\tau$. On peut v\'erifier que 
\begin{eqnarray} \label{chp de Jacobi1}
X(t) & = & (\cosh \lambda t) X_t 
\end{eqnarray}
est un champ de Jacobi le long de $\tau$ v\'erifiant (\ref{conditions initiales}). Soit 
$L_0 \in T_{\tau_0 }(X_V )^{\perp} $ un vecteur propre de $R(.,Y)Y$ pour la valeur propre 
$-\lambda^2$ ($\lambda \geq 0$). Soit $L_t$ le transport parall\`ele de $L_0$ le long de $\tau$. On peut v\'erifier que 
\begin{eqnarray} \label{chp de Jacobi2}
L(t) & = & \left\{ \begin{array}{cc}
                   (\sinh \lambda t) L_t & \mbox{ si } \lambda \neq 0, \\
                   t L_t                 & \mbox{ sinon}
                   \end{array} \right.
\end{eqnarray}
est un champ de Jacobi le long de $\tau$ v\'erifiant (\ref{conditions initiales}). L'espace des champs de Jacobi v\'erifiant 
(\ref{conditions initiales}) a pour dimension $pq$. Cet espace est engendr\'e par les champs de Jacobi construits 
ci-dessus.

\bigskip

L'espace  $X - X_V$ se d\'ecompose en un produit :
$${\cal F} \times ]0, +\infty [ $$
o\`u ${\cal F}$ est l'hypersurface de $X$ constitu\'ee des points \`a distance $1$ de 
$X_V$ et o\`u nous identifions un point $(W,t) \in {\cal F} \times ]0, +\infty [$ avec le point 
$Z \in X$ \`a distance $t$ de $X_V$ et tel que la g\'eod\'esique passant par 
$Z$ et $W$ soit perpendiculaire \`a $X_V$.

Si $S$ est un sous-ensemble mesurable de ${\cal F}$, nous noterons $\omega (t,S)$ le volume
de $\{ Z \in {\cal F} \times \{ t \} \; : \; Z_1 \in S \}$. D'apr\`es (\ref{action}), $\omega (t,S) = \omega (t ,gS)$ pour 
tout $g\in G_V$. On peut donc voir $\omega (t,S)$, pour chaque $t$, comme une mesure invariante sur 
$X_V$. Il existe alors une fonction $f(t)$ telle que $\omega (t,S) = f(t) \mbox{vol} (S)$.
 
\begin{lem} \label{croissance du volume}
Il existe une constante $c$ telle que~:
\begin{enumerate}
\item si $r=1$, 
$$\omega (t,S) = c {\rm vol}(S) (\sinh t)^{p-1} (\cosh t)^{q}   ,$$
\item en g\'en\'eral,
$$\omega (t,S) \leq c {\rm vol}(S) (1+t^{p(q+r)} ) e^{(p+q+r-1)\sqrt{m} t)}  ,$$
o\`u $m=\min \{ r,p\}$.
\end{enumerate}
\end{lem}
{\it D\'emonstration.} Si l'on \'ecrit $\omega (t,S) = f(t) \mbox{vol} (S)$, il nous faut estimer $f(t)$, par exemple en la comparant
\`a la constante $f(1)$. Soit $x_s$ une courbe dans ${\cal F}$. On note $x_s^t$ le point $(x_s , t)$, $x_{(s)}^t$ la 
courbe \`a $s$ fix\'e et $x_s^{(t)}$ la courbe \`a $t$ fix\'e. Les courbes $x_{(s)}^t$ sont des g\'eod\'esiques et 
$\dot{x}_s^{(t)}$ est un champ de Jacobi le long de cette g\'eod\'esique qui v\'erifie (\ref{conditions initiales}).
Mais d'apr\`es le Lemme \ref{champs de Jacobi} et (\ref{chp de Jacobi1}), (\ref{chp de Jacobi2}), 
l'espace $T_{x_s} ({\cal  F})$ admet un base orthonorm\'ee r\'eelle
$$X_1 , \ldots , X_{q} , X_{q+1} , \ldots , X_{2q}, \ldots \ldots , X_{qp} ,$$
$$Y_1 , \ldots , Y_{rp-1}$$
et il existe des r\'eels 
$$\lambda_1 \geq \ldots \geq \lambda_l \geq 0 , \; l=\max \{ r, p\}$$
v\'erifiant~:
\begin{enumerate}
\item $\lambda_1^2 +\ldots + \lambda_l^2 =1$,
\item $\lambda_i =0$ pour tout $i> \min \{ r,p\}$,
\item au point $(x_s , t)$,
\begin{eqnarray*}
\begin{array}{l}
||X_1 ||= \ldots = ||X_{q} || =\frac{\cosh \lambda_1 t}{\cosh \lambda_1} , \\
\  \  \  \vdots \\
||X_{q(p-1)+1}|| = \ldots = ||X_{qp}|| = \frac{\cosh \lambda_p t}{\cosh \lambda_p} ,
\end{array}
\end{eqnarray*}
\item au point $(x_s , t)$, l'ensemble des $||Y_j ||$ pour $1\leq j \leq rp-1$ (compt\'ees avec multiplicit\'es) co\"{\i}ncide,
\`a une permutation pr\`es, avec l'ensemble des $a_{ij}$ pour $1\leq i \leq r$, $1\leq j \leq p$ et $(i,j) \neq (1,1)$ tels que
$$b_{ij} = \left\{
\begin{array}{lc}
\frac{\sinh |\lambda_i -\lambda_j |t}{\sinh |\lambda_i -\lambda_j |}, & \mbox{ si } \lambda_i \neq \lambda_j , \\
t, & \mbox{ si } \lambda_i = \lambda_j . 
\end{array} \right.$$
\end{enumerate}
On en d\'eduit alors facilement que 
$$f(t) = f(1) \frac{\sinh^{p-1} t \cosh^{q} t}{\sinh^{p-1} 1 \cosh^{q} 1} \  \   \mbox{  si  }  r=1$$
et en g\'en\'eral qu'il existe une constante $c$ telle que
$$f(t) \leq c (1+t^{p(q+r)}) e^{(p+q+r-1)\sqrt{m} t} , $$
o\`u $m = \min \{ r,p\}$.~$\Box$

\bigskip

Il r\'esulte du paragraphe pr\'ec\'edent que la fonction $\frac{B}{A}$ est fortement reli\'ee \`a (et plus naturelle que) 
la fonction distance $d(.,X_V )$. Nous \'etudions maintenant la croissance du volume \`a l'aide de la fonction $\frac{B}{A}$.

Notons d'abord que si $g\in G$ et $Z\in X$,
$$d(gZ) = l (g,Z) dZ j(g,Z)^{-1} .$$
On sait que 
$$\det ( l (g,Z )) = \det (j(g,Z))^{-1} .$$
Donc si l'on pose
$$\{ dZ \} = \prod_{i=1}^{q+r} \prod_{j=1}^p dZ_{ij} $$
et si $g\in G$,
$$\{ dgZ \} = \det (j(g,Z))^{-(p+q+r)} \{ dZ \} .$$
Puis d'apr\`es le Lemme \ref{h}, 
$$A (gZ) = \det (j(g,Z))^{-2} A(Z) .$$
La forme volume invariante $dv_X$ sur $X$ s'\'ecrit donc 
\begin{eqnarray} \label{forme volume de D}
dv_X = A^{-(p+q+r)/2} \{ dZ \} .
\end{eqnarray}

Si $\left( 
\begin{array}{c}
Z_1 \\
0
\end{array} \right) \in X_V$, soit $F_{Z_1}$ la fibre au-dessus de ce point dans le fibr\'e $X \rightarrow 
X_V$. Autrement dit
$$F_{Z_1} = \{ Z\in X \; : \; Z_1 \mbox{ fix\'e} \} .$$
Soit $g\in G_V$ l'\'el\'ement 
$$g = \left( 
\begin{array}{ccc}
(1_{q} - Z_1 {}^tZ_1 )^{-1/2} & 0 & -( 1_{q} - Z_1 {}^tZ_1 )^{-1/2} Z_1 \\
0 & 1_r & 0 \\
-(1_p - {}^t Z_1  Z_1 )^{-1/2} {}^t Z_1 & 0 & (1_p - {}^t Z_1 Z_1 )^{-1/2}
\end{array} \right) .$$
Alors $g$ envoie $F_{Z_1}$ isom\'etriquement sur $F_0$, et 
$$ g \left( 
\begin{array}{c}
Z_1 \\
Z_2
\end{array} \right) = \left( 
\begin{array}{c}
0 \\
Z_2 (1_p - {}^tZ_1 Z_1 )^{-1/2} 
\end{array} \right) .$$
Sur $F_0$, l'\'el\'ement de volume est $\det (1_p - {}^t Z_2 Z_2 )^{-(r+p)/2} \{ dZ_2 \}$, l'\'el\'ement de volume
sur $F_{Z_1}$ est donc
$$dv_F =  A^{-r/2} \left( \frac{B}{A} \right)^{p/2} \{ dZ_2 \} .$$
D'o\`u il d\'ecoule que 
\begin{eqnarray} \label{decomposition du volume}
dv_{{\cal D}} = \left( \frac{B}{A} \right)^{q/2} dv_{X_V} dv_F ,
\end{eqnarray}
o\`u $dv_{X_V} = B^{-(p+q)/2} \{ dZ_1 \}$ est la forme volume invariante sur $X_V$. 

\begin{lem} \label{formules d'integration}
On a les formules d'int\'egration~:
\begin{enumerate}
\item
$$\int_{X} A^{s/2} \{ dZ \} = \pi^{p(q+r)/2} \prod_{i=1}^{q+r} \frac{\Gamma ((s+i+1)/2)}{\Gamma ((s+p+i+1)/2)},$$
d\`es que Re$(s) >-2$; et 
\item si $\Gamma_V$ est un sous-groupe discret sans torsion et cocompact dans $G_V$,
$$\int_{\Gamma_V \backslash X} \left( \frac{A}{B} \right)^{s/2} dv_{X} = \pi^{rp/2} 
\prod_{i=1}^r \frac{\Gamma ((s-p-q-r+i+1)/2)}{\Gamma ((s-q-r+i+1)/2)} {\rm vol} (\Gamma_V \backslash X_V ) ,$$
d\`es que ${\rm Re} (s) >p+q+r-2$.
\end{enumerate}
\end{lem}
{\it D\'emonstration.} On introduit tout d'abord $f(s,q+r,p) = \int_{X} A^s \{ dZ \}$. 
On d\'eduit de (\ref{decomposition du volume}), avec $r=1$, la relation de r\'ecurrence
$$f(s,q+r,p) = f(s+1, q+r-1 , p)f(s,1,p) .$$
Mais,
\begin{eqnarray*}
f(s,1,p) & = & \int_{\sum_i x_i^2 \leq 1} (1-(x_1^2 + \ldots +x_p^2 ))^{s/2} dx_1  \ldots dx_p  \\
& = & \frac{\pi^{p/2}}{\Gamma (p/2)} \int_0^1 (1-t)^{s/2} t^{p/2-1} dt \\
& = & \pi^{p/2} \frac{\Gamma (s/2 +1)}{\Gamma (s/2+p/2+1)} \; \mbox{ d\`es que Re}(s)>-2 . 
\end{eqnarray*}
Alors le premier point du Lemme \ref{formules d'integration} d\'ecoule d'une simple r\'ecurrence.

Concernant le deuxi\`eme point, il d\'ecoule de (\ref{decomposition du volume}) que l'int\'egrale vaut
$$\int_{\Gamma_V \backslash X} \left( \frac{A}{B} \right)^{(s-q)/2} dv_F dv_{X_V} .$$
Puisque $\frac{A}{B}$ et $dv_F$ sont $G_V$-invariants, l'int\'egrale
$$\int_{F_{Z_1}} \left( \frac{A}{B} \right)^{(s-q)/2} dv_F$$
est ind\'ependante de $Z_1$. En $Z_1 =0$, sa valeur est 
$$\int_{X_2} \det (1_p - {}^t Z_2 Z_2 )^{(s-p-q-r)/2} \{ dZ_2 \} ,$$
o\`u $X_2 = \{ Z_2 \; : \; {}^t Z_2  Z_2 < 1_p \}$. Le deuxi\`eme point du Lemme \ref{formules d'integration}
d\'ecoule donc directement du premier point.~$\Box$

\medskip

\subsection{Fonction distance \`a l'hypersurface}

Dans l'optique de calculer la cohomologie $L^2$ des quotients $\Gamma_V \backslash X$ nous aurons besoin de d\'eterminer le 
hessien de la fonction distance g\'eod\'esique \`a la sous-vari\'et\'e $X_V$ dans $X$.
Rappelons que le {\it hessien} d'une fonction $C^2$ $F$ de $X$ dans ${\Bbb R}$ est la seconde d\'eriv\'ee 
covariante $\nabla^2 F$ de $F$, {\it i.e.} 
$$\nabla^2 F (U,V) = U(VF) - (\nabla_U V)F ,$$
pour n'importe quels champs de vecteurs $U$, $V$ sur $X$ et o\`u $\nabla$ est la connexion de 
Levi-Civit\`a induite par la structure riemannienne de $X$. Le hessien $\nabla^2 F$ d\'efinit donc  un
tenseur sym\'etrique de type $(0,2)$. Nous appelons  {\it valeurs propres du hessien} les fonctions qui \`a 
chaque point $Z$ de $X$ associent les valeurs propres de la matrice associ\'ee dans n'importe quelle base 
orthonorm\'ee de l'espace tangent \`a $X$ au point $Z$. 

\bigskip

Soit $F$ la fonction distance g\'eod\'esique \`a la sous-vari\'et\'e $X_V$. La fonction $Z \mapsto F(Z)$ est bien \'evidemment
lisse pour $Z\in X-X_V$. 
  
\begin{prop} \label{fonction distance a l'hypersurface}
Supposons $\mathbf{r=1}$. Notons $\{ \gamma_i (Z) \}_{1 \leq i \leq p(q+1)}$ les valeurs propres du hessien $\nabla^{2} F$ en un point $Z \in X$. 
Alors, quitte \`a r\'eordonner les $\gamma_i (Z)$, 
$$\gamma_1 (Z) = \tanh F(Z) , \ldots , \gamma_{q} (Z) = \tanh F(Z) ,$$
$$\gamma_{j} (Z) = 0 \mbox{  pour } j= q+1 , \ldots , pq , $$
et les $\gamma_k (Z)$ pour $pq < k \leq p(q+1)$, sont (\`a permutations pr\`es) 
$0$ et  $\coth F(Z)$ avec multiplicit\'e $p-1$.
\end{prop}
{\it D\'emonstration.} Soit toujours $\tau$ une g\'eod\'esique perpendiculaire \`a $X_V$ avec 
$\tau =\tau_{t}$ o\`u $t$ est la longueur d'arc de $X_V$ \`a $\tau_t$.  Soit $Y= \dot{\tau}$. 
D'apr\`es le Lemme \ref{champs de Jacobi} et puisque $r=1$, il existe  
un champs de bases orthonorm\'ees le long de $\tau$~: $\{ e_{j}, f_{k}  \; : \; 1\leq j \leq pq \mbox{ et } 1\leq k \leq p \}$ tel
que pour tout entier $1 \leq j \leq pq$ le vecteur 
$e_j (0 ) \in T_{\tau_0 } (X_V )$ et soit un vecteur $- 1$-propre si $1 \leq j \leq q$ (resp. $0$-propre si $q+1 \leq j \leq pq$) de $R(.,Y)Y$ et que pour tout entier
$1\leq k \leq p$ le vecteur $f_{k} (0) \in T_{\tau_0} (X_V)^{\perp}$ et soit
un vecteur $0$-propre si $k=1$ (resp. $-1$-propre si $k\geq 2$) de $R(.,Y)Y$.
Nous supposerons de plus (ce que l'on peut bien \'evidemment faire) que le vecteur $Y$ est \'egal au vecteur 
$f_{1}  (0)$. 

Alors d'apr\`es (\ref{chp de Jacobi1}) et pour tout entier $1 \leq j \leq q$ (resp. $ q+1\leq j  \leq pq$),  les champs de vecteurs~:
\begin{eqnarray} \label{base1}
\begin{array}{l}
v_{j} (t) = \cosh t e_{j} (t) \\
({\rm resp. } \   v_j (t) = e_j (t) .)
\end{array}  
\end{eqnarray}
sont des champs de Jacobi le long de $\tau$ v\'erifiant (\ref{conditions initiales}).
Puis, d'apr\`es (\ref{chp de Jacobi2}) et pour tout entier $1 \leq k \leq p$ les champs de vecteurs~:
\begin{eqnarray} \label{base3} 
w_{k}  (t) = \left\{
\begin{array}{ll}
\sinh t f_{k}  (t), &  \mbox{ si } k\geq 2 \\
t f_{k}  (t), & \mbox{ si } k=1  ,\\
\end{array} \right.
\end{eqnarray}
sont des champs de Jacobi le long de $\tau$ v\'erifiant (\ref{conditions initiales}).
De plus, nous avons vu que les champs de vecteurs (\ref{base1}) et (\ref{base3})
forment une base orthogonale de l'espace des champs de Jacobi le long de $\tau$ v\'erifiant (\ref{conditions initiales}). 
La formule de la variation seconde \cite{Sakai} nous dit alors que le Hessien 
$\nabla^2 t (= \nabla^2 F)$ se diagonalise dans la base 
$\{ e_{j} \}_{1\leq j \leq pq} \cup \{ f_{k}  \}_{1\leq k \leq p }$. Et plus pr\'ecisemment permet de calculer par exemple
$$\nabla^2 t(y) (e_{j} ,e_{j} ) = \frac{d^2}{ds^2}_{|s=0} L(\tau^s ) , $$
o\`u si $\tau$ va de $x:= \tau_0 \in X_V$ \`a $y:= \tau_{t(y)}$, $\tau^s$ d\'esigne la g\'eod\'esique minimisante
joignant $X_V$ au point $\exp_y (s e_{j} )$ et $L( \tau^s )$ sa longueur. Or, si $j$ est par exemple compris entre $1$ et 
$q$, le champ de vecteur $\hat{v}_{j} = \frac{v_j}{\sinh t(y)}$ est un champ de Jacobi le long de $\tau$, perpendiculaire \`a $\dot{\tau}$ et 
v\'erifiant~: $\hat{v}_{j} (t(y)) =e_{j} (t(y))$ et (\ref{conditions initiales}). La formule de la variation seconde implique alors~:
\begin{eqnarray*}
\nabla^2 t(y) (e_{j} ,e_{j} ) & = & \langle \nabla \hat{v}_{j} (t(y)), \hat{v}_{j} (t(y)) \rangle , \\
& = & \frac{\sinh t(y)}{\cosh t(y)} .
\end{eqnarray*}
De la m\^eme mani\`ere, si $j$ est un entier v\'erifiant $q+1 \leq j \leq pq$, on obtient~:
$$\nabla^2 t(y) (e_{j} ,e_{j} ) = 0$$
et si $k$ est un entier v\'erifiant $1\leq k \leq p$,
$$
\nabla^2 t(y) (f_{k} ,f_{k} )  = \left\{
\begin{array}{ll} 
\frac{\cosh t(y)}{\sinh t(y)}, & \mbox{ si } k \geq 2,  \\
0, & \mbox{ si } k=1. 
\end{array} \right. $$
Ce qui conclut la d\'emonstration de la Proposition \ref{fonction distance a l'hypersurface}.~$\Box$

\bigskip

Lorsque $r>1$, comme dans les paragraphes pr\'ec\'edents, plut\^ot
que la fonction distance g\'eod\'esique \`a $X_V$ il est plus naturel de consid\'erer la fonction $\log \left( \frac{B}{A} \right)$.

\begin{prop} \label{asymptotique des vp}
Les valeurs propres du hessien $\nabla^2 \log \left( \frac{B}{A} \right)$ en un point $Z \in X$ sont toutes 
positives, inf\'erieures (ou \'egales) \`a $1$, et parmi celles-ci au moins $q+pr-1$ tendent vers $1$ lorsque $\frac{B}{A} (Z)$ tend vers l'infini.
\end{prop}
{\it D\'emonstration.} Nous d\'emontrons par r\'ecurrence sur un entier $k\geq 1$ que si 
$V$ est un sous-espace de ${\Bbb R}^n$ de dimension $k$, les valeurs propres du hessien $\nabla^2 \log \left( \frac{B}{A} \right)$
en un point $Z \in X$ sont toutes 
positives, inf\'erieures (ou \'egales) \`a $1$, et parmi celles-ci au moins $q+r -k +pk-1$ tendent vers $1$ lorsque $\frac{B}{A} (Z)$ tend vers l'infini. La Proposition correspond donc au cas $k=r$. 
Lorsque $k=1$, la fonction $\log \left( \frac{A}{B} \right)$ co\"{\i}ncide avec la fonction 
distance \`a $X_V$, et la Proposition \ref{asymptotique des vp} d\'ecoule de la Proposition \ref{fonction distance a l'hypersurface}. 
Supposons donc la Proposition d\'emontr\'ee au rang
$k\geq 1$. Notons $V' \subset V$ deux sous-espaces de ${\Bbb R}^n$ de dimensions respectives $1$ et $k+1$. On a alors $X_V \subset X_{V '} \subset X$. Notons $B$ et $B'$ les
fonctions correspondantes \`a la fonction $B$ plus haut pour $X_V$ et $X_{V'}$ respectivement. D'apr\`es (\ref{action}) la fonction $\frac{B}{A}$ est $G_V$-invariante, il nous suffit
donc de d\'eterminer le hessien $\nabla^2 \log \left( \frac{B}{A} \right)$ aux points $Z$ tels que $Z_1 =0$. Mais,
$\nabla^2 \log \left( \frac{B}{A} \right) = \nabla^2 \log \left( \frac{B '}{A} \right) + \nabla^2 \log \left( \frac{B}{B' } \right)$, et la d\'emonstration de la Proposition \ref{fonction distance a l'hypersurface} montre que les valeurs propres sont toutes positives, 
inf\'erieures (ou \'egales) \`a $1$ et que parmi celles-ci, au moins 
$$q+r-(k+1) + p(k+1) -1 = (q+r-1-k +pk -1)+ p$$
tendent vers $1$ lorsque $\frac{B}{A} (Z)$ tend vers l'infini. Ce qui conclut la d\'emonstration de la Proposition \ref{asymptotique des vp}.~$\Box$

\bigskip

\subsection{S\'eries de Poincar\'e}

Soit $\phi$ une forme diff\'erentielle de degr\'e $l$ sur $X$. Nous notons $|| \phi ||$ (resp. $|| \phi ||_0 $) la 
norme ponctuelle induite par la m\'etrique $g$ (resp. la m\'etrique euclidienne).

\begin{lem} \label{comparaison des normes}
On a les in\'egalit\'es suivantes~:
$$||\phi ||_0 \geq ||\phi || \geq A^l || \phi ||_0 ,$$
o\`u $A=\det (1_p - {}^t Z Z )$.
\end{lem}
{\it D\'emonstration.} La m\'etrique riemannienne de $X$, s'\'ecrit $ds^2 = {\rm tr} ((1_{q+r} -Z {}^t Z )^{-1} dZ (1_p -{}^t Z
Z)^{-1} d {}^t Z )$. Il est donc imm\'ediat que 
$${\rm tr} (dZ d{}^t Z ) \leq ds^2 \leq A^{-2} {\rm tr} ( dZd{}^t Z ).$$
Le Lemme \ref{comparaison des normes} d\'ecoule trivialement de ces derni\`eres in\'egalit\'es.~$\Box$

\bigskip

\begin{cor} \label{comparaison}
Soit $\phi$ une forme diff\'erentielle $G_V$-invariante de degr\'e $l$. Supposons que chaque coefficient 
de $\theta_1 \wedge \ldots \wedge \theta_l$, avec $\theta_1 , \ldots , \theta_l \in \{ dZ_{ij}  \; : \; 
1\leq i \leq q+r , \; 1\leq j \leq p \}$, soit born\'e en $\left(
\begin{array}{c}
0 \\
Z_2 
\end{array} \right)$. Il existe alors deux constantes $C_1$, $C_2 >0$ telles que 
$$C_1 \geq || \phi || \geq C_2 \left( \frac{A}{B} \right)^l .$$
\end{cor}
{\it D\'emonstration.} D'apr\`es (\ref{action}), il suffit de le v\'erifier en $\left(
\begin{array}{c}
0 \\
Z_2 
\end{array} \right)$. Mais en $\left(
\begin{array}{c}
0 \\
Z_2 
\end{array} \right)$, $B=1$ et le Corollaire \ref{comparaison} d\'ecoule alors du Lemme \ref{comparaison des normes}.~$\Box$

\bigskip

Soit $\Gamma \subset G$ un sous-groupe discret sans torsion de type fini. Alors 
$M= \Gamma \backslash X $ est une vari\'et\'e riemannienne compl\`ete de dimension $p(q+r)$ localement model\'ee sur $X_{p,q+r}$. 
Soit $\Gamma_V = \Gamma \cap G_V$, et soit 
$C_V =\Gamma_V \backslash X_V $. On obtient alors le diagramme commutatif suivant~:
$$\begin{array}{ccccc}
   & X_V          & \hookrightarrow           & X         & \\
   & \downarrow                         &                           & \downarrow                    & \\
C_V = & \Gamma_V \backslash X_V  & \stackrel{i}{\rightarrow} & \Gamma \backslash X & =M
\end{array}$$
o\`u l'application $i$ est induite par l'inclusion de $X_V$ dans $X$. 
En g\'en\'eral, le groupe $\Gamma_V$ est r\'eduit \`a l'identit\'e. Dans la suite nous supposons que $C_V$ est de 
volume fini (plus loin nous supposerons m\^eme que $C_V$ est compacte). Soit 
$M_V =\Gamma_V \backslash X$.
Remarquons que la fibration naturelle
$$\pi : \left\{ 
\begin{array}{rcl}
X & \rightarrow & X_V \\
Z & \mapsto & \left( 
\begin{array}{c}
Z_1 \\
0
\end{array} \right) 
\end{array} \right. $$ 
induit une fibration, nous la notons \'egalement $\pi$: $M_V = \Gamma_V \backslash X 
\rightarrow \Gamma_V \backslash X_V =C_V$.

\bigskip

La Proposition 3 de \cite{MargulisSoifer} (ou le Lemme principal de \cite{EnseignMath}) implique(nt) le lemme suivant.

\begin{lem} \label{effeuillage}
Il existe une suite $\{ \Gamma_m \}$ de sous-groupes d'indices finis dans $\Gamma$, d\'ecroissante pour l'inclusion, 
telle que 
$$\Gamma_V = \bigcap_{m\in {\Bbb N}} \Gamma_m \mbox{  et  } \Gamma_0 =\Gamma .$$
Si de plus $\Gamma$ est un sous-groupe de congruence, on peut choisir les $\Gamma_m$ de congruence.
\end{lem}

Le Lemme \ref{effeuillage} implique que lorsque $\Gamma$ est de type fini, la vari\'et\'e $M$ admet une suite 
croissante $\{ M_m \}$ de rev\^etements finis telle que la suite $\{ M_m \}$ converge uniform\'ement sur tout compact 
vers la vari\'et\'e $M_V$ (il suffit de poser $M_m =\Gamma_m \backslash X$). Nous appelons une telle 
suite de rev\^etements finis, une {\it tour d'effeuillage autour de $C_V$}. Dans la suite, nous supposons que $M$ poss\`ede 
une telle tour et notons $\Gamma_m$ le groupe fondamental de $M_m$.

Nous allons travailler avec des formes diff\'erentielles sur $X$, $M_V$ ou 
$M_m$. Il sera plus commode de consid\'erer toutes ces formes diff\'erentielles comme d\'efinies sur 
$X$ et invariantes sous l'action des groupes $\{ e \}$, $\Gamma_V$ ou $\Gamma_m$. \'Etant donn\'e 
un entier $m_0$, un \'el\'ement $\gamma \in \Gamma_{m_0}$ et une forme diff\'erentielle $\omega$ sur $M_m$ (avec 
$m\in {\Bbb N} \cup \{ \infty \}$, $m\geq m_0$), nous pourrons notamment parler de la forme diff\'erentielle 
$\gamma^* \omega$. 

Passons donc \`a l'\'etude des s\'eries de Poincar\'e.

Soient $Z_1 ,Z_2 \in X$, $t \in {\Bbb R}$, $t>0$. On introduit~: 
\begin{eqnarray} 
\nu  (Z_1 ,Z_2 ,t) := | \{ \gamma \in \Gamma \; : \; d(Z_1 , \gamma Z_2 ) \leq t \} |,
\end{eqnarray}
et
\begin{eqnarray} \label{denombrement}
N (Z, t) := |\{ \gamma \in \Gamma_V \backslash \Gamma \; : \; d(\gamma Z , X_V ) \leq t \} |.
\end{eqnarray}

\begin{lem} \label{comptage}
Il existe une constante $c_1 (Z) >0$ (qui d\'epend de $\Gamma$) telle que pour tout $t>0$ on ait~:
$$N (Z,t)  \leq c_1 (Z)  \int_0^{t+1} (1+t^{p(q+r)}) e^{(p+q+r-1) \sqrt{m} t} dt .$$
De plus on peut choisir $c_1 (Z)$ de mani\`ere \`a ce qu'elle soit born\'ee sur les compacts de ${\cal D}$.
\end{lem}
{\it D\'emonstration.} Soit $\varepsilon$ un nombre r\'eel strictement compris entre $0$ et $1$ et suffisamment petit pour que 
$$B(Z ,\varepsilon ) \cap B(\gamma Z,\varepsilon ) \neq \emptyset \Rightarrow \gamma = e ,$$
o\`u $B(Z, \varepsilon)$ d\'esigne la boule de rayon $\varepsilon$ autour du point $Z$ et 
$e$ d\'esigne l'\'el\'ement neutre du groupe $\Gamma$. 
Dans la suite \'etant donn\'ee une sous-vari\'et\'e ${\cal V}$ de $X$, nous noterons $B({\cal V}, \rho )$ l'ensemble des points de $X$
\`a distance plus petite que $\rho$ de ${\cal V}$.
On a alors : 
\begin{eqnarray*}
N (Z,t)            & \leq & | \{ [\gamma] \in \Gamma_V \backslash \Gamma \; : \; \gamma (B(Z,\varepsilon )) \subset B(X_V ,t+\varepsilon ) \} |.
\end{eqnarray*}
Mais, d'apr\`es (\ref{action}), si $\gamma \in \Gamma$ v\'erifie que $\gamma (B(Z, t+\varepsilon )) \subset B(X_V , t+\varepsilon )$
quitte \`a translater $\gamma$ par un \'el\'ement de $\Gamma_V$, on peut supposer que $\gamma (B(Z, \varepsilon ))
\subset B(S , \varepsilon )$, o\`u $S$ est un domaine fondamental mesurable pour l'action de $\Gamma_V$ sur 
$X_V$. On d\'eduit alors du Lemme \ref{croissance du volume}~:
\begin{eqnarray*}
N(Z,t)            & \leq & \frac{\mbox{vol}(B(S,t+\varepsilon ))}{\mbox{vol}(B(Z,\varepsilon ))} \\
                      & \leq & \frac{c}{\mbox{vol}(B(P,\varepsilon ))} \int_0^{t+\varepsilon} (1+t^{p(q+r)}) e^{(p+q+r-1) \sqrt{m}t} dt. \\
\end{eqnarray*}  
Ce qui ach\`eve la d\'emonstration du Lemme \ref{comptage} \footnote{Contrairement \`a ce que pourrait laisser croire la
d\'emonstration, le Lemme \ref{comptage} reste valable lorsque $C_V$ est de volume fini mais non compact, cf. \cite{TongWang2}}.~$\Box$

\bigskip

Remarquons que pour $q=0$, la d\'emonstration du Lemme \ref{comptage} permet d'estimer $\nu (Z_1, Z_2 ,t)$
uniform\'ement par rapport \`a $Z_2$. On obtient, en effet, que
pour tout $Z_2 \in X$ et pour $t>0$, 
\begin{eqnarray} \label{nu}
\nu (Z_1 , Z_2 , t) \leq c_1 (Z_1 ) \int_0^{t+1} (1+t^{p(q+r)} )e^{(p+q+r-1) \sqrt{p} t} dt .
\end{eqnarray}
On en d\'eduit la proposition suivante.

\begin{prop} \label{serie de Poincare}
Soit $K$ un compact de $X$. Alors il existe une constante $c_2 (K)$
(qui d\'epend de $\Gamma$) telle que pour tout point $Z_1 \in K$, tout point $Z_2 \in X$
et $t \geq 0$, on ait~:
$$\sum_{\begin{array}{c}
\gamma \in \Gamma \\
d(Z_1 , \gamma Z_2) \leq t
\end{array}} e^{-(p+q+r-1+s) d(Z_1 ,Z_2 ) } \leq c_2 (K)\left( 1+ \frac{1}{s} \left( 1 + \frac{1}{s^{p(q+r)}} \right) \right) , $$
pour tout $s>0$.
\end{prop}
{\it D\'emonstration.} D'apr\`es (\ref{nu}), il existe une constante $c_1 (K)$ telle que 
$$d\nu (Z_1 , Z_2 , t ) \leq c_1 (K) (1+(t+1)^{p(q+r)} )e^{(p+q+r-1) \sqrt{p} t} dt ,$$
pour tout $Z_1 \in K$, $Z_2 \in X$ et $t>0$. On a donc~: 
\begin{eqnarray*}
\sum_{\begin{array}{c}
\gamma \in \Gamma \\
d(Z_1 , \gamma Z_2) \leq t
\end{array}} e^{-((p+q+r-1)\sqrt{p} + s) d(Z_1 ,Z_2 ) } & = & \int_0^t e^{-((p+q+r-1)\sqrt{p} +s)t}  d \nu (Z_1 ,Z_2 ,t)  \\ 
                                                                                   & = & c_1 (K) \int_0^t (1+(t+1)^{p(q+r)} ) e^{-st} dt .  
\end{eqnarray*}  
Et la Proposition \ref{serie de Poincare} d\'ecoule d'un calcul simple et d'approximations grossi\`eres.~$\Box$                                                      

\bigskip

De mani\`ere analogue on d\'emontre la proposition suivante.

\begin{prop} \label{serie de formes}
Soit $\phi$ une forme diff\'erentielle $\Gamma_V$-invariante de degr\'e $l$ sur $X$.
Si $||\phi || \leq c \left( \frac{A}{B} \right)^{(p+q+r-1) \sqrt{m}/2 + \varepsilon }$, $m= \min \{ r,p \}$
pour un r\'eel strictement positif $\varepsilon >0$, alors la s\'erie
$$\sum_{\Gamma_V \backslash \Gamma} \gamma^* \phi $$
converge uniform\'ement sur les compacts de $X$.
\end{prop}
{\it D\'emonstration.} Commen\c{c}ons par remarquer que la norme $||\gamma^* \phi ||$ au point $Z$ est 
\'egale \`a la norme $|| \phi ||$ au point $\gamma Z$. D'apr\`es l'hypoth\`ese faite sur la norme de $\phi$,
\begin{eqnarray*}
\sum_{\Gamma_V \backslash \Gamma} ||\gamma^* \phi || & \leq & c \sum_{\Gamma_V \backslash \Gamma} \left( \frac{A}{B} (\gamma Z) \right)^{(p+q+r-1)\sqrt{m}/2 +\varepsilon} \\
                                                                                             & \leq & \frac{c}{4^m} \sum_{\Gamma_V \backslash \Gamma} e^{-((p+q+r-1)\sqrt{m} +\varepsilon )d(\gamma Z , X_V )},
\end{eqnarray*}
d'apr\`es la Proposition \ref{distance2}. On conclut alors facilement comme pour la 
Proposition \ref{serie de Poincare}.~$\Box$

\subsection{Tours de rev\^etements finis}

Dans \cite{Wang} Wang d\'efinit une famille de formes lisses, ferm\'ees et $G_V$-invariantes $\omega_s$ pour $s \in {\Bbb C}$ telles que 
\begin{eqnarray} \label{norme de omegas}
|| \omega_s || \prec \left( \frac{B}{A} \right)^{r/2+rp -({\rm Re}(s) -qr/2)} ,
\end{eqnarray}
(o\`u le signe $\prec$ signifie que l'on a une in\'egalit\'e $\leq$ \`a une constante positive pr\`es) 
et qui, pour Re$(s)>>0$, peut se voir comme la forme duale \`a $\Gamma_V \backslash X_V$ dans 
$\Gamma_V \backslash X$. Plus pr\'ecisemment, consid\'erons une forme $\phi$ sur $C_V$ de degr\'e $pq$. Supposons
\begin{eqnarray} \label{condition sur phi}
||\phi || \prec \left( \frac{A}{B} \right)^N 
\end{eqnarray}
pour un certain entier $N$. Remarquons que la condition (\ref{condition sur phi}) est v\'erifi\'ee par toute forme born\'ee pour $N=0$. D'apr\`es 
(\ref{norme de omegas}) et le Lemme \ref{formules d'integration}, l'int\'egrale $\int_{C_V} \omega_s \wedge \phi$ est absolument convergente pour 
Re$(s)>>0$, la constante ne d\'ependant que de $N$. Wang d\'emontre alors le th\'eor\`eme suivant.

\begin{thm} \label{dualite}
Soit $\phi$ une forme ferm\'ee, lisse, de degr\'e $pq$ sur $\Gamma_V \backslash X$ et v\'erifiant la 
condition (\ref{condition sur phi}). Alors,
$$\int_{\Gamma_V \backslash X} \omega_s \wedge \phi = 
\int_{C_V}  \phi \ \ (\mbox{Re}(s)>>0).$$
\end{thm}

\bigskip

Consid\'erons maintenant $\mu$ une $k$-forme {\bf harmonique} sur $C_V = \Gamma_V \backslash X_V$. 
Nous notons~:
\begin{eqnarray} \label{Omega1}
\Omega_{\mu} (s) = \omega_s \wedge \pi^* \mu .
\end{eqnarray}

D'apr\`es la Proposition \ref{serie de formes} et (\ref{norme de omegas}), et pour Re$(s) >>0$, la s\'erie  
\begin{eqnarray} \label{omegasm}
\Omega_{\mu}^m (s)    =  \sum_{\gamma \in \Gamma_V \backslash \Gamma_m  } \gamma^* \Omega_{\mu} (s) 
\end{eqnarray}
converge uniform\'ement sur tout compact de $X$ et d\'efinit une $(k+pr)$-forme ferm\'ee sur 
$M_m$. Le Th\'eor\`eme \ref{dualite} implique que pour toute $(pq-k)$-forme ferm\'ee $\eta $ sur $M_m$ on a~:
$$\int_{M_m} \Omega_{\mu}^m (s) \wedge \eta = \int_{M_V} \Omega_{\mu} (s) \wedge \eta = \int_{C_V} \mu \wedge \eta .$$

On obtient donc le th\'eor\`eme suivant.

\begin{thm} 
L'application $\mu \mapsto \Omega_{\mu}^m (s)$ (pour Re$(s)>>0$) induit en cohomologie, l'application naturelle
$$H^k (C_V ) \rightarrow H^k (M_m )$$
``cup-produit avec $[C_V]$''.
\end{thm}

Nous cherchons dans ce paragraphe \`a comprendre comment \'evolue cette application lorsque $m$ tend vers l'infini.

\bigskip

Soit $\Delta = \delta d +d \delta$, o\`u $\delta$ est l'adjoint de $d$, l'op\'erateur laplacien que l'on \'etend en un 
op\'erateur, toujours not\'e $\Delta$, agissant sur l'espace $L^2 \Omega^{k} (X )$ des 
$k$-formes de carr\'e int\'egrable sur $X$ de fa\c{c}on essentiellement auto-adjointe. Alors le 
Th\'eor\`eme spectral s'applique et il existe une famille spectrale $\{ P_{\lambda} \; : \; \lambda \in [0, +\infty [\}$ associ\'ee 
\`a $\Delta$.

Notons $P_{\lambda} (x,y)$ le noyau de Schwartz de $P_{\lambda}$. 
On a $\Delta = \int_0^{+\infty} \lambda dP_{\lambda}$. \`A toute fonction $f\in C_0 ([0,+\infty [)$, on associe l'op\'erateur 
$$f(\Delta ) = \int_0^{+\infty} f(\lambda )dP_{\lambda }.$$
Le laplacien est un op\'erateur elliptique. Soit $\omega \in L^2 \Omega^{n-p} (X)$.
On a :
$$\Delta (f(\Delta ) \omega ) = F(\Delta )\omega $$
o\`u $F$ est la fonction qui \`a $x$ associe $xf(x)$. D'apr\`es le Th\'eor\`eme de r\'egularit\'e sur les op\'erateurs 
elliptiques, la forme $f(\Delta )\omega$ est lisse. De plus, pour tout $x \in X$, il existe une 
constante $C(x,f,X )$ telle que~:
$$|f(\Delta )\omega |(x) \leq C(x,f,X ) ||\omega ||_{L^2 (X )} .$$
En particulier, l'application 
$$\left\{
\begin{array}{ccc}
L^2 \Omega^{k} (X ) & \rightarrow & L^2 \Omega^{k}_x (X ) \\
\omega & \mapsto & f(\Delta )\omega (x) 
\end{array}
\right. $$
est continue. D'apr\`es le Th\'eor\`eme de Riesz, il existe donc 
$f(\Delta )(x,.) \in L^2 \Omega^{k} (X )$ tel que 
$$f(\Delta )\omega (x) = \int_{X} f(\Delta )(x,y) \omega (y) dy .$$
De plus, pour tout compact $K$ de $X$, 
$\int_{K \times X} ||f(\Delta ) (x,y)||^{2} dxdy \leq \int_K C(x,f,X)^2 dx $. 
Donc $f(\Delta )(.,.) \in L^2_{loc} (X \times X)$. Or 
$$(\Delta_x + \Delta_y )f(\Delta )(x,y) = 2F(\Delta ) (x,y)$$
et l'op\'erateur $\Delta_x + \Delta_y$ est elliptique. Le Th\'eor\`eme de r\'egularit\'e elliptique implique donc que 
$f(\Delta ) (x,y)$ est une fonction $C^{\infty}$ en $x$ et $y$.

\bigskip

Sur les vari\'et\'es $M_m$ et $M_{V}$, nous notons le laplacien respectivement $\Delta_m$ et $\Delta_{\infty}$; 
de m\^eme nous notons respectivement $P_{\lambda}^m$ et $P_{\lambda}^{\infty}$ les familles spectrales associ\'ees. 
Si l'on d\'esigne, de mani\`ere coh\'erente avec les notations pr\'ec\'edentes,
le noyau de la chaleur sur les $k$-formes de $X$ par $e^{-t\Delta} (x,y)$, il est connu \cite{Donnelly} 
que pour tout $T>0$, il existe une constante $\alpha >0$ et une constante $C_T$ (d\'ependante de $T$) telles que~:
\begin{eqnarray} \label{noyau de la chaleur}
|e^{-t\Delta} (x,y)| & \leq & C_T e^{-\alpha d(x,y)^2 /t } ,
\end{eqnarray}
pour tout $t\in ]0,T]$. 

\begin{lem} \label{noyau de la chaleur 2}
Pour $t>0$ fix\'e, la s\'erie 
$$\sum_{\gamma \in \Gamma} |e^{-t\Delta} (x,\gamma y)|$$
converge uniform\'ement pour $x \in X$ et $y$ dans un compact.
\end{lem}
{\it D\'emonstration.} Soit $K$ un compact de $X$ et soit $t$ un r\'eel strictement positif. D'apr\`es 
le Lemme \ref{comptage}, il existe une constante $c_1 (K)$ telle que 
$$\nu(x,y,R ) := |\{ \gamma \in \Gamma \; : \;  d(x,\gamma y) \leq R \} |\leq c_1 (K) \int_0^{R+1} (1+u^{p(q+r)} ) e^{(p+q+r-1)\sqrt{p} u} du$$
pour tout $x \in X$ et $y \in K$. 

Soient $C=C_t$ et $\beta = \alpha /t$. Alors, d'apr\`es l'in\'egalit\'e (\ref{noyau de la chaleur}), pour 
$x \in X$ et $y \in K$, on a~:
\begin{eqnarray*}
\sum_{ 
\begin{array}{c}
\gamma \in \Gamma \\
d(x,\gamma y) \leq R
\end{array} 
}
|e^{-t\Delta } (x,\gamma y )| & \leq & 
C \sum_{ 
\begin{array}{c}
\gamma \in \Gamma \\
d(x,\gamma y) \leq R
\end{array} 
}
e^{-\beta d(x,\gamma y )^2 } \\
& \leq & C \left( \int_0^R e^{-\beta r^2} d\nu (x,y,r) \right) \\
& \leq & C \left( [e^{-\beta r^2 } \nu (x,y,r) ]_0^R + 2 \beta \int_0^R re^{-\beta r^2} \nu (x,y,r) dr \right) \\
& \leq & C c_2 (K) \left( e^{-\beta R^2} \int_0^{R+1} (1+u^{p(q+r)} )e^{(p+q+r-1) \sqrt{p} u} du  \right. \\
&        & \left. +2\beta \int_0^R re^{-\beta r^2} \int_0^{r+1} (1+u^{p(q+r)} ) e^{(p+q+r-1) \sqrt{p} u} du dr \right) .
\end{eqnarray*}
Le Lemme d\'ecoule de ces in\'egalit\'es en faisant tendre $R$ vers l'infini.~$\Box$  

\bigskip

On obtient alors que pour tout $t>0$,
$$e^{-t\Delta_m } (x,y) = \sum_{\gamma \in \Gamma_m} (\gamma_y )^* e^{-t\Delta } (x,y) $$
$$e^{-t\Delta_{\infty}} (x,y) = \sum_{\gamma \in \Gamma_V } (\gamma_y )^* e^{-t\Delta } (x,y) .$$
Et la convergence est absolue et uniforme pour $x$, $y$ dans un compact. En particulier, le noyau de la chaleur 
$e^{-t\Delta_m } (x,y)$ est $\Gamma_m$-invariant. De plus, puisque d'apr\`es le Lemme \ref{effeuillage}, 
$\cap_m \Gamma_m = \Gamma_V$ et $\Gamma_{m+1} \subset \Gamma_m$, on en d\'eduit que si $t>0$ est fix\'e, 
$$e^{-t\Delta_{\infty}} (x,y) = \lim_{m\rightarrow +\infty} e^{-t\Delta_m} (x,y) $$
uniform\'ement pour $x$ et $y$ dans un compact de $X$. 

\bigskip

Le lemme suivant est une cons\'equence du Th\'eor\`eme d'approximation de Weierstrass.

\begin{lem}[cf. \cite{Donnelly2}]  \label{approx par des polynomes}   
Soit $f \in C_0 ([0, +\infty [)$. Alors $f$ peut \^etre uniform\'ement approch\'ee sur $[0, +\infty[$ par une combinaison     
lin\'eaire finie d'exponentielles $e^{-tx}$, $t>0$.
\end{lem}

Comme Donnelly dans \cite{Donnelly} (cf. aussi \cite{BergeronClozel}) on peut d\'eduire du Lemme \ref{approx par des polynomes} le lemme suivant \footnote{Ici l'hypoth\`ese de compacit\'e sur $M$ peut sembler n\'ecessaire, il n'en est rien le finitude du volume
est suffisante \`a condition de consid\'erer une fonction $f$ dont le support \'evite le spectre continue de $M$.}. 

\begin{lem} \label{convergence des noyaux}
Pour tout $f\in C_0 ([0,+\infty [)$, la suite $ f(\Delta_m )(x,y)$ converge vers $f(\Delta_{\infty} )(x,y)$ 
uniform\'ement pour $x$ et $y$ dans un compact. Et l'expression
$$| f(\Delta_m )(x, y)) -f(\Delta_{\infty} )(x,y) |$$
est uniform\'ement born\'ee (ind\'ependamment de $m$) pour $x\in X$ et $y$ dans un compact.
\end{lem}

\bigskip

Remarquons que le noyau $f(\Delta_m ) (.,.)$ est $\Gamma_m$-bi-invariant.

\begin{prop} \label{convergence des noyaux 2}
Soient $f \in C_0 ([0,+\infty [)$ et $s\in {\Bbb C}$, ${\rm Re}(s) >>0$. Alors, la suite 
$f (\Delta_{m} ) \Omega^m_{\mu} (s) $ converge vers $f(\Delta_{\infty } ) \Omega_{\mu} (s)$ uniform\'ement sur les compacts.
\end{prop}
{\it D\'emonstration.} Soit $s\in {\Bbb C}$, Re$(s) >>0$. Si ${\cal F}_m$ est un domaine fondamental pour 
l'action de $\Gamma_m$ sur $X$, on a~:
\begin{eqnarray*}    
f(\Delta_m )\Omega_{\mu}^m (s) (.) & = & \int_{M_m} \Omega_{\mu}^m (s) (x) \wedge     *f(\Delta_m ) (x,.) dx \\
    & = & \int_{M_m} \left( \sum_{\gamma \in \Gamma_V \backslash     \Gamma_m } \gamma^* \Omega_{\mu} (s) (x) \right) \wedge *f(\Delta_m     )(x,.) dx \\
    & = & \sum_{\gamma \in \Gamma_V \backslash     \Gamma_m } \int_{{\cal F}_m}  \gamma^* \Omega_{\mu} (s) (x) \wedge *f(\Delta_m     )(x,.) dx \\
      & = & \sum_{\gamma \in \Gamma_V \backslash     \Gamma_m } \int_{\gamma {\cal F}_m} \Omega_{\mu} (s) (x) \wedge *f(\Delta_m     )(x,.) dx \\
    & & \mbox{(car $f(\Delta_m )(.,.)$ est $\Gamma_m$-bi-invariant)} \\
    & = & \int_{M_{V}} \Omega_{\mu} (s) (x) \wedge *f(\Delta_m     )(x,.) dx .
\end{eqnarray*}

On obtient donc sur $X$~:
\begin{eqnarray*}    
f(\Delta_m )\Omega_{\mu}^m (s) (.) - f(\Delta_{\infty } ) \Omega_{\mu} (s) (.) & = & \int_{M_{V}} \Omega_{\mu} (s) (x) \wedge *f(\Delta_m ) (x,.) dx \\
&  & - \int_{M_{V}} \Omega_{\mu} (s) (x) \wedge *f(\Delta_{\infty} ) (x,.) dx \\
& = & \int_{M_{V}} \Omega_{\mu} (s) (x) \wedge *(f(\Delta_m ) (x,.) -f(\Delta_{\infty} ) (x,.)) dx .
\end{eqnarray*}

De plus d'apr\`es le Lemme \ref{convergence des noyaux}, l'expression 
$$|f(\Delta_m ) (x,y) -f(\Delta_{\infty} ) (x,y)|$$
est uniform\'ement born\'ee (ind\'ependamment de $m$) pour $x\in M_{V}$ et $y$ dans un compact.
Puisque la forme $\Omega_{\mu} (s)$ est dans $L^1$, le Th\'eor\`eme de convergence domin\'ee et le Lemme \ref{convergence des noyaux} impliquent 
que pour tout r\'eel $\varepsilon >0$ et pour tout compact $K$, il existe un entier $m_0$ tel que pour tout $m\geq m_0$, 
les applications $f(\Delta_m )\Omega_{\mu}^m (s)$ et $f(\Delta_{\infty } ) \Omega_{\mu} (s)$ sont $\varepsilon$-proches sur $K$. 
Ce qui ach\`eve la d\'emonstration de la Proposition \ref{convergence des noyaux 2}.~$\Box$

\bigskip

\section{Calcul de la cohomologie $L^2$}

Nous conservons dans cette section les notations pr\'ec\'edentes. Soient donc toujours $G = O(p,q+r)$ ($p\leq q+r$), $X$
l'espace sym\'etrique associ\'e, $X_V$ le sous-espace totalement g\'eod\'esique de $X$ associ\'e \`a un 
sous-espace vectoriel de dimension $r$ de ${\Bbb R}^{p+q+r}$, $G_V$ le sous-groupe de $G$ pr\'eservant 
$X_V$ et $\Gamma_V$ un sous-groupe discret sans torsion et cocompact dans $G_V$. Nous notons $C_V = \Gamma_V 
\backslash X_V$ et $(M= )M_V = \Gamma_V \backslash X$. Dans la suite $H^k_2 (M_V )$ d\'esigne toujours le 
groupe de cohomologie $L^2$ r\'eduite de degr\'e $k$ de $M_V$.
Le but de cette section est la d\'emonstration du th\'eor\`eme suivant.

\begin{thm} \label{cohom l2}
Pour tout entier, $k < (q+pr-1)/2$, l'application naturelle ``cup-produit avec $[C_V]$''
$$H^{k-pr} (C_V ) \rightarrow H_{2}^{k} (M_V )  $$
est un {\bf isomorphisme}.
Si de plus, $p=1$ et $q+r$ est pair, l'espace $H_{2}^{(q+r)/2} (M_V )$ est de dimension infini et 
l'application ci-dessus reste {\bf injective} pour $k=(q+r)/2$. 
\end{thm}

Remarquons que le cas $p=1$ d\'ecoule d'un Th\'eor\`eme de Mazzeo et Philips \cite{MazzeoPhillips}.

Avec les notations du \S 2.6 (Groupes orthogonaux), le Th\'eor\`eme \ref{cohom l2} implique que, pour tout entier $k\leq (q- pr -1)/2$, l'application ``cup-produit avec $[F]$'' 
\begin{eqnarray} \label{c}
H^{k} (F ) \rightarrow H^{k+ pr}_2 (M)
\end{eqnarray}
est {\bf injective}. Ce qui est un cas particulier de la Conjecture \ref{conjl2O}.

\subsection{Une proposition de Donnelly et Xavier}

L'ingr\'edient principal de la d\'emonstration du Th\'eor\`eme \ref{cohom l2} est la proposition suivante due \`a Donnelly et Xavier \cite{DonnellyXavier}.
 
\begin{prop} \label{DX}
Soit $\phi$ une forme diff\'erentielle dans $C_0^{\infty} (\bigwedge^{k} T^* M)$.
Soient $F$ une fonction r\'elle $C^2$ sur le support de $\phi$ et
$\gamma_1 , \ldots , \gamma_{p(q+r)}$ les valeurs propres r\'eelles du hessien de $F$. Si $|dF| \leq 1$, alors~:
$$[||d\phi ||_{2} + ||\delta \phi ||_{2} ] ||\phi ||_{2} \geq \frac12 \int_M \left[ \sum \gamma_i - 2k \max_i (\gamma_i ) \right] |\phi |^2 .$$
\end{prop}

La vertu essentielle de la Proposition \ref{DX} est de contr\^oler le spectre essentiel de $M_V$.

Soit $\psi$ la fonction r\'elle sur $X$ d\'efinie par 
$$\psi (Z) = \log \left( \frac{B}{A} \right). $$
D'apr\`es la Proposition \ref{asymptotique des vp}, le hessien 
$\nabla^2 \psi$ de $\psi$ a pour valeurs propres des r\'eels
$\gamma_1 (Z) , \ldots , \gamma_{p(q+r)} (Z)$ tels qu'il existe une constante $c_0 >0$ telle que 
pour tout $Z$ tel que $\psi (Z) > c_0$, 
\begin{eqnarray} \label{1/10}
\sum_i \gamma_i (Z) - 2k \max_i (\gamma_i (Z)) \geq 1/10 ,
\end{eqnarray}
pour $k < (p+qr -1)/2$. 

Remarquons que la fonction $\psi $ est $G_V$-invariante et descend donc en une fonction sur $M_V = \Gamma_V 
\backslash X$. Dans la suite, $M= M_V $ et 
$M_c = \{ Z \in M \; : \; \psi (Z )\leq c\} $ pour tout r\'eel $c>0$. 

\begin{lem} \label{spectre isole de 0}
Pour tout degr\'e $k$ tel que $k <(p+qr-1)/2$, $0$ n'est pas dans le spectre essentiel du laplacien sur les 
formes de degr\'e $k$ sur $M$.
\end{lem}
{\it D\'emonstration.} Il est bien connu que le spectre essentiel du laplacien ne 
d\'epend que de la g\'eom\'etrie \`a l'infini. 
Or, la Proposition \ref{DX} et l'in\'egalit\'e (\ref{1/10}) impliquent que si $\omega$ est une forme diff\'erentielle de 
degr\'e $k$, $k < (p+qr-1)/2$, \`a support dans le compl\'ementaire de $M_{c_0}$ alors~:
\begin{eqnarray*}
[||d\omega ||_2 + ||\delta \omega ||_2 ] ||\omega ||_2 & \geq &\int_M \left[ \sum \gamma_i -k \max_i (\gamma_i ) \right] |\omega |^2 \\
& \geq & \frac{1}{10} ||\omega ||_2^2 .
\end{eqnarray*}
Le spectre du laplacien \`a l'infini (et donc le spectre essentiel) est donc isol\'e de $0$.
Le Lemme \ref{spectre isole de 0} est d\'emontr\'e.~$\Box$

\bigskip

Il est naturel de se demander (comme au \S 2.6) si dans l'\'enonc\'e du Th\'eor\`eme \ref{cohom l2} ou du Lemme \ref{spectre isole de 0}
le nombre $(p+qr-1)/2$ est optimal. En g\'en\'eral on ne sait pas r\'epondre \`a cette question. Peut-\^etre la formule de Plancherel
pour les espaces sym\'etriques pseudoriemannien $G/H$ peut-elle apporter une r\'eponse. Remarquons
n\'eanmoins le lemme suivant (qui explique la deuxi\`eme partie du Th\'eor\`eme \ref{cohom l2}).

\begin{lem} \label{q=1}
Si $p=1$ et $q+r$ est pair, $0$ n'est pas dans le spectre essentiel du laplacien sur les formes de degr\'e $(q+r)/2$ sur $M$.
\end{lem}
{\it D\'emonstration.} C'est \'evident puisque, le laplacien commutant 
aux op\'erateurs $d$ et $\delta$, son spectre sur les formes de degr\'e $k$ est 
contenu dans la r\'eunion du spectre du laplacien sur les formes de degr\'e $k-1$, du spectre du laplacien sur 
les formes de degr\'e $k+1$ et (\'eventuellement) de la valeur propre $0$. Mais par dualit\'e de Hodge le 
spectre du laplacien sur les formes de degr\'e $k+1$ co\"{\i}ncide avec le spectre du laplacien sur les formes
de degr\'e $k-1$. Le Lemme \ref{q=1} d\'ecoule donc du Lemme \ref{spectre isole de 0}.~$\Box$

\bigskip

Nous aurons besoin d'une version \`a bord de la Proposition \ref{DX}. Soit $D$ une sous-vari\'et\'e compacte \`a bord contenu dans $M$ et soit $M_0 =M-D$.
Alors $M_0$ est une vari\'et\'e riemannienne \`a bord compact. Nous cherchons \`a \'etudier le cas o\`u la fonction $F$ n'est d\'efinie que dans un voisinage de $M_0$ 
dans $M$.

Soit $\phi$ une forme diff\'erentielle d\'efinie sur $M_0$. Le long du bord $\partial M_0$ de $M_0$ on peut \'ecrire $\phi = \phi_{\rm tan} + \phi_{\rm norm} \wedge \nu$ o\`u
$\nu$ est un vecteur normal pointant vers l'int\'erieur de $M_0$. 

\begin{prop} \label{DX2}
Soit $\phi$ une $k$-forme diff\'erentielle \`a support compact dans $M_0$ et v\'erifiant la condition au bord $\phi_{\rm norm} =0$ le long de $\partial M_0$.
Soit $F$ une fonction r\'elle comme dans la Proposition \ref{DX}. Supposons de plus que $dF_{|\partial M_0} = 0$ et $dF (\nu ) \geq 0$ le long de $\partial M_0$.
Si $|dF| \leq 1$, alors~:
$$[||d\phi ||_{2} + ||\delta \phi ||_{2} ] ||\phi ||_{2} \geq \frac12 \int_M \left[ \sum \gamma_i - 2k \max_i (\gamma_i ) \right] |\phi |^2 .$$
\end{prop}
{\it D\'emonstration.} La Proposition \ref{DX2} d\'ecoule essentiellement de la d\'emonstration de \cite[Theorem 2.2]{DonnellyXavier}. 
L'\'etape initiale de celle-ci consiste en effet en une int\'egration par
parties, il s'agit ici d'utiliser les hypoth\`eses $\phi_{\rm norm}=0$ et $dF (\nu ) \geq 0$ pour s'assurer de la positivit\'e de l'int\'egrale le long du bord $\partial M_0$.  
La suite de la d\'emonstration est identique.~$\Box$

\bigskip

\subsection{Cohomologie $L^2$ relative}

Si $M_0$ est une vari\'et\'e riemannienne \`a bord compacte et m\'etriquement compl\`ete, on peut d\'efinir des espaces de cohomologie $L^2$ absolue et relative. Notons
$C_b^{\infty} (\bigwedge^k T^* M_0 )$ l'espace des $k$-formes lisses \`a support born\'e dans $M_0$, le support pouvant rencontrer le bord (contrairement \`a ce qui 
se passe pour les \'el\'ements de $C_0^{\infty}$). L'espace de cohomologie $L^2$ absolue est alors d\'efinie par~:
$$H_2^k (M_0 ) = ( \delta C_0^{\infty} (\bigwedge^{k+1} T^* M_0 ))^{\perp} / \overline{dC_b^{\infty} (\bigwedge^{k-1} T^* M_0 )} .$$
(L'orthogonal et l'adh\'erence sont, ici encore, pris dans l'espace $L^2 (\bigwedge^k T^* M_0 )$.)
Cet espace est isomorphe \`a un espace de formes harmoniques v\'erifiant la condition absolue sur le bord $\partial M_0$~: 
$$H_2 ^k (M_0 ) \cong {\cal H}_A^k (M_0 ) := \{ \alpha \in L^2 (\bigwedge^k T^* M_0 ) \; : \; d \alpha = \delta \alpha = 0, \; \alpha_{\rm norm} = 0 \} .$$

La cohomologie $L^2$ relative, quant \`a elle, est d\'efinie par 
$$H_2^k (M_0 , \partial M_0 ) = (\delta C_b^{\infty} (\bigwedge^{k+1} T^* M_0 ))^{\perp} / \overline{dC_0^{\infty} (\bigwedge^{k-1} T^* M_0 )} .$$
Elle est isomorphe \`a un espace de formes harmoniques v\'erifiant la condition relative au bord~: 
$$H_2 ^k (M_0 , \partial M_0 ) \cong {\cal H}_R^k (M_0 ) := \{ \alpha \in L^2 (\bigwedge^k T^* M_0 ) \; : \; d \alpha = \delta \alpha = 0, \; \alpha_{\rm tan} = 0 \} .$$

\bigskip

Bien que la d\'emonstration du Th\'eor\`eme \ref{cohom l2} ne n\'ecessite que des r\'esultats plus anciens de Cheeger, nous utiliserons le r\'esultat se suite
exacte commode suivant d\^u \`a Yeganefar \cite{Yeganefar}.

\begin{thm} \label{yeganefar}
Soient $M$ une vari\'et\'e riemannienne compl\`ete, $D$ une sous-vari\'et\'e \`a bord compact et r\'egulier de $M$ et $M_0 = M-D$. On suppose que pour un certain entier $k$, $0$
n'est pas dans le spectre essentiel du laplacien sur les $k$-formes. Alors nous avons la suite exacte
$$
\begin{array}{l}
H_2^{k-1} (M ) \stackrel{r}{\rightarrow} H_2^{k-1} (M_0 ) \stackrel{b}{\rightarrow} H_2^k (D , \partial D ) \stackrel{e}{\rightarrow} H_2^k (M ) \\
\stackrel{r}{\rightarrow} H_2^k (M_0 ) \stackrel{b}{\rightarrow} H_2^{k+1} (D , \partial D ) \stackrel{e}{\rightarrow} H_2^{k+1} (M) \stackrel{r}{\rightarrow} H_2^{k+1} (M_0 ) .
\end{array}
$$
\end{thm}

Pr\'ecisons les applications $r$, $e$ et $b$ qui interviennent dans la suite exacte.
\begin{itemize}
\item $r:H_2^* (M) \rightarrow H_2^* (M_0 )$ est l'application de restriction induite en cohomologie $L^2$ par l'inclusion $M_0 \hookrightarrow M$.
\item $e$ est l'extension par z\'ero~: si $[\alpha ] \in H_2^* (D, \partial D)$, de repr\'esentant $\alpha$, on d\'efinit $\tilde{\alpha}$ sur $M$ par $\tilde{\alpha} = \alpha$
sur $D$ et $\tilde{\alpha} = 0$ sur $M_0$. On v\'erifie que $\tilde{\alpha}$ est faiblement ferm\'ee et on pose $e([\alpha ]) = [\tilde{\alpha}]$. Ceci est ind\'ependant du 
repr\'esentant choisi.
\item $b$ se d\'efinit comme l'homomorphisme cobord ordinaire en cohomologie de de Rham~: si $[\alpha ]$ est une classe dans $H_2^k (M_0 )$, on peut choisir un 
repr\'esentant $\alpha$ qui soit ferm\'e et lisse. Il existe alors une $k$-forme $\overline{\alpha}$ sur $M$ qui co\"{\i}ncide avec $\alpha$ sur $M_0$ et qui v\'erifie $d\overline{\alpha} =0$
sur un voisinage de $M_0$. On impose de plus que $\overline{\alpha}$ et $d\overline{\alpha}$ soient de carr\'e int\'egrable (c'est toujours possible car le bord $\partial D$
\'etant compact, on peut choisir $\overline{\alpha}$ telle que $\overline{\alpha}_{|D}$ soit \`a support born\'e). On v\'erifie que la classe de $d\overline{\alpha}$ dans 
$H_2^{k+1} (D, \partial D)$ est ind\'ependante du repr\'esentant lisse $\alpha$ choisi, ainsi que du prolongement $\overline{\alpha}$, et on pose $b([\alpha]) = [d\overline{\alpha}]$.
\end{itemize}

\bigskip

Remarquons maintenant qu'en utilisant la Proposition \ref{DX2} plut\^ot que la Proposition \ref{DX}, la d\'emonstration du Lemme \ref{spectre isole de 0} implique le lemme suivant.

\begin{lem} \label{5.1}
Soit $c$ un r\'eel strictement positif suffisamment grand ($\geq c_0$). Alors, pour tout entier $k < (p+qr-1)/2$, l'espace 
$$H_2^k (M-M_c ) \cong {\cal H}_A^k (M-M_c ) = \{ 0 \} .$$
\end{lem}

\subsection{D\'emonstration du Th\'eor\`eme \ref{cohom l2}}

Fixons $c$ comme dans le Lemme \ref{5.1} et notons $D= M_c$ et $M_0 = M-D$. Fixons enfin un entier $k < (p+qr-1)/2$. D'apr\`es le  
Lemme \ref{spectre isole de 0}, $0$ n'est pas dans le spectre essentiel du laplacien sur les 
formes de degr\'e $k$ sur $M$. Nous pouvons donc appliquer le Th\'eor\`eme \ref{yeganefar}, on a en particulier la suite exacte
\begin{eqnarray} \label{se}
H_2^{k-1} (M_0 ) \rightarrow H_2^k (D , \partial D ) \rightarrow H_2^k (M ) \rightarrow H_2^k (M_0 ).
\end{eqnarray}

Remarquons maintenant que puisque $k-1 < k < (p+qr-1)/2$, le Lemme \ref{5.1} implique que 
$$H_2^{k-1} (M_0 ) \cong H_2^{k} (M_0 ) \cong \{ 0 \} .$$
La suite exacte (\ref{se}) se r\'eduit donc \`a l'isomorphisme
$$H_2^k (M ) \cong H_2^k (D, \partial D ) .$$
Mais le domaine $D$ est compact, on a donc 
$$H_2^k (D , \partial D) \cong H^k (D , \partial D) \cong H_c^k (M).$$
La vari\'et\'e $M$ \'etant par ailleurs hom\'eomorphe au produit $C_V \times {\Bbb R}^{qr}$, la formule de K\"unneth implique que 
$$H^{k-pr} (C_V ) \cong H_c^k (M) .$$
Et la premi\`ere partie de Th\'eor\`eme \ref{cohom l2} est d\'emontr\'ee. 

Supposons maintenant $p=1$ et $q+r$ pair. On a bien s\^ur toujours l'isomorphisme $H_2^{(q+r)/2} (M) \cong {\cal H}_2^{(q+r)/2} (M)$. Le groupe ${\cal H}_2^{(q+r)/2} (M)$
($(q+r)/2$ est ici la dimension r\'eelle moiti\'ee de $M$) ne d\'epend que de la structure conforme de $M$. Il est donc facile de v\'erifier que l'espace ${\cal H}_2^{(q+r)/2} (M)$
est de dimension infini. 

Enfin d'apr\`es le Lemme \ref{spectre isole de 0}, $0$ n'est pas dans le spectre essentiel du laplacien sur les 
formes de degr\'e $(q+r)/2 -1$ sur $M$. Nous pouvons donc appliquer le Th\'eor\`eme \ref{yeganefar} pour $k=(q+r)/2 -1$, on a en particulier la suite 
exacte
\begin{eqnarray} \label{se2}
H_2^{(q+r)/2 -1} (M_0 ) \rightarrow H_2^{(q+r)/2} (D , \partial D ) \rightarrow H_2^{(q+r)/2} (M).
\end{eqnarray}
Mais puisque $(q+r)/2-1 < (q+r-1)/2$ le Lemme \ref{5.1} implique que $H_2^{(q+r)/2 -1} (M_0 ) \cong \{ 0\}$ et la suite exacte (\ref{se2}) se r\'eduit \`a l'injection 
de $H_2^{(q+r)/2} (D , \partial D )$ dans $H_2^{(q+r)/2} (M)$. La d\'emonstration de la derni\`ere partie du Th\'eor\`eme \ref{cohom l2} suit alors celle de la premi\`ere partie.~$\Box$

\bigskip

\section{D\'emonstration des principaux r\'esultats}

\subsection{Autour d'un Th\'eor\`eme de Burger et Sarnak}

Soient $G$ un groupe r\'eductif, connexe et anisotrope sur ${\Bbb Q}$ et $H$ un sous-groupe rationnel de $G$, r\'eductif et connexe. 
Soit $\Gamma$ un sous-groupe de congruence de $G$. D'apr\`es la formule de Matsushima, si $\pi \in \widehat{G}_0^{\rm nc}$
intervient discr\`etement dans $L^2 (\Gamma \backslash G_0^{\rm nc} )$ avec multiplicit\'e $n_{\Gamma} (\pi )$ alors
$$H^* (\pi : \Gamma ) \simeq n_{\Gamma} (\pi ) H^* (\mathfrak{g} , K ; \pi ) .$$

Fixons $\pi$ est une repr\'esentation cohomologique de $G_0^{\rm nc}$ de degr\'e fortement primitif $R$, $\varphi$ un morphisme non triviale dans 
${\rm Hom}_{\mathfrak{g}, K} (\pi , C^{\infty} ( \Gamma \backslash G_0^{\rm nc} ))$ et $\omega : \bigwedge^R \mathfrak{p} \rightarrow \pi$ une $K$-application
non nulle. Ceci d\'efinit une classe $\omega_{\varphi} \in H^R (\pi : \Gamma )$. D\'emontrer l'injectivit\'e de la restriction \`a $H^R (\pi : \Gamma )$ de l'application de restriction
virtuelle ${\rm Res}_H^G$ revient \`a d\'emontrer qu'il existe un \'el\'ement $g \in G({\Bbb Q})$ tel que la classe $j_g^* \omega_{\varphi} \in H^R ((H \cap g^{-1} \Gamma g)
\backslash X_H )$ soit non nulle.

Harris et Li ont remarqu\'e dans \cite{HarrisLi} que pour d\'emontrer ceci il suffit de d\'emontrer que 
\begin{enumerate}
\item il existe une repr\'esentation cohomologique $\sigma$ de degr\'e (fortement primitif) $R$ de $H_0^{\rm nc}$ qui intervient discr\`etement dans la restriction de $\pi$
\`a $H_0^{\rm nc}$;
\item l'application de restriction naturelle
\begin{eqnarray} \label{3.9}
H^R (\mathfrak{g} , K ; \pi ) \rightarrow H^R (\mathfrak{h} , K^H ; \sigma )
\end{eqnarray}
est injective;
\item il existe un \'el\'ement $g\in G({\Bbb Q})$ tel que l'application $H_0^{\rm nc}$-\'equivariante
\begin{eqnarray} \label{3.10}
\psi_g : \left\{
\begin{array}{rcl}
\pi & \rightarrow & C^{\infty} ((H \cap g^{-1} \Gamma g) \backslash X_H ) , \\
\nu & \mapsto & \varphi (\nu ) (g.) 
\end{array} \right.   
\end{eqnarray}
soit non nulle.
\end{enumerate}

On peut penser aux deux premiers points ci-dessus comme \`a un analogue local (r\'eduit \`a de l'alg\`ebre lin\'eaire) de la conclusion.

\bigskip

Harris et Li ont montr\'e comment paraphraser un Th\'eor\`eme de Burger et Sarnak \cite{BurgerSarnak} pour obtenir un crit\`ere de non nullit\'e de (\ref{3.10}). La tr\`es 
l\'eg\`ere modification suivante de ce crit\`ere est d\'emontr\'ee dans \cite{IRMN}.

\begin{prop} \label{P3.2}
Soit $\pi$ (resp. $\sigma$) une repr\'esentation cohomologique, de degr\'e fortement primitif $R$, du groupe $G_0^{\rm nc}$ (resp. $H_0^{\rm nc}$).
Supposons que 
\begin{enumerate}
\item la repr\'esentation $\sigma$ appara\^{\i}t discr\`etement dans la restriction de $\pi$ \`a $H_0^{\rm nc}$; 
\item la repr\'esentation $\sigma$ est isol\'ee dans le dual automorphe de $H$ sous la condition $d=0$.
\end{enumerate}
Il existe alors un \'el\'ement $g \in G({\Bbb Q})$ tel que l'application $\psi_g$ (\ref{3.10}) est non nulle.
 \end{prop}

\bigskip

Consid\'erons maintenant deux repr\'esentations cohomologiques $\pi_1$ et $\pi_2$ de $G_0^{\rm nc}$ de degr\'es fortement primitifs respectifs $R_1$ et $R_2$ et 
deux classes de cohomologie $\alpha_{\varphi_1} \in H^{R_1} (\pi : \Gamma )$ et $\beta_{\varphi_2} \in H^{R_2} (\pi : \Gamma )$ avec comme au-dessus
$\varphi_i \in {\rm Hom}_{\mathfrak{g}, K} (\pi_i , C^{\infty} ( \Gamma \backslash G_0^{\rm nc} ))$ non nul ($i=1,2$) et $\alpha : \bigwedge^R \mathfrak{p} \rightarrow \pi_1$
(resp.  $\beta : \bigwedge^R \mathfrak{p} \rightarrow \pi_2$) une $K$-application non nulle.

De la m\^eme mani\`ere que pour la restriction, 
remarquons que si $g$ est un \'el\'ement de $G({\Bbb Q})$, pour d\'emontrer que le cup-produit $g(\alpha_{\varphi_1} ) \wedge \beta_{\varphi_2}$ est non nul il suffit 
de d\'emontrer que
\begin{enumerate}
\item il existe une repr\'esentation cohomologique $\sigma$ de degr\'e (fortement primitif) $R=R_1 +R_2$ de $G_0^{\rm nc}$ qui intervient discr\`etement dans la restriction de 
$\pi_1 \otimes \pi_2$ \`a $G_0^{\rm nc}$;
\item l'application de restriction naturelle
\begin{eqnarray} \label{3.9'}
H^{R_1} (\mathfrak{g} , K ; \pi_1 ) \otimes H^{R_2} (\mathfrak{g} , K ; \pi_2 ) \rightarrow H^R (\mathfrak{g} , K ; \sigma )
\end{eqnarray}
est injective;
\item l'application $G_0^{\rm nc}$-\'equivariante
\begin{eqnarray} \label{3.10'}
\psi_g : \left\{
\begin{array}{rcl}
\pi_1 \otimes \pi_2 & \rightarrow & C^{\infty} ((\Gamma \cap g^{-1} \Gamma g) \backslash X_G ) , \\
\nu_1 \otimes \nu_2 & \mapsto & \varphi_1 (\nu_1 ) (g.) \varphi_2 (\nu_2) (.) 
\end{array} \right.   
\end{eqnarray}
est non nulle.
\end{enumerate}

\begin{prop} \label{P3.2'}
Soient $\pi_1$, $\pi_2$ et $\sigma$ trois repr\'esentations cohomologiques, de degr\'es fortement primitifs respectifs $R_1$, $R_2$ et $R=R_1 +R_2$, du groupe $G_0^{\rm nc}$.
Supposons que 
\begin{enumerate}
\item la repr\'esentation $\sigma$ appara\^{\i}t discr\`etement dans la restriction de $\pi_1 \otimes \pi_2$ \`a $G_0^{\rm nc}$; 
\item la repr\'esentation $\sigma$ est isol\'ee dans le dual automorphe de $G$ sous la condition $d=0$.
\end{enumerate}
Il existe alors un \'el\'ement $g \in G({\Bbb Q})$ tel que l'application $\psi_g$ (\ref{3.10'}) est non nulle. 
\end{prop}
{\it D\'emonstration.} Fixons $f_1$ et $f_2  \in C^{\infty} (\Gamma \backslash G_0^{\rm nc})$ dans les images respectives de $\varphi_1$ et $\varphi_2$ et telles que 
la fonction $f_1 \otimes f_2$ engendre la repr\'esentation $\sigma$ de $G_0^{\rm nc}$ (c'est possible d'apr\`es notre premi\`ere hypoth\`ese). 
Le coefficient matriciel de $\sigma$ associ\'e \`a ce vecteur est \'egal \`a 
\begin{eqnarray} \label{3.11}
\psi (g) = \int_{\Gamma \backslash G_0^{\rm nc}} \int_{\Gamma \backslash G_0^{\rm nc}} f_1 (g_1 ) f_2 (g_2 ) \overline{f_1 (g_1 g) f_2 (g_2 g)} dg_1 dg_2 \  (g \in G_0^{\rm nc} ).
\end{eqnarray}
Dans \cite{BurgerSarnak}, Burger et Sarnak montrent que $\psi$ est la limite, uniforme sur les compactes, de sommes finis de coefficients matriciels de la forme
\begin{eqnarray} \label{3.12}
\int_{(\Gamma \cap \delta^{-1} \Gamma \delta) \backslash G_0^{\rm nc}} f_1 (\delta g') f_2 ( g') \overline{f_1 (\delta g' g) f_2 (g' g)} dg' , %a verifier
\end{eqnarray} 
o\`u $\delta \in G({\Bbb Q})$. \footnote{Dans le cas isotrope, le Th\'eor\`eme de Burger et Sarnak cit\'e ci-dessus est encore vrai tant que $f$ est uniform\'ement continue.}

Le coefficient (\ref{3.12}) est un coefficient matriciel associ\'e \`a un vecteur de l'espace $\oplus_{g \in G({\Bbb Q})} L^2 ((\Gamma \cap g^{-1} \Gamma g) \backslash G_0^{\rm nc})$.
La repr\'esentation $\sigma$ est donc faiblement contenue dans le $G_0^{\rm nc}$-spectre de $\oplus_{g \in G({\Bbb Q})} L^2 ((\Gamma \cap g^{-1} \Gamma g) \backslash G_0^{\rm nc})$, elle appartient donc au dual automorphe de $G$.

La d\'emonstration de \cite[Fact 3.3]{IRMN} montre plus pr\'ecisemment que $\sigma$ est faiblement contenue dans $\{ \rho \in \widehat{G}_{\rm Aut} \; : \; d(C^R (\rho )) = 0 \}$.
Puisque nous avons suppos\'e $\sigma$ isol\'ee dans le dual automorphe de $G$ sous la condition $d=0$, la Proposition \ref{P3.2'} est finalement d\'emontr\'ee.~$\Box$

\subsection{Restriction et cup-produit virtuels}

Dans le cas des groupes orthogonaux et de l'application de restriction nous d\'eduisons de tout ceci le th\'eor\`eme suivant.

\begin{thm} \label{anaO}
Soit $G$ un groupe alg\'ebrique r\'eductif, connexe et anisotrope sur ${\Bbb Q}$ tel que $G^{\rm nc} = O (p,q)$. Soit $Sh^0 H \subset Sh^0 G$ avec $H^{\rm nc} = O(p,q-r)$ plong\'e
de mani\`ere standard dans $G^{\rm nc}$ avec $p,q \geq 2$. Soit $i$ un entier $\leq (q-r-2)/2$ tel que $p+q-r-2i \geq 5$. Alors, l'application
$$H^{(i^p)} (Sh^0 G) \rightarrow \prod_{g \in G({\Bbb Q})} H^{*}_{{\rm prim} +} (Sh^0 H) $$
obtenue en composant l'application ${\rm Res}_H^G$ et la projection sur la composante fortement primitive de la cohomologie de $Sh^0 H$
est {\bf injective}. Son image est contenue dans $\prod_{g\in G({\Bbb Q})} H^{(i^p)} (Sh^0 H)$.
\end{thm}
{\it D\'emonstration.} D'apr\`es le Th\'eor\`eme \ref{res cohomO} et la Proposition \ref{P3.2}, il suffit de v\'erifier que la repr\'esentation $A((i^p))_H$ est isol\'ee dans le dual 
automorphe de $H$, mais d'apr\`es la Proposition \ref{Oisol} et pour $i$ v\'erifiant les conditions du Th\'eor\`eme, on a beaucoup plus puisque la repr\'esentation 
$A((i^p))_H$ est en fait isol\'ee dans le dual unitaire de $H_0^{\rm nc}$. Enfin la partie sur la composante fortement primitive de la cohomologie r\'esulte de la Remarque 
qui suit le Lemme \ref{resO}.~$\Box$

\bigskip

\noindent
{\bf Remarques.} Le Th\'eor\`eme \ref{anaO} (et sa d\'emonstration) implique(nt) 
le premier point du Th\'eor\`eme \ref{opq}~: il est imm\'ediat (cf. par exemple \cite[Fait 32]{BergeronTentative})
qu'une classe de cohomologie de degr\'e $\leq p+q-4$ appartient \`a $H^* ( A((i^p)) : Sh^0 G)$ ou \`a $H^* (A((q^j)) : Sh^0 G)$ (avec $i,j$ entiers). On peut donc
restreindre l'\'etude de Res$_H^G$ au sous-espace $H^* ( A((i^p)) : Sh^0 G)$, puis au sous-espace $H^{(i^p)} (Sh^0 G)$. Enfin le Corollaire \ref{mino} implique
que toute repr\'esentation cohomologique de $SO_0 (p,q-1)$ de degr\'e fortement primitif $< p+q-4$ est isol\'ee dans le dual unitaire de $SO_0 (p,q-1)$, ce qui 
implique imm\'ediatement que toute repr\'esentation cohomologique de $SO_0 (p,q-1)$ de degr\'e fortement primitif $\leq p+q-4$ est isol\'ee dans le dual unitaire
sous la  condition $d=0$. Le Th\'eor\`eme \ref{anaO} (et sa d\'emonstration) implique(nt) alors (avec $r=1$) que la restriction de Res$_H^G$ au sous-espace $H^{(i^p)} (Sh^0 G)$
est injective. Ce qui d\'emontre le premier point du Th\'eor\`eme \ref{opq}. Le Th\'eor\`eme \ref{anaO} n'implique qu'une version plus faible du premier point du Th\'eor\`eme \ref{o2n}. 
Celui-ci est n\'eanmoins d\'emontr\'e par Venkataramana \cite{Venky}.

Le Th\'eor\`eme \ref{anaO} est un analogue (bien que sa d\'emonstration soit compl\`etement diff\'erente) du Th\'eor\`eme \ref{analogue}. 
Plus g\'en\'eralement et au vu de la Conjecture \ref{CU2} nous conjecturons le r\'esultat suivant.

\begin{conj} \label{CanaO}
Soit $G$ un groupe alg\'ebrique r\'eductif, connexe et anisotrope sur ${\Bbb Q}$ tel que $G^{\rm nc} = O (p,q)$. Soit $Sh^0 H \subset Sh^0 G$ avec $H^{\rm nc} = O(p,q-r)$ plong\'e
de mani\`ere standard dans $G^{\rm nc}$, $1\leq p ,q$ et $1 \leq r <q$. Soit $\lambda$ une partition orthogonale dans $p\times q$.
Alors, l'application 
$$H^{\lambda} (Sh^0 G)_{\pm_1}^{\pm_2} \rightarrow \prod_{g \in G({\Bbb Q})} H^{|\lambda |}_{{\rm prim} +} (Sh^0 H) $$
obtenue en composant l'application ${\rm Res}_H^G$ et la projection sur la composante fortement primitive de la cohomologie de $Sh^0 H$
est {\bf injective} si et seulement si la partition $(r^p)$ s'inscrit dans le diagramme gauche $\hat{\lambda} / \lambda$. 
Son image est alors contenue dans $\prod_{g \in G({\Bbb Q})} H^{\lambda} (Sh^0 H)_{\pm_1 '}^{\pm_2 '}$.
\end{conj}

\bigskip

Concernant la restriction des groupes unitaires vers les groupes orthogonaux nous obtenons le r\'esultat suivant.

\begin{thm} \label{anaUO}
Soit $G$ un groupe alg\'ebrique r\'eductif, connexe et anisotrope sur ${\Bbb Q}$ tel que $G^{\rm nc} = U (p,q)$. 
Soit $Sh^0 H \subset Sh^0 G$ avec $H^{\rm nc} = O(p,q)$ plong\'e
de mani\`ere standard dans $G^{\rm nc}$ avec $p,q \geq 2$. Soit $i$ un entier $\leq (q-2)/2$ tel que $p+q-2i \geq 5$. 
Alors, l'application
$$H^{(i^p)} (Sh^0 G) \rightarrow \prod_{g \in G({\Bbb Q})} H^{*}_{{\rm prim} +} (Sh^0 H) $$
obtenue en composant l'application ${\rm Res}_H^G$ et la projection sur la composante fortement primitive de la cohomologie de $Sh^0 H$
est {\bf injective}. Son image est contenue dans $\prod_{g\in G({\Bbb Q})} H^{(i^p)} (Sh^0 H)$.
\end{thm}
{\it D\'emonstration.} Elle est identique \`a celle du Th\'eor\`eme \ref{anaO} \`a condition de remplacer le 
Th\'eor\`eme \ref{res cohomO} par le  Th\'eor\`eme \ref{res cohomUO} et la Remarque 
suivant le Lemme \ref{resO} par la Remarque 
qui suit le Lemme \ref{resUO}.~$\Box$

\bigskip

\noindent
{\bf Remarques.} Le Th\'eor\`eme \ref{anaO} (et sa d\'emonstration) implique(nt) 
les premiers points des Th\'eor\`emes \ref{u et o} et \ref{u et o2}.
Les deuxi\`emes points de ces Th\'eor\`emes proviennent de ce qu'\`a un niveau fini on a un plongement totalement r\'eel d'une 
vari\'et\'e r\'eelle de dimension paire dans une vari\'et\'e complexe. La multiplication par $J$, donn\'e par la structure complexe de
la vari\'et\'e ambiante, \'echange donc le fibr\'e tangent et le fibr\'e normal de la sous-vari\'et\'e totalement r\'eelle. Le
nombre d'Euler de son fibr\'e normal dans la vari\'et\'e complexe ambiante est donc \'egal \`a sa
caract\'eristique d'Euler qui est non nulle.

Le Th\'eor\`eme \ref{anaUO} (et sa d\'emonstration) implique(nt) le Corollaire \ref{CORUO} de l'Introduction, il faut juste 
remplacer la Proposition \ref{Oisol} par le Corollaire \ref{isolq}, lorsque $p=2$ et $q=3$, 
par la Proposition \ref{isolq1}, lorsque $p=1$, et utiliser 
\cite[Th\'eor\`eme 5.1]{STSG} lorsque $p=q=2$. 

Plus g\'en\'eralement et au vu de la Conjecture \ref{CU2} nous conjecturons le r\'esultat suivant.

\begin{conj} \label{CanaUO}
Soit $G$ un groupe alg\'ebrique r\'eductif, connexe et anisotrope sur ${\Bbb Q}$ tel que $G^{\rm nc} = U (p,q)$. 
Soit $Sh^0 H \subset Sh^0 G$ avec $H^{\rm nc} = O(p,q)$ plong\'e
de mani\`ere standard dans $G^{\rm nc}$, $1\leq p ,q$. Soient $\lambda$ et $\mu$ deux partitions incluses dans $p\times q$
formant un couple compatible. Si l'application 
\begin{eqnarray} \label{applUO}
H^{\lambda , \mu} (Sh^0 G)_{\pm_1}^{\pm_2} \rightarrow \prod_{g \in G({\Bbb Q})} H^{|\lambda |}_{{\rm prim} +} (Sh^0 H) 
\end{eqnarray}
obtenue en composant l'application ${\rm Res}_H^G$ et la projection sur la composante fortement primitive de la cohomologie de $Sh^0 H$
est non nulle alors $\lambda =0$ ou $\mu = p\times q$. Supposons par exemple $\mu = p\times q$. Alors, l'application 
(\ref{applUO}) est {\bf injective} si et seulement si la partition $\lambda$ est orthogonale.
Son image est alors contenue dans $\prod_{g \in G({\Bbb Q})} H^{\lambda} (Sh^0 H)_{\pm_1 '}^{\pm_2 '}$.
\end{conj}

\bigskip

Concernant le cup-produit nous obtenons le r\'esultat suivant.

\begin{thm} \label{cupO}
Soit $G$ un groupe alg\'ebrique r\'eductif, connexe et anisotrope sur ${\Bbb Q}$ tel que $G^{\rm nc} = O (p,q)$, avec $p,q \geq 2$. 
Soient $\alpha$ et $\beta$
deux classes de cohomologie appartenant respectivement \`a $H^{(k^p)} (Sh^0 G)$ et  $H^{(l^p)} (Sh^0 G)$ avec $k+l \leq (q-2)/2$ et $p+q-2(k+l) \geq 5$.
Il existe alors un \'el\'ement $g \in G({\Bbb Q})$ tel que le projet\'e de 
$$g(\alpha ) \wedge \beta \neq 0 $$
dans $H^{(k+l)^p} (Sh^0 G)$ soit non nul.
\end{thm}
{\it D\'emonstration.} D'apr\`es le Th\'eor\`eme \ref{pdt cohomO} et la Proposition \ref{P3.2'}, il suffit de v\'erifier que la repr\'esentation $A(((k+l)^p))$ est isol\'ee dans le dual 
automorphe de $H$, mais d'apr\`es la Proposition \ref{Oisol} et pour $k$, $l$ v\'erifiant les conditions du Th\'eor\`eme, on a beaucoup plus puisque la repr\'esentation 
$A(((k+l)^p))_H$ est en fait isol\'ee dans le dual unitaire de $G_0^{\rm nc}$.~$\Box$

\bigskip

\noindent
{\bf Remarques.} L\`a encore, le Th\'eor\`eme \ref{cupO} implique le Th\'eor\`eme \ref{cup o}, lorsque $p\geq 3$. Dans le cas $p=2$, il n'implique qu'une version plus faible. 
Le reste du Th\'eor\`eme \ref{cup o} est n\'eanmoins d\'emontr\'e par Venkataramana \cite{Venky}.
  
Le Th\'eor\`eme \ref{cupO} est un analogue (bien que la d\'emonstration soit compl\`etement diff\'erente) du Th\'eor\`eme \ref{cup u}. 
Plus g\'en\'eralement et au vu du Th\'eor\`eme \ref{pdt cohomO} nous 
conjecturons le r\'esultat suivant (qui g\'en\'eralise la Conjecture \ref{cup hyp}).

\begin{conj} \label{CanaO}
Soit $G$ un groupe alg\'ebrique r\'eductif, connexe et anisotrope sur ${\Bbb Q}$ tel que $G^{\rm nc} = O (p,q)$, avec $p,q \geq 1$. 
Soient $\alpha$ et $\beta$
deux classes de cohomologie appartenant respectivement \`a $H^{(k^p)} (Sh^0 G)$ et  $H^{(l^p)} (Sh^0 G)$..
Alors, il existe un \'el\'ement $g \in G({\Bbb Q})$ tel que le projet\'e de 
$$g(\alpha ) \wedge \beta \neq 0 $$
dans la partie fortement primitive de la cohomologie de $Sh^0 G$ soit non nul si et seulement si $k+l \leq q/2$. Le projet\'e appartient alors 
\`a $H^{(k+l)^p} (Sh^0 G)$. 
\end{conj}

\bigskip

\subsection{L'application ``cup-produit avec $[Sh^0 H]$''}

Concernant l'application ``cup-produit avec $[Sh^0 H]$'' nous commen\c{c}ons par d\'emontrer le r\'esultat g\'en\'eral suivant.

\begin{thm}    \label{sur l'homologie}
Soit $G$ un groupe alg\'ebrique r\'eductif, connexe et anisotrope sur ${\Bbb Q}$ tel que $G^{\rm nc} = O (p,q+r)$. Soit $Sh^0 H \subset Sh^0 G$ avec $H^{\rm nc} = O(p,q)$ plong\'e
de mani\`ere standard dans $G^{\rm nc}$ avec $p,q ,r\geq 1$. Soit $\lambda$ une partition orthogonale dans $p\times q$ telle que la partition $(r^p)$ s'inscrive
dans le diagramme $\hat{\lambda}/\lambda$. Supposons~:
\begin{enumerate}
\item que $\lambda$ v\'erifie la Conjecture \ref{conjl2O}, et
\item que $A(\lambda+(r^p))_{\pm_1}^{\pm_2}$ est isol\'ee dans le dual automorphe de $G$ sous la condition $d=0$.
\end{enumerate} 
Alors, l'application 
$$H^{\lambda} (Sh^0 H) \rightarrow H^{\lambda +(r^p)} (Sh^0 G),$$
obtenue en composant l'application ``cup-produit avec $[Sh^0 H]$'' et la projection sur la composante $H^{\lambda +(r^p)} (Sh^0 G) \subset H^* (Sh^0 G)$,
est {\bf injective}.

Si de plus $\alpha$ est une classe non triviale dans $H^{\lambda} (Sh^0 H)$, l'espace vectoriel engendr\'e par les translat\'es de Hecke de l'image de $\alpha$ dans
$H^{\lambda +(r^p)} (Sh^0 G)$ est de dimension infini.
\end{thm}
{\it D\'emonstration.} Soit $\alpha$ une classe non triviale dans $H^{\lambda} (Sh^0 H)$. La classe $\alpha$ est repr\'esent\'ee par un classe non triviale (que nous notons toujours 
$\alpha$) d\'efinie sur un niveau fini. On peut donc supposer qu'il existe un sous-groupe de congruence $\Gamma$ dans $G$ et un sous-espace 
$V$ de dimension $r$ dans ${\Bbb R}^{p+q+r}$ tels que $H^{\rm nc} = G^{\rm nc}_V$ et $\alpha \in H^{\lambda} (\Gamma_V \backslash X_V )$, o\`u
$\Gamma_V = \Gamma \cap G^{\rm nc}_V$. Repr\'esentons la classe $\alpha$ par une forme harmonique que nous notons $\mu$. 
(Nous utiliserons tout au long de la d\'emonstration les notations du \S 5.)
 
Notons $M= \Gamma \backslash X$ et fixons $\{ M_m = \Gamma_m \backslash X \}$ une tour d'effeuillage (fournie par le Lemme \ref{effeuillage}) constitu\'ee 
de rev\^etements de congruence.

\bigskip

Notons $\delta$ l'unique $K$-type minimal de la repr\'esentation $A(\lambda + (r^p))$ de $G$ et $P_0^{m, \delta }$ ($m\in {\Bbb N} \cup \{ \infty \}$) 
l'application obtenue en composant la projection $P_0^m$ (sur les formes harmoniques) par la projection sur le sous-espace $\delta$-isotypique $H^* (M_m)_{\delta}$ de la cohomologie
(\'eventuellement $L^2$). 

\begin{lem} \label{fait1} 
Soit $s \in {\Bbb C}$, Re$(s)>>0$. Alors, la suite $P_0^{m , \delta} \Omega_{\mu}^m (s)$ converge uniform\'ement 
sur tout compact de $X$ vers $P_0^{\infty , \delta} \Omega_{\mu} (s)$.
\end{lem}
{\it D\'emonstration.} Puisque par hypoth\`ese la repr\'esentation $A(\lambda+(r^p))_{\pm_1}^{\pm_2}$ est isol\'ee dans le dual automorphe de $G$ sous la condition $d=0$,
il existe un r\'eel $\nu$ strictement positif tel que la premi\`ere valeur propre non nulle du laplacien sur les 
$(k+pr)$-formes ferm\'ees de $M_m$ (resp. $M_{\infty}$) engendrant le $K$-type $\delta$ soit strictement sup\'erieure \`a $\nu$. 
Notons alors $h_{\nu}$ la fonction $\in C_0 ([0, +\infty [)$ qui vaut $1$ sur 
l'intervalle $[0, \frac{\nu }{2} ]$, $0$ sur l'intervalle $[\nu , +\infty [$ et qui d\'ecroit lin\'eairement sur 
$[\frac{\nu}{2} ,\lambda ]$. Puisque~:
\begin{enumerate}
\item la seule valeur propre du laplacien sur les $(k+pr)$-formes ferm\'ees de $M_m$ (resp. $M_{\infty}$) engendrant le $K$-type $\delta$
et strictement inf\'erieure \`a $\nu$ est $0$,
\item l'espace des formes ferm\'ees engendrant le $K$-type $\delta$ est ferm\'e, et
\item les formes $\Omega_{\mu}^m (s)$ et $\Omega_{\mu} (s)$ sont ferm\'ees,
\end{enumerate}
la projection sur la composante $\delta$-isotypique de $h_{\nu} (\Delta_m ) \Omega_{\mu}^m (s)$ (resp. $h_{\nu} (\Delta_{\infty} ) 
\Omega_{\mu} (s)$)
 est \'egale \`a $P_0^{m ,\delta} \Omega_{\mu}^m (s)$ (resp. $P_0^{\infty, \delta} \Omega_{\mu} (s)$).

Or, d'apr\`es la Proposition \ref{convergence des noyaux 2} et pour $s \in {\Bbb C}$, Re$(s)>>0$, la suite $h_{\nu} (\Delta_m )
\Omega_{\mu}^m (s)$ 
converge vers $h_{\nu } (\Delta_{\infty} )\Omega_{\mu} (s)$ uniform\'ement sur les compacts.
Ce qui conclut la d\'emonstration du Lemme \ref{fait1}.~$\Box$

\bigskip

D'apr\`es le Th\'eor\`eme \ref{dualite}, l'application $\mu \mapsto \Omega_{\mu} (s)$ (pour Re$(s)>>0$) induit l'application naturelle 
$$H^k (C_V) \mapsto H^{k+pr}_2 (M_V )$$
``cup-produit avec $[C_V]$''. Puisque par hypoth\`ese $\lambda$ v\'erifie la la Conjecture \ref{conjl2O}, la projection $P_0^{\infty , \delta} \Omega_{\mu} (s)$ doit \^etre non
nulle (et ind\'ependante de $s$, Re$(s)>>0$).

On peut maintenant d\'emontrer la premi\`ere partie du Th\'eor\`eme \ref{sur l'homologie}. 
L'application $H^{\lambda} (C_V ) \rightarrow H^{\lambda + (r^p)} (M_m )$ correspond en effet
\`a l'application $\mu \mapsto P_0^{m , \delta} \Omega_{\mu}^m (s)$. Mais, d'apr\`es le Lemme \ref{fait1}, cette derni\`ere converge simplement vers 
l'application
$$H^{\lambda} (C_V ) \rightarrow H_2^{\lambda + (r^p)} (M_V ),$$
injective par hypoth\`ese. Or $H^{\lambda} (C_V )$ est de dimension finie donc la convergence est uniforme et pour 
$m$ grand, l'application $H^{\lambda} (C_V ) \rightarrow H^{\lambda +(r^p)} (M_m )$ est injective. Comme $M_m$ est un rev\^etement de
congruence de 
$M$, la premi\`ere partie du Th\'eor\`eme \ref{sur l'homologie} est d\'emontr\'ee.

\bigskip

Pour d\'emontrer la deuxi\`eme partie du Th\'eor\`eme \ref{sur l'homologie}, nous allons d'abord d\'eduire des faits pr\'ec\'edents la proposition suivante.

\begin{prop} \label{translate}   
On se place sous les hypoth\`eses du Th\'eor\`eme \ref{sur l'homologie}.  Soient $s \in {\Bbb C}$, Re$(s)>>0$, et $\mu$ une forme harmonique repr\'esentant une classe non nulle dans
$H^{\lambda} (C_V)$.  Supposons $P_0^{0, \delta} \Omega_{\mu}^0 (s) \neq 0$. Il existe alors un entier $m \geq 0$ et un \'el\'ement $\gamma \in \Gamma$ tels que les formes 
harmoniques $P_0^{m , \delta} \Omega_{\mu}^m (s)$ et $\gamma^* P_0^{m , \delta} \Omega_{\mu}^m (s)$ soient lin\'eairement ind\'ependantes.
\end{prop}
{\it D\'emonstration.} Nous montrons d'abord par l'absurde qu'il existe un entier $m \geq 0$ tel que la forme 
$P_0^{m , \delta}  \Omega_{\mu}^m$ ne soit pas invariante sous l'action de $\Gamma$. Soit $m$ un entier $\geq 0$. Supposons que 
la forme $P_0^{m , \delta} \Omega_{\mu}^m (s)$  soit invariante sous l'action de $\Gamma$. Alors,
$$P_0^{0 , \delta} \Omega_{\mu}^{0} (s) = [\Gamma  :\Gamma_m ] P_0^{m , \delta} \Omega_{\mu}^m (s) .$$
Or $[\Gamma :\Gamma_m ]$ tend vers l'infini avec $m$, donc $P_0^{m , \delta} \Omega_{\mu}^m (s)$ tend vers $0$ avec $m$ ce qui 
contredit la premi\`ere partie du Th\'eor\`eme \ref{sur l'homologie} d\'emontr\'e ci-dessus. Il existe donc un entier $m \geq 0$ et un \'el\'ement $\gamma \in \Gamma$ tels que les formes 
$P_0^{m , \delta} \Omega_{\mu}^m (s)$ et $\gamma^* P_0^{m ,\delta} \Omega_{\mu}^m (s)$ soient distinctes. Puisque 
$$\sum_{g \in \Gamma_m \backslash \Gamma } g^* P_0^{m , \delta} \Omega_{\mu}^m (s) = \sum_{g\in \Gamma_m \backslash \Gamma} g^* (\gamma^* P_0^{m, \delta} \Omega_{\mu}^m (s) )=P_0^{0, \delta} \Omega_{\mu}^{0} (s) \neq 0, $$
les formes $P_0^{m , \delta} \Omega_{\mu}^m (s)$ et $\gamma^* P_0^{m , \delta} \Omega_{\mu}^m (s)$ sont en fait n\'ecessairement lin\'eairement 
ind\'ependantes. Ce qui ach\`eve la d\'emonstration de la Proposition \ref{translate}.~$\Box$

\bigskip

\`A l'aide de la Proposition \ref{translate}, nous pouvons maintenant conclure la d\'emonstration du Th\'eor\`eme \ref{sur l'homologie}.

Soit $\mu$ une forme harmonique sur $C_V$ repr\'esentant une classe non nulle dans $H^{\lambda} (C_V)$. Nous allons montrer par r\'ecurrence sur $N\geq 1$ qu'il 
existe un rev\^etement fini de congruence $M_N$ de $M$ et $N$ translat\'es par $\Gamma$ de l'image de $\mu$ dans $H^{\lambda + (r^p)} (M_N)$ qui soient lin\'eairement 
ind\'ependants. Le Th\'eor\`eme \ref{sur l'homologie} en d\'ecoule imm\'ediatement. 

Supposons qu'il existe un tel rev\^etement $M_N$ pour un certain $N\geq 1$.
Notons $\omega_1 , \ldots , \omega_N$ les $N$ formes harmoniques ind\'ependantes obtenues.
On peut supposer que la forme $\omega_1$ est la forme harmonique associ\'ee \`a $\mu$ (qui est une forme 
harmonique sur une pr\'eimage de $C_V$ dans $M_N$). La Proposition \ref{translate} implique qu'il existe un 
rev\^etement fini $M_{N+1}$ de $M_N$, une forme harmonique $\hat{\omega}_1$ sur $M_{N+1}$ et un \'el\'ement 
$\gamma \in \pi_1 M_N$ tels que les formes harmoniques $\hat{\omega}_1$ et $\gamma^* \hat{\omega}_1$ soient 
lin\'eairement ind\'ependantes. Supposons qu'il existe $N+1$ r\'eels $\alpha_0 , \alpha_1,\ldots ,\alpha_N$ tels que 
$$\alpha_0 \hat{\omega}_1 + \alpha_1 \gamma^* \hat{\omega}_1 +\alpha_2 \omega_2 +\cdots +\alpha_N \omega_N =0.$$

Alors en moyennant par $\pi_1 M_{N+1} \backslash \pi_1 M_{N}$, on obtient~:
$$(\alpha_1 +\alpha_2) \omega_1 + [\pi_1 M_N  : \pi_1 M_{N+1} ] \{ \alpha_2 \omega_2 +\cdots +\alpha_N \omega_N \} =0.$$
L'hypoth\`ese de r\'ecurrence implique donc~:
$$\alpha_0 +\alpha_1 =\alpha_2 =\ldots =\alpha_N =0.$$
Et puisque les formes harmoniques $\hat{\omega}_1$ et $\gamma^* \hat{\omega}_1$ sont lin\'eairement 
ind\'ependantes, on obtient finalement que~:
$$\alpha_0 =\alpha_1 =0.$$
Enfin puisque le rev\^etement de $M_{N+1}$ sur $M$ est fini, on peut le supposer galoisien, la forme 
$\gamma^* \hat{\omega}_1$ est alors d\'efinie sur $M_{N+1}$. Ce qui ach\`eve la r\'ecurrence et la d\'emonstration du 
Th\'eor\`eme \ref{sur l'homologie}.~$\Box$

\bigskip

Remarquons maintenant qu'en vertu de la Proposition \ref{Oisol} la deuxi\`eme hypoth\`ese du Th\'eor\`eme \ref{sur l'homologie} est v\'erifi\'ee par tout diagramme
orthogonal $\lambda$ de poids $| \lambda | \leq p+q+r-rp-3$, si $p,q \geq 3$ et de poids $|\lambda |\leq [(q+r)/2]-2r$, si $p=2 \leq q+r$. 
Il d\'ecoule par ailleurs du Th\'eor\`eme \ref{cohom l2} que la premi\`ere hypoth\`ese du Th\'eor\`eme \ref{sur l'homologie}
est v\'erifi\'ee pour tout diagramme orthogonal $\lambda$ de poids $| \lambda | \leq (q-pr)/2 -1$. Le Th\'eor\`eme \ref{sur l'homologie} implique donc le th\'eor\`eme suivant qui 
\`a son tour implique imm\'ediatement les deuxi\`emes points des Th\'eor\`emes \ref{opq} et \ref{o2n}.

\begin{thm}
Supposons fix\'ee une donn\'ee $Sh^0 H \subset Sh^0 G$ avec $H^{\rm nc} = O(p,q)$ plong\'e de mani\`ere standard dans $G^{\rm nc} = O(p,q+r)$, avec $p,q \geq 2$. 
Alors, pour tout degr\'e 
$k \leq \min (p+q+r-rp-3 , (q-pr)/2 -1)$,  l'application ``cup-produit avec $[Sh^0 H]$'' 
$$\bigwedge_H^G : H^k (Sh^0 H) \rightarrow H^{k+rp} (Sh^0 G) $$
est {\bf injective}.
\end{thm}

Plus g\'en\'eralement, et au vu de la Conjecture \ref{conjl2O}, nous conjecturons le r\'esultat suivant.

\begin{conj} \label{C100}
Supposons fix\'ee une donn\'ee $Sh^0 H \subset Sh^0 G$ avec $H^{\rm nc} = O(p,q)$ plong\'e de mani\`ere standard dans $G^{\rm nc} = O(p,q+r)$. Si
$\lambda$ est une partition incluse dans $p\times q$, alors l'application 
$$H^{\lambda} (Sh^0 H) \rightarrow H^{|\lambda| + rp}_{{\rm prim} +} (Sh^0 G)$$
obtenue en composant l'application ``cup-produit avec $[Sh^0 H]$'' et la projection sur la composante fortement primitive de la cohomologie de $Sh^0 G$ est 
{\bf injective} si et seulement si la partition $(r^p)$ s'inscrit dans le diagramme $\hat{\lambda} / \lambda$. Son image est alors contenue dans 
$H^{\lambda + (r^p)} (Sh^0 G)$.
\end{conj} 

\bigskip

\noindent
{\bf Remarque.} L'analogue de la partie de l'\'enonc\'e de la Conjecture \ref{C100} concernant la composante fortement primitive devrait pouvoir \^etre d\'emontr\'ee
sous les hypoth\`ese du Th\'eor\`eme \ref{sur l'homologie}. Nous nous contenterons de traiter (plus loin) le cas des symboles modulaires.

\bigskip

\subsection{Applications}

Commen\c{c}ons par d\'eduire des Th\'eor\`emes \ref{opq} et \ref{o2n} le corollaire suivant.

\begin{cor} \label{bp}
Soit $G$ un groupe alg\'ebrique et anisotrope sur ${\Bbb Q}$ obtenu par restriction des scalaires \`a partir d'un groupe orthogonal sur un corps de nombres totalement r\'eel 
et tel que $G^{\rm nc} = O(p,q)$ avec $q\geq p \geq 2$. Alors , 
$$H^p (Sh^0 G) \neq 0.$$
\end{cor}

\bigskip

\noindent
{\bf Remarques.} Le Corollaire \ref{bp} reste vrai pour $p, q \geq 1$, cela d\'ecoule en particulier du Corollaire \ref{C1} sur lequel nous revenons au paragraphe pr\'ec\'edent. 
Ce r\'esultat est un Th\'eor\`eme de Millson et Raghunathan qui montrent plus g\'en\'eralement que tous les groupes $H^{kp} (Sh^0 G)$, $k= 0, \ldots ,q$, sont non triviaux. 
Nous d\'eduisons ici directement le Corollaire \ref{bp} des Th\'eor\`emes \ref{opq} et \ref{o2n} pour illustrer comment ceux-ci s'appliquent dans un cas relativement simple.

\bigskip

\noindent
{\it D\'emonstration du Corollaire \ref{bp}.} Puisque le groupe $G$ provient d'un groupe orthogonal sur un corps de nombre, il est imm\'ediat qu'il existe une donn\'ee 
$Sh^0 H \subset Sh^0 G$ avec $H^{\rm nc} = O(p,q-1)$. Remarquons d'abord que pour $q$ grand par rapport \`a $p$ (plus 
pr\'ecisemment pour $q \geq p+3$), les points 2. (avec $k=0$) des Th\'eor\`emes \ref{opq} et \ref{o2n} impliquent la conclusion du Corollaire \ref{bp}. Il nous reste \`a montrer
comment faire diminuer $q$, autrement dit d\'eduire du r\'esultat pour $O(p,q+1)$ le r\'esultat pour $O(p,q)$, pour tout $q\geq p \geq 2$.

Mais, et toujours puisque le groupe $G$ provient d'un groupe orthogonal sur un corps de nombre, il existe une donn\'ee 
$Sh^0 G \subset Sh^0 G_1$ avec $G_1^{\rm nc} = O(p,q+1)$. Remarquons maintenant que si $q\geq p \geq 3$, alors $p \leq p+ (q+1) -4$ et que si $p=2$, alors 
$2 \leq (q+1)-1$. On peut donc appliquer les points 1., avec $k=p$, des Th\'eor\`emes \ref{opq} et \ref{o2n} tant que $q\geq p \geq2$. D'o\`u l'on d\'eduit que si 
$H^p (Sh^0 G_1)$ est non nul alors $H^p (Sh^0 G) \neq 0$. Et le Corollaire \ref{bp} d\'ecoule d'un simple argument de descente.~$\Box$

\bigskip

Les Th\'eor\`emes \ref{opq} et \ref{o2n} sont bien s\^ur bien plus g\'en\'eraux et deviennent plus surprenant lorsqu'on les applique \`a des classes de cohomologie plus profondes.
Dans les vari\'et\'es kaehl\'eriennes le Th\'eor\`eme de Lefschetz fort implique que la multiplication par la forme de Kaehler propage injectivement les classes de cohomologie. 
De mani\`ere surprenante un ph\'enom\`ene analogue \`a lieu dans la cohomologie $H^* (Sh^0 O(p,q))$. Explicitons par exemple ce qu'il se passe \`a partir 
d'une classe de degr\'e $p$.

\begin{cor} \label{fort}
Supposons fix\'ees des donn\'ees $Sh^0 H \subset Sh^0 G$ avec $H^{\rm nc} = O(p,q-1)$, $G^{\rm nc} = O(p,q)$ et $p,q\geq 2$. Alors, pour tout 
$k\leq q/4$, l'application ``cup-produit avec $[Sh^0 H]^k$''
$$H^p (Sh^0 G ) \rightarrow \prod_{g \in G({\Bbb Q})} H^{p(k+1)} (Sh^0 G)$$ 
est {\bf injective}.
\end{cor} 

\bigskip

Le Corollaire \ref{CORUO} permet de d\'emontrer un r\'esultat bien plus g\'en\'eral que le Corollaire \ref{bp} \`a savoir, 
le Corollaire \ref{application} de l'Introduction. Raghunathan
et Venkataramana \cite{RaghunathanVenky} montrent en effet que si $G$ est un groupe alg\'ebrique obtenu par restriction des scalaires \`a
partir d'un groupe de type $D_n$, $\neq {}^{3,6} D_4$ et $n>2$, sur un corps de nombres totalement r\'eel de telle mani\`ere que 
$G^{{\rm nc}}=O(p,q)$, avec $1\leq p \leq q$, il existe alors un groupe $G_1$ obtenu par restriction des scalaires \`a partir
d'un groupe unitaire sur un corps de nombres totalement r\'eel et tel que $G_1^{{\rm nc}} = U(p,q)$ et 
$Sh^0 G \subset Sh^0 G_1$. Mais dans \cite{Wallach} Wallach montre (\`a l'aide de la correspondance th\'eta globale) que   
$H^p (Sh^0 G_1 )  \neq 0$. Le Corollaire \ref{CORUO} implique donc imm\'ediatement le Corollaire \ref{application} de 
l'introduction.

\subsection{Sur la classe de cohomologie des symboles modulaires}

Concernant la classe de cohomologie des symboles modulaires, nos m\'ethodes permettent de d\'emontrer le Th\'eor\`eme suivant.

\begin{thm} \label{symbmodul}
Supposons fix\'ees des donn\'ees $Sh^0 H \subset Sh^0 G$ avec $G$ anisotrope, $H^{\rm nc} = O(p,q)$, $G^{\rm nc} = O(p,q+r)$ et $p,q \geq 2$.
Supposons $\mathbf{q \geq r+2}$ {\bf et} $\mathbf{p+q-r \geq 5}$. Alors, la classe de $[Sh^0 H]$ est non triviale dans $H^* (Sh^0 G)$. 
Plus pr\'ecis\'ement, la projection de $[Sh^0 H]$ dans la composante fortement primitive de la cohomologie de $Sh^0 G$ est non triviale et appartient \`a
$H^{(r^p)} (Sh^0 G)$.
\end{thm}
{\it D\'emonstration.} Le Corollaire \ref{cohoml2O} implique l'hypoth\`ese 1. du Th\'eor\`eme \ref{sur l'homologie} pour la partition nulle. 
La Proposition \ref{Oisol} implique que sous les conditions $q\geq r+2$ et $p+q-r \geq 5$, l'hypoth\`ese 2. du Th\'eor\`eme \ref{sur l'homologie} est v\'erifi\'ee par la partition nulle.
Le Th\'eor\`eme \ref{sur l'homologie} implique alors que la projection de $[Sh^0 H]$ dans $H^{(r^p)} (Sh^0 G)$ est non triviale. Ce qui implique le premier point du Th\'eor\`eme
\ref{symbmodul}. La trivialit\'e de la projection de $[Sh^0 H]$ dans les autre composantes fortement primitives de la cohomologie de $Sh^0 G$, provient de ce que parmi les 
repr\'esentations cohomologiques appartenant \`a la s\'erie discr\`ete de $G^{\rm nc} / H^{\rm nc}$ (qui sont connues, cf. \S 2.6), 
$A((r^p))$ est la seule de degr\'e fortement primitif $pr$.~$\Box$

\bigskip

\noindent
{\bf Remarque.} Dans ce cas on peut \'egalement obtenir des r\'esultats sur la composante non fortement primitive de la cohomologie de $[Sh^0 H]$~: 
la projection de la classe $[Sh^0 H]$ dans $H^{pr} ( A(\lambda ) : Sh^0 G)$ est triviale pour toute partition non nulle $\lambda \subset [p/2] \times (q+r)$. Cela d\'ecoule en effet
d'un crit\`ere \cite[Theorem 2.8]{KobayashiOda} de Kobayashi et Oda, et du Th\'eor\`eme \ref{kobaO}.

Enfin, remarquons que le Corollaire \ref{isolq}, la Proposition \ref{isolq1} et la d\'emonstration du Th\'eor\`eme \ref{symbmodul} impliquent le Corollaire \ref{C1} de l'Introduction,

\bigskip

Dans le cas unitaire, l'analoque du premier point du Th\'eor\`eme \ref{symbmodul} est bien s\^ur imm\'ediat (puisque les symboles modulaires sont alors des
sous-vari\'et\'es complexes d'une vari\'et\'e kaehl\'erienne). Il n'est n\'eanmoins pas clair que la classe de cohomologie d'un symbole modulaire ne se projette pas  
trivialement sur la composante fortement primitive (ou m\^eme non triviale) de la cohomologie.
Le th\'eor\`eme suivant se d\'emontre exactement de la m\^eme mani\`ere que dans le cas orthogonal.

\begin{thm} \label{symbmodulU}
Supposons fix\'ees des donn\'ees $Sh^0 H \subset Sh^0 G$ avec $G$ anisotrope, $H^{\rm nc} = U(p,q)$, $G^{\rm nc} = U(p,q+r)$ et $p,q \geq 2$.
Supposons $\mathbf{q \geq r+2}$. Alors, la projection de $[Sh^0 H]$ dans la composante fortement primitive de la cohomologie de $Sh^0 G$ est non triviale 
et appartient \`a $H^{(r^p), (q^p)} (Sh^0 G)$. 
\end{thm}

\bigskip

\noindent
{\bf Remarque.} Dans ce cas le crit\`ere de Kobayashi et Oda montre (d'apr\`es le Th\'eor\`eme \ref{kobaU}) que la projection de $[Sh^0 H]$ dans 
$H^{2pr} ( A(\lambda ,\mu ) : Sh^0 G)$ est triviale pour tout couple $(\lambda , \mu )\neq (0 , p\times (q+r))$ de partitions telles que deux \'el\'ements de 
$\lambda$ et $\hat{\mu}$ ne soient jamais align\'es. On retrouve en particulier que la classe de cohomologie $[Sh^0 H]$ n'a pas de composante holomorphe. 

\section{G\'en\'eralisations et Perspectives}

\subsection{G\'en\'eralisations}

Les r\'esultats de la section pr\'ec\'edente se g\'en\'eralisent de fa\c{c}on assez naturelle dans deux directions~: (1) lorsque les syst\`emes de coefficients (de la cohomologie)
sont non triviaux, et (2) lorsque les groupes sont isotropes.

\bigskip

Le cas (1) est relativement imm\'ediat. Soit $E$ une repr\'esentation de dimension finie de $G$. La repr\'esentation $E$ d\'efinit un syst\`eme local ${\cal E}$ sur tous les quotients
$\Gamma \backslash X_G$ consid\'er\'es dans cet article. Supposons pour simplifier $E$ irr\'eductible. Toujours d'apr\`es la th\'eorie de Vogan et Zuckerman, l'alg\`ebre 
gradu\'ee $H^* (Sh^0 G , {\cal E} )$ peut \^etre d\'ecompos\'ee selon des repr\'esentations $A_{\mathfrak{q}} (E)$ associ\'ees aux sous-alg\`ebres paraboliques de $\mathfrak{g}$
(voir \cite{VoganZuckerman}). Les r\'esultats des sections pr\'ec\'edentes se traduisent alors mot \`a mot. Dans le cas de l'application ``cup-produit avec $[Sh^0 H]$'' il faut quand 
m\^eme penser \`a tordre la classe fondamentale $[Sh^0 H]$. Plus pr\'ecisemment et en se pla\c{c}ant \`a un niveau fini $\Gamma$, on consid\`ere toujours la sous-vari\'et\'e
$(\Gamma \cap H) \backslash X_H \rightarrow \Gamma \backslash X_G$, que l'on note $C_{\Gamma}^H$. Il s'agit alors de construire une section $s$ du fibr\'e 
${\cal E}_{|C_{\Gamma}^H}$. Ceci est fait par Tong et Wang dans \cite{TongWang}, dont on peut \'egalement extraire (comme au \S 2.4) la construction de la classe duale
\`a $(C_{\Gamma}^H , s)$ dans $H_2^{d_G - d_H } (\Gamma \backslash X_G , {\cal E})$. La reste se g\'en\'eralise imm\'ediatement. Notons m\^eme que
concernant les symboles modulaires les d\'emonstrations 
se simplifient dans nombre de cas o\`u la forme duale construite par Tong et Wang est $L^1$. On peut en effet alors former directement une s\'erie de Poincar\'e convergente
(sans avoir recours \`a un poids et donc au param\`etre $s$ du \S 5.6) et se passer de l'hypoth\`ese 2. dans le Th\'eor\`eme \ref{sur l'homologie}. Ceci explique les
constructions par Tong et Wang de classes de cohomologie non triviales (pour des syst\`emes de coefficients non triviaux). Nous avons pr\'ef\'er\'e n\'eglig\'e ici les 
syst\`emes de  coefficients, ils nous paraissent en effet un peu hors sujet en ce qui concerne les propri\'et\'es de Lefschetz.

\bigskip

Le cas (2) est plus d\'elicat et tr\`es int\'eressant. Consid\'erons donc maintenant un groupe $G$ isotrope sur ${\Bbb Q}$. La cohomologie de $Sh^0 G$ n'est plus 
naturellement reli\'ee \`a la th\'eorie des repr\'esentations, il n'y a plus de d\'ecomposition de Hodge ou de Vogan-Zuckerman. Il est dans ce contexte plus naturel de 
consid\'erer la cohomologie $L^2$, $H^*_2 (Sh^0 G)$, ou encore la cohomologie cuspidale $H^*_{{\rm cusp}} (Sh^0 G)$. L'espace $H^*_{{\rm cusp}}  (Sh^0 G)$ est un 
sous-espace de $H^* (Sh^0 G)$ comme de $H^*_2 (Sh^0 G)$, celui des classes repr\'esent\'ees par des formes cuspidales (qui sont born\'ees et donc $L^2$). Les 
th\'eories de Matsushima et de Vogan-Zuckerman s'appliquent \`a l'espace $H^*_2 (Sh^0 G)$. Le sous-espace $H^*_{{\rm cusp}} (Sh^0 G)$ h\'erite alors de la d\'ecomposition
de Vogan-Zuckerman.

Commen\c{c}ons par consid\'erer l'application de restriction stable de $G$ \`a $H$ et adoptons les notations du \S 8.1. Soit donc $\omega_{\varphi} \in H^R_2 (\pi : \Gamma )$.
Remarquons que d'apr\`es le Lemme \ref{borne} la forme harmonique $\omega_{\varphi}$ est born\'ee sur $\Gamma \backslash X_G$. La restriction de la
forme $\omega_{\varphi}$, via les applications $j_g$ ($g \in G({\Bbb Q})$), aux sous-vari\'et\'es $(H \cap g^{-1} \Gamma g) \backslash X_H$ sont born\'ees et 
donc dans $L^2$. L'application 
$${\rm Res}_{2, \; H}^G : H_2^* (Sh^0 G) \rightarrow \prod_{g \in G({\Bbb Q})} H^*_2 (Sh^0 H)$$
est bien d\'efinie.
Si de plus  $\omega_{\varphi} \in H_{\rm cusp}^R (\pi : \Gamma)$, un r\'esultat de Clozel et Venkataramana \cite[Lemma 2.9]{ClozelVenky} affirme que la forme 
$j_g^* \omega_{\varphi}$ est rapidement d\'ecroissante le long de $(H \cap g^{-1} \Gamma g) \backslash X_H$, et d\'efini en particulier une forme cuspidale. 
L'application 
$${\rm Res}_{{\rm cusp} , \; H}^G : H_{\rm cusp}^* (Sh^0 G) \rightarrow \prod_{g \in G({\Bbb Q})} H^*_{\rm cusp} (Sh^0 H)$$
est bien d\'efinie, elle est induite par la restriction de ${\rm Res}_{2 , \; H}^G$ au sous-espace $H_{\rm cusp}^* (Sh^0 G)$. 

La question de l'injectivit\'e des applications ${\rm Res}_{2, \; H}^G$ et ${\rm Res}_{{\rm cusp} , \; H}^G$ est donc bien pos\'ee. Les techniques de \cite{BergeronTentative}
ne se g\'en\'eralisent pas au cas isotrope (non-compact), le Th\'eor\`eme \ref{analogue} avec ${\rm Res}_H^G$ remplac\'e par ${\rm Res}_{2, \; H}^G$ ou 
${\rm Res}_{{\rm cusp} , \; H}^G$ est ouvert. Suivons donc la m\^eme approche (\`a savoir la d\'emonstration du Th\'eor\`eme \ref{anaO})
pour les groupes orthogonaux et unitaires dans le cas isotrope. 
La seule \'etape de la d\'emonstration du Th\'eor\`eme \ref{anaO} qui utilise l'hypoth\`ese $G$ anisotrope sur ${\Bbb Q}$ est la Proposition \ref{P3.2}, c'est \`a dire de \cite[Proposition 3.2]{IRMN}. Comme le 
remarquent Harris et Li dans \cite{HarrisLi} la d\'emonstration de la Proposition \ref{P3.2} reste valable dans le cas isotrope \`a condition de v\'erifier 
que chaque forme $f$ dans l'image de $\varphi$ est uniform\'ement continue (cf. la note de bas page no. 12 pour une remarque analogue). Mais ceci (et m\^eme $f$ uniform\'ement
H\"older) peut \^etre d\'eduit du Lemme \ref{borne} en utilisant une ``in\'egalit\'e de Harnak'' comme dans \cite[Proposition 9.4 (Chapitre 14)]{Taylor}.
(Le cas cuspidal est \'evidemment plus facile \`a traiter puisque les formes $j_g^* \omega_{\varphi}$ sont rapidement d\'ecroissantes 
le long de $(H \cap g^{-1} \Gamma g) \backslash X_H$.)

La d\'emonstration du Th\'eor\`eme \ref{anaO} implique alors les th\'eor\`emes suivants dans le cas isotrope.

\begin{thm} \label{anaO'}
Soit $G$ un groupe alg\'ebrique r\'eductif et connexe sur ${\Bbb Q}$ tel que $G^{\rm nc} = O (p,q)$. Soit $Sh^0 H \subset Sh^0 G$ avec $H^{\rm nc} = O(p,q-r)$ plong\'e
de mani\`ere standard dans $G^{\rm nc}$ avec $p,q \geq 2$. Soit $i$ un entier $\leq (q-r-2)/2$ tel que $p+q-r-2i \geq 5$. Alors, l'application
$$H^{(i^p)}_2 (Sh^0 G) \rightarrow \prod_{g \in G({\Bbb Q})} H^{*}_{2, \; {\rm prim} +} (Sh^0 H) $$
obtenue en composant l'application ${\rm Res}_{2 , \; H}^G$ et la projection sur la composante fortement primitive de la cohomologie $L^2$ de $Sh^0 H$
est {\bf injective}. Son image est contenue dans $\prod_{g\in G({\Bbb Q})} H_2^{(i^p)} (Sh^0 H)$.
\end{thm}

\begin{thm} \label{anaU'}
Soit $G$ un groupe alg\'ebrique r\'eductif et connexe sur ${\Bbb Q}$ tel que $G^{\rm nc} = U (p,q)$. Soit $Sh^0 H \subset Sh^0 G$ avec $H^{\rm nc} = U(p,q-r)$ plong\'e
de mani\`ere standard dans $G^{\rm nc}$ avec $p,q \geq 2$. Soient $i$ et $j$ deux entiers naturels de somme $i+j \leq q-r-2$. 
Alors, l'application
$$H^{(i^p),((q-j)^p)}_2 (Sh^0 G) \rightarrow \prod_{g \in G({\Bbb Q})} H^{*}_{2, \; {\rm prim} +} (Sh^0 H) $$
obtenue en composant l'application ${\rm Res}_{2 , \; H}^G$ et la projection sur la composante fortement primitive de la cohomologie $L^2$ de $Sh^0 H$
est {\bf injective}. Son image est contenue dans $\prod_{g\in G({\Bbb Q})} H_2^{(i^p) , ((q-r-j)^p)} (Sh^0 H)$.
\end{thm}

\bigskip

\noindent
{\bf Remarque.} Dans le cas unitaire, et en ce qui concerne la cohomologie holomorphe l'application ${\rm Res}_{2, \; H}^G$ est bien comprise, d'apr\`es les 
travaux de Clozel et Venkataramana \cite{ClozelVenky}. 

\bigskip

Concernant la restriction des groupes unitaires vers les groupes orthogonaux nous obtenons de la m\^eme mani\`ere le r\'esultat suivant.

\begin{thm} \label{anaUO'}
Soit $G$ un groupe alg\'ebrique r\'eductif et connexe sur ${\Bbb Q}$ tel que $G^{\rm nc} = U (p,q)$. Soit $Sh^0 H \subset Sh^0 G$ 
avec $H^{\rm nc} = O(p,q)$ plong\'e de mani\`ere standard dans $G^{\rm nc}$ avec $p,q \geq 2$. Soit $i$ un entier $\leq (q-2)/2$ 
tel que $p+q-2i \geq 5$. Alors, l'application
$$H^{(i^p)}_2 (Sh^0 G) \rightarrow \prod_{g \in G({\Bbb Q})} H^{*}_{2, \; {\rm prim} +} (Sh^0 H) $$
obtenue en composant l'application ${\rm Res}_{2 , \; H}^G$ et la projection sur la composante fortement primitive 
de la cohomologie $L^2$ de $Sh^0 H$
est {\bf injective}. Son image est contenue dans $\prod_{g\in G({\Bbb Q})} H_2^{(i^p)} (Sh^0 H)$.
\end{thm}

\bigskip

\noindent
{\bf Remarque.} Le Th\'eor\`eme \ref{anaUO'} (et sa d\'emonstration) implique(nt), comme annonc\'ee dans l'Introduction, que le
Corollaire \ref{application} reste vrai dans le cas isotrope avec la m\^eme d\'emonstration.

\bigskip

Concernant le cup-produit nous obtenons de mani\`ere analogue les r\'esultats suivants.

\begin{thm} \label{cupO'}
Soit $G$ un groupe alg\'ebrique r\'eductif et connexe sur ${\Bbb Q}$ tel que $G^{\rm nc} = O (p,q)$, avec $p,q \geq 2$. 
Soient $\alpha$ et $\beta$
deux classes de cohomologie $L^2$ appartenant respectivement \`a $H_2^{(k^p)} (Sh^0 G)$ et  $H_2^{(l^p)} (Sh^0 G)$ avec $k+l \leq (q-2)/2$ et $p+q-2(k+l) \geq 5$.
Il existe alors un \'el\'ement $g \in G({\Bbb Q})$ tel que le projet\'e de 
$$g(\alpha ) \wedge \beta \neq 0 $$
dans $H_2^{((k+l)^p)} (Sh^0 G)$ soit non nul.
\end{thm}

\begin{thm} \label{cupU'}
Soit $G$ un groupe alg\'ebrique r\'eductif et connexe sur ${\Bbb Q}$ tel que $G^{\rm nc} = U (p,q)$, avec $p,q \geq 2$. 
Soient $\alpha$ et $\beta$
deux classes de cohomologie $L^2$ appartenant respectivement \`a $H_2^{(i^p) , (q-j)^p} (Sh^0 G)$ et  $H_2^{(k^p) , (q-l)^p} (Sh^0 G)$ avec $i+j+k+l \leq q-2$.
Il existe alors un \'el\'ement $g \in G({\Bbb Q})$ tel que le projet\'e de 
$$g(\alpha ) \wedge \beta \neq 0 $$
dans $H_2^{((i+k)^p) , ((q-j-l)^p)} (Sh^0 G)$ soit non nul.
\end{thm}

\bigskip

Consid\'erons maintenant l'application ``cup-produit avec $[Sh^0 H]$''. La seule \'etape de la d\'emonstration du Th\'eor\`eme \ref{sur l'homologie} 
qui utilise de mani\`ere essentielle le fait que $G$ est anisotrope est la d\'emonstration du Th\'eor\`eme \ref{cohom l2}. La g\'en\'eralisation de ce r\'esultat
au cas isotrope semble n\'ecessiter des id\'ees nouvelles ou \`a tout le moins une description plus fine de la g\'eom\'etrie \`a l'infini de $M_V$,
o\`u $\Lambda$ est un r\'eseau non cocompact de $H$. N\'eanmoins il est facile de v\'erifier que la non isotropie de $G$ n'est pas n\'ecessaire dans les d\'emonstrations 
des Corollaires \ref{conjl2} et \ref{conjl2O}, ni dans dans le \S 5.6 (la construction de Wang est valable pour $\Gamma_V \backslash X_V$ de volume fini). On peut donc
au moins appliquer notre m\'ethode aux symboles modulaires pour obtenir les Th\'eor\`emes suivants.

\begin{thm} \label{symbmodul'}
Supposons fix\'ees des donn\'ees $Sh^0 H \subset Sh^0 G$ avec $H^{\rm nc} = O(p,q)$, $G^{\rm nc} = O(p,q+r)$ et $p,q \geq 2$.
Supposons $\mathbf{q \geq r+2}$ {\bf et} $\mathbf{p+q-r \geq 5}$. Alors, la classe $[Sh^0 H]$ est non triviale dans $H_2^* (Sh^0 G)$. 
Plus pr\'ecis\'ement, la projection de $[Sh^0 H]$ dans la composante fortement primitive de la cohomologie $L^2$ de $Sh^0 G$ est non triviale et appartient \`a
$H_2^{(r^p)} (Sh^0 G)$.
\end{thm}

\begin{thm} \label{symbmodulU}
Supposons fix\'ees des donn\'ees $Sh^0 H \subset Sh^0 G$ avec $H^{\rm nc} = U(p,q)$, $G^{\rm nc} = U(p,q+r)$ et $p,q \geq 2$.
Supposons $\mathbf{q \geq r+2}$. Alors, la classe $[Sh^0 H]$ est non triviale dans $H_2^* (Sh^0 G)$. Plus 
pr\'ecis\'ement, la projection de $[Sh^0 H]$ dans la composante fortement primitive de la cohomologie $L^2$ de $Sh^0 G$ est non triviale 
et appartient \`a $H_2^{(r^p), (q^p)} (Sh^0 G)$. 
\end{thm}

\bigskip

Dans un preprint r\'ecent \cite{SpehVenky}, Speh et Venkataramana d\'emontrent que si $Sh^0 H \subset Sh^0 G$ avec $H^{\rm nc} =
U(1,q)$, $G^{\rm nc} = U(1,q+r)$ et $r=1$ ou $2$, alors la classe $[Sh^0 H]$ est non triviale dans $H^* (Sh^0 G)$ et engendre
sous l'action des op\'erateurs de Hecke un espace dimension infini. 
Ils montrent en fait la non trivialit\'e de la projection de la classe
$[Sh^0 H]$ dans la cohomologie triviale. De mani\`ere compl\'ementaire notre m\'ethode permet de d\'emontrer (au moins lorsque $r=1$) 
la non trivialit\'e de la projection de $[Sh^0 H]$ dans la cohomologie fortement primitive. Nous montrons plus pr\'ecis\'ement
le th\'eor\`eme suivant.

\begin{thm}
Supposons fix\'ees des donn\'ees $Sh^0 H \subset Sh^0 G$ avec $H^{\rm nc} = O(1,q)$ (resp. $=U(1,q)$), $G^{\rm nc} = O(1,q+1)$
(resp. $U(1,q+1)$) et $q\geq 2$ (resp. $q\geq 1$). Alors, la projection de la classe $[Sh^0 H]$ dans la composante 
fortement primitive de $H_2^* (Sh^0 G)$ est non triviale et le sous-espace engendr\'e par ses translat\'es de Hecke est de dimension 
infini. En particulier, la classe $[Sh^0 H]$ est non triviale dans $H^* (Sh^0 G)$ en engendre   
sous l'action des op\'erateurs de Hecke un espace de dimension infini.
\end{thm}
{\it D\'emonstration.} La d\'emonstration est identique \`a celles des Th\'eor\`emes \ref{symbmodul'} et \ref{symbmodulU}, sauf
que l'on ne peut plus appliquer les Propositions \ref{Oisol} et \ref{Uisol}. Celles-ci peuvent n\'eanmoins 
\^etre respectivement remplac\'ees par la Proposition \ref{isolq1} et \ref{isolq2}.
L'assertion sur les translat\'es de Hecke d\'ecoule de la d\'emonstration du Th\'eor\`eme \ref{sur l'homologie}.
On passe finalement de la cohomologie $L^2$ \`a la cohomologie usuelle \`a l'aide d'un Th\'eor\`eme de Zucker \cite{Zucker}.~$\Box$

\bigskip

\subsection{Perspectives}

Il n'y a \'evidemment aucune raison autre que technique pour restreindre l'\'etude des propri\'et\'es de Lefschetz automorphes
au cas des groupes unitaires ou orthogonaux. Conluons alors cet article par deux conjectures (peut-\^etre un peu 
optimistes et ambitieuses) qui d\'ecrivent
les propri\'et\'es de Lefschetz automorphes auquelles on s'attend pour des groupes g\'en\'eraux.

La premi\`ere conjecture concerne l'application de restriction stable. Elle est motiv\'ee par les Th\'eor\`emes \ref{analogue} et 
\ref{anaO}, l'article \cite{BergeronTentative} 
et (surtout) par un r\'esultat g\'en\'eral de Venkataramana \cite[Theorem 6]{Venky} dans le cas hermitien.
Nous avons besoin de quelques pr\'eliminaires pour \'enoncer cette conjecture.

Condid\'erons $G$ un groupe alg\'ebrique r\'eductif et connexe sur ${\Bbb Q}$, presque simple sur ${\Bbb Q}$ modulo 
son centre. Nous supposons de plus que que le groupe
$G({\Bbb R})$ de ses points r\'eels est semi-simple et non compact modulo un centre compact. Soit 
$H \subset G$ un sous-groupe r\'eductif connexe d\'efini sur ${\Bbb Q}$ dont le groupe $H({\Bbb R})$ intersecte 
un compact maximal $K$ de $G({\Bbb R})$ selon un sous-groupe compact maximal de $H({\Bbb R})$. 
Soient $\mathfrak{g}=\mathfrak{k} \oplus \mathfrak{p}$ la d\'ecomposition de Cartan du complexifi\'e de l'alg\`ebre de Lie de 
$G({\Bbb R})$ et $\mathfrak{k}_H \oplus \mathfrak{p}_H$ la d\'ecomposition correspondante pour $H$. 
Notons $\mathfrak{s}$ le suppl\'ementaire orthogonal (pour la forme de Killing) de $\mathfrak{h}$ dans $\mathfrak{g}$.
Soit maintenant $T\subset K$ un tore maximal,  
$\mathfrak{t}_0 = {\rm Lie} (T)$ et $\mathfrak{t}= \mathfrak{t}_0 \otimes {\Bbb C}$. Fixons $\Delta^+ (\mathfrak{k} , \mathfrak{t})$
un syst\`eme positif de racines dans $\Delta (\mathfrak{k} , \mathfrak{t})$. 
Posons alors
\begin{eqnarray} \label{e+}
e_H = \bigwedge^{{\rm dim}_{{\Bbb C}} (\mathfrak{s} \cap \mathfrak{p} )} (\mathfrak{s} \cap \mathfrak{p} )
\end{eqnarray}
et consid\'erons $V_H$ le plus petit sous-espace $K$-stable de $\bigwedge \mathfrak{p}$ 
(cet espace n'est en g\'en\'eral par irr\'eductible). 

\'Etant donn\'e un \'el\'ement $X \in i \mathfrak{t}_0$ tel que $\alpha (X) \geq 0$ pour toute racine $\alpha 
\in \Delta^+ (\mathfrak{k} , \mathfrak{t})$, on pose
$$\mathfrak{q} = \mathfrak{q} (X) = \mathfrak{l} \oplus \mathfrak{u}, \: \mathfrak{l}=\mathfrak{g}^X , \; \mathfrak{u} = \oplus_{
\alpha (X) >0} \mathfrak{g}_{\alpha} .$$
L'alg\`ebre $\mathfrak{q}$ est une sous-alg\`ebre parabolique $\theta$-stable de $\mathfrak{g}$. 
Notons $E(G,L)$ le sous-espace de 
$\bigwedge \mathfrak{p}$ engendr\'e par les translat\'es par $K$ du sous-espace $\bigwedge (\mathfrak{p} \cap \mathfrak{l})$ et 
toujours $R={\rm dim} (\mathfrak{u} \cap \mathfrak{p})$.

\begin{conj} \label{CONJ1}
L'application 
$$H_2^R (A_{\mathfrak{q}} : Sh^0 G) \rightarrow \prod_{g\in G({\Bbb Q})} H_{2, \; {\rm prim}+}^R (Sh^0 H)$$
obtenue en composant l'application ${\rm Res}_{2, \; H}^G$ et la projection sur la composante fortement primitive
de la cohomologie $L^2$ de $Sh^0 H$ est {\bf injective} si et seulement si l'intersection
$$V_H \cap E(G,L) \neq 0 .$$
\end{conj}

\bigskip

Consid\'erons maintenant une sous-alg\`ebre parabolique $\theta$-stable de $\mathfrak{h}$, que nous notons toujours $\mathfrak{q}$.
La sous-alg\`ebre $\mathfrak{l}$ de $\mathfrak{h}$ d\'efinit bien \'evidemment une sous-alg\`ebre de $\mathfrak{g}$. On peut
donc toujours parl\'e de $E(G,L)$.
 
\begin{conj} \label{CONJ2}
La projection de la classe $[Sh^0 H]$ dans la cohomologie $L^2$ fortement primitive $H_{2, \; {\rm prim}+}^{*} (Sh^0 G )$ est 
non nulle si et seulement si 
$${\rm rang}_{{\Bbb C}} (G / H) = {\rm rang}_{{\Bbb C}} (K / (K\cap H)) .$$
Dans ce cas, l'application
$$H_2^R ( A_{\mathfrak{q}} : Sh^0 H) \rightarrow H_{2 , \; {\rm prim}+}^{R+d} (Sh^0 G )$$
obtenue en composant l'application ``cup-produit avec $[Sh^0 H]$'' et la projection sur la composante fortement primitive
de la cohomologie $L^2$ de $Sh^0 G$ est {\bf injective} si et seulement si l'intersection
$$V_H \cap E(G,L) \neq 0 .$$
\end{conj} 

\bigskip

\noindent
{\bf Remarques.} Il n'est m\^eme pas \'evident que les Conjectures \ref{CONJ1} et \ref{CONJ2} impliquent les
r\'esultats et conjectures \'enonc\'es plus haut dans le cas des groupes unitaires et orthogonaux. Cela semble
n\'eanmoins raisonnable au vu de \cite{BergeronTentative}. Enfin on peut \'evidemment formuler une conjecture 
g\'en\'erale analogue concernant l'application de cup-produit.

\bibliography{bibliographie}

\def\cprime{$'$} \def\cprime{$'$} \def\cprime{$'$} \def\cprime{$'$}
  \def\cprime{$'$} \def\cprime{$'$}
\begin{thebibliography}{10}

\bibitem{IRMN}
N.~Bergeron.
\newblock Lefschetz properties for arithmetic real and complex hyperbolic
  manifolds.
\newblock {\em Int. Math. Res. Not.}, (20):1089--1122, 2003.

\bibitem{BergeronTentative}
Nicolas Bergeron.
\newblock Tentative d'\'epuisement de la cohomologie d'une vari\'et\'e de
  shimura par restriction \`a ses sous-vari\'et\'es.
\newblock arXiv:math.NT/0403407 v1 24 Mar 2004.

\bibitem{EnseignMath}
Nicolas Bergeron.
\newblock Premier nombre de {B}etti et spectre du laplacien de certaines
  vari\'et\'es hyperboliques.
\newblock {\em Enseign. Math. (2)}, 46(1-2):109--137, 2000.

\bibitem{STSG}
Nicolas Bergeron.
\newblock Propri\'et\'es de {L}efschetz dans la cohomologie de certaines
  vari\'et\'es arithm\'etiques: le cas des surfaces modulaires de {H}ilbert.
\newblock In {\em S\'eminaire de Th\'eorie Spectrale et G\'eom\'etrie. Vol. 21.
  Ann\'ee 2002--2003}, volume~21 of {\em S\'emin. Th\'eor. Spectr. G\'eom.},
  pages 75--101. Univ. Grenoble I, Saint, 2003.

\bibitem{Produit}
Nicolas Bergeron.
\newblock Produits dans la cohomologie des vari\'et\'es de shimura~: quelques
  calculs.
\newblock {\em C. R. Math. Acad. Sci. Paris}, 339(11):751--756, 2004.

\bibitem{BergeronClozel}
Nicolas Bergeron and Laurent Clozel.
\newblock {\em Spectre automorphe des vari\'et\'es hyperboliques et
  applications topologiques}.
\newblock arXiv:math.NT/0411385 v1 17 Nov 2004.

\bibitem{BorelHarishChandra}
Armand Borel and Harish-Chandra.
\newblock Arithmetic subgroups of algebraic groups.
\newblock {\em Ann. of Math. (2)}, 75:485--535, 1962.

\bibitem{BorelWallach}
Armand Borel and Nolan~R. Wallach.
\newblock {\em Continuous cohomology, discrete subgroups, and representations
  of reductive groups}, volume~94 of {\em Annals of Mathematics Studies}.
\newblock Princeton University Press, Princeton, N.J., 1980.

\bibitem{BurgerLiSarnak}
M.~Burger, J.-S. Li, and P.~Sarnak.
\newblock Ramanujan duals and automorphic spectrum.
\newblock {\em Bull. Amer. Math. Soc. (N.S.)}, 26(2):253--257, 1992.

\bibitem{BurgerSarnak}
M.~Burger and P.~Sarnak.
\newblock Ramanujan duals. {II}.
\newblock {\em Invent. Math.}, 106(1):1--11, 1991.

\bibitem{ClozelVenky}
L.~Clozel and T.~N. Venkataramana.
\newblock Restriction of the holomorphic cohomology of a {S}himura variety to a
  smaller {S}himura variety.
\newblock {\em Duke Math. J.}, 95(1):51--106, 1998.

\bibitem{Clozel}
Laurent Clozel.
\newblock D\'emonstration de la conjecture {$\tau$}.
\newblock {\em Invent. Math.}, 151(2):297--328, 2003.

\bibitem{Deligne}
Pierre Deligne.
\newblock Travaux de {S}himura.
\newblock In {\em S\'eminaire Bourbaki, 23\`eme ann\'ee (1970/71), Exp. No.
  389}, pages 123--165. Lecture Notes in Math., Vol. 244. Springer, Berlin,
  1971.

\bibitem{Donnelly2}
Harold Donnelly.
\newblock On the spectrum of towers.
\newblock {\em Proc. Amer. Math. Soc.}, 87(2):322--329, 1983.

\bibitem{Donnelly}
Harold Donnelly.
\newblock Elliptic operators and covers of {R}iemannian manifolds.
\newblock {\em Math. Z.}, 223(2):303--308, 1996.

\bibitem{DonnellyXavier}
Harold Donnelly and Frederico Xavier.
\newblock On the differential form spectrum of negatively curved {R}iemannian
  manifolds.
\newblock {\em Amer. J. Math.}, 106(1):169--185, 1984.

\bibitem{FlenstedJensen}
Mogens Flensted-Jensen.
\newblock Discrete series for semisimple symmetric spaces.
\newblock {\em Ann. of Math. (2)}, 111(2):253--311, 1980.

\bibitem{FlenstedJensen2}
Mogens Flensted-Jensen.
\newblock {\em Analysis on non-{R}iemannian symmetric spaces}, volume~61 of
  {\em CBMS Regional Conference Series in Mathematics}.
\newblock Published for the Conference Board of the Mathematical Sciences,
  Washington, DC, 1986.

\bibitem{Fulton}
William Fulton.
\newblock {\em Young tableaux}, volume~35 of {\em London Mathematical Society
  Student Texts}.
\newblock Cambridge University Press, Cambridge, 1997.
\newblock With applications to representation theory and geometry.

\bibitem{FultonHarris}
William Fulton and Joe Harris.
\newblock {\em Representation theory}, volume 129 of {\em Graduate Texts in
  Mathematics}.
\newblock Springer-Verlag, New York, 1991.
\newblock A first course, Readings in Mathematics.

\bibitem{GGPS}
I.~M. Gel{\cprime}fand, M.~I. Graev, and I.~I. Pyatetskii-Shapiro.
\newblock {\em Representation theory and automorphic functions}, volume~6 of
  {\em Generalized Functions}.
\newblock Academic Press Inc., Boston, MA, 1990.
\newblock Translated from the Russian by K. A. Hirsch, Reprint of the 1969
  edition.

\bibitem{HarishChandra}
Harish-Chandra.
\newblock Representations of a semisimple {L}ie group on a {B}anach space. {I}.
\newblock {\em Trans. Amer. Math. Soc.}, 75:185--243, 1953.

\bibitem{HarrisLi}
Michael Harris and Jian-Shu Li.
\newblock A {L}efschetz property for subvarieties of {S}himura varieties.
\newblock {\em J. Algebraic Geom.}, 7(1):77--122, 1998.

\bibitem{Helgason}
Sigurdur Helgason.
\newblock {\em Differential geometry, {L}ie groups, and symmetric spaces},
  volume~34 of {\em Graduate Studies in Mathematics}.
\newblock American Mathematical Society, Providence, RI, 2001.
\newblock Corrected reprint of the 1978 original.

\bibitem{Howe}
R.~Howe.
\newblock {$\theta $}-series and invariant theory.
\newblock In {\em Automorphic forms, representations and $L$-functions (Proc.
  Sympos. Pure Math., Oregon State Univ., Corvallis, Ore., 1977), Part 1},
  Proc. Sympos. Pure Math., XXXIII, pages 275--285. Amer. Math. Soc.,
  Providence, R.I., 1979.

\bibitem{Howe3}
Roger Howe.
\newblock Reciprocity laws in the theory of dual pairs.
\newblock In {\em Representation theory of reductive groups (Park City, Utah,
  1982)}, volume~40 of {\em Progr. Math.}, pages 159--175. Birkh\"auser Boston,
  Boston, MA, 1983.

\bibitem{Howe2}
Roger Howe.
\newblock Transcending classical invariant theory.
\newblock {\em J. Amer. Math. Soc.}, 2(3):535--552, 1989.

\bibitem{Knapp}
Anthony~W. Knapp.
\newblock {\em Representation theory of semisimple groups}.
\newblock Princeton Landmarks in Mathematics. Princeton University Press,
  Princeton, NJ, 2001.
\newblock An overview based on examples, Reprint of the 1986 original.

\bibitem{KobayashiNomizu}
Shoshichi Kobayashi and Katsumi Nomizu.
\newblock {\em Foundations of differential geometry. {V}ol. {II}}.
\newblock Wiley Classics Library. John Wiley \& Sons Inc., New York, 1996.
\newblock Reprint of the 1969 original, A Wiley-Interscience Publication.

\bibitem{Kobayashi2}
Toshiyuki Kobayashi.
\newblock Discrete decomposability of the restriction of {$A_{\mathfrak{q}}
  (\lambda)$} with respect to reductive subgroups and its applications.
\newblock {\em Invent. Math.}, 117(2):181--205, 1994.

\bibitem{Kobayashi}
Toshiyuki Kobayashi.
\newblock Discrete decomposability of the restriction of {$A_{\mathfrak{q}}
  (\lambda)$} with respect to reductive subgroups. {III}. {R}estriction of
  {H}arish-{C}handra modules and associated varieties.
\newblock {\em Invent. Math.}, 131(2):229--256, 1998.

\bibitem{KobayashiOda}
Toshiyuki Kobayashi and Takayuki Oda.
\newblock A vanishing theorem for modular symbols on locally symmetric spaces.
\newblock {\em Comment. Math. Helv.}, 73(1):45--70, 1998.

\bibitem{Kudla}
Stephen~S. Kudla.
\newblock Seesaw dual reductive pairs.
\newblock In {\em Automorphic forms of several variables (Katata, 1983)},
  volume~46 of {\em Progr. Math.}, pages 244--268. Birkh\"auser Boston, Boston,
  MA, 1984.

\bibitem{Kumaresan}
S.~Kumaresan.
\newblock On the canonical {$k$}-types in the irreducible unitary {$g$}-modules
  with nonzero relative cohomology.
\newblock {\em Invent. Math.}, 59(1):1--11, 1980.

\bibitem{Li}
Jian-Shu Li.
\newblock Theta lifting for unitary representations with nonzero cohomology.
\newblock {\em Duke Math. J.}, 61(3):913--937, 1990.

\bibitem{Li2}
Jian-Shu Li.
\newblock Nonvanishing theorems for the cohomology of certain arithmetic
  quotients.
\newblock {\em J. Reine Angew. Math.}, 428:177--217, 1992.

\bibitem{LiMillson}
Jian-Shu Li and John~J. Millson.
\newblock On the first {B}etti number of a hyperbolic manifold with an
  arithmetic fundamental group.
\newblock {\em Duke Math. J.}, 71(2):365--401, 1993.

\bibitem{Littlewood}
D.~E. Littlewood.
\newblock Invariant theory, tensors and group characters.
\newblock {\em Philos. Trans. Roy. Soc. London. Ser. A.}, 239:305--365, 1944.

\bibitem{MargulisSoifer}
G.~A. Margulis and G.~A. So{\u\i}fer.
\newblock Maximal subgroups of infinite index in finitely generated linear
  groups.
\newblock {\em J. Algebra}, 69(1):1--23, 1981.

\bibitem{Matsushima}
Yoz{\^o} Matsushima.
\newblock A formula for the {B}etti numbers of compact locally symmetric
  {R}iemannian manifolds.
\newblock {\em J. Differential Geometry}, 1:99--109, 1967.

\bibitem{MazzeoPhillips}
Rafe Mazzeo and Ralph~S. Phillips.
\newblock Hodge theory on hyperbolic manifolds.
\newblock {\em Duke Math. J.}, 60(2):509--559, 1990.

\bibitem{MillsonRaghunathan}
John~J. Millson and M.~S. Raghunathan.
\newblock Geometric construction of cohomology for arithmetic groups. {I}.
\newblock In {\em Geometry and analysis}, pages 103--123. Indian Acad. Sci.,
  Bangalore, 1980.

\bibitem{OshimaMatsuki}
Toshio {\=O}shima and Toshihiko Matsuki.
\newblock A description of discrete series for semisimple symmetric spaces.
\newblock In {\em Group representations and systems of differential equations
  (Tokyo, 1982)}, volume~4 of {\em Adv. Stud. Pure Math.}, pages 331--390.
  North-Holland, Amsterdam, 1984.

\bibitem{Parthasarathy}
R.~Parthasarathy.
\newblock Criteria for the unitarizability of some highest weight modules.
\newblock {\em Proc. Indian Acad. Sci. Sect. A Math. Sci.}, 89(1):1--24, 1980.

\bibitem{RaghunathanVenky}
M.~S. Raghunathan and T.~N. Venkataramana.
\newblock The first {B}etti number of arithmetic groups and the congruence
  subgroup problem.
\newblock In {\em Linear algebraic groups and their representations (Los
  Angeles, CA, 1992)}, volume 153 of {\em Contemp. Math.}, pages 95--107. Amer.
  Math. Soc., Providence, RI, 1993.

\bibitem{Rohlfs}
J.~Rohlfs.
\newblock Projective limits of locally symmetric spaces and cohomology.
\newblock {\em J. Reine Angew. Math.}, 479:149--182, 1996.

\bibitem{Sakai}
Takashi Sakai.
\newblock {\em Riemannian geometry}, volume 149 of {\em Translations of
  Mathematical Monographs}.
\newblock American Mathematical Society, Providence, RI, 1996.
\newblock Translated from the 1992 Japanese original by the author.

\bibitem{Schlichtkrull}
Henrik Schlichtkrull.
\newblock The {L}anglands parameters of {F}lensted-{J}ensen's discrete series
  for semisimple symmetric spaces.
\newblock {\em J. Funct. Anal.}, 50(2):133--150, 1983.

\bibitem{SpehVenky}
B.~Speh and T.N. Venkataramana.
\newblock Construction of some generalised modular symbols.
\newblock arXiv:math.GR/0409376 v1 21 sept 2004.

\bibitem{Taylor}
Michael~E. Taylor.
\newblock {\em Partial differential equations. {III}}, volume 117 of {\em
  Applied Mathematical Sciences}.
\newblock Springer-Verlag, New York, 1997.
\newblock Nonlinear equations, Corrected reprint of the 1996 original.

\bibitem{TongWang2}
Y.~L. Tong and S.~P. Wang.
\newblock Harmonic forms dual to geodesic cycles in quotients of {${\rm
  SU}(p,\,1)$}.
\newblock {\em Math. Ann.}, 258(3):289--318, 1981/82.

\bibitem{TongWang}
Y.~L. Tong and S.~P. Wang.
\newblock Geometric realization of discrete series for semisimple symmetric
  spaces.
\newblock {\em Invent. Math.}, 96(2):425--458, 1989.

\bibitem{Venky}
T.~N. Venkataramana.
\newblock Cohomology of compact locally symmetric spaces.
\newblock {\em Compositio Math.}, 125(2):221--253, 2001.

\bibitem{Vogan}
D.~A. Vogan.
\newblock Isolated unitary representations.
\newblock to appear in the 2002 Park City summer school volume.

\bibitem{Vogan3}
David~A. Vogan, Jr.
\newblock The algebraic structure of the representation of semisimple {L}ie
  groups. {I}.
\newblock {\em Ann. of Math. (2)}, 109(1):1--60, 1979.

\bibitem{VoganZuckerman}
David~A. Vogan, Jr. and Gregg~J. Zuckerman.
\newblock Unitary representations with nonzero cohomology.
\newblock {\em Compositio Math.}, 53(1):51--90, 1984.

\bibitem{Wallach}
Nolan~R. Wallach.
\newblock {\em Harmonic analysis on homogeneous spaces}.
\newblock Marcel Dekker Inc., New York, 1973.
\newblock Pure and Applied Mathematics, No. 19.

\bibitem{Wang}
S.~P. Wang.
\newblock Correspondence of modular forms to cycles associated to {${\rm
  O}(p,q)$}.
\newblock {\em J. Differential Geom.}, 22(2):151--213, 1985.

\bibitem{Yeganefar}
Nader Yeganefar.
\newblock Sur la {$L\sp 2$}-cohomologie des vari\'et\'es \`a courbure
  n\'egative.
\newblock {\em Duke Math. J.}, 122(1):145--180, 2004.

\bibitem{Yeung}
Sai-Kee Yeung.
\newblock Betti numbers on a tower of coverings.
\newblock {\em Duke Math. J.}, 73(1):201--226, 1994.

\bibitem{Zucker}
Steven Zucker.
\newblock {$L\sb{2}$} cohomology of warped products and arithmetic groups.
\newblock {\em Invent. Math.}, 70(2):169--218, 1982/83.

\end{thebibliography}

\bibliographystyle{plain}

\bigskip

\noindent
Unit\'e Mixte de Recherche 8628 du CNRS, \\
Laboratoire de Math\'ematiques, B\^at. 425, \\
Universit\'e Paris-Sud, 91405 Orsay Cedex, France \\
{\it adresse electronique :} \texttt{Nicolas.Bergeron@math.u-psud.fr}

\end{document}